\documentclass[11pt,a4paper]{article}
\usepackage{amsmath}
\usepackage{amssymb}
\usepackage{theorem}
\usepackage{euscript}
\usepackage{fancyhdr}
\usepackage[dvipsnames]{xcolor}
\usepackage{enumitem}
\usepackage{hyperref}
\hypersetup{colorlinks=true,citecolor=Blue,linkcolor=Blue}

\usepackage{tikz}
\usetikzlibrary{arrows}

\renewcommand{\title}[1]{
  \addvspace{3\baselineskip}  
  \begin{center} \LARGE \bf #1
  \end{center}
  \addvspace{2\baselineskip}}   

\renewcommand{\author}[1]{
  \addvspace{-1\baselineskip}  
  \begin{center} \large \sc #1
  \end{center}
  \addvspace{2\baselineskip}}   

\makeatletter

\def\section{%
        \@startsection{section}{1}{\z@}%
        {8ex plus 6ex minus 3ex}{\baselineskip}%
        {\normalfont\large\scshape\centering}%
        }
\def\@seccntformat#1{\csname the#1\endcsname.\hspace*{0.5em}}

\renewcommand{\paragraph}[1]{{\par\removelastskip\vskip.5\baselineskip%
         \indent{\bfseries{#1}}{\ifperiod.\else\global\periodtrue\fi}%
         \rm \ignorespaces}}

\let\goth=\mathfrak
\let\calligraphy=\mathcal

\def\CC{{\mathbb C}}

\def\NN{{\mathbb N}}

\def\PP{{\mathbb P}}
\def\QQ{{\mathbb Q}}
\def\RR{{\mathbb R}}

\def\Cc{{\calligraphy C}}

\def\Ee{{\calligraphy E}}
\def\Ff{{\calligraphy F}}
\def\Gg{{\calligraphy G}}

\def\Ii{{\calligraphy I}}
\def\Jj{{\calligraphy J}}

\def\Ll{{\calligraphy L}}
\def\Mm{{\calligraphy M}}
\def\Nn{{\calligraphy N}}
\def\Oo{{\calligraphy O}}
\def\Pp{{\calligraphy P}}
\def\Qq{{\calligraphy Q}}

\def\Ss{{\calligraphy S}}
\def\Tt{{\calligraphy T}}
\def\Uu{{\calligraphy U}}

\def\aaa{{\goth a}}

\def\DDD{{\goth D}}

\def\FFF{{\goth F}}

\def\III{{\goth I}}

\def\nnn{{\goth n}}

\def\PPP{{\goth P}}

\def\hbar{{\,\overline{\!h}}}

\def\htilde{{\,\widetilde{\!h}}}

\def\stilde{{\,\widetilde{\!s}}}

\def\phitilde{{\,\widetilde{\!\phi}}}
\def\Alb{\operatorname{Alb}}

\def\coker{\operatorname{coker}}

\def\End{\operatorname{End}}

\def\Ext{\operatorname{Ext}}

\def\Gr{\operatorname{Gr}}

\def\im{\operatorname{im}}

\def\Jac{\operatorname{Jac}}

\newcommand{\NS}{\operatorname{NS}}

\def\Pic{\operatorname{Pic}}

\def\rank{\operatorname{rank}}

\def\Span{\operatorname{Span}}

\def\Sing{\operatorname{Sing}}

\def\Sym{\operatorname{Sym}}

\def\Bigwedge{%
  \tikz[baseline=.3ex]{\draw (0ex,0ex) -- (.8ex,2.15ex) -- (1.6ex,0ex);}
}
\def\contraction{%
  \tikz[baseline=.35ex]{\draw (0ex,0ex) -- (.9ex,0) -- (.9ex,2ex);}\hspace{.4ex}
}

\let\lra=\longrightarrow

\def\ie{{\it i.e.,}~}
\def\eg{{\it e.g.,}~}
\def\cf{{\it cf.}~}

\def\inv{^{-1}}
\def\dual{^\vee}

\let\phi=\varphi
\let\epsilon=\varepsilon

\newcommand{\textq}[1]{\quad\text{#1}\quad}

\newcommand{\listspace}{\setlength{\itemsep}{-.1\baselineskip}}

\newcounter{icounter}

\def\theoname{Theorem}
\def\lemmaname{Lemma}
\def\propositionname{Proposition}
\def\notationname{Notation}
\def\corollaryname{Corollary}
\def\conjecturename{Conjecture}
\def\remarkname{Remark}
\def\remarksname{Remarks}
\def\examplename{Example}
\def\examplesname{Examples}
\def\definitionname{Definition}
\def\definitionsname{Definitions}
\def\notationname{Notation}

\def\proofname{Proof}

\def\Dquad{\hskip 0.6em plus .02em minus .2em}  
\def\Dpar{\belowdisplayskip=0pt\belowdisplayshortskip=0pt\par}

\def\bigpenalty{\interlinepenalty=\@M}
\def\smallpenalty{\interlinepenalty=100}

\newif\ifperiod \periodtrue 

\def\D@makemargins{%
  \labelsep=0pt
  \itemindent=0pt
  \labelwidth=0pt}

\def\D@restoremargins{%
  \labelsep=5pt
  \itemindent=0pt
  \leftmargin=5mm  
  \labelwidth=\leftmargin \advance\labelwidth by -\labelsep}

\def\th@Dindent{\hspace\parindent}
\def\th@Dheadingshape{\scshape}

\gdef\th@DthAndSuchtheo{%
  \D@makemargins%
  \def\@begintheorem##1##2{%
  \item[]\th@Dindent{\th@Dheadingshape ##1~\rm ##2.}\Dquad         
        \D@restoremargins}%
  \def\@opargbegintheorem##1##2##3{\def\next{##3}%
  \item[]\th@Dindent{\th@Dheadingshape ##1~\rm ##2\ifx\next\empty
  \else\ {\normalfont(##3)}\fi.}         
        \D@restoremargins}}

\gdef\th@DthAndSuchtheostar{%
  \D@makemargins%
  \def\@begintheorem##1##2{%
  \item[]\th@Dindent{\th@Dheadingshape ##1.}\Dquad     
        \D@restoremargins}%
  \def\@opargbegintheorem##1##2##3{\def\next{##3}%
  \item[]\th@Dindent{\th@Dheadingshape ##1\ifx\next\empty
  \else\ ##3\fi.}\Dquad         
        \D@restoremargins}}

\gdef\th@DthAndSuchliketheo{
  \D@makemargins%
  \def\@begintheorem##1##2{%
    \@latex@error{likethm: You must provide an argument in square brackets,
    though it may be empty [] !}%
    }%
  \def\@opargbegintheorem##1##2##3{%
        \def\next{##3}\ifx\next\empty\item[\th@Dindent]\else
        \item[]\th@Dindent{\th@Dheadingshape \next.}\Dquad\fi
        \D@restoremargins}}

\gdef\th@DdefAndSuch{%
  \D@makemargins%
  \def\@begintheorem##1##2{%
  \item[]\th@Dindent{\def@Dheadingshape ##1~\rm ##2.}\Dquad         
        \D@restoremargins}%
  \def\@opargbegintheorem##1##2##3{\def\next{##3}%
  \item[]\th@Dindent{\def@Dheadingshape ##1~\rm ##2\ifx\next\empty
  \else\ {\normalfont(##3)}\fi.}         
        \D@restoremargins}}

\gdef\th@DdefAndSuchStar{%
  \D@makemargins%
  \def\@begintheorem##1##2{%
  \item[]\th@Dindent{\def@Dheadingshape ##1.}\Dquad     
        \D@restoremargins}%
  \def\@opargbegintheorem##1##2##3{\def\next{##3}%
  \item[]\th@Dindent{\th@Dheadingshape ##1\ifx\next\empty
  \else\ ##3\fi.}\Dquad         
        \D@restoremargins}}

\def\th@Dheadingshape{\scshape}
\def\def@Dheadingshape{\itshape}

\theoremstyle{DthAndSuchliketheo}
\theorembodyfont{\bigpenalty\itshape}   
\newtheorem{likethm}{}
\theorembodyfont{\rmfamily}  

\theoremstyle{DthAndSuchtheostar}
\theorembodyfont{\bigpenalty\itshape}
\newtheorem{thm*}{\theoname}
\newtheorem{lem*}{\lemmaname}
\newtheorem{pro*}{\propositionname}
\newtheorem{cor*}{\corollaryname}
\newtheorem{conjecture*}{\conjecturename}
\theorembodyfont{\smallpenalty\rmfamily}
\newtheorem{notation*}{\notationname}
\newtheorem{exa*}{\examplename}
\newtheorem{examples*}{\examplesname}
\theoremstyle{DdefAndSuchStar}
\theorembodyfont{\smallpenalty\rmfamily}
\newtheorem{definition*}{\definitionname}
\newtheorem{definitions*}{\definitionsname}
\newtheorem{rem*}{\remarkname}
\newtheorem{remarks*}{\remarksname}
\theoremstyle{DthAndSuchtheo}
\theorembodyfont{\bigpenalty\itshape}
\newtheorem{thm}{\theoname}[section]
\newtheorem{lem}[thm]{\lemmaname}
\newtheorem{pro}[thm]{\propositionname}
\newtheorem{cor}[thm]{\corollaryname}

\theorembodyfont{\smallpenalty\rmfamily}
\newtheorem{exa}[thm]{\examplename}

\theoremstyle{DdefAndSuch}
\theorembodyfont{\smallpenalty\rmfamily}

\newtheorem{rem}[thm]{\remarkname}

\newtheorem{remarks}[thm]{\remarksname}

\theoremstyle{DthAndSuchtheo}               
\theorembodyfont{\itshape}

\newcommand{\proof}[1][]{{\par\removelastskip\vskip.6\baselineskip   
    \noindent\th@Dindent\def\next{#1}%
    {\itshape\proofname\ifx\next\empty\else\next\fi\ifperiod.%
      \else\global\periodtrue\fi\Dquad}%
    \clubpenalty=5000\rm\ignorespaces}\setcounter{step}{0}\setcounter{case}{0}}

\newcounter{step}
\newcommand{\step}{\stepcounter{step}
  \par\indent{\itshape Step \thestep .\hspace{1ex}}}

\newcounter{case}

\def\case{\stepcounter{case}%
  \def\mynext{\thecase}
  \ifnum\mynext=1{\bfseries First case.}
  \else
    \ifnum\mynext=2{\bfseries Second case.}
    \else
      \ifnum\mynext=3{\bfseries Third case.}
      \else
        \ifnum\mynext=4{\bfseries Fourth case.}
        \fi      
      \fi
    \fi
  \fi
  \hspace{1ex}}

\newcommand{\likeproof}[1][]{{\par\removelastskip\vskip.6\baselineskip
    \noindent\th@Dindent\def\next{#1}%
    {\itshape\ifx\next\empty\else\next\fi\ifperiod.%
      \else\global\periodtrue\fi\Dquad}%
    \clubpenalty=5000\rm\ignorespaces}\hspace{-2pt}\setcounter{step}{0}}

\def\qed{{\ifmmode\hskip 6mm plus 1mm minus 3mm{$\square$}
    \else
    \hfil\hskip1em\null\nobreak\hfil
    {\clubpenalty=50000\bigpenalty\hfill
      $\square$\parfillskip=0pt\finalhyphendemerits=0
      \let\par=\endgraf\par}
    \fi
    \Dpar\penalty-150\vskip.6\normalbaselineskip}}

\makeatother

\begin{document}
\numberwithin{equation}{section}
\titlepage
\title{Surfaces in $\PP^4$ lying on small degree hypersurfaces}

\author{Daniel Naie and Igor Reider}

\begin{abstract}
Since the work of Ellingsrud and Peskine at the end of 1980s, it has
been known that smooth compact complex surfaces in $\PP^4$ with
prescribed Chern classes, with the exception of a finite number of
families, must lie on hypersurfaces of degree $m\leq 5$.  Hence the
motivation for the present work: to study smooth surfaces contained in
a hypersurface of degree $m\leq 5$ (the meaning of `small degree' in
the title). There are two main issues considered in the paper:
\begin{enumerate}[label={\rm (\Roman*)}]
  \listspace
\item
  an analogue of the Hartshorne-Lichtenbaum finiteness results for
  smooth surfaces of general type contained in a small degree
  hypersurface in $\PP^4$, 
\item
  a study of the irregularity of smooth surfaces contained in a small
  degree hypersurface in $\PP^4$.
\end{enumerate}

For (I) we show that for $m\leq 4$, the number of families is
controlled by a function depending on the {\it ratio} 
$\displaystyle{\alpha=\frac{K^2}{\chi}}$ of the Chern invariants 
$(K^2, \chi)$ of surfaces. The same result holds for $m=5$, with a 
possible exception of $\alpha=6$. 

For (II) we determine all irregular surfaces contained in a
hypersurface of degree $m\leq 3$.  We do the same in case $m=4$, under
the additional assumption that a quartic hypersurface has only
isolated double points.  In general, we show that the Albanese
dimension of surfaces contained in quartic hypersurfaces is at most
$1$.

For $m=5$, we show that minimal surfaces of Albanese dimension $2$ have
the irregularity at most $3$ and describe the hypothetical surfaces
with irregularity $3$.

Conceptually, the main idea underlying the above results as well as the 
whole approach of our paper can be termed as a representation of 
various geometric and cohomological entities attached to a surface 
in $\PP^4$ in the category of coherent sheaves on that surface.
\end{abstract}


\tableofcontents

\section{Introduction} 

Smooth compact complex surfaces in $\PP^4$ constitute an interesting
and important part of the study of subvarieties of projective spaces.
They are naturally situated at the cross-road of the theories of
surfaces, vector bundles, and algebraic cycles.  It is well-known that
every smooth projective surface can be embedded into $\PP^5$.  To fit
into $\PP^4$, a surface must satisfy an obstruction, known as double
point formula,
\begin{equation}
  \label{eq:dPFormula}
  d^2 -5d -10(g-1) +(c_2-K_X^2) = 0,
\end{equation}
which ties together the degree $d$, the sectional genus $g$, and the
basic topological invariants (the Chern numbers) $K^2$ and $c_2$ of
the surface.  So one of the basic goals, which is still out of reach,
is to find all the surfaces that can be embedded into $\PP^4$.

Another line of inquiry into the geometry of surfaces in $\PP^4$ has
been motivated by a conjecture of Hartshorne and Lichtenbaum stating
that rational surfaces in $\PP^4$ form a finite number of families.
The work of Ellingsrud and Peskine, \cite{ElPe}, solved a more general
problem of finiteness of families of surfaces not of general type in
$\PP^4$.  A key observation of \cite{ElPe} is that a surface $X$ lying
in a hypersurface of degree $m$ in $\PP^4$ has the holomorphic Euler
characteristic $\chi(\Oo_X)$ bounded from below by a certain cubic
polynomial $P_m(d)$ in the degree $d$ of $X$.  Since then, this
result has been greatly improved and clarified.  Notably, Decker and
Schreyer's work, \cite{DeSch}, gives a precise expression for
$P_m(d)$,
\begin{equation}
  \label{eq:polynomialPm}
  P_m(d)
  = m\binom{\frac{d}{m} +\frac{m-3}{2}}{3}
  -\frac{(m-1)^2}{2m}\,d(d-3) -\binom{m-1}{4} +1,
\end{equation}
where $m=m_X$ is the smallest degree of a hypersurface in $\PP^4$
containing $X$,
\begin{equation}
  \label{eq:mX}
  m_X := \min \{k \in \NN \mid h^0(\Jj_X(k))>0 \}.
\end{equation}
Then, in \cite{DeSch}, it is shown that 
\begin{equation}
\label{eq:chiPm}
  \chi(\Oo_X) \geq P_m(d),
\end{equation}
provided $d \geq (m-1)^2+2$.  This result together with the bound on
the sectional genus of \cite{ElPe} implies

\begin{thm}[Ellingsrud and Peskine; Decker and Schreyer]
  \label{th:ff} 
Given integers $\chi$ and $m\geq 2$, the surfaces in $\PP^4$ having
holomorphic Euler characteristic $\chi$ and lying on a hypersurface of
degree $m$ and not on one of a smaller degree form at most a finite
number of families.
\end{thm}

Therefore, the ``world'' of surfaces in $\PP^4$ is governed by the
pairs of integers $(\chi,m)$ as in Theorem~\ref{th:ff} and emphasis is
placed on the understanding of surfaces that can lie on a hypersurface
of given degree.  From this point of view, the study of surfaces lying
on a small degree hypersurface in $\PP^4$---{\it small} meaning 
$m\in \{2,3,4,5\}$---would appear, and could perhaps be justified, as
a way of obtaining empirical data leading to a better conceptual
understanding of surfaces in $\PP^4$.  But in fact, taking the ideas
of \cite{ElPe} and \cite{DeSch} a bit further, one can argue that
these small degrees really matter precisely for their conceptual
significance.  To explain this point we fix a pair of integers 
$(K^2,\chi)$ and observe

\begin{pro}[\cite{NaRe}]
  \label{p:Chern-d} 
There exists a number $d(K^2,\chi)$, depending only on $K^2$ and
$\chi$, such that all surfaces in $\PP^4$ with Chern numbers 
$(K^2,\chi)$ and degree $d>d(K^2,\chi)$ lie on a hypersurface of
degree $\leq 5$.
\end{pro}

This result tells us that, with a possible exception of finitely many
families, the study of surfaces in $\PP^4$ having {\it prescribed}
Chern invariants comes down to understanding surfaces lying on
hypersurfaces of degree $m\in\{2,3,4,5\}$.  Extrapolating further
Proposition \ref{p:Chern-d}, we suggest

\begin{likethm}[Metha-principle]
  \label{th:principle}
An understanding of a property $\Pp$ for surfaces in $\PP^4$, with the
exception of a finite number of families, comes down to studying the
property $\Pp$ for surfaces contained in hypersurfaces of degree 
$m\leq5$.
\end{likethm}

After clarifying the origins and motivations for studying surfaces in
$\PP^4$ lying on hypersurfaces of small degree let us give an overview
of the main results of this paper.

\paragraph{Main results of the paper} 
There are two main issues considered in this work:
\begin{enumerate}[label={\rm (\Roman*)}]
  \listspace
\item
  an analogue of Hartshorne-Lichtenbaum finiteness results for smooth
  surfaces of general type contained in a small degree hypersurface in
  $\PP^4$, 
\item
  a study of the irregularity of smooth surfaces contained in a small
  degree hypersurface in $\PP^4$.
\end{enumerate}

We approach (I) as students of the theory of surfaces of general type.
To be more precise, let us recall that one of the main problems of
that theory is the ``geography'' problem: characterize the pairs of
integers $(k,c)$ which are respectively $K^2_X$ and $\chi(\Oo_X)$ of
some minimal surface $X$ of general type.  The integers
$\chi(\Oo_X)$ and $K^2_X$ are often referred to as the
Chern numbers of $X$---terminology\footnote{In view of the Noether
formula $12\chi=K_X^2+c_2$, either $(K_X^2,c_2)$ or $(K_X^2,\chi)$
will be referred to as Chern numbers of $X$.} that we adopt in this
paper---and their ratio $\alpha_X := K^2_X/\chi(\Oo_X)$ (called the
slope in the sequel) provides many important dividing lines in this
$2$-dimensional world of Chern invariants of surfaces of general type.
In this notation we often omit the reference to $X$, if no ambiguity
is likely.

It is well known that for a surfaces $X$ of general type
$\chi(\Oo_X)>0$, so it makes sense to speak about the slope $\alpha_X$
of the Chern numbers even when the surface $X$ is not minimal.  This
is what we do for surfaces of general type in $\PP^4$.  It turns out
that the number of families of such surfaces contained in a small
degree hypersurface can be controlled only by the value of the slope
$\alpha$.  With the notation \eqref{eq:mX} in mind, we can formulate a
sample result concerning part (I).

\begin{thm}
  \label{th:ISample}
For all surfaces $X$ of general type in $\PP^4$ with $m_X \leq 4$ the
following assertions hold.
\begin{enumerate}[label={\rm \arabic*)}]
  \listspace
\item 
  The slope $\alpha$ of the Chern numbers is smaller than $6$. 
\item
  For every rational number $\alpha<6$, there exists an integer
  $d(\alpha)$ such that every surface of slope $\alpha$ has the degree
  $\leq d(\alpha)$.
\item
  For every rational number $\alpha<6$ and integer $d\leq d(\alpha)$, 
  there exists $\chi(\alpha,d)$ such that every surface of slope 
  $\alpha$ and degree $d$ has the holomorphic Euler characteristic 
  $\leq\chi(\alpha,d)$.
\end{enumerate}
\end{thm}

A similar but somewhat more involved statement holds for surfaces $X$
with $m_X=5$, see Proposition \ref{p:d-chi-finite-5}.  It should be
also pointed out that the expressions $d(\alpha)$ and $\chi(\alpha,d)$
in the above theorem are effectively computable.  For example,
$\chi(\alpha,d)$ is an explicit rational function of $\alpha$ and $d$.

The above results, in essence, are obtained by a combination of two
ingredients: the bound \eqref{eq:polynomialPm} of Decker and Schreyer,
and inequalities of the form
\begin{equation}
  \label{eq:ineqOnChernNumbers}
  c_2 -K^2 \geq a\, H\cdot K_X +bd,
\end{equation}
where $\Oo_X(H)$ is the line bundle embedding $X$ into $\PP^4$, $d$ is
the degree of $X$, and $a$ and $b$ are positive rational numbers
(depending on $m_X$ and explicitly determined in the main body of the
paper, see Theorem~\ref{th:mainTechnical}, Theorem~\ref{th:m=2} and
Theorem~\ref{th:m=3}).  So our contribution to the
Hartshorne-Lichtenbaum problem for surfaces of general type are the
inequalities \eqref{eq:ineqOnChernNumbers} above.  Of course, one most
certainly wonders where those inequalities come from.  This and other
results of the paper will be explained shortly.  For now, let us just
say that the existence of these inequalities is an {\it a priori}
consequence of our approach toward the study of surfaces in $\PP^4$
contained in a hypersurface of a small degree.

We now turn to the results concerning the issue (II), the irregularity
of smooth surfaces contained in a small degree hypersurface in
$\PP^4$.

\begin{thm}
  \label{th:irreg}
Let $X\subset\PP^4$ be a smooth surface and $m_X$ be the smallest
degree of a hypersurface containing $X$.  
\begin{enumerate}[label={\rm \arabic*)}]
  \listspace
\item 
  If $m_X=2$, then $X$ is regular.
\item
  If $m_X=3$ and $X$ is irregular, then $X$ is an elliptic scroll of
  degree $d=5$.  Furthermore, a general cubic hypersurface containing
  $X$ is a Segre cubic\footnote{Such a cubic has ten nodes, the
    maximal possible number of nodes for a cubic hypersurface in
    $\PP^4$.} and the surface $X$
  must pass through the ten nodes of every Segre cubic containing it.  
\item
  If $m_X=4$, then the Albanese dimension of $X$ is at most $1$.
  Furthermore, if $X$ lies on a quartic hypersurface with only
  ordinary double points, then $X$ is regular, with a possible
  exception of $X$ being an elliptic conic bundle of degree
  $d=-K^2_X=8$. 
\item
  If $m_X=5$, $X$ is minimal, and its Albanese dimension is $2$, then
  the irregularity of $X$ is $2$ or $3$.
\end{enumerate}
\end{thm}

The above statements illustrate several objectives of our inquiry
about the irregularity of surfaces in $\PP^4$:
\begin{enumerate}[label={\rm \alph*)}]
  \listspace
\item 
  determine all irregular surfaces for a given value of $m_X$,
\item
  for a given $m_X$, determine all possible values of the Albanese
  dimension of irregular surfaces,
\item
  determine the upper bound on the irregularity for every value of the
  Albanese dimension that may occur.
\end{enumerate}
Statements 1) and 2) of Theorem \ref{th:irreg} are examples of a),
while 3) and 4) are partial answers to b) and c).

\paragraph{The outline of our approach}
In the rest of the introduction we discuss our approach to the study
of smooth surfaces contained in a small degree hypersurface in
$\PP^4$.  The main idea consists of interpreting the extrinsic datum
of a small degree hypersurface as an intrinsic one.  This is done
first, by thinking of a hypersurface of the minimal degree $m=m_X$
containing a surface $X\subset\PP^4$ as a nonzero global section of
$\Nn_X(-K_X-(5-m)H)$, the normal bundle $\Nn_X=\Nn_{X/\PP^4}$ tensored
with $\Oo_X(-K_X-(5-m)H)$, and second, by attaching to that global
section a cohomology class, call it $\xi$, in
$H^1(\Theta_X(-K_X-(5-m)H))$, where $\Theta_X$ is the holomorphic
tangent bundle of $X$.  The last step is achieved via the coboundary
homomorphism
\[
  H^0(\Nn_X(-K_X-(5-m)H)) \lra H^1(\Theta_X(-K_X-(5-m)H))
\]
coming from the normal exact sequence of $X\subset\PP^4$ tensored with
$\Oo_X(-K_X-(5-m)H)$.  Next, via the natural identification
$H^1(\Theta_X(-K_X-(5-m)H))\cong\Ext^1(\Omega_X,\Oo_X(-K_X-(5-m)H))$,
we interpret the cohomology class $\xi$ as the corresponding extension
\begin{equation}
  \label{eq:extensionSequence}
  0 \lra  \Oo_X(-K_X-(5-m)H) \lra \Tt_\xi \lra \Omega_X \lra 0.
\end{equation}
The above extension sequence can be viewed as a reincarnation of a
hypersurface of degree $m$ containing $X$ in the category of complexes
of coherent sheaves on $X$.  It is clear that it can be viewed as an
independent entity and the understanding of its properties constitutes
an important part of our approach.  Namely, our strategy is to extract
as much information as possible from \eqref{eq:extensionSequence}
independently of the fact that $X$ is embedded into $\PP^4$ and then
to use the acquired data to gain an additional insight into the
embedding $X\subset\PP^4$.  As an example, let us take up the question
of the (semi)stability of the sheaf $\Tt_\xi$.  This is something {\it
a priori} independent of the fact that $X$ lies on a hypersurface of
degree $m$ in $\PP^4$, and it gives, by the Bogomolov-Gieseker
inequality, the constraint on the Chern invariants of $\Tt_\xi$,
\[
  3c_2(\Tt_\xi) \geq c^2_1(\Tt_\xi).
\]
This together with the Chern invariants of $\Tt_\xi$, determined
from the defining sequence \eqref{eq:extensionSequence}, provide a
prototypical example for the inequalities in
\eqref{eq:ineqOnChernNumbers}, one of the ingredients to prove Theorem
\ref{th:ISample}.  Of course, there is no reason for $\Tt_\xi$ to be
semistable.  However, even in the unstable case one has a sufficient
control of the destabilizing filtration of $\Tt_\xi$ to recapture the
spirit of those inequalities.  Though this constitutes a somewhat
technical part of our considerations, the main point is quite
transparent: the control of the properties of the destabilizing
filtration of $\Tt_\xi$ is enabled by
\begin{itemize}
  \listspace
\item 
  the fact that $\Tt_\xi$ is a part of the extension sequence
  \eqref{eq:extensionSequence} and in particular, that the cotangent
  sheaf $\Omega_X$ is a quotient bundle of $\Tt_\xi$, and
\item
  the geometric origin of the extension sequence
  \eqref{eq:extensionSequence} which allows one to relate certain
  properties of the destabilizing filtration of $\Tt_\xi$ to the
  embedding of $X$ in $\PP^4$.
\end{itemize}
These remarks indicate that, though the semistable case provides the
strongest form of the inequalities \eqref{eq:ineqOnChernNumbers}, it
is the unstable case that is more interesting.  Not only a numerical
constraint in the form of \eqref{eq:ineqOnChernNumbers} is obtained,
but a certain amount of geometric data encoded in the destabilizing
filtration of $\Tt_\xi$ is gained.

The results of Theorem \ref{th:ISample} and Proposition
\ref{p:d-chi-finite-5} exploit only the numerical part of the study of
\eqref{eq:extensionSequence}.  In this respect, the problem of the
irregularity considered in (II) allows to reveal some more geometric
aspects of our approach.  Namely, the destabilizing subsheaf of
$\Tt_\xi$, call it $\Gg$, is defined as the saturation of the subsheaf
of $\Tt_\xi$ generated by its global sections.  These sections are
connected to the homomorphism
\[
  H^0(\Tt_\xi) \lra H^0(\Omega_X) 
\]
induced by the epimorphism in \eqref{eq:extensionSequence}.  It takes
no effort to work out conditions for the above homomorphism to be
surjective in the case $m\leq 4$. When $m=5$, we need the hypothesis
of minimality for $X$ for our approach to go through.  This
consideration, for example, immediately implies that for $m\leq 4$
the rank of $\Gg$ is at most $2$ and the proof of the statement 3) of
Theorem \ref{th:irreg} comes down to ruling out the possibility
$\rank(\Gg)=2$.

The geometric destabilization of $\Tt_\xi$ just outlined works well as
long as the cotangent bundle $\Omega_X$ is generically generated by
its global sections, \ie when the Albanese dimension of $X$ is $2$.
In particular, it gives the result on the possible values of the
irregularity in Theorem \ref{th:irreg} as well as establishes a short
list of hypothetical surfaces with irregularity $3$, see Theorem
\ref{th:V5Alb2}.

However, for surfaces of Albanese dimension $1$, the above approach
fails.  This brings us to the second way of associating an extension
sequence to a reduced irreducible hypersurface of degree $m$
containing a surface $X\subset\PP^4$.  Here again we think of such a
hypersurface as a nonzero global section of
\hspace{1ex}$\Nn_X(-K_X-(5-m)H)$ and then take the Koszul sequence
associated to it to obtain the extension
\begin{equation}
  \label{intro:KoszVm}
  0 \lra \Oo_X(K_X +(5-m)H) \lra \Nn_X \lra \Jj_Z(mH) \lra 0.
\end{equation}
The sequence, in general, is not exact, but let us ignore this for now
and assume it is\footnote{The interested reader will find a more
  detailed discussion of this approach in the introduction of \S
  \ref{sec:irreg}.}. 
Its pertinence to the problem of irregularity of $X$ comes from the
identifications 
\[
  H^1(\Nn_X(K_X))
  \cong H^1(\Nn_X^\ast)^\ast \cong H^0(\Omega_X)^\ast,
\]
where the first isomorphism is the Serre duality and the second is a
general fact valid for any smooth subvariety of dimension at least $2$
in a projective space. From this and \eqref{intro:KoszVm} tensored
with $\Oo_X(K_X)$ it follows that the the irregularity of $X$ is
controlled by two cohomology groups $H^1(\Oo_X(2K_X +(5-m)H)$ and
$H^1(\Jj_Z(K_X+mH))$.  For $m\leq 5$, the nonvanishing of the first
group gives rather strong restrictions on $X$.  The understanding of
the nonvanishing of the second group depends largely on the knowledge
of the subscheme $Z$ which could be related to the singular locus of
the hypersurface we started with.  It is clear that this approach can
only work when a good understanding of the singular locus of the
hypersurface in question is available.  This is the case for $m=3$,
\ie for cubic hypersurfaces, when everything can be analyzed
completely leading to the elliptic scroll as the only irregular
surface with $m=3$; see Theorem \ref{th:irreg}, 2).

\paragraph{Relation to other works} 
The subject of surfaces in $\PP^4$ goes back to the classical
algebraic geometry, see \cite{Ro} and the references therein.  Most of
the results obtained in the subject in the last $30$ years are based
on the methods of syzygies and of construction of bundles on $\PP^4$.
Our approach of interpreting hypersurfaces containing a surface in
$\PP^4$ as certain extensions of sheaves on the surface itself seems
to be a relative newcomer in the subject.  It was initiated in our
previous work \cite{NaRe} with an eye toward the problem of bounding
the irregularity of surfaces in $\PP^4$.  Here we enlarge its scope by
addressing the Hartshorne-Lichtenbaum problem as well as the problem
of classification of irregular surfaces in $\PP^4$. If for the first
problem our contribution is largely tributary to the works
\cite{ElPe} and \cite{DeSch}, it is with the second problem that we
have tried to be as self-contained in our treatment as possible.  In
particular, in deriving Theorem \ref{th:irreg}, 2), 3), we have
avoided to call upon the results on classification of surfaces of
small degree in $\PP^4$. This is motivated (and hopefully justified)
by our objective to show/explore various aspects of using the
extension constructions to gain an insight into the geometry of
surfaces.  More importantly, we wanted to see (and show to the reader)
how the extension construction `pins down' (hypothetical) irregular
surfaces in $\PP^4$.  The proofs of Theorem \ref{th:irreg}, 2) - 4)
and Theorem \ref{th:V5Alb2} provide a substantial evidence that our
approach is useful for classifying surfaces in $\PP^4$.

From the conceptual point of view, our approach could be termed as
representing various geometric or cohomological entities by (short
exact) complexes of coherent sheaves on a surface in question.  The
complexes or, better, distinguished triangles in the derived category
(of coherent sheaves) we are using turn out to be unstable either for
numerical (Bogomolov instability) or for more subtle geometrical
reasons.  This instability gives rise to a new distinguished triangle
which carries more geometry than the initial one.  It is in relating
the two triangles that one is able to obtain new geometric insights.
This is very much in line with more recent developments of methods of
derived categories in algebraic geometry such as Bridgeland's
stability conditions, see \eg \cite{Br}, \cite{ArBe}.

\paragraph{Organization of the paper}  
In \S \ref{sec:prelimin}, preliminary
material is gathered.  We start by recalling some facts about the
Bogomolov instability and then go on explaining how to relate
hypersurfaces in $\PP^4$ containing a surface to the extension
sequences of sheaves on that surface. One of the technical results
used throughout the paper is Lemma \ref{l:mainTechnical}.

The sections \S\S \ref{s:m=4} - \ref{subsecm235} are devoted to the
Hartshorne-Lichtenbaum problem for surfaces of general type in
$\PP^4$.  The main results here are Proposition \ref{p:d-chi-finite}
and Proposition \ref{p:d-chi-finite-5}, see also Theorem
\ref{th:ISample} in the introduction.

The rest of the paper is devoted to the problem of the irregularity of
surfaces in $\PP^4$.  In \S5 we explain how the main ideas of our
approach are connected with this problem and in \S6 we illustrate some
of these ideas in the case of surfaces lying on a quadric
hypersurface, see Theorem \ref{th:q=0V2}.

In \S7 the surfaces on a cubic hypersurface are treated.  Theorem
\ref{th:onV3} is the main result of this section. 
\S8 is an interlude about elliptic scrolls in $\PP^4$.  The subject is
well-known, see \cite{Hu,ADHPR, ADHPR2}, but we approach it from the 
point of view of the (twisted) conormal bundle of a scroll.  The main 
result is Theorem \ref{th:cubicsSections}. 

In \S9 we treat the case of surfaces on a quartic hypersurfaces with
ordinary double points.  The main result here is Theorem
\ref{th:irrsurfV4}.

In \S10 the Albanese dimension of surfaces contained in a quartic
(rep. quintic) hypersurface is considered: Theorem
\ref{th:albaneseDimension} and Theorem \ref{th:V5Alb2} are the main
results of this section.

The Appendix of the paper returns again to the case of an elliptic
scroll.  The main objective here is to show how the classical
configuration of $10$ nodes of a Segre cubic hypersurface in $\PP^4$
is related to the geometry of an elliptic scroll contained in it.  In
particular, we show how the extension construction lifts the famous
configuration $(10_4, 15_6)$ of Segre to the category of (short exact
complexes of) coherent sheaves on a scroll and we suggest that this
should lead to a categorification of $(10_4, 15_6)$ configuration of
Segre.

\section{Notation and preliminaries}
\label{sec:prelimin}

\subsection{Bogomolov instability}
\label{sec:Bogomolov}

Let $X$ be a smooth complex projective surface%
\footnote{These hypotheses are assumed throughout the paper whenerver
  we speak about surfaces.}.  We denote by $\NS(X)$ the N\'eron-Severi
group of $X$.  Its rank $\rho$ is called the Picard number of $X$ and
the intersection product defines an integral quadratic form on
$\NS(X)$, whose real extension to $N(X):=\NS_\RR(X)$ is of type
$(1,\rho-1)$, by the Hodge Index Theorem.  The {\it positive cone} of
$X$ is the open cone
\[
  N^+(X) = \{D\in N(X)  \mid D^2>0, H\cdot D>0,
  \text{ for some (hence any)  ample divisor class $H$ on $X$}\}
.
\]
Note that $N^+(X)$ contains the ample cone and is contained in the 
cone of effective divisors.

Let $\Ff$ be a coherent sheaf on $X$ of rank $r=r_\Ff$.  The
discriminant of $\Ff$ is the expression 
\[
  \Delta(\Ff) = 2r\,c_2(\Ff)-(r-1)\,c_1^2(\Ff).
\]
A more geometric way to think about $\Delta(\Ff)$ for sheaves of rank
$r\geq 1$ is to observe that
\[
  \frac{\Delta(\Ff)}{2r} = c_2(\Ff)-\frac{r-1}{2r}\,c_1^2(\Ff) 
  = c_2\Big(\Ff \otimes \Oo_X\big(-\frac{1}{r} c_1(\Ff)\big)\Big).
\]
The next result is due to Bogomolov and it is used constantly in the
sequel.

\begin{likethm}[Bogomolov Theorem]
Let $\Ff$ be a torsion free coherent sheaf on a surface $X$.  If
$\Delta(\Ff)<0$, then there exists a maximal non-trivial saturated
subsheaf $\Ff'$ such that
\begin{itemize}
  \listspace
\item 
  $\Delta(\Ff')\geq0$,
\item
  $\dfrac{c_1(\Ff')}{r_{\Ff'}}-\dfrac{c_1(\Ff)}{r_{\Ff}} \in N^+(X)$
  and 
  $\bigg( c_1(\Ff')-\dfrac{r_{\Ff'}}{r_\Ff}\,c_1(\Ff) \bigg)^{2}
  \geq -\dfrac{\Delta(\Ff)}{2r_\Ff}$.
\end{itemize}
In particular, if $\Ff$ is
$D$-semistable\footnote{The sheaf $\Ff$ is $D$-semistable if
  $\dfrac{c_1(\Ff')}{r_{\Ff'}}\cdot D \leq
  \dfrac{c_1(\Ff)}{r_{\Ff}}\cdot D$, for any nonzero subsheaf
  $\Ff'\subset \Ff$.} with respect to an ample  
divisor $D$, then $\Delta(\Ff)\geq0$.  
\end{likethm}

A torsion free sheaf is called Bogomolov unstable if
$\Delta(\Ff)<0$ and Bogomolov semistable if $\Delta(\Ff)\geq0$.  The
theorem asserts that a torsion free Bogomolov unstable sheaf contains
a maximal Bogomolov semistable subsheaf which destabilizes it with
respect to every polarization.  Such a subsheaf is called a {\it
  maximal Bogomolov destabilizing} subsheaf of the given sheaf.

\begin{lem}
  \label{l:filtration}
Let $\Ff$ be a locally free sheaf on the surface $X$.  There exists a
unique Bogomolov filtration of $\Ff$,
\[
  0=\Ff_0 \subset \Ff_1 \subset\cdots\subset\Ff_m=\Ff
\]
such that for each $1\leq i\leq m$, $\Ff_i/\Ff_{i-1}$ is the maximal
Bogomolov destabilizing subsheaf of $\Ff_j/\Ff_{i-1}$ for every $j>i$.
\end{lem}

\proof
One can argue by induction on the rank $r=\rank(\Ff)$. For $r=1$ the
statement is obvious, since by definition locally free sheaves of rank
$1$ are Bogomolov semistable. So we assume $r\geq 2$ and suppose that
the theorem holds for all locally free sheaves of inferior
rank. Furthermore, we can assume that $\Ff$ is Bogomolov unstable
(since otherwise there is nothing to prove).

Let $\Ff_1$ be a maximal Bogomolov destabilizing subsheaf of $\Ff$.
By assumption, $\Ff_1\neq\Ff$.  Since $\Ff_1$ is saturated, the
quotient $\Ff/\Ff_1$ is torsion free, and therefore $\Ff_1$ is
reflexive (\cf \cite[Proposition 5.22]{Ko}), hence locally free, since
$X$ is a surface.  Now, if the quotient $\Ff/\Ff_1$ is Bogomolov
stable, the filtration reduces to
$0=\Ff_0\subset\Ff_1\subset\Ff_2=\Ff$ and we are done.  If not, the
quotient $\Ff/\Ff_1$ has the rank strictly smaller than $r$ and hence
the theorem holds for (the reflexive hull or the double dual of)
$\Ff/\Ff_1$.  Hence $\left(\Ff/\Ff_1\right)^{\ast\ast}$ admits a
unique Bogomolov filtration. Lifting this filtration to $\Ff$ gives
the desired filtration of $\Ff$.  It is enough to describe the
procedure for the lifting of the maximal Bogomolov destabilizing
subsheaf, call it $\Gg'$, of $\Ff/\Ff_1$ and then apply it inductively
for other pieces of the Bogomolov filtration of
$\left(\Ff/\Ff_1\right)^{\ast\ast}$.

Let $\Gg''$ be the quotient of the inclusion $\Gg' \subset \Ff/\Ff_1$.
We have the diagram
\[
  \begin{tikzpicture}[every node/.style={draw=none},
    ->,inner sep=1.1ex]
    \matrix [draw=none,row sep=3.75ex,column sep=4ex]
    {
      && \node (02) {$0$};
      & \node (03) {$0$}; \\
      \node (10) {$0$};
      & \node (11) {$\Ff_1$};
      & \node (12) {$\Ff_2$};
      & \node (13) {$\Gg'$};
      & \node (14) {$0$}; \\
      \node (20) {$0$};
      & \node (21) {$\Ff_1$};
      & \node (22) {$\Ff$};
      & \node (23) {$\Ff/\Ff_1$};
      & \node (24) {$0$}; \\
      && \node (32) {$\Gg''$};
      & \node (33) {$\Gg''$}; \\
      && \node (42) {$0$};
      & \node (43) {$0$}; \\
    };
    \path
      (10) edge[dashed] (11) 
      (11) edge[dashed] (12) 
      (12) edge[dashed] (13)
      (13) edge[dashed] (14)
      (02) edge (12)
      (03) edge (13)
      (11) edge[-,double distance=.33ex] (21)
      (12) edge (22)
      (13) edge (23)
      (22) edge (32)
      (23) edge (33)
      (32) edge (42)
      (33) edge (43)

      (20) edge (21)
      (21) edge (22)
      (22) edge (23)
      (23) edge (24)

      (32) edge[-,double distance=.33ex] (33)
    ;
  \end{tikzpicture}
\]
where $\Ff_2$ is the kernel of the epimorphism $\Ff\to\Gg''$.  As
before, in this short exact sequence $\Ff$ is locally free and $\Gg''$
is torsion free, hence $\Ff_2$ is locally free.  Clearly
$\Ff_1\subset\Ff_2$ and $\Gg'\cong\Ff_2/\Ff_1$.  We must show that
$\Ff_2$ is Bogomolov unstable and that $\Ff_1$ is a maximal Bogomolov
destabilizing subsheaf of $\Ff_2$. 

Set $r=\rank(\Ff)$, $r_j=\rank(\Ff_j)$, and $r_{\Gg'}=\rank(\Gg')$.
Since
\[
\begin{split}
  \frac{c_1(\Ff_2)}{r_2}&-\frac{c_1(\Ff)}{r}
  = \bigg(\frac{c_1(\Ff_1)}{r_1}-\frac{c_1(\Ff)}{r}\bigg)
  -\bigg(\frac{c_1(\Ff_1)}{r_1}-\frac{c_1(\Ff_2)}{r_2}\bigg)\\
  &= \bigg(\frac{c_1(\Ff_1)}{r_1}-\frac{c_1(\Ff)}{r}\bigg)
  -\frac{r_{\Gg'}}{r_2}\bigg(
  \frac{c_1(\Ff_1)}{r_1}-\frac{c_1(\Gg')}{r_{\Gg'}}
  \bigg)\\
  &= \bigg(\frac{c_1(\Ff_1)}{r_1}-\frac{c_1(\Ff)}{r}\bigg)
  -\frac{r_{\Gg'}}{r_2}\bigg(
  \frac{c_1(\Ff_1)}{r_1}-\frac{c_1(\Ff/\Ff_1)}{r-r_1}
  \bigg)
  +\frac{r_{\Gg'}}{r_2}\bigg(
  \frac{c_1(\Gg')}{r_{\Gg'}}-\frac{c_1(\Ff/\Ff_1)}{r-r_1}
  \bigg)\\
  &= \bigg(1-\frac{r_{\Gg'}}{r_2}\frac{r}{r-r_1}\bigg)
  \bigg(\frac{c_1(\Ff_1)}{r_1}-\frac{c_1(\Ff)}{r}\bigg)
  +\frac{r_{\Gg'}}{r_2}\bigg(
  \frac{c_1(\Gg')}{r_{\Gg'}}-\frac{c_1(\Ff/\Ff_1)}{r-r_1}
  \bigg)\\
  &= \frac{r(r-r_2)}{r(r-r_1)}
  \bigg(\frac{c_1(\Ff_1)}{r_1}-\frac{c_1(\Ff)}{r}\bigg)
  +\frac{r_{\Gg'}}{r_2}\bigg(
  \frac{c_1(\Gg')}{r_{\Gg'}}-\frac{c_1(\Ff/\Ff_1)}{r-r_1}
  \bigg),
\end{split}
\]
we see that $c_1(\Ff_2)/r_2-c_1(\Ff)/r\in N^+(X)$.  Hence
$\Delta(\Ff_2)<0$, since otherwise $\Ff_2$ would be a Bogomolov
destabilizing subsheaf of $\Ff$ and this contradicts the maximality of
$\Ff_1$.

Thus we now have constructed a Bogomolov unstable subsheaf $\Ff_2$ of
$\Ff$ and we claim that $\Ff_1$ is its maximal Bogomolov destabilizing
subsheaf.  Indeed, if $\Ff_1$ is not a maximal Bogomolov destabilizing
subsheaf of $\Ff_2$, then there exists a Bogomolov semistable (locally
free) subsheaf $\Ff'$ such that $\Ff_1\subset\Ff'\subset\Ff_2$ and
such that $c_1(\Ff')/r'-c_1(\Ff_2)/r_2\in N^+(X)$.  But then, by the
previous argument,
\[
  \frac{c_1(\Ff')}{r'}-\frac{c_1(\Ff)}{r}
  = \bigg(\frac{c_1(\Ff')}{r'}-\frac{c_1(\Ff_2)}{r_2}\bigg)
  +\bigg(\frac{c_1(\Ff_2)}{r_2}-\frac{c_1(\Ff)}{r}\bigg) \in N^+(X),
\]
contradicting the maximality of $\Ff_1$. 
\qed

\subsection{From hypersurfaces  of small degree to 
extension classes}
\label{subs:hypext}
Let $X\subset\PP^4$ be a smooth surface.  In what follows we denote by
$\Nn_X=\Nn_{X/\PP^4}$ the normal bundle of $X$ in $\PP^4$ and by $H$ a
hyperplane section of $X$.  The normal bundle $\Nn_X$ is of rank $2$
on $X$ with determinant $\det(\Nn_X)=\Bigwedge^2\Nn_X=\Oo_X(K_X+5H)$.
The conormal bundle satisfies
\begin{equation}
  \label{eq:nor-conor}
  \Jj_X/\Jj_X^2=\Nn^\ast_X \cong \det(\Nn^\ast)\otimes\Nn_X 
  = \Nn_X(-K_X-5H),
\end{equation}
where $\Jj_X$ is the ideal sheaf of $X$ in $\PP^4$ and the second
identification comes from the fact that the rank of $\Nn^\ast_X$ is
$2$.

We now assume that $X$ lies on a hypersurface of degree $m\leq 4$ and
not on one of a smaller degree, \ie $m=m_X$.  From the first
equality in \eqref{eq:nor-conor}, it follows that
$H^0(\Nn^\ast_X(mH))$ does not vanish.  Hence, by two other
identifications in \eqref{eq:nor-conor}, we obtain
\begin{equation}
  \label{eq:nonVaninshing4}
  H^0(\Nn_X(-K_X-(5-m)H)) = H^0(\Nn^\ast_X(mH)) \neq 0.
\end{equation}
Set $t=5-m$ and observe that $t \in \{1,2,3\}$.  We wish to interpret
nonzero sections of $\Nn_X(-K_X-tH)$ cohomologically.  For this,
consider the normal sequence of $X$ in $\PP^4$ tensored with
$\Oo_X(-K_X-tH)$,
\[
  0 \lra \Theta_X(-K_X-tH) \lra \Theta_{\PP^4}|_X(-K_X-tH) 
  \lra \Nn_X(-K_X-tH) \lra 0.
\]
This implies that $H^0(\Nn_X(-K_X-tH))$ fits into the following
exact sequence of cohomology groups
\begin{equation}
  \label{eq:coh-seq}
  H^0(\Theta_{\PP^4}|_X(-K_X-tH)) \lra
  H^0(\Nn_X(-K_X-tH)) \stackrel{\delta_X}{\lra}
  H^1(\Theta_X(-K_X-tH)).
\end{equation}
The following result improves a part of Lemma~5.3 in \cite{NaRe}.

\begin{lem}
  \label{l:EP-coh}
Let $s$ be a nonzero global section of $\Nn_X(-K_X-tH)$, 
$1\leq t\leq3$, corresponding to a $3$-fold of degree $m=5-t$
containing $X$.  If the Kodaira dimension of $X$ is non-negative, then
the cohomology class $\delta_X(s)\neq0$.
\end{lem}

\proof
We only prove the case $t=1$, \ie $X$ is contained in a quartic 
hypersurface, since the other cases are much easier. 

Assume $\delta_X(s)=0$ in $H^1(\Theta_X(-K_X-H))$, then $s$ is the
image of a nonzero global section $\stilde$ of
$\Theta_{\PP^4}|_X(-K_X-H))$.  From the Euler sequence of
$\Theta_{\PP^4}$, we deduce that either $H^0(\Oo_X(-K_X))\neq0$ or
$\ker(H^1(-K_X-H)\to H^0(H)^\ast\otimes H^1(-K_X))\neq0$.  In both
cases we see that $H^0(\Oo_C(-K_X))\neq0$ for every $C\in|H|$. Hence
$H\cdot(-K_X)\geq0$.

If $H\cdot(-K_X)>0$, then
$h^0(\Oo_X(mK_X))=0$ for every positive integer $m$, hence $X$
rational or irrationally ruled.  If $H\cdot(-K_X)=0$, then
$\Oo_C(-K_X)=\Oo_C$, for every $C\in|H|$ This tells us that
$H^0(\Oo_X(-K_X))\neq 0$ and hence
 $K_X =0$.  Therefore, $X$ is minimal and it is
either a K3 or an abelian surface.

If $X$ is a minimal K3 surface, then the exact
sequence 
\[
0 \lra H^0(\Oo_X) \lra H^0(\Oo_X(H)) \lra H^0(\Oo_C(H)) \lra 0,
\]
with $C$ a general curve in $|H|$, implies that
$g(C)=h^0(\Oo_C(H))=h^0(\Oo_X(H))-1=4$.  Hence $d=H^2=2g(C)-2=6$.
But then the sequence  
\[
  0 \lra \Jj_X(2) \lra \Oo_{\PP^4}(2) \lra \Oo_X(2H) \lra 0
\]
implies 
\[
  h^0(\Jj_X(2) \geq h^0(\Oo_{\PP^4}(2))-h^0(\Oo_X(2H))
  = \binom{6}{2}-\frac{(2H)^2}{2}-2 = 1,
\]
\ie there is a quadric passing through $X$ and this is contrary to our
assumption.

Thus we are left with the second case: $X$ a minimal abelian surface.
From the double point formula \eqref{eq:dPFormula} it follows that
$d=10$.  We show that a minimal abelian surface of degree $10$ can not
lie on a hypersurface of degree $4$. This follows from the following
observation.

\paragraph{Claim}
\label{page:claim}
A quartic hypersurface $Q\in\PP^4$ containing a minimal abelian
surface $X$ must be a cone over a quartic surface $S\subset\PP^3$ with
at most isolated singularities.

\medskip

Let us assume the claim and derive a contradiction.  We view $\PP^4$
as $\PP(H^0(\Oo_X(H))^\ast)$ and denote by $[v]$, 
$v\in H^0(\Oo_X(H))^\ast$, the vertex of $Q$.  Under the projection
from $[v]$ the surface $X$ becomes a finite covering of a quartic
surface $S$ lying in some $\PP^3$, complementary to $[v]$.  Let $m$ be
the degree of this covering.  It is related to the intersection of $X$
with a ruling $l$ of the cone $Q$ as follows.
\[
  (X\cdot l)_Q =
  \begin{cases}
    m &\text{if }\, [v]\notin X \\
    m-1 &\text{if }\, [v]\in X.
  \end{cases}
\] 
But for a general plane $\Lambda$ passing through $[v]$, 
\[
  10 = d = (X\cdot \Lambda)_{\PP^4} =
  \begin{cases}
    4m &\text{if }\, [v]\notin X \\
    4(m-1)+1 &\text{if }\, [v]\in X
  \end{cases}  
\]
and neither case is possible. 

\medskip

We now turn to the proof of the claim.  Let us recall the situation:
$s$ is a nonzero global section in $\Nn_X(-H)$ corresponding to a
quartic hypersurface $Q$ containing $X$, and it is the image of a
global section $\stilde$ of $\Theta_{\PP^4}|_X(-H)$.  The argument is
divided into two steps.

\step
We claim that the scheme of zeros $Z_s=\{s=0\}$ of the global section
$s$ is $0$-dimensional.  Indeed, let us assume that this is not the
case and let $\Gamma$ be a reduced, irreducible curve in $Z_s$.  This
means that $s=\gamma s'$, where $\gamma\in H^0(\Oo_X(\Gamma))$ is a
global section defining $\Gamma$ and $s'\in H^0(\Nn_X(-H-\Gamma))$.
From the commutative diagram
\[
  \begin{tikzpicture}[every node/.style={draw=none},
    ->,inner sep=1.1ex]
    \matrix [draw=none,row sep=4.5ex,column sep=4.2ex]
    {
      \node (11) {$H^0(\Theta_{\PP^4}|_X(-H-\Gamma))$};
      & \node (12) {$H^0(\Nn_X(-H-\Gamma))$};
      & \node (13) {$H^1(\Theta_X(-H-\Gamma))$}; \\
      \node (21) {$H^0(\Theta_{\PP^4}|_X(-H))$};
      & \node (22) {$H^0(\Nn_X(-H))$}; \\
    };
    \path
      (11) edge (12)
      (12) edge (13)
      (11) edge node[left=-.1ex,scale=.75] {$\gamma\cdot$} (21)
      (12) edge node[right=-.1ex,scale=.75] {$\gamma\cdot$} (22)
      (21) edge (22)
    ;
  \end{tikzpicture}
\]
it follows that either $H^0(\Theta_{\PP^4}|_X(-H-\Gamma))\neq0$ or 
$H^1(\Theta_X(-H-\Gamma))\neq0$.  The latter possibility leads to
$H^1(\Oo_X(-H-\Gamma))\neq0$, since
$H^1(\Theta_X(-H-\Gamma))\cong\oplus H^1(\Oo_X(-H-\Gamma))$.  But
then, $H^0(\Oo_\Gamma(-H))\cong H^1(\Oo_X(-H-\Gamma))\neq0$ which is
impossible.  The former possibility leads, using the Euler
sequence for $\Theta_{\PP^4}$, to 
\[
  \ker \left(
    H^1(\Oo_X(-H-\Gamma))
    \to H^0(\Oo_X(H))^\ast\otimes H^1(\Oo_X(-\Gamma)
  \right) \neq 0,
\]
which means that $H^0(\Oo_C(-\Gamma))\neq0$ for every $C\in|H|$, which
is not possible either.

\step
$Q$ is a cone over a quartic surface.  Indeed, on the one hand we
consider the diagram
\vspace{6ex}
\begin{equation}
  \label{v-diagram}
  \begin{tikzpicture}[overlay,every node/.style={draw=none},
    ->,inner sep=1.1ex]
    \matrix [draw=none,row sep=4ex,column sep=4.2ex]
    {
      \node (11) {$H^0(\Oo_X(H))^\ast$}; \\
      \node (21) {$H^0(\Theta_{\PP^4}|_X(-H))$};
      & \node (22) {$H^0(\Nn_X(-H))$}; \\
    };
    \path
      (11) edge (21)
      (11) edge (22)      
      (21) edge (22)
    ;
  \end{tikzpicture}
  \vspace{6ex}
\end{equation}
where the horizontal arrow comes from the normal sequence of $X$ in
$\PP^4$ and the vertical one is part of the Euler sequence of
$\Theta_{\PP^4}$ (tensored with $\Oo_{\PP^4}(-1)$) restricted to $X$.
Both maps are isomorphisms.  This implies that the section $s$ is the
image of a unique element $v\in H^0(\Oo_X(H))^\ast$.

On the other hand, we have the Koszul sequence associated to
$s$,
\[
  0 \lra \Oo_X \stackrel{s}\lra \Nn_X(-H) \stackrel{s\wedge}\lra
  \Jj_{Z_s}(3H) \lra 0.
\]
which is exact by {\it Step 1}.  Combining this with the slanted arrow
in \eqref{v-diagram} gives the following diagram
\vspace{6ex}
\begin{equation}
  \label{s-diagram}
  \begin{tikzpicture}[overlay,every node/.style={draw=none},
    ->,inner sep=1.1ex]
    \matrix [draw=none,row sep=4.2ex,column sep=4.2ex]
    {
      && \node (12) {$H^0(\Oo_X(H))^\ast\otimes\Oo_X$}; \\
      \node (20) {$0$};
      & \node (21) {$\Oo_X$};
      & \node (22) {$\Nn_X(-H)$};
      & \node (23) {$\Jj_{Z_s}(3H)$};
      & \node (24) {$0$}; \\
    };
    \path
      (12) edge node[left=-.2ex,scale=.75] {$e$} (22)
      (12) edge (23)
      
      (20) edge (21)
      (21) edge node[above=-.2ex,scale=.75] {$s$} (22)
      (22) edge node[above=-.2ex,scale=.75] {$s\wedge$} (23)
      (23) edge (24)
    ;
  \end{tikzpicture}
  \vspace{6ex}
\end{equation}
Furthermore, the homomorphism
\begin{equation}
  \label{diff}
  \partial: H^0(\Oo_X(H))^\ast \longrightarrow H^0(\Jj_{Z_s}(3H))
\end{equation}
induced from the above diagram on the level of global sections is the
partial differentiation. Namely, let $F$ be a homogeneous polynomial
defining $Q$, then the homomorphism $\partial$ is given by 
\[
  H^0(\Oo_X(H))^\ast \ni u \mapsto
  \partial(u) = \partial_uF|_X \in H^0(\Jj_{Z_s}(3H)).
\]
By construction
$\partial(v)$ factors via  $H^0(\Nn_X(-H))$, \ie we have
\[
  \partial_v F|_X = \partial(v) =s\wedge (e(v)) =s\wedge s =0.
\]
This means that the homogeneous polynomial of degree $3$,
$\partial_vF\in\Sym^3H^0(\Oo_X(H))$, vanishes on $X$.  But since the
surface $X$ is not contained in any hypersurface of degree less than
$4$, we conclude that $\partial_vF=0$.  Equivalently,
$F\in\Sym^4(\ker(v))$, \ie the $3$-fold of degree four $Q$ is the cone
in $\PP(H^0(\Oo_X(H))^\ast)$ with vertex $[v]$ and base the quartic
surface defined by $F$ in $\PP^3 =\PP(\ker(v)^\ast)$.

The last assertion, stating that the quartic surface $S$ defined by
$F$ in $\PP(\ker(v)^\ast)$ has at most isolated singularities, follows
from the observation that a curve, call it $\Gamma$, in the singular
locus of $S$ produces a surface $\Sigma$ in $Q$---the cone over
$\Gamma$ with vertex at $[v]$---and this surface lies in the singular
locus of the quartic $Q$. But then the surfaces $X$ intersects
$\Sigma$ along a curve which is part of the zero-locus of the section
$s$. This contradicts {\it Step 1}.
\qed

From now on we assume that the cohomology class 
$\delta_X(s) \in H^1(\Theta_X(-K_X-tH))$ in Lemma \ref{l:EP-coh} is
nonzero.  The identification
\[
  H^1(\Theta_X(-K_X-tH)) \cong \Ext^1(\Omega_X,\Oo_X(-K_X-tH))
\]
allows us to interpret a cohomology class on the left as an
extension sequence of sheaves on $X$.  The following result
constitutes one of the main technical ingredients of this study.

\begin{lem}
  \label{l:mainTechnical}
Let $X$ be a smooth projective surface and let $M$ be a divisor on
$X$. Let $\xi\in H^1(\Theta_X(-K_X-M))$ be a nonzero cohomology class
and let
\begin{equation}
  \label{eq:extSeq0}
  0 \lra \Oo_X(-K_X-M) \lra \Tt_\xi \lra \Omega_X \lra 0.
\end{equation}
be the corresponding extension sequence.  Assume that $\Tt_\xi$
contains a subsheaf $\Ff$ of rank $2$ such that the induced morphism
$\Ff\to\Omega_X$ is generically an isomorphism. Then the following
holds.
\begin{enumerate}[label={\rm \arabic*)}]
  \listspace
\item
  The canonical divisor of $X$ decomposes as $K_X=L+E$, where
  $L=c_1(\Ff)$ and $E$ is the support of the cokernel
  $\coker(\Ff\to\Omega_X)$, a nonzero effective divisor on $X$.
\item
  If $e$ is a section of $\Oo_X(E)$ defining $E$, \ie 
  $E=(e=0)$, then the cohomology class $\xi$ is annihilated by $e$,
  \ie
  \[
    e\,\xi=0 \textq{in} H^1(\Theta_X(E-K_X-M)).
  \]
\item
  If, in addition, $X\subset\PP^4$ lies on a $3$-fold of degree
  $m\leq 5$ and $\xi=\delta_X(s)$, where $\delta_X$ is the coboundary
  map in \eqref{eq:coh-seq}, then 
  \[
    H^0(\Theta_{\PP^4}|_X(E-K_X-(5-m)H))
    = H^0(\Theta_{\PP^4}|_X(-L-(5-m)H))
    \neq 0. 
  \]
  In particular, $H\cdot L \leq (m-4)H^2$.
\end{enumerate}
\end{lem}

\proof
Set
\[
  \phi_\xi : \Ff \lra \Omega_X
\]
to be the morphism defined by the composition of the inclusion
$\Ff{\hookrightarrow}\Tt_\xi$ together with the epimorphism of the
extension sequence \eqref{eq:extSeq0}. By assumption ${\phi_\xi}$ is
generically an isomorphism. This implies that the support of 
$coker ({\phi_\xi})$ is a nonzero effective divisor, since otherwise
${\phi_\xi}$ is an isomorphism and the exact sequence
\eqref{eq:extSeq0} splits or, equivalently, $\xi=0$.

Writing out the exact sequence
\[
  0 \lra \Ff \stackrel{\phi_\xi}{\lra}
  \Omega_X \lra \coker ({\phi_\xi}) \lra 0
\]
we obtain the decomposition of the canonical divisor asserted in 1) of
the lemma. 

To prove 2) we consider the diagram
\vspace{16ex}
\begin{equation}
  \label{d:unstable}
  \begin{tikzpicture}[overlay,every node/.style={draw=none},
    ->,inner sep=1.1ex]
    \matrix [draw=none,row sep=3.5ex,column sep=4ex]
    {
      && \node (02) {$0$}; \\
      && \node (12) {$\Ff$}; \\
      \node (20) {$0$};
      & \node (21) {$\Oo_X(-K_X-M)$};
      & \node (22) {$\Tt_\xi$};
      & \node (23) {$\Omega_X$};
      & \node (24) {$0$}; \\
      && \node (32) {$\Jj_Z(-L-M)$}; \\      
      && \node (42) {$0$}; \\      
    };
    \path
      (02) edge (12) 
      (12) edge (22)
      (22) edge (32)
      (32) edge (42)

      (12) edge node[above right=-.3ex,scale=.75] {$\phi_\xi$} (23)
      (21.-11.5) edge node[below left=-.3ex,scale=.75] {$e$} (32)
      
      (20) edge (21)
      (21) edge (22)
      (22) edge (23)
      (23) edge (24)
    ;
  \end{tikzpicture}
  \vspace{16.5ex}
\end{equation}
where the slanted arrow in the lower part of the diagram is the
morphism given by multiplication with the section $e$.
Dualizing and ten\-so\-ring the diagram with $\Oo_X(-L-M)$ we arrive
at
\vspace{6ex}
\begin{equation}
  \label{d:unstableD}
  \begin{tikzpicture}[overlay,every node/.style={draw=none},
    ->,inner sep=1.1ex]
    \matrix [draw=none,row sep=4ex,column sep=4.2ex]
    {
      && \node (12) {$\Oo_X$}; \\
      \node (20) {$0$};
      & \node (21) {$\Theta_X(-L-M)$};
      & \node (22) {$\Tt_\xi^\ast(-L-M)$};
      & \node (23) {$\Oo_X(K_X-L)$};
      & \node (24) {$0$}; \\
    };
    \path
      (12) edge (22)
      (12) edge  node[above right=-.3ex,scale=.75] {$e$} (23)
      
      (20) edge (21)
      (21) edge (22)
      (22) edge (23)
      (23) edge (24)
    ;
  \end{tikzpicture}
  \vspace{5.5ex}
\end{equation}
Since the coboundary map 
\[
  H^0(\Oo_X(K_X-L)) 
  = H^0(\Oo_X(E)) \lra H^1(\Theta_X(-L-M))
\]
in long exact sequence of cohomology groups of the horizontal
sequence in \eqref{d:unstableD} is given by the cup-product with the
class $\xi$ and since the section $e$ is in its kernel, we deduce
$e\,\xi=0$ in $H^1(\Theta_X(-L-M))$.

For 3), we use the fact that $\xi=\delta_X(s)$, where $\delta_X$ is as
in \eqref{eq:coh-seq}.  Set $\Delta=K_X+(5-m)H$ and consider the
commutative diagram
\vspace{-1ex}
\[
  \begin{tikzpicture}[every node/.style={draw=none},
    ->,inner sep=1.1ex]
    \matrix [draw=none,row sep=4.2ex,column sep=4.2ex]
    {
      & \node (12) {$H^0(\Nn_X(-\Delta))$};
      & \node (13) {$H^1(\Theta_X(-\Delta))$}; \\
      \node (21) {$H^0(\Theta_{\PP^4}|_X(E-\Delta))$};
      & \node (22) {$H^0(\Nn_X(E-\Delta))$};
      & \node (23) {$H^1(\Theta_X(E-\Delta))$}; \\
    };
    \path
      (12) edge node[above=-.2ex,scale=.75] {$\delta_X$} (13)
      (12) edge node[left=-.2ex,scale=.75] {$e$} (22)
      (13) edge node[right=-.2ex,scale=.75] {$e$} (23)
      
      (21) edge (22)
      (22) edge node[above=-.2ex,scale=.75] {$\delta'_X$} (23)
    ;
  \end{tikzpicture}
  \vspace{-2ex}
\]
From this and 2) of the lemma, we obtain
\[
  \delta^{\prime}_X(es)=e\,\delta_X(s)=e\,\xi=0.
\]
Hence the global section $es \in H^0(\Nn_X(E-K_X-(5-m)H))$, being
obviously nonzero, comes from a nonzero section in
$H^0(\Theta_{\PP^4}|_X(E-K_X-(5-m)H))$.  This proves the first
assertion of 3).

To see the inequality $H\cdot L\leq (m-4)H^2$, we restrict the Euler
sequence for $\PP^4$ to $X$ and tensor it with
$\Oo_X(E-K_X-(5-m)H)=\Oo_X(-L-(5-m)H)$ to arrive at
\[
  0 \to
  \Oo_X(-L-(5-m)H) \to
  H^0(\Oo_X(H))^\ast\otimes\Oo_X(-L-(4-m)H) \to
  \Theta_{\PP^4}|_X(-L-(5-m)H) \to 0.
\]
Since $H^0(\Theta_{\PP^4}|_X(E-K_X-(5-m)H))\neq0$, by the first part
of 3), the above sequence implies that either
$H^0(\Oo_X(-L-(4-m)H))\neq0$, or
\[
  \ker\Big(H^1(\Oo_X(-L-(5-m)H))\lra
  H^0(\Oo_X(H))^\ast\otimes H^1(\Oo_X(-L-(4-m)H))\Big) \neq 0.
\]
The first possibility immediately gives the assertion 
$H\cdot L\leq(m-4)H^2$.  The second one implies that the homomorphism
$H^1(\Oo_X(-L-(5-m)H))\stackrel{h}{\to}H^1(\Oo_X(-L-(4-m)H))$ given by
the multiplication by any global section $h \in H^0(\Oo_X(H))$ has a
nonzero kernel.  Since that kernel comes from
$H^0(\Oo_{C_h}(-L-(4-m)H))$, where $C_h =(h=0)$, we deduce that
$H^0(\Oo_{C}(-L-(4-m)H))\neq 0$, for any divisor $C$ in the linear
system $|H|$.  This implies $H\cdot(-L-(4-m)H)\geq 0$ and hence the
assertion $H\cdot L\leq (m-4)H^2$.
\qed

\vspace{2\baselineskip}

\part*{\Large  (I)\hspace{.75ex}
  Hartshorne-Lichtenbaum for surfaces
  of general type}

\section{Numerical invariants for surfaces on degree $4$
  hypersurfaces}
\label{s:m=4} 

In this section we prove the inequalities
\eqref{eq:ineqOnChernNumbers} of the form 
$c_2-K_X^2\geq a\,H\cdot K_X+bd$, stated in the introduction.  Let $X$
be a smooth surface in $\PP^4$ lying on a $3$-fold $V=V_4$ of degree
four and not on any of a smaller degree.  Unless stated otherwise, we
assume that the Kodaira dimension of $X$ is non-negative. This
assumption, according to Lemma~\ref{l:EP-coh}, gives a nonzero
cohomology class $\delta_X(s)\in H^1(\Theta_X(-K_X-H))$ (see
Lemma~\ref{l:EP-coh} for notation) which we denote by $\xi$.  As we
already explained, this cohomology class is used to build the
extension \eqref{eq:extSeq0} and the focus of study becomes the vector
bundle $\Tt_\xi$ sitting in the middle of that sequence.  The
inequalities we are after are a consequence of Bogomolov semistability
or instability of $\Tt_\xi$.

\subsection{The inequalities \eqref{eq:ineqOnChernNumbers}}

\begin{thm}
  \label{th:mainTechnical}
Let $\Tt_\xi$ be the sheaf in the middle of the extension sequence
\eqref{eq:extSeq0} associated to $\xi=\delta_X(s)$.  
\begin{enumerate}[label={\rm \arabic*)}]
  \listspace
\item
If $\Tt_\xi$ is Bogomolov semistable, then 
$c_2-K_X^2\geq H\cdot K_X+\frac{1}{3}\,d$.
\item
If $\Tt_\xi$ is Bogomolov unstable, then
$c_2-K_X^2\geq 
\min\Big(
\frac{3}{4}H\cdot K_X,\frac{1}{2}\,H\cdot K_X+\frac{1}{4}\,d
\Big)$.
\end{enumerate}
\end{thm}

\proof
From the exact sequence \eqref{eq:extSeq0} it follows that $\Tt_\xi$
is a locally free sheaf of rank $3$ with the Chern invariants
$c_1(\Tt_\xi)=-H$ and $c_2(\Tt_\xi)=c_2-K_X^2-H\cdot K_X$. Therefore
the Bogomolov semistability condition for $\Tt_\xi$ reads as follows
\[
  6c_2(\Tt_\xi)-2c_1^2(\Tt_\xi)
  = 6(c_2-K_X^2)-6H\cdot K_X-2H^2 \geq 0
\]
and this is equivalent to the inequality in 1) of the theorem.

We now turn to the case when $\Tt_\xi$ is Bogomolov unstable. To
analyse the situation we use the Bogomolov filtration of $\Tt_\xi$,
see Lemma~\ref{l:filtration}.  In particular, according to the shape
of that filtration, we obtain the following 
inequalities\footnote{The inequality of the last line in
  \eqref{eq:lbound-filt} is strict, see Lemma \ref{l:lastCase} for
  details.}: 
\begin{equation}\label{eq:lbound-filt}
  c_2-K_X^2 \geq 
\begin{cases}
  H\cdot K_X, &
  \textq{if} 0\subset\Ff_1\subset\Ff_2=\Tt_\xi
  \textq{and} \rank(\Ff_1)=2 \\
  \frac{1}{2}\,H\cdot K_X+\frac{1}{4}\,d, &
  \textq{if} 0\subset\Ff_1\subset\Ff_2=\Tt_\xi
  \textq{and} \rank(\Ff_1)=1 \\
  \frac{3}{4}\,H\cdot K_X, &
  \textq{if} 0\subset\Ff_1\subset\Ff_2\subset\Ff_3=\Tt_\xi.
\end{cases}
\end{equation}
This implies 2) of the theorem.
\qed

Before we proceed with the proof of \eqref{eq:lbound-filt}, we would
like to provide the reader with the conducting line of the proofs of
the lemmas below.  The basic idea is to use the Bogomolov filtration
of $\Tt_\xi$ for writing down a ``good'' lower bound for the second
Chern number of $\Tt_\xi$.  ``Good'' here means that a sought after
estimate should imply a lower bound for $c_2-K_X^2$ as a positive
function of $d$ or/and $H\cdot K_X$. This is possible in view of the
following special features of $\Tt_\xi$:
\begin{itemize}\listspace
\item 
The subsheaves of rank $2$ involved in the filtration of $\Tt_\xi$
satisfy the hypotheses of the technical Lemma \ref{l:mainTechnical}. 
\item
The subsheaves of rank $1$ involved in the filtration of $\Tt_\xi$
having positive degree (with respect to some polarization of $X$)
inject into $\Omega_X$ and hence must be of Iitaka dimension at most
one (Bogomolov lemma); furthermore, the generic semi-positivity of
$\Omega_X$ insures that the quotient sheaf must be of non-negative
degree.
\end{itemize}

With these remarks in mind, we now consider all possible
filtrations of the Bogomolov unstable vector bundle $\Tt_\xi$.

\begin{lem}
  \label{l:rank2ForF}
If the Bogomolov filtration of $\Tt_\xi$ is
$0=\Ff_0\subset\Ff_1\subset\Ff_2=\Tt_\xi$ with $\Ff_1$ of rank $2$,
then it gives rise to a divisor $B_1$ in the positive cone of $N(X)$
and to an effective nonzero divisor $E$ such that the following hold: 
\begin{enumerate}[label={\rm \arabic*)}]
  \listspace
\item 
$c_1(\Ff_1) =\frac{1}{3}\,(B_1 -2H)$ and $B_1\cdot H \leq 2d$,
\item
$K_X = \frac{1}{3}\,(B_1 -2H) +E$,
\item
$c_2-K_X^2\geq H\cdot K_X$.
\end{enumerate}
\end{lem}

\proof
The bundle $\Tt_\xi$ is the middle term of two exact sequences, as in
diagram \eqref{d:unstable}, where the maximal Bogomolov destabilizing
subsheaf $\Ff_1$ takes the place of $\Ff$.

We set $L_1=c_1(\Ff_1)$ and
$L_2=c_1(\Tt_\xi/\Ff_1)$. Using the equality 
\[
  -H = c_1(\Tt_\xi) = L_1+L_2,
\]
the Bogomolov destabilizing condition for $\Ff_1$ tells us that the
$\QQ$-divisor
\[
  \frac{c_1(\Ff)}{2}-\frac{c_1(\Tt_\xi)}{3} 
  = \frac{L_1}{2}-\frac{-H}{3} = \frac{1}{6}\,(L_1-2L_2),
\]
lies in the positive cone $N^+(X)$ of $N(X)$. Thus the divisor
$B_1=L_1-2L_2$ lies in the positive cone $N^+(X)$ and we write $L_i$,
$i=1,2$, as a linear combination of $H$ and $B_1$ as follows:
\begin{equation}
  \label{eq:LandL}
  \begin{aligned}
     L_1 &= -\frac{2}{3}\,H+\frac{1}{3}\,B_1 \\
     L_2 &= -\frac{1}{3}\,H-\frac{1}{3}\,B_1.
  \end{aligned}
\end{equation}

Recall that $c_2(\Tt_\xi) =c_2-K_X^2-H\cdot K_X $. Computing that
Chern class using the vertical sequence in \eqref{d:unstable} with
$\Ff=\Ff_1$ and the quotient sheaf $\Tt_\xi/\Ff_1 =\Ii_Z(L_2)$, we
obtain
\begin{equation}
  \label{eq:c2}
  c_2-K_X^2-H\cdot K_X = c_2(\Ff_1)+L_1\cdot L_2+\deg(Z) 
  \geq \frac{1}{4}\,L_1^2+L_1\cdot L_2
  =\frac{1}{12}(4H^2-B_1^2).
\end{equation}
where the inequality uses $c_2(\Ff_1)\geq\frac{1}{4}\,L_1^2$, the
Bogomolov semistability of $\Ff_1$, and the last equality comes from
substituting the expressions from \eqref{eq:LandL}.

Next we claim that the slanted arrow $\phi_\xi:\Ff_1\to\Omega_X$ in
the diagram \eqref{d:unstable} is generically an isomorphism.  Indeed,
if $\phi_\xi$ drops its rank everywhere, then we obtain the
commutative diagram
\vspace{-1ex}
\[
  \begin{tikzpicture}[every node/.style={draw=none},
    ->,inner sep=1.1ex]
    \matrix [draw=none,row sep=4ex,column sep=4.2ex]
    {
      && \node (02) {$0$}; \\
      \node (10) {$0$};
      & \node (11) {$\Oo_X(-K_X-H)$};
      & \node (12) {$\Ff_1$};
      & \node (13) {$\im(\phi_\xi)$};
      & \node (14) {$0$}; \\
      \node (20) {$0$};
      & \node (21) {$\Oo_X(-K_X-H)$};
      & \node (22) {$\Tt_\xi$};
      & \node (23) {$\Omega_X$};
      & \node (24) {$0$}; \\
    };
    \path
      (10) edge (11) 
      (11) edge (12) 
      (12) edge (13)
      (13) edge (14)
      (02) edge (12)
      (11) edge (21)
      (12) edge (22)
      (13) edge (23)
      
      (20) edge (21)
      (21) edge (22)
      (22) edge (23)
      (23) edge (24)
    ;
  \end{tikzpicture}
  \vspace{-1ex}
\]
where the sheaf $\im(\phi_\xi)$ is the image of $\phi_\xi$. It is a
torsion free subsheaf of rank $1$ of $\Omega_X$ with the first Chern
class
\[
  c_1(\im(\phi_\xi)) = c_1(\Ff_1)+(K_X+H) 
  = -\frac{2}{3}\,H+\frac{1}{3}\,B_1+K_X+H
  = K_X+\frac{1}{3}\,H+\frac{1}{3}\,B_1.
\]
But this means that $\Omega_X$ contains a rank $1$ subsheaf of Iitaka
dimension $2$ which is impossible in view of Bogomolov Lemma.

Once we know that $\phi_\xi$ is generically of maximal rank,
Lemma~\ref{l:mainTechnical} can be applied.  In particular, we obtain
the decomposition asserted in 2) of that lemma, where the divisor $E$
is the support of the cokernel of $\phi_\xi$, and, from the part 3) of
Lemma~\ref{l:mainTechnical}, we deduce that $H\cdot L_1\leq0$.  This
together with the formula for $L_1$ in \eqref{eq:LandL} implies
\[
  0 \geq H\cdot L_1 = H\cdot(-\frac{2}{3}\,H+\frac{1}{3}\,B_1) 
  = \frac{1}{3}(B_1\cdot H -2d).
\]
Hence $B_1\cdot H \leq 2d$ as asserted in 1) of the lemma.

The above inequality and the Hodge Index
Theorem give $B_1^2\leq4d$.  Substituting this inequality in
\eqref{eq:c2}, we deduce assertion 3) of the lemma.
\qed

\begin{lem}
  \label{l:Qstable}
If the Bogomolov filtration of $\Tt_\xi$ is
$0=\Ff_0\subset\Ff_1\subset\Ff_2=\Tt_\xi$ with $\Ff_1$ of rank $1$,
then
\[
  c_2-K_X^2 \geq 
  \begin{cases}
    H\cdot K_X+\frac{1}{4}\,d\,,& \text{if $H\cdot c_1(\Ff_1)\leq0$} \\
    \frac{1}{2}\,H\cdot K_X+\frac{1}{4}\,d\,,& 
    \text{if $H\cdot c_1(\Ff_1)>0$}. 
  \end{cases}
\]
\end{lem}

\proof
Let $\Qq=\Tt_\xi/\Ff_1$ and set $L_1=c_1(\Ff_1)$ and
$L_2=c_1(\Qq)$. 
The bundle $\Tt_\xi$ becomes the middle term of two exact sequences.
\vspace{16ex}
\begin{equation}
  \label{d:forQuotientQ}
  \begin{tikzpicture}[overlay,every node/.style={draw=none},
    ->,inner sep=1.1ex]
    \matrix [draw=none,row sep=3.5ex,column sep=4ex]
    {
      && \node (02) {$0$}; \\
      && \node (12) {$\Oo_X(L_1)$}; \\
      \node (20) {$0$};
      & \node (21) {$\Oo_X(-K_X-H)$};
      & \node (22) {$\Tt_\xi$};
      & \node (23) {$\Omega_X$};
      & \node (24) {$0$}; \\
      && \node (32) {$\Qq$}; \\      
      && \node (42) {$0$}; \\      
    };
    \path
      (02) edge (12) 
      (12) edge (22)
      (22) edge (32)
      (32) edge (42)

      (12) edge (23)
      (21.-12.5) edge (32)
      
      (20) edge (21)
      (21) edge (22)
      (22) edge (23)
      (23) edge (24)
    ;
  \end{tikzpicture}
  \vspace{16.5ex}
\end{equation}
The condition that $\Oo_X(L_1)$ is a Bogomolov destabilizing subsheaf
of $\Tt_\xi$ and the equality $(-H)=L_1+L_2$ imply that the
$\QQ$-divisor
\[
  L_1-\frac{c_1(\Tt_\xi)}{3}
  = L_1-\frac{-H}{3} = L_1-\frac{L_1+L_2}{3} 
  = \frac{2}{3}\Big(L_1-\frac{1}{2}\,L_2\Big)
\]
lies in the positive cone of $N(X)$. Setting 
\begin{equation*}
  \label{eq:deltaAgain}
  B_1 = L_1-\frac{1}{2}\,L_2 \in N^+(X),
\end{equation*}
we express $L_1$ and $L_2$ as linear combinations of
$H$ and $B_1$:
\begin{equation}
  \label{eq:L1andL2}
  \begin{aligned}
     L_1 &= -\frac{1}{3}\,H+\frac{2}{3}\,B_1 \\
     L_2 &= -\frac{2}{3}\,H-\frac{2}{3}\,B_1.
  \end{aligned}
\end{equation}
From here on we use the same argument as in Lemma \ref{l:rank2ForF}.
Namely, we use the vertical sequence in \eqref{d:forQuotientQ} to
estimate the second Chern class of $\Tt_\xi$:
\begin{equation}
  \label{eq:ChernCl}
  c_2-K_X^2-H\cdot K_X = c_2(\Qq)+L_1\cdot L_2
  \geq \frac{1}{4}\,L_2^2+L_1\cdot L_2 =\frac{1}{3}(d-B_1^2),
\end{equation}
where the inequality uses the condition that the quotient 
$\Qq =\Ff_2/\Ff_1$ is Bogomolov se\-mi\-stable, \ie
$c_2(\Qq)\geq \frac{1}{4}\,L_2^2$, and the last equality comes from
substituting the expressions from \eqref{eq:L1andL2}.
 
To conclude the argument we need an appropriate upper bound on the
self-intersection $B_1^2$.  We argue according to the sign of 
$H\cdot L_1$.

\medskip 
\case 
If $H\cdot L_1\leq0$, then we are essentially in the same situation as
in the proof of Lemma \ref{l:rank2ForF}.  Namely, we have
\[
  0\geq H\cdot L_1 = \frac{1}{3}H\cdot(-H+2B_1) 
  = \frac{1}{3} \left(2B_1\cdot H -d \right),
\]
where the first equality uses the formula for $L_1$ in
\eqref{eq:L1andL2}.  Thus we obtain $H\cdot B_1\leq\frac{1}{2}\,d$
and hence, by the Hodge index, the upper bound
$B_1^2\leq\frac{1}{4}\,d$.  The inequality \eqref{eq:ChernCl} then
becomes
\[
  c_2-K_X^2-H\cdot K_X
  \geq \frac{1}{3}\,(d-B_1^2) 
  \geq \frac{1}{3}\Big(d-\frac{1}{4}\,d\Big) 
  = \frac{1}{4}\,d
\]
and this is equivalent to the first inequality of the lemma.

\medskip
\case 
If $H\cdot L_1>0$, then the morphism $\Oo_X(L_1)\to\Omega_X$ given by
slanted arrow in the diagram \eqref{d:forQuotientQ}  is
nonzero.%
\footnote{Otherwise  the divisor $-(H+K_X+L_1)$ 
must be effective; however substituting the formula for $L_1$ from
\eqref{eq:L1andL2}, one obtains 
$H+K_X+L_1=\frac{3}{4} H +K_X +\frac{2}{3}B_1$ which is of the Iitaka
dimension $2$.} 
Hence $\Oo_X(L_1)$ injects into $\Omega_X$ and by the Bogomolov Lemma,
$L_1^2\leq0$.  This inequality and the formula \eqref{eq:L1andL2} for
$L_1$ give 
\[
0\geq(-H+2B_1)^2=H^2-4H\cdot B_1+4B_1^2,
\]
or equivalently, 
\begin{equation}
  \label{eq:Delta2}
  B_1^2 \leq H\cdot B_1-\frac{1}{4}\,d.
\end{equation}

On the other hand, the generic semi-positivity of $\Omega_X$
stipulates that the quotient sheaf $\Omega_X/\Oo_X(L_1)$ has
non-negative degree with respect to any ample divisor on $X$. In
particular, we deduce that
\begin{equation}
  \label{eq:gen-semipos}
  0 
  \leq H\cdot c_1(\Omega_X/\Oo_X(L_1)) 
  = H\cdot K_X -H\cdot L_1.
\end{equation}    
The above inequality and the formula for $L_1$ in \eqref{eq:L1andL2}
imply
\begin{equation}
  \label{eq:DeltaH}
  H\cdot B_1 \leq \frac{1}{2}\,(3H\cdot K_X+d).
\end{equation}
Combining  \eqref{eq:Delta2} and \eqref{eq:DeltaH} we obtain
\[
 B_1^2 \leq \frac{3}{2}H\cdot K_X +\frac{1}{4}\,d.
\]
This inequality together with \eqref{eq:ChernCl} gives the estimate
\[
  c_2-K_X^2 -H\cdot K_X \geq 
  -\frac{1}{2}\,H\cdot K_X +\frac{1}{4}\,d,
\]
which is equivalent to the second inequality of the lemma. 
\qed

\begin{lem}
  \label{l:lastCase}
If the Bogomolov filtration of $\Tt_\xi$ is
$0=\Ff_0\subset\Ff_1\subset\Ff_2\subset\Ff_3=\Tt_\xi$, then it
determines two divisors $B_i$, $i=1,2$, in the positive cone of $N(X)$
and an effective nonzero divisor $E$ such that the following
conditions hold:
\begin{enumerate}[label={\rm \arabic*)}]
  \listspace
\item 
  $c_1(\Ff_1) = \frac{1}{3} \left(2B_1+B_2-H \right)$,
\item
  $K_X = \frac{1}{3}\,(B_1+2B_2-2H) +E$,
\item
  $c_2-K_X^2 \geq
  \begin{cases}
    H\cdot K_X+B^2_1 +B_1\cdot B_2, &
    \text{if}\quad  H\cdot c_1(\Ff_1)\leq0 \\ \\
    \dfrac{3}{4}\,H\cdot K_X
    +\dfrac{1}{8}\,\dfrac{(H\cdot B_1)(H\cdot B_2)}{d}\, & 
    \text{if}\quad  H\cdot c_1(\Ff_1)>0. 
  \end{cases}$
\end{enumerate}
\end{lem}

\proof
Set $L_i=c_1(\Ff_i/\Ff_{i-1})$, $1\leq i\leq3$. Hence
$c_1(\Ff_i)=L_1+\cdots+L_i$, for $i=1,2,3$. In particular, 
\[
  -H = c_1(\Tt_\xi) = c_1(\Ff_3) = L_1+L_2+L_3.
\]
Since $\Ff_1$ is the maximal destabilizing subsheaf of $\Ff_2$
\begin{equation}
  \label{eq:B1positive}
  B_1 = L_1-L_2 = 2\bigg(L_1-\frac{c_1(\Ff_2)}{2}\bigg) \in N^+(X),
\end{equation}
and similarly, since $\Ff_2/\Ff_1$ is the maximal destabilizing
subsheaf of $\Ff_3/\Ff_1$,
\begin{equation}
  \label{eq:B2positive}
  B_2 = L_2-L_3 =
  2\bigg(c_1(\Ff_2/\Ff_1)-\frac{c_1(\Ff_3/\Ff_1)}{2}\bigg) \in N^+(X).
\end{equation}
Combining \eqref{eq:B2positive} and \eqref{eq:B1positive}
together with the decomposition $H=-L_1-L_2-L_3$, we obtain the
following formulas 
\begin{equation}
  \label{eq:theOtherBasis}
\begin{aligned}
  L_1 &= -\frac{1}{3}\,H+\frac{2}{3}\,B_1+\frac{1}{3}\,B_2 \\
  L_2 &= -\frac{1}{3}\,H-\frac{1}{3}\,B_1+\frac{1}{3}\,B_2 \\
  L_3 &= -\frac{1}{3}\,H-\frac{1}{3}\,B_1-\frac{2}{3}\,B_2.
\end{aligned}
\end{equation}

As in the two previous lemmas, we use the Bogomolov filtration of
$\Tt_\xi$ to estimate its second Chern number
$c_2(\Tt_\xi)=c_2-K_X^2-H\cdot K_X$. Since the filtration is a
maximal ladder, it yields
\[
c_2-K_X^2-H\cdot K_X \geq L_1\cdot L_2 +L_1\cdot L_3 +L_2\cdot L_3.
\]
Substituting the formulas from \eqref{eq:theOtherBasis} leads to
\begin{equation}
  \label{eq:conclusion}
  c_2-K_X^2-H\cdot K_X \geq \frac{1}{3}\,(d-B_1^2-B_2^2-B_1\cdot B_2) 
  = \frac{1}{3}
  \left( d-\frac{1}{4}(2B_1 +B_2)^2 -\frac{3}{4}B^2_2 \right) .
\end{equation}
The argument continues, as in in the proof of Lemma \ref{l:Qstable},
according to the sign of the intersection $H\cdot L_1$.

\medskip 
\case 
If $H\cdot L_1\leq 0$, then the formula for $L_1$ in
\eqref{eq:theOtherBasis} gives 
$d\geq2H\cdot B_1+H\cdot B_2= H\cdot(2B_1 +B_2)$. This and the Hodge
Index Theorem imply 
\begin{equation}
  \label{eq:si-d}
 (2B_1+B_2)^2 \leq d.
\end{equation}
As a consequence, the inequality \eqref{eq:conclusion} becomes
\[
  c_2-K_X^2-H\cdot K_X \geq 
  \frac{1}{4}\,(d-B_2^2) \geq B^2_1 +B_1\cdot B_2,
\]
where the last inequality is obtained by substituting the upper bound
for $B_2^2$ from \eqref{eq:si-d}.  Hence the first inequality in part
3) of the lemma.

\medskip 
\case 
If $H\cdot L_1>0$, then as in the proof of Lemma \ref{l:Qstable}, we
obtain $L^2_1 \leq 0$. This and the formula for $L_1$ in
\eqref{eq:theOtherBasis} imply
\begin{equation}
  \label{eq:L1SquareNonPositive}
  2H\cdot(2B_1+B_2) \geq d+(2B_1+B_2)^2.
\end{equation}

Next we exploit the subsheaf $\Ff_2$ of the Bogomolov filtration of
$\Tt_\xi$.  Namely, combining the inclusion $\Ff_2\subset\Tt_\xi$
with the epimorphism of the extension sequence \eqref{eq:extSeq0}
gives rise to a nonzero morphism
\[
  u: \Ff_2 \longrightarrow \Omega_X.
\]
Arguing exactly as in the proof of Lemma \ref{l:rank2ForF} one shows
that this morphism is generically an isomorphism.  Thus we can apply
Lemma \ref{l:mainTechnical} to obtain the decomposition
\[
  K_X = c_1(\Ff_2)+E = L_1+L_2+E = -H-L_3+E,
\]
where $E$ is the support of the cokernel of $u$.  This is the
decomposition asserted in part 2) of the lemma.

Furthermore, the last assertion in Lemma \ref{l:mainTechnical} tells
us that $(-H-L_3)\cdot H \leq 0$.  This inequality combined with the
formula for $L_3$ in \eqref{eq:theOtherBasis} gives the following
inequality:
\[
  0 \leq H\cdot(H+L_3) =
  d-\frac{1}{3}\,H\cdot(H+B_1+2B_2)
  = \frac{1}{3}\,(2d-H\cdot B_1-2H\cdot B_2).
\]
Thus
\begin{equation}
  \label{eq:B1HDelta2H}
  H\cdot B_1+2H\cdot B_2 \leq 2d.
\end{equation}
Rewriting this inequality as 
$H\cdot B_2\leq d-\frac{1}{2}\,H\cdot B_1$, implies that 
$H\cdot B_2 < d$ and hence, by the Hodge Index Theorem, we obtain
\begin{equation}
  \label{eq:B2Square}
  B_2^2 \leq \frac{H\cdot B_2}{H^2}\,H\cdot B_2  
  \leq \frac{H\cdot B_2}{d}\bigg(d-\frac{1}{2}\,H\cdot B_1\bigg) 
  = H\cdot B_2 -\frac{1}{2d}\,(H\cdot B_1)(H\cdot B_2). 
\end{equation}

We now return to the inequality \eqref{eq:conclusion}.  Substituting
the upper bounds for the second (resp. third) term from
\eqref{eq:L1SquareNonPositive} (resp. \eqref{eq:B2Square}), 
we obtain 
\[
\begin{aligned}
  c_2-&K_X^2-\,H\cdot K_X 
  \geq \frac{1}{12}\,\left( 5d -4H\cdot B_1 -5H\cdot B_2 \right)
  +\frac{1}{8d}\,(H\cdot B_1)(H\cdot B_2) \\
  &= \frac{1}{12} \left( 4d -2H\cdot B_1 -4H\cdot B_2 \right)
  +\frac{1}{12} \left(d -2H\cdot B_1 -H\cdot B_2 \right)
  +\frac{1}{8d}\,(H\cdot B_1)(H\cdot B_2) \\
  &\geq -\frac{1}{4}\,H\cdot L_1 + \frac{1}{8d}\,(H\cdot B_1)(H\cdot B_2),
\end{aligned}
\]
where the last inequality is obtained by using \eqref{eq:B1HDelta2H}
for the first parenthesis, and the formula for $L_1$ in
\eqref{eq:theOtherBasis} for the second one.  Thus we obtain
\[
  c_2-K_X^2 
  \geq \frac{3}{4}\,H\cdot K_X +\frac{1}{4}\,(K_X -L_1)\cdot H 
  +\frac{1}{8d}\,(H\cdot B_1)(H\cdot B_2) 
  \geq \frac{3}{4}\,H\cdot K_X +\frac{1}{8d}(H\cdot B_1)(H\cdot B_2),
\] 
where the last inequality uses $(K_X -L_1)\cdot H \geq 0$ coming
from the generic semi-positivity of $\Omega_X$, see the discussion
just above \eqref{eq:gen-semipos}.  This completes the proof of the
second inequality in part 3) of the lemma.
\qed

As a corollary of Theorem \ref{th:mainTechnical} we deduce the
following property of surfaces lying on a quartic hypersurface in
$\PP^4$.

\begin{thm}
  \label{th:main}
If $X\subset\PP^4$ is a smooth surface lying on a hypersurface of
degree four and not on one of a smaller degree, then 
$K_X^2 < 6\chi(\Oo_X)$.
\end{thm}

\proof
Assume $K_X^2\geq 6\chi(\Oo_X)$.  We claim that $X$ must be rational
or irrationally ruled.  Indeed, if this is not the case, Theorem
\ref{th:mainTechnical} tells us that $K_X^2=6\chi(\Oo_X)$ and 
$H\cdot K_X=0$.  The second identity implies that $K_X=0$ in $N(X)$.
Hence $K_X^2=\chi(\Oo_X)=0$.  In particular, $X$ must be irregular
with irregularity $q=1$ or $2$, and of degree $d=10$ (this is obtained
from the double point formula).

Furthermore, Theorem \ref{th:mainTechnical} tells us that the vector
bundle $\Tt_\xi$ is Bogomolov unstable with the filtration treated in
Lemma \ref{l:rank2ForF}. From the proof of that lemma we deduce the
isomorphisms
\[
  H^0(\Omega_X) \cong H^0(\Tt_\xi) \cong H^0(\Ff_1).
\]
This implies that we must have $q=1$ (otherwise $X$ is an abelian
surface and the above isomorphisms imply $\Ff_1\cong\Omega_X$, hence
a splitting of the extension sequence \eqref{eq:extSeq0}).  Thus the
Albanese variety $\Alb(X)$ is an elliptic curve, call it $A$, and the
Albanese map $a:X\to A$ is an elliptic fibration.  In addition, since
$c_2(X)=0$, we also see that every {\it reduced} fibre is smooth.  In
particular, there are no smooth rational curves on $X$.  This fact, in
turn, implies that any divisor in the positive cone of $N(X)$ is
ample.

Another aspect of the filtration in Lemma \ref{l:rank2ForF} is the
decomposition of the canonical divisor
\[
  0 = K_X = L_1 +E = \frac{1}{3}\,(B_1-2H)+E,
\]
where $L_1 =c_1(\Ff_1) =\dfrac{1}{3}\,(B_1-2H)$, see
\eqref{eq:lbound-filt}, and $E$ is an effective nonzero divisor.  The
above equation can be rewritten as follows 
\begin{equation}
  \label{eq:HEB1}
  2H = 3E +B_1,
\end{equation}
with $B_1$ in the positive cone.  Hence, as remarked above, the
divisor $B_1$ is ample.

Next we investigate the divisor $E$.  Using formula \eqref{eq:HEB1},
we deduce that
\begin{equation*}
  \label{compE}
  C\cdot E +\frac{1}{3}\,C\cdot B_1 =  \frac{2}{3}\,H\cdot C,
\end{equation*}
for any reduced irreducible component $C$ of $E$.  On the other hand,
we know that for every $C$ as above the subsheaf $\Ff_1\subset\Tt_\xi$
gives rise to a nonzero morphism (see \cite{NaRe})
\begin{equation}
  \label{eq:morph-C-E}
  \Oo_C(L_1+H) \lra \Theta_X \otimes \Oo_C.
\end{equation}
Using again the expression for $L_1$ from \eqref{eq:LandL}, we obtain
\begin{equation}
  \label{eq:Cdeg-lb}
  (L_1+H)\cdot C  = \frac{1}{3}\,(H+B_1)\cdot C > 0. 
\end{equation}
Hence the composition of
the morphism \eqref{eq:morph-C-E} with the morphism
$\Theta_X\otimes\Oo_C\to\Oo_C(C)$ coming from the normal sequence of
$C\subset X$ yields a nonzero morphism\footnote{The fact that $X$
  contains no smooth rational curve is used here again.}
\[
  \Oo_C(L_1+H) \longrightarrow \Oo_C(C).
\]
This morphism together with formula \eqref{eq:Cdeg-lb} leads to
\begin{equation}
  \label{eq:C2}
  C^2 \geq (L_1 +H)\cdot C = \frac{1}{3}\,(H+B_1)\cdot C > 0,
\end{equation}
implying that $E$ lies in the positive cone of $N(X)$ and hence it is
ample.

From \eqref{eq:HEB1}, it follows that
\[
  20 = 2H^2 = 3E\cdot H +B_1\cdot H \geq 3E\cdot H +4.
\]
Thus $E\cdot H\leq 5$ and by the Hodge Index Theorem, we obtain
$E^2\leq 2$.  Since the intersection form on $\NS(X)$ is even and $E$
is ample, we deduce the equality $E^2 =2$.  Hence $E$ is reduced and
irreducible.  Thus replacing $C$ by $E$ in \eqref{eq:C2}, we obtain
\[
  6 = 3E^2 \geq H\cdot E +B_1\cdot E.
\] 
But the Hodge index tells us that $H\cdot E=5$ and that 
$B_1\cdot E\geq 2$ which contradicts the above inequality.

\medskip

Now we know that $X$ is either rational or irrationally ruled.  In the
latter case $\chi(\Oo_X)\leq 0$ and 
$8\chi(\Oo_X)\geq K_X^2\geq 6\chi(\Oo_X)$, implying
$K_X^2=\chi(\Oo_X)=0$.  This identifies $X$ as the projectivization
$\PP(\Ee)$ of a rank $2$ bundle $\Ee$ over an elliptic curve. But then
it is well-known that $X$ is an elliptic scroll of degree 
$d=-H\cdot K_X =5$, see Lemma \ref{l:scrollsInP4}.  It is easy to see
(and again well-known, see Theorem \ref{th:cubicsSections}) that such an
$X$ is contained in cubic hypersurfaces.  This of course is contrary
to our assumption that $4$ is the smallest degree of a hypersurface
containing $X$.

We now turn to the remaining possibility: $X$ is rational.  In this
case $\chi(\Oo_X)=1$ and hence, $K_X^2\geq6\chi(\Oo_X)=6$.  By
Riemann-Roch applied to $\Oo_X(-K_X)$ it follows that
$h^0(\Oo_X(-K_X))\geq 7$. Therefore, $(-K_X)$ is an effective nonzero
divisor.

Next we wish to have an upper bound for $(-H\cdot K_X)$.  This can be
done by observing that
\[
  h^0(\Oo_X(H)) \geq \chi(\Oo_X(H)) = \frac{H\cdot(H-K_X)}{2}+1,
\] 
which yields
\[
  d -H\cdot K_X \leq 2h^0(\Oo_X(H)) -2.
\]
Furthermore, it is well-known that $h^0(\Oo_X(H))=5$, unless $X$ is
the projection of the Veronese surface from a general point in 
$ \PP^5$. But such a projection lies on a hypersurface of degree $3$
and this is contrary to our assumption that the smallest degree of a
hypersurface containing $X$ is $4$.  Thus $h^0(\Oo_X(H))=5$ and we
obtain the upper bound
\begin{equation}
  \label{eq:-KH}
  -H\cdot K_X \leq 8-d.
\end{equation}
However, one observes the following.

\paragraph{Claim} $H^0(\Oo_X(-K_X -H))\neq 0$.

\medskip This concludes the argument, since 
$H^0(\Oo_X(-K_X -H))\neq 0$ together with the earlier estimate
$h^0(\Oo_X(-K_X))\geq 7$ imply $(-H\cdot K_X)>d$.  From this
inequality and \eqref{eq:-KH}, one obtains $d=3$.  It follows
immediately that $X$ is the projection of the Veronese surface from a
point on it.  But then, an easy dimension count implies that $X$
lies on a hypersurface of degree $2$, contradicting our hypothesis.

We now return to the proof of the claim.  It is equivalent to showing
that the homomorphism $H^0(\Oo_X(-K_X))\to H^0(\Oo_C(-K_X))$ induced
by the restriction of sections to a smooth irreducible curve $C\in|H|$
is not injective.  Assume on the contrary that it is injective. Then
\[
  7 \leq h^0(\Oo_X(-K_X)) \leq h^0(\Oo_C(-K_X)).
\]
If $\Oo_C(-K_X)$ is non-special, then we compute the right hand side
in the above inequality from the Riemann-Roch and obtain
\[
  7 \leq h^0(\Oo_X(-K_X)) 
  \leq h^0(\Oo_C(-K_X))
  = \deg(-K_X|_C)+1-g(C) = \frac{-3H\cdot K_X-d}{2}.
\]
Equivalently, $-H\cdot K_X \geq \dfrac{14+d}{3}$, which together with
\eqref{eq:-KH} tells us that $d\leq 2$ which is impossible since $d$
must be at least $3$.

If $\Oo_C(-K_X)$ is special, then by the Clifford inequality
\[
  14 \leq 2h^0(\Oo_X(-K_X)) 
  \leq 2h^0(\Oo_C(-K_X)) 
  \leq -H\cdot K_X+2.
\]
Hence $(-H\cdot K_X)\geq 12$, which again is impossible in
view of the upper bound \eqref{eq:-KH}. 
\qed

\subsection{The proof of Theorem \ref{th:ISample} when $m_X=4$} 
\label{subsec:m4proof}

For the reader's convenience we restate the theorem in this case as
follows.

\begin{pro}
  \label{p:d-chi-finite}
For all surfaces of general type contained in a quartic hypersurface 
in $\PP^4$ and not in one of a smaller degree the following assertions
hold.
\begin{enumerate}[label={\rm \arabic*)}]
  \listspace
\item 
The slope of their Chern numbers 
$\alpha=\frac{K^2}{\chi}<6$. 
\item
For every rational number $\alpha<6$, there exists
$d(\alpha)\in\NN$ such that every surface of slope $\alpha$
has the degree $d\leq d(\alpha)$.
\item
For every rational number $\alpha<6$, all surfaces of the slope
$\alpha$ and degree $d$ have the holomorphic Euler characteristic 
$\chi\leq\chi(\alpha,d)$, where
\[
\chi(\alpha,d) = \frac{d^2 +20d}{8(6-\alpha)} .
\]  
\end{enumerate}
\end{pro}

\proof
Let $X$ be a smooth surface in $\PP^4$ with $m_X =4$, see
\eqref{eq:mX} for notation.  The first assertion of the proposition,
$\alpha_X<6$, is a reformulation of Theorem \ref{th:main}.

For the second assertion, we use the inequality \eqref{eq:chiPm} for
$m=m_X=4$ to obtain
\begin{equation}
  \label{eq:DeSchBound}
  \chi:=\chi(\Oo_X) \geq
  \frac{1}{96}\,d^3-\frac{19}{16}\,d^2+\frac{10}{3}\,d+\frac{5}{4},
\end{equation}
provided $d\geq 11$.

On the other hand the Chern numbers of $X$ are related to the degree
$d$ via the double point formula which we write in the form 
\[
  12\chi -2K^2_X = c_2 -K^2_X = 5H\cdot K_X +10d -d^2.
\]
Setting $K^2_X =\alpha \chi$ and substituting into the above equation
give
\begin{equation}
  \label{eq:1}
  5H\cdot K_X +10d -d^2
  = 12\chi -2K^2_X = 12\chi -2\alpha \chi
  = 2(6-\alpha)\chi.
\end{equation}
From \cite{ElPe} the sectional genus of $X$ has the following upper
bound 
\begin{equation}
\label{eq:ElPe4}
  2g(H)-2=H^2+H\cdot K_X\leq \frac{1}{4}\,d^2.
\end{equation}
Using this in the equation \eqref{eq:1} we obtain
\begin{equation}
  \label{eq:upperBound}
  2(6-\alpha)\chi
  = 5H\cdot K_X +10d -d^2
  = 5(H\cdot K_X +d) +5d -d^2
  \leq \frac{1}{4}\,d^2 +5d.
\end{equation}
This inequality combined with the first assertion of the proposition
and the lower bound for $\chi$ in \eqref{eq:DeSchBound} give
\begin{equation}
  \label{eq:aldeg}
  (6-\alpha)
  \bigg( 
    \frac{1}{12}\,d^3-\frac{19}{2}\,d^2+\frac{80}{3}\,d+10
  \bigg) 
  \leq d^2+20d,
\end{equation}
for all $d\geq 11$.  From this it follows that one can explicitly
determine a positive integer $d_0(\alpha)$ depending on $\alpha$ only,
which is an upper bound for the solutions of the above inequality.
Setting $d(\alpha)=\max(10,d_0(\alpha))$, one deduces the second
assertion of the proposition.

The upper bound for $\chi$ in the third assertion of the proposition
is obtained from \eqref{eq:upperBound} and the first assertion of the
proposition.
\qed

As a consequence, we deduce the following finiteness result.

\begin{cor}
  \label{c:slopefinite}
The number of families of surfaces in $\PP^4$ of general type with
fixed slope $\alpha$ and contained in a quartic hypersurface
(and not in one of a smaller degree) is at most finite.
\end{cor}

\proof
The components of the Hilbert scheme of surfaces in $\PP^4$ are
labeled by Hilbert polynomials.  So the proof of the statement comes
down to checking that there is at most a finite number of such
polynomials for the surfaces subject to the hypotheses of the
corollary.  Since the Hilbert polynomial of a surface $X\subset \PP^4$
of degree $d$ has the form
\[
  P_X(t) 
  = \frac{d}{2}\,t^2 - \frac{H\cdot K_X}{2}\,t +\chi(\Oo_X),
\]
one needs to see that there is only a finite number of possibilities
for the triples \linebreak
$\left( d,\, H\cdot K_X,\, \chi(\Oo_X) \right)$.  This is exactly what
Proposition \ref{p:d-chi-finite} tells us: parts 2) and 3) give a
finite number of possibilities for $d$ and $\chi(\Oo_X)$ respectively.
Once $d$ takes a finite number of values, the inequality
\eqref{eq:ElPe4} insures that there is only a finite number of values
for $H\cdot K_X$ as well.
\qed

\begin{remarks*}
1) If $\alpha\leq 5$, the inequality \eqref{eq:aldeg} gives rise to
the relation
\[
d^3 -126d^2 +80d +120 \leq 0,
\]
provided the polynomial on the left side of \eqref{eq:aldeg} is
non-negative\footnote{This is the case for $d\geq 112$.}.  From this
it follows that all surfaces $X$ of general type with slope
$\alpha\leq 5$ and $m_X =4$ have degree $d\leq 125$.

\smallskip

2) For complete intersections $(4,a)$ in $\PP^4$ with
$a=\frac{d}{4}\geq2$, we have $K_X=(a-1)H$ and
$\chi(\Oo_X)=p_g(X)+1=\binom{a+3}{4}-\binom{a-1}{4}$.  Hence
\[
\begin{aligned}
  K_X^2 &= 4a(a-1)^2 \\
  \chi(\Oo_X) &= \frac{1}{3}(2a^3-3a^2+7a-3) \\
  \alpha_X &= f(a)=\frac{12a(a-1)^2}{2a^3-3a^2+7a-3} \lra 6
\end{aligned}
\]
as $a\lra \infty$.  It is easily checked that $f(a)$ is an increasing
function of $a$ in the interval $[2, +\infty[$.  In particular,
$\alpha_X=f(a)>5$, for all $a\geq7$.  These examples suggest the
following questions: besides complete intersections, 
\begin{enumerate}[label={\rm \alph*)}]
  \listspace
\item 
  are there surfaces of general type in $\PP^4$ which are contained in
  a quartic $3$-fold (and not on one of a smaller degree) and whose
  slope $\alpha=K^2/\chi>5$?
\item
  is there an infinite sequence $(X_n)_n$ of surfaces of general type
  in $\PP^4$ with $m_{X_n} =4$ such that the slopes
  $\alpha_{X_n}=K^2_{X_n}/\chi(\Oo_{X_n})$ converge to $6$? 
\end{enumerate}
\end{remarks*}

\section{Numerical invariants for surfaces with $m_X =2,3,5$ } 
\label{subsecm235}

This section treats the remaining values of the smallest degree $m_X$
of a hypersurface containing $X\subset\PP^4$.  As in Section
\ref{s:m=4}, we assume that the Kodaira dimension of $X$ is
non-negative.  Consider the extension sequence
\eqref{eq:extensionSequence} corresponding to 
$\xi=\delta_X(s)\in H^1(\Theta_X(-K_X-(5-m)H))$ provided by Lemma
\ref{l:EP-coh}.

We begin with cases $m_X=2$ or $3$.

\begin{thm}
  \label{th:m=2} 
If $m_X=2$, then
\[
c_2-K_X^2\geq
\begin{cases}
  3H\cdot K_X+3\,d,
  & \text{if $\Tt_\xi$ is Bogomolov semistable} \\
  \dfrac{3}{2}\,H\cdot K_X +\dfrac{9}{4}\,d,
  & \text{if $\Tt_\xi$ is Bogomolov unstable}.
\end{cases}
\]
\end{thm}

\begin{thm}
  \label{th:m=3} 
If $m_X=3$, then
\[
c_2-K_X^2\geq
\begin{cases}
  2H\cdot K_X+\,\dfrac{4}{3}\,d,
  & \text{if $\Tt_\xi$ is Bogomolov semistable} \\
  \min \bigg(
  K\cdot H +d,\,\, \dfrac{3}{2}\,K\cdot H +\dfrac{1}{3}\,d
  \bigg),
  & \text{if $\Tt_\xi$ is Bogomolov unstable}.
\end{cases}
\]
\end{thm}

We omit the proofs since they follow exactly the same pattern as the
one of Theorem \ref{th:mainTechnical}.  We only mention, that in the
case $m_X=2$ the only possibility for the Bogomolov filtration is as
in Lemma \ref{l:Qstable}.

We now turn to the case $m_X =5$.  The sequence \eqref{eq:coh-seq}
takes the form
\[
  H^0(\Theta_{\PP^4}|_X(-K_X)) \lra
  H^0(\Nn_X(-K_X)) \stackrel{\delta_X}{\lra}
  H^1(\Theta_X(-K_X))
\]
and we seek an analogue of Lemma \ref{l:EP-coh}.

\begin{lem}
  \label{l:s5}
Let $s$ be a global section of $\Nn_X(-K_X)$ corresponding to a
quintic hypersurface containing $X$.  If $\delta_X(s)=0$ in
$H^1(\Theta_X(-K_X))$, then $d\leq 14$ and $K_X^2\leq6\chi(\Oo_X)$.
Furthermore, if $X$ is of general type, then $p_g\leq2$.
\end{lem}

\proof
The vanishing of $\delta_X(s)$ implies that the section $s$ of
$\Nn_X(-K_X)$ comes from a nonzero section of
$\Theta_{\PP^4}\otimes\Oo_X(-K_X)$.  Thus 
$H^0(\Theta_{\PP^4}\otimes\Oo_X(-K_X))\neq0$.  This together with
the Euler sequence of $\Theta_{\PP^4}$ leads to two possibilities:
\begin{enumerate}
  \listspace
\item 
  $H^0(\Oo_X(H-K_X)) \neq 0$,
\item
  $\ker \left( H^1(\Oo_X(-K_X))\to 
    H^0(\Oo_X(H))^\ast\otimes H^1(\Oo_X(H-K_X)) \right) \neq 0$.
\end{enumerate}
Both of them tell us that $H\cdot K_X\leq H^2=d$. Furthermore, the
equality holds if and only if $\Oo_X(K_X)=\Oo_X(H)$ and then one
knows that $X$ must be a complete intersection which is impossible due
to the condition $m_X =5$.  Thus we have 
\begin{equation}
  \label{eq:HKd}
  H\cdot K_X < d.
\end{equation} 
This and the double point formula give
\begin{equation}
  \label{eq:dp-0}
  2(K_X^2-6\chi) = d^2 -10d -5H\cdot K_X > d^2 -15d.
\end{equation}
In particular, one obtains that $d\leq 14$, provided 
$K_X^2 \leq 6\chi$.  Thus the first assertion of the lemma is a
consequence of the inequality $K_X^2 \leq 6\chi$.
 
Assume $K_X^2>6\chi$.  This and the assumption that the Kodaira
dimension of $X$ is non-negative implies $\chi\geq 1$.  Hence
$K_X^2\geq 7$ and $H\cdot K_X >0$ tell us that $X$ is of general
type.

Observe that the inequality $H\cdot K_X<d$ and the Hodge index give
also the upper bound $K_X^2<d$.  Substituting into \eqref{eq:dp-0}, we
deduce
\[
  0 > d^2-17d+12, 
\]
and hence $d\leq 16$.  This and the inequalities $6\chi<K_X^2<d$ imply 
\[
  K_X^2 \leq 15  \textq{and} \chi\leq2.
\]
Furthermore, in the case $K_X^2=15$ the degree $d$ must be $16$ and
the Hodge index $(H\cdot K_X)^2 \geq H^2 K^2\geq 16\cdot 15$ implies
$H\cdot K_X \geq 16 \geq d$ contrary to \eqref{eq:HKd}.  Thus one has
$K^2_X \leq 14$.

We now examine the remaining possibilities according to two possible
values of $\chi=1$ or $2$.
If $\chi=1$, then the double point formula reads
\begin{equation}
  \label{dp-5chi1}
  d^2 -10d -5H\cdot K_X = 2K_X^2 -12.
\end{equation}
By Bogomolov-Miyaoka-Yau inequality, $K_X^2 \leq 9\chi =9$.  Thus
$K_X^2\in\{7,8,9\}$, implying that
\begin{equation}
  \label{ineq:dp-5chi1}
  d^2 -10d \geq 2+5H\cdot K_X \geq 32,
\end{equation}
where the last inequality comes from the Hodge index 
$(H\cdot K_X)^2\geq H^2 K_X^2=d\,K_X^2\geq 5\cdot7=35$ and where the
last inequality uses $d\geq 5$, coming from the fact that $X$ is of
general type.  Hence $d\geq 13$.  Using this lower bound in the Hodge
Index estimate of $H\cdot K_X$ once again, we obtain 
$H\cdot K_X \geq 10$.  Substituting this in \eqref{ineq:dp-5chi1},
implies $d\geq 14$.  With the previously derived upper bound 
$d\leq 16$, we obtain $d= 14, 15,16$ are the only possible values.  A
direct check shows that $d=15$ is incompatible with the double point
formula \eqref{dp-5chi1}, while for $d=14$ and $16$, that formula
forces $K_X^2 =9$ and hence, $H\cdot K_X =10$ and $18$, respectively.
The latter value contradicts the inequality \eqref{eq:HKd}, while the
former, the Hodge Index inequality.

If $\chi=2$, then the double point formula becomes
\begin{equation}
  \label{dp-5chi2}
  d^2 -10d -5H\cdot K_X = 2K_X^2 -24,
\end{equation}
while $13,14$ are the only possible values for $K_X^2$.  Arguing as in
the previous case we obtain the lower bound $H\cdot K_X \geq 9$ which
together with \eqref{dp-5chi2} gives the possibilities
$d=14,15,16$ for the degree of $X$.  But then, going back to the lower
bound estimate for $H\cdot K_X$, we obtain 
$H\cdot K_X \geq 14$. Substituting into \eqref{dp-5chi2} gives
$d\geq 15$ and hence, $d=15$ or $16$.  The value $d=15$ is again
incompatible with the double point formula \eqref{dp-5chi2}, while for
$d=16$ that formula tells us that $K_X^2 =14$ and $H\cdot K_X =16$.
But this contradicts the inequality \eqref{eq:HKd}.

The last statement saying that $p_g \leq 2$, for $X$ of general type,
can be seen as follows.  Let $X_0$ be the minimal model of $X$ and set
$\sigma: X\to X_0$ to be a sequence of blow-down maps.  Then the
canonical divisor $K_X $ can be written as follows
\[
  K_X = \sigma^\ast K_{X_0} +D,
\]
where $D$ is an effective divisor composed of curves contracted by
$\sigma$.  From the nonvanishing of
$H^0(\Theta_X(-K_X))=H^0(\Theta_X(-\sigma^\ast K_{X_0}-D)$ it follows
that $H^0(\Theta_X(-\sigma^\ast K_{X_0}))\neq0$ as well. Running the
argument of the first paragraph of the proof for that group and using
the fact that $H^1(\Oo_X(-\sigma^\ast K_{X_0}))=0$, we obtain
$H^0(\Oo_X(H-\sigma^\ast K_{X_0}))\neq 0$.

Next we observe that $(H-\sigma^\ast K_{X_0})\neq 0$, since otherwise
$H=\sigma^\ast K_{X_0}=K_X$ and, as it was argued above, this is
incompatible with the condition $m_X =5$.

Once $(H-\sigma^\ast K_{X_0})\neq 0$, we take a divisor 
$\Gamma\in| H-\sigma^\ast K_{X_0} |$ and identify $p_g$ as follows
\[
  p_g = h^0(\Oo_X(\sigma^\ast K_{X_0})) = h^0(\Oo_X(H-\Gamma)), 
\]
\ie $p_g$ is the dimension of the space of hyperplanes in $\PP^4$
containing $\Gamma$.  From this it follows that $p_g \leq 3$ and the
equality holds if and only if $\Gamma$ is a line in $\PP^4$.  We claim
that this is impossible.  Indeed, if $\Gamma$ is a line then the
identity $H = \Gamma+\sigma^\ast K_{X_0}$ implies
\begin{equation}
  \label{line}
  1 = \Gamma\cdot H
  = \Gamma^2 +\Gamma\cdot\sigma^\ast K_{X_0}.
\end{equation}
Hence $\Gamma\cdot \sigma^\ast K_{X_0} >0$ and hence $\Gamma$ is not
in the exceptional divisor $D$.  Therefore $\Gamma\cdot D\geq0$.
But then the identity \eqref{line} can be rewritten as 
\[
  1 = \Gamma^2 +\Gamma\cdot \sigma^\ast K_{X_0}
  = \Gamma^2 +\Gamma\cdot(K_{X}-D) 
  = (\Gamma^2 +\Gamma\cdot  K_{X}) -\Gamma\cdot D 
  = -2 -\Gamma\cdot D, 
\]
which is absurd.  Thus $p_g\leq 2$ as asserted.
\qed

\paragraph{Question} Can one enumerate all surfaces in $\PP^4$ with
$m_X=5$ and $\delta_X(s)=0$?

\medskip 

From now on we set $\xi=\delta_X(s)\in H^1(\Theta_X(-K_X))$ and assume
it to be nonzero.  The corresponding extension sequence has the form
\begin{equation}
  \label{eq:extSeq2}
  0 \lra \Oo_X(-K_X) \lra \Tt_\xi \lra \Omega_X \lra 0.
\end{equation} 
The Bogomolov semistability/instability considerations of the sheaf
$\Tt_\xi$ above give the following.

\begin{thm}
  \label{th:m5}
Let $X$ be a surface in $\PP^4$ with $m_X =5$.  Then, either
$K_X^2\leq c_2$, or $X$ is a surface of general type subject to the
following properties:
\begin{enumerate}[label={\rm \roman*)}]
  \listspace
\item
  The canonical divisor $K_X$ admits a distinguished decomposition
  $K_X =L+E$, such that $L$ is in the positive cone of $N(X)$ and $E$
  is an effective nonzero divisor.
\item
  $0 < K_X^2-c_2 \leq \dfrac{2}{3}\,L^2 \leq \dfrac{2}{3}\,d$.
\end{enumerate}
\end{thm}

\proof
We may assume that $K_X^2 > c_2$, since otherwise there is nothing to
prove.  This assumption and Lemma \ref{l:s5} insure the nonvanishing
of the cohomology class $\xi=\delta_X(s)$ defined by a section 
$s\in H^0(\Nn_X(-K_X))$ arising from a hypersurface of degree $5$
containing $X$.  We associate to $\xi$ the extension sequence
\eqref{eq:extSeq2} and observe that the assumption $K_X^2>c_2$ implies
that the sheaf $\Tt_\xi$ in the middle of that sequence is Bogomolov
unstable.  Let $\Ff$ be the maximal Bogomolov destabilizing subsheaf
of $\Tt_\xi$.  Observe that the Bogomolov destabilizing property
implies that $\Oo_X(L)$, the determinant of $\Ff$, is in the positive
cone of $N(X)$.

Putting the inclusion $\Ff\subset \Tt_\xi$ together with the extension
sequence \eqref{eq:extSeq2} we obtain the diagram
\vspace{7ex}
\begin{equation}
  \label{d:Fseq-xi}
  \begin{tikzpicture}[overlay,every node/.style={draw=none},
    ->,inner sep=1.1ex]
    \matrix [draw=none,row sep=3.8ex,column sep=4.2ex]
    {
      && \node (02) {$0$}; \\
      && \node (12) {$\Ff$}; \\
      \node (20) {$0$};
      & \node (21) {$\Oo_X(-K_X)$};
      & \node (22) {$\Tt_\xi$};
      & \node (23) {$\Omega_X$};
      & \node (24) {$0$}; \\
    };
    \path
      (02) edge (12)
      (12) edge (22)
      (12) edge node[above right=-.3ex,scale=.75] {$\phi$} (23)
      (20) edge (21)
      (21) edge (22)
      (22) edge (23)
      (23) edge (24)
    ;
  \end{tikzpicture}
  \vspace{8.5ex}
\end{equation}
where the slanted arrow is the composition of the inclusion together
with the epimorphism in \eqref{eq:extSeq2}.  We claim that $\Ff$ must
be of rank $2$ and the morphism $\phi$ in the above diagram is
generically an isomorphism.  Indeed, if the rank of $\Ff$ is one, then
it is $\Oo_X(L)$ and $\phi$ must be zero, since by the Bogomolov Lemma,
$\Omega_X$ admits no rank $1$ subsheaves of Iitaka dimension $2$.  But
the vanishing of $\phi$ implies $\Oo_X(L)=\Oo_X(-K_X)$.  Hence $X$
is rational with $K_X^2>6\chi=6$.  This together with the Hodge index
gives a lower bound on the degree of $-K_X$:
\begin{equation}
  \label{lbound5}
  -K_X\cdot H \geq 5.
\end{equation}
On the other hand, the Riemann-Roch for $\Oo_X(H)$ yields
\[
  h^0(\Oo_X(H)) \geq \frac{d-K_X\cdot H}{2} +1.
\]
Furthermore, $h^0(\Oo_X(H))=5$, since $X$ is not the projection of
the Veronese surface from a point outside its secant variety%
\footnote{Such a surface is contained in cubic hypersurfaces
  and this is contrary to our assumption that the minimal degree of a 
  hypersurface containing $X$ is $5$.}. 
Substituting this in the above inequality we obtain 
\[
  d-K_X\cdot H \leq 8.
\]
This and the lower bound \eqref{lbound5} imply $d=3$ and 
$-K_X\cdot H=5$.  Therefore, a general hyperplane section, call it $C$, 
is a rational normal cubic in $\PP^3$ and hence $C$ is cut out by 
quadrics in $\PP^3$.  But then, in view of the isomorphism
\[
  H^0(\Jj_X(2)) \cong H^0(\Jj_C (2)),
\]
one obtains quadric hypersurfaces containing $X$ contradicting the
assumption $m_X =5$.

We know now that the rank of $\Ff$ is $2$ and we need to check that
the morphism $\phi$ in the diagram \eqref{d:Fseq-xi} is generically an
isomorphism.  Assuming that this is not the case and taking the second
exterior power of the diagram \eqref{d:Fseq-xi}, we deduce a nonzero
morphism $\Oo_X(L)\to\Theta_X$.  This implies once again that $X$ must
be rational with $K_X^2>6\chi=6$.  Thus we are back in the situation
considered in the preceding paragraph implying that $\phi$ in
\eqref{d:Fseq-xi} is generically an isomorphism.  Hence we are in the
position to apply Lemma \ref{l:mainTechnical}; in particular, we
obtain the decomposition of the canonical divisor $K_X$ asserted in i)
of the theorem.  Furthermore, Lemma \ref{l:mainTechnical}, 3) gives us
the upper bound
\[
  L\cdot H \leq H^2 = d.
\]
This, by the Hodge Index Theorem, implies 
\begin{equation}
  \label{L2bound}
  L^2 \leq d.
\end{equation}  
The part ii) of the theorem is obtained by estimating the second Chern
number of $\Tt_\xi$ from the vertical sequence in the diagram
\eqref{d:Fseq-xi}.  Namely, we complete it to the exact sequence
\[ 
  0 \lra \Ff \lra \Tt_\xi \lra \Jj_Z(-L) \lra 0
\]
from which we obtain
\[
  c_2-K_X^2 
  = c_2(\Tt_\xi) 
  = c_2(\Ff) -L^2 +\deg(Z) 
  \geq c_2(\Ff) -L^2 
  \geq \frac{1}{3}\,L^2 -L^2 
  = -\frac{2}{3}\,L^2, 
\]
where the second inequality is a result of Miyaoka, \cite[Remark
4.18]{Miy}.  Hence the first upper bound
\[
  K_X^2 -c_2 \leq \frac{2}{3}\,L^2
\]
asserted in the part  ii) of the theorem.  This together with
\eqref{L2bound} gives the second upper bound in ii). 
\qed

\begin{pro}
  \label{p:d-chi-finite-5}
For all surfaces of general type lying on a quintic hypersurface in
$\PP^4$ and not on one of a smaller degree, the following assertions
hold:
\begin{enumerate}[label={\rm \arabic*)}]
  \listspace
\item
  For every rational number $\alpha \neq 6$, there exists
  $d(\alpha)\in\NN$ such that every surface of slope $\alpha$,
  has the degree $d\leq d(\alpha)$.
\item
  For every rational number $\alpha \neq 6$, all surfaces with slope
  $\alpha$ and degree $d\leq d(\alpha)$ have the holomorphic Euler
  characteristic subject to the inequalities 
\[
\chi \leq \begin{cases}  
\frac{5d}{6-\alpha},&\mbox{ if $\alpha<6$}, \\ 
  \frac{d}{3(\alpha-6)},& \mbox{ if $\alpha>6$}.
\end{cases}
\]
\item
  The slope of Chern numbers satisfies
  $\alpha=\frac{K^2}{\chi}<6+\frac{1}{2}$, with possible exceptions of
  surfaces having degree $\leq37$ and holomorphic Euler characteristic
  $\leq24$.
\end{enumerate}
\end{pro}

\proof
For $\alpha <6$ the argument is analogous to the one in the proof of
Proposition \ref{p:d-chi-finite}.  Namely, the Decker-Schreyer
polynomial in \eqref{eq:polynomialPm} for $m=5$ has the form
\[
  P_5(d) = \frac{1}{25}\,d(d^2 -40d +95)
\]
and their result gives the lower bound 
\begin{equation}
  \label{D-S-5}
  \chi \geq \frac{1}{25}\,d(d^2 -40d +95)
\end{equation}
for $\chi$, the holomorphic Euler characteristic, for every surface of
degree $d\geq 18$ contained in a quintic hypersurface and not in one
of a smaller degree.  Combining this with the double point formula, we
obtain
\begin{equation}
  \label{chi-d-l6}
  \frac{2}{25}\,(6-\alpha)d(d^2 -40d +95) 
  \leq 2(6-\alpha)\chi 
  = 12\chi -K^2 
  = 5H\cdot K_X +10d-d^2.
\end{equation}
From \cite{ElPe} we also know that the genus $g_H$ of a (smooth)
hyperplane section of $X$ is subject to the inequality
\begin{equation}
  \label{E-P-5}
  d+H\cdot K_X = 2g_H -2 \leq \frac{d^2+5d}{5}.
\end{equation} 
From this and \eqref{chi-d-l6}, we deduce 
\[
  \frac{1}{25}(6-\alpha)(d^2 -40d +95) \leq 5d.
\]
Hence one can explicitly determine a positive integer $d_0(\alpha)$
depending only on $\alpha$, which is an upper bound for the integer
solutions of the above inequality.  Setting 
$d(\alpha):= \max(17,d_0(\alpha))$, we obtain assertion 1) of the
proposition in the range $\alpha<6$.

If $\alpha>6$, then Theorem \ref{th:m5} ii) tells us that
$(\alpha-6)\chi=K^2-6\chi\leq\frac{1}{3}d$.  Combining this with
\eqref{D-S-5} gives the inequality
\begin{equation}
  \label{chi-d-g6}
  \frac{1}{25}\,(\alpha-6)(d^2 -40d +95) 
  \leq \frac{(\alpha-6)\chi}{d}  
  \leq \frac{1}{3},
\end{equation}
which provides the assertion 1) of the proposition in the range
$\alpha >6$. 

The second inequality for $\chi$ in 2) is just a restatement of
Theorem \ref{th:m5}, ii), while the first one is obtained by combining
\eqref{chi-d-l6} and the upper bound $5H\cdot K_X\leq d^2$ from
\eqref{E-P-5}.

For the third part of the proposition which asserts the upper bound
for the slope $\alpha$, we may assume $\alpha>6$ and go back to the
inequality \eqref{chi-d-g6}.  For $d\geq 38$, that inequality implies
\[
  \alpha-6 \leq \frac{25}{3\cdot 19} < \frac{1}{2}.
\]
For $d\leq 37$, we use the inequality in Theorem \ref{th:m5}, ii), to
deduce
\[
  (\alpha-6)\chi = K^2 -6\chi 
  \leq \frac{1}{3}\,d \leq \frac{37}{3}.
\] 
Hence $(\alpha-6)\chi\leq 12$ and $\alpha<6+\frac{1}{2}$, unless
$\chi\leq 24$.
\qed

As a consequence, we obtain the following finiteness result.

\begin{cor}
  \label{c:slopefinite5}
The number of families of surfaces in $\PP^4$ of general type with
fixed slope $\alpha\neq 6$ and contained in a quintic hypersurface
(and not in one of a smaller degree) is at most finite.
\end{cor}

\proof
See the proof of Corollary \ref{c:slopefinite}.
\qed

\vspace{1\baselineskip} 

\part*{\Large (II)\hspace{.75ex}
  Irregularity of surfaces in $\PP^4$}

\section{The irregularity of surfaces lying on a small 
  degree hypersurface} \label{sec:irreg}

One of the outstanding problems about surfaces in $\PP^4$ is the
control of their irregularity.  Since the beautiful work of Horrocks
and Mumford, \cite{Ho-Mu}, which contains a construction of an abelian
surface of degree $10$ in $\PP^4$, the irregularity $2$ remains the
maximal known value for surfaces in $\PP^4$.  Indeed, it is
conjectured that no surface of irregularity greater than $2$ can be
embedded into $\PP^4$.  To our knowledge there is no conceptual reason
for this phenomenon.

In our previous work, \cite{NaRe}, we were able to show that the
irregularity of surfaces in $\PP^4$ is bounded by $3$ under a certain
precise set of conditions, see \cite[Theorem 5.1]{NaRe}.  What is
perhaps more interesting, is that we have uncovered how to use the
cohomological condition $H^1(\Theta_X(-K_X))\neq 0$ to bound the
irregularity of a surface $X\subset\PP^4$.  In the previous sections
we have seen that this non-vanishing condition arises naturally
whenever $X$ is contained in a hypersurface of a degree $m\leq5$.
This phenomenon together with the metha-principle formulated in the
introduction suggests that the irregularity of surfaces in $\PP^4$,
with the exception of a finite number of families, is bounded by the
irregularity of surfaces contained in hypersurfaces of small degree.
Guided by this heuristic principle, this section begins the
investigation of the irregularity of surfaces in $\PP^4$ contained in
a small degree hypersurface.

Our study of irregular surfaces in $\PP^4$ lying on a small degree
hypersurfaces has the same unifying theme as before: the extension
construction.  We consistently interpret a small degree hypersurface
containing our surface as an extension sequence of sheaves on $X$.
Let $m$ be the smallest degree of a hypersurface containing a surface
$X\subset\PP^4$ and let $V_m$ be such a hypersurface.  We recall that
$\Nn_X$ (resp. $\Nn^\ast_X$) denotes the normal (resp. conormal)
bundle of $X$ in $\PP^4$ and one has the following identifications
\begin{equation*}
  \label{eq:nor-conorAgain}
  \Jj_X/\Jj_X^2 = \Nn^\ast_X \cong \det(\Nn^\ast)\otimes\Nn_X 
  = \Nn_X(-K_X-5H),
\end{equation*}
where $ \Jj_X $ is the ideal sheaf of $X$ in $\PP^4$ and the second
identification is due to the rank of $\Nn^\ast_X$ being two.  This
leads to the non-vanishing
\begin{equation*}
  \label{eq:nonVaninshing4Again}
  H^0(\Nn_X(-K_X-(5-m)H)) = H^0(\Nn^\ast_X(mH)) \neq 0
\end{equation*}
already encountered in \eqref{eq:nonVaninshing4}.
Thus, we associate to $V_m$ a nonzero global section, denoted by $s$,
of the twisted normal bundle $\Nn_X(-K_X-(5-m)H)$.  Our approach
consists of using this section to build up an appropriate extension
sequence.

In the context of the irregularity, there are two lines of thinking.
The first one is to use the coboundary map $\delta_X$ in the normal
exact sequence for $X\subset\PP^4$ to produce the cohomology class
$\xi=\delta_X(s)\in H^1(\Theta_X(-K_X-(5-m)H))$ and then, to view it,
via the natural identification
\[
  H^1(\Theta_X(-K_X-(5-m)H)) \cong 
  \Ext^1(\Omega_X,\Oo_X(-K_X-(5-m)H)),
\]
as an extension
\[
  0\lra \Oo_X(-K_X-(5-m)H) \lra \Tt_\xi \lra \Omega_X \lra 0.
\]
This was the idea exploited in \cite{NaRe} and it produces
satisfactory results provided that $X$ is of Albanese dimension
$2$, \ie the image of the Albanese morphism of $X$ is of dimension
$2$.  However, if $X$ fibers over a curve $B$ of genus $g(B)=q(X)$,
the method fails.

This brings us to the second way of associating an extension sequence
to the section $s\in H^0(\Nn_X(-K_X-(5-m)H))$.  Namely, we simply take
the Koszul sequence associated to $s$ to obtain
\begin{equation}
  \label{eq:koszulForS}
  0 \lra \Oo_X(K_X+(5-m)H) \stackrel{s}{\lra}
  \Nn_{X}\stackrel{\wedge s}{\lra} \Jj_Z(mH) \lra 0,
\end{equation}
where $Z\subset X$ is the scheme of zeros of $s$ and $\Jj_Z$ its ideal
sheaf.  This is of course very classical and yet efficient in
addressing the irregularity problem, provided we have a good control
of the subscheme $Z$.  Let us explain the main ingredients of this
approach as well as set up the stage for more technical considerations
in the subsequent sections.

The extension \eqref{eq:koszulForS} fails to be exact precisely when
$Z=(s=0)$ has divisorial part.  If this is the case, let $Z_1$ be the
divisorial part and $Z_0$ be the residual part of $Z_1$ in $Z$.  If
$s_1$ is a section of $\Oo_X(Z_1)$ defining $Z_1$, then $s=s_1s_0$,
where $s_0$ is a section of $\Nn_X(-K_X-(5-m)H-Z_1)$ whose zero-locus
is $Z_0$, a $0$-dimensional subscheme of $X$.  We now have the short
exact sequence
\begin{equation}
  \label{eq:koszulForS0}
  0 \lra \Oo_X(K_X+(5-m)H +Z_1) \stackrel{s_0}{\lra}
  \Nn_{X} \stackrel{\wedge s_0}{\lra} \Jj_{Z_0}(mH-Z_1) \lra 0,
\end{equation}
where $\Jj_{Z_0}$ is the ideal sheaf of $Z_0$.

At this stage we bring in the irregularity.  The main point in
relating the irregularity of $X$ to the extension sequence
\eqref{eq:koszulForS0} is the following general fact.

\begin{lem}
  \label{l:cotang-conorm}
If $X \subset \PP^n$ is a complex projective manifold of dimension
bigger than $1$, then 
\[
  H^0(\Omega_X) \cong H^1(\Nn_X^\ast),
\]
where $\Omega_X$(resp. $\Nn_X^\ast$) is the cotangent
(resp. conormal) bundle of $X$. 
\end{lem}

\proof
From the conormal exact sequence 
\[
  0 \lra \Nn_X^\ast \lra
  \Omega_{\PP^n}|_X \lra \Omega_X 
  \lra 0 
\]
of $X\subset\PP^n$ we obtain
\[
  0 \lra H^0(\Omega_X) \lra H^1(\Nn_X^\ast) \lra
  H^1(\Omega_{\PP^n}|_X) \stackrel{r}{\lra} H^1(\Omega_X) 
\]
where the injectivity on the left is the vanishing of
$H^0(\Omega_{\PP^n}\big|_X)$.  The asserted isomorphism follows from
the injectivity of the homomorphism 
$r:H^1(\Omega_{\PP^n}|_X)\to H^1(\Omega_X)$.  To establish it we use
the dual of the Euler sequence and the assumption $\dim(X)>1$ to
deduce $H^1(\Omega_{\PP^n}|_X)\cong H^0(X,\Oo_X)\cong\CC$.  At the
same time, we have the linear map
\[
  \CC \cong H^1(\Omega_{\PP^n})
  \stackrel{i^\ast}{\longrightarrow} H^1(\Omega_X)
\]
defined by the pullback $i^\ast$, where $i: X\hookrightarrow\PP^n$
is the inclusion morphism.  This linear map factors through
$H^1(\Omega_{\PP^n}|_X)$ to give rise to the commutative diagram
\[
  \begin{tikzpicture}[every node/.style={draw=none},
    ->,inner sep=1.1ex]
    \matrix [draw=none,row sep=4.5ex,column sep=4ex]
    {
      \node (11) {$H^1(\Omega_{\PP^n})$};
      && \node (13) {$H^1(\Omega_X)$}; \\
      & \node (22) {$H^1(\Omega_{\PP^n}|_X)$}; \\
    };
    \path
      (11) edge node[above=-.1ex,scale=.75] {$i^\ast$} (13) 
      (11) edge (22)
      (22) edge node[above left=-.3ex,scale=.75] {$r$} (13)
    ;
  \end{tikzpicture}
  \vspace{-2.5ex}
\]
and the injectivity of $r$ follows from the injectivity of $i^\ast$.
The latter is injective, since it sends the generator
$c_1(\Oo_{\PP^n}(1))\in H^1(\Omega_{\PP^n})$ to the class of a
hyperplane section of $X$.
\qed

In case $X$ is a surface, the isomorphism in Lemma
\ref{l:cotang-conorm} and the Serre duality yield the identification
\begin{equation*}
  \label{cot-norm}
  H^0(\Omega_X) \cong H^1(\Nn_X^\ast)\cong H^1(\Nn_X(K_X))^\ast.
\end{equation*}
Now we can see that the extension sequence \eqref{eq:koszulForS0}
tensored with $\Oo_X(K_X)$ is tied to the irregularity of $X$ via the
exact sequence
\begin{equation}
  \label{eq:koszulForS-H1}
  H^1(\Oo_X(2K_X+(5-m)H+Z_1))
  \stackrel{s_0}{\lra} H^1(\Nn_{X}(K_X))
  \stackrel{\wedge s_0}{\lra} H^1(\Jj_{Z_0}(K_X +mH-Z_1)). 
\end{equation}
Thus the problem of computing or bounding the irregularity comes down
to the understanding of the cohomology groups
$H^1(\Oo_X(2K_X+(5-m)H+Z_1))$ and $H^1(\Jj_{Z_0}(K_X+mH-Z_1)$).  This
in turn depends on controlling the subscheme $Z$ and the decomposition
\begin{equation*}
  \label{eq:ZDecomp}
  Z = Z_1 +Z_0.
\end{equation*}
The subscheme $Z$ is related to the singular locus
of the hypersurface $V_m$ which was used to define the section $s$.
Let us spell out this relationship.

The normal sequence of $X\subset\PP^4$ and the Euler sequence of
$\Theta_{\PP^4}$ give rise to the surjective morphism
\[
  H^0(\Oo_X(H))^\ast \otimes \Oo_X(H)\longrightarrow \Nn_X.
\]
This together with the section 
$s\in H^0(\Nn^\ast_X(mH))$---corresponding to $V_m$---yields
the commutative diagram
\vspace{6ex}
\begin{equation}
  \label{JZjac}
  \begin{tikzpicture}[overlay,every node/.style={draw=none},
    ->,inner sep=1.1ex]
    \matrix [draw=none,row sep=4.2ex,column sep=4ex]
    {
      \node (11) {$H^0(\Oo_X(H))^\ast\otimes\Oo_X(H)$}; \\
      \node (21) {$\Nn_X$};
      & \node (22) {$\Jj_Z(mH)$}; \\
    };
    \path
      (11) edge (21) 
      (11) edge (22)
      (21) edge node[above =-.2ex,scale=.75] {$s\wedge$} (22)
    ;
  \end{tikzpicture}
  \vspace{6ex}
\end{equation}
where the composite (slanted) arrow is given by the partial
differentiation of a homogeneous polynomial defining $V_m$.  Thus,
denoting by $\III_{V_m}$ the sheaf of ideals defined by the Jacobian
ideal of $V_m$, the diagram \eqref{JZjac} tells us that we have the
inclusion 
\begin{equation*}
  \label{JZjac1}
  \Jj_Z \supset \III_{V_m} \otimes \Oo_X
\end{equation*}
or, equivalently, that $Z$ is a subscheme of the scheme-theoretic
intersection of the singular locus $\Sing(V_m)$ with $X$.  So it is
reasonable to expect that one could control $Z$ and hence, the
cohomology groups in \eqref{eq:koszulForS-H1}, for small values of
$m$.  As an easy illustrative example of the above ideas, let us work
out the case $m=2$.

\section{Irregular surfaces on hypersurfaces of degree 2}
\label{s:V2}

The following theorem, no doubt, is well-known to experts, see \eg
\cite{Ro}, p.~152.

\begin{thm}
  \label{th:q=0V2}
If $X\subset\PP^4$ is a smooth surface lying on a quadric hypersurface
$V_2$, then the irregularity of $X$ vanishes.
\end{thm}

\proof
Suppose that $X\subset V_2$.  We may assume that $V_2$ is singular,
since otherwise $X$ is a complete intersection and we are done by
Lefschetz hyperplane theorem.  The singular locus $\Sing(V_2)$ is
either a single point, the vertex $p$ of $V_2$, or a line, $L$.

\paragraph{The case $\Sing(V_2)=\{p\}$}
The extension sequence \eqref{eq:koszulForS} takes the form
\[
  0 \lra \Oo_X(2K_X+3H) \stackrel{s}{\lra}
  \Nn_{X}(K_X) \lra \Jj_Z(K_X+2H) \lra 0
\]  
where $Z$ is either empty or $p$.  Thus the sequence of cohomology
groups \eqref{eq:koszulForS-H1} becomes
\[
  H^1(\Oo_X(2K_X+3H))
  \stackrel{s}{\lra}
  H^1(\Nn_{X}(K_X)) \lra  H^1(\Jj_Z(K_X+2H)).
\]  
Thus the irregularity is bounded as follows.
\[
\begin{gathered}
  q(X) = h^1(\Nn_{X}(K_X)) 
  \leq h^1(\Oo_X(2K_X+3H)) +h^1(\Jj_Z(K_X+2H)).
\end{gathered}
\]
From \cite{Re} it follows that the divisor $K_X+3H$ is very ample.
Hence, by the Kodaira vanishing and the Serre duality, we obtain
$h^1(\Oo_X(2K_X+3H))=0$.

We turn now toward $h^1(\Jj_Z(K_X+2H))$.  If $Z=0$, then
\[
  H^1(\Jj_Z(K_X+2H)) = H^1(\Oo_X(K_X+2H)) \stackrel{(SD)}{\cong}
  H^1(\Oo_X(-2H))^\ast = 0,
\] 
and hence, $q(X)=0$ in this case.  If $Z=p$, then the non-vanishing of
$H^1(\Jj_p(K_X+2H))$ means that $p$ is a base point of
$\Oo_X(K_X+2H)$, but this is ruled out by \cite[Theorem 1,(i)]{Re}.
Hence $H^1(\Jj_p(K_X+2H))=0$.  This completes the proof of the theorem
when $\Sing(V_2)$ is a point.

\paragraph{The case $\Sing(V_2) =L$} 
In this case $V_2$ is a singular rational scroll over a smooth conic
$C$ lying in a plane $\Pi$ complementary to the line $L$.  It is ruled
by the one parameter family of planes $\{P_t\}_{t\in C}$, where the
plane $P_t$ is the span of $t$ and $L$.  There are two cases to
consider according to whether or not the line $L$ is contained
in $X$.

\medskip

\noindent$\bullet$ 
If $L\subset X$, then the projection from $L$ defines the morphism
\[
  \phi : X \lra C \cong \PP^1
\] 
with the fibre $F_t$ over a point $t\in C$ being the component of the
intersection $P_t\cdot X$ complementary to $L$.  This implies that the
geometric ingredients $Z_1$ and $Z_0$ in the cohomology sequence
\eqref{eq:koszulForS-H1} are the line $L$ and the empty set,
respectively.  Hence that sequence takes the form
\[
  H^1(\Oo_X(2K_X+3H+L)) 
  \lra H^1(\Nn_{X}(K_X)) 
  \lra  H^1(\Oo_X(K_X+2H-L)). 
\]
We claim that the cohomology groups $H^1(\Oo_X(2K_X+3H+L))$ and
$H^1(\Oo_X(K_X+2H-L))$ vanish, since the divisors $2H-L$ and 
$K_X+3H +L$ are both ample.  Indeed, writing $2H-L=H+(H-L)$ and using
the fact that the linear system $|H-L|$ is base point free (the linear
system defines the morphism $\phi$), we deduce that $2H-L$ is very
ample.  The ampleness of $K_X+3H+L$ is checked easily using the very
ampleness of $K_X+3H$ and the Nakai-Moishezon criterion, see \cite{Ha}.

\medskip

\noindent$\bullet$ 
If $L\not\subset X$, then the sequence \eqref{eq:koszulForS}
associated to $V_2$ is exact and has the form
\begin{equation}
  \label{eq:seqV2line}
  0 \lra \Oo_X(K_X+3H) \lra \Nn_X \lra \Jj_Z(2H) \lra 0,
\end{equation}
where $Z$ is the $0$-dimensional subscheme of $X$, the
scheme-theoretic intersection of $X$ and the line $L$.  
This gives the relations
\begin{equation}
  \label{eq:ndpQHypersurface}
  d^2 -4d-\deg(Z) = 4(g_H-1)\,\,\text{and}\,\,\deg(Z)<d.
\end{equation}
We will now calculate $\deg(Z)$ using the geometry of the singular
scroll $V_2$ containing $X$.  Namely, we go back to the plane $\Pi$
containing the conic $C$ and complementary to $L$ and consider a
hyperplane $\Span(L,\Lambda)$ in $\PP^4$ spanned by $L$ and
a general line in $\Lambda\subset\Pi$.

The hyperplane $\Span(L,\Lambda)$ intersects $V_2$ in the union of two
planes $P_t$, where $t\in \Lambda \cdot C$.  Hence 
$H=\Span(L,\Lambda)\cdot X=2F$, where $F$ is the class of the divisor
$P_{t'}\cdot X$, for a closed point $t'\in C$.  From this it follows
that $d=4F^2$, and, by adjunction,
\[
  2(g_H-1) = 2F\cdot K_X + 4F^2
  = 2\deg(\omega_F)+2F^2
  = 2\deg(\omega_F)+\frac{d}{2},
\]
where $\omega_F$ is the dualizing sheaf of $F$.  On the other hand,
since $F$ is a plane divisor, its dualizing sheaf is subject to 
\[
  \deg(\omega_F) = d_F(d_F-3)
  = \frac{d}{2} \bigg(\frac{d}{2}-3 \bigg)
  = \frac{1}{4}\,d(d-6).
\]
Substituting into the previous identity, we have
\[
  2(g_H-1) = \frac{1}{2}\,d(d-6) + \frac{d}{2}
  = \frac{1}{2}\,d(d-5).
\]
Putting it together with the identity in \eqref{eq:ndpQHypersurface},
we obtain $\deg(Z)=d$ which contradicts the inequality in
\eqref{eq:ndpQHypersurface}.
\qed

We close this `warm up' section with a simple general observation
concerning the sheaf $\Jj_{Z_0}(mH-Z_1)$ in the sequence
\eqref{eq:koszulForS0}.  This observation will be important in
studying irregular surfaces lying on cubic hypersurfaces.

\begin{lem}
  \label{l:glGenAndConeCondition}
The sheaf $\Jj_{Z_0}((m-1)H-Z_1)$ (see the sequence
\eqref{eq:koszulForS0} for notation), is generated by global sections
outside the $0$-dimensional subscheme $Z_0$. Furthermore,
$h^0(\Jj_{Z_0}((m-1)H-Z_1))\geq 5$, unless the hypersurface $V_m$ is
a cone over a surface in $\PP^3$.
\end{lem}

\proof
From the sequence \eqref{eq:koszulForS0} tensored with $\Oo_X(-H)$ it
follows that $\Jj_{Z_0}((m-1)H-Z_1)$ is the quotient of $\Nn_X(-H)$
which is globally generated.  Hence the first statement of the
lemma.

For the second statement, we use the diagram \eqref{JZjac} tensored
with $\Oo_X(-H)$.  The slanted arrow is the morphism
\[
  H^0(\Oo_X(H))^\ast \otimes \Oo_X \lra \Jj_Z((m-1)H)
\]
which factors through $\Jj_{Z_0}((m-1)H-Z_1)$ and induces, at the
level of global sections, the linear map
\begin{equation}
  \label{partderiv}
  \partial(f)|_X:H^0(\Oo_X(H))^\ast \lra H^0(\Jj_{Z_0}((m-1)H-Z_1)),
\end{equation}
where $f$ is a homogeneous polynomial defining $V_m$.  The map
$\partial(f)|_X$ takes a vector $v\in H^0(\Oo_X(H))^\ast$ to the
partial derivative $\partial_v (f)|_X$ restricted to $X$ and then
divides it by a section defining $Z_1$.  From this it follows that the
kernel of $\partial(f)|_X$ consists of vectors 
$v\in H^0(\Oo_X(H))^\ast$ for which $\partial_v (f)$ is a homogeneous
polynomial of degree $m-1$ vanishing on $X$.  Since by definition, $m$
is the least degree of such polynomials, we deduce that 
$\partial_v(f)=0$ in $\Sym^{m-1}(H^0(\Oo_X(H)))$ or, equivalently,
\[
  f \in \Sym^{m-1}(v^{\perp}),
\]
where $v^{\perp}=\{ l\in H^0(\Oo_X(H)) \mid l(v)=0 \}$.  In
particular, for $v\neq 0$, the above means that the point 
$[v] \in \PP(H^0(\Oo_X(H))^\ast) =\PP^4$ is the vertex of the cone
over the surface in $\PP((v^{\perp})^\ast) \cong \PP^3$ defined by
$f$, viewed as a homogeneous polynomial on
$\PP((v^{\perp})^\ast)$. Hence, unless $V_m$ is a cone over a surface
in $\PP^3$, the operator $\partial(f)|_X$ in \eqref{partderiv} is
injective and, therefore,
\[
h^0(\Jj_{Z_0}((m-1)H-Z_1)) \geq h^0(\Oo_X(H)) \geq 5
\]
as asserted in the lemma. 
\qed

\section{Irregular surfaces on hypersurfaces of degree $3$}
\label{s:V3WithIsolatedSingularities}

The main result of this section is the following characterization of
irregular surfaces contained in a cubic hypersurface in $\PP^4$.
Results concerning surfaces contained in a cubic hypersurface in
$\PP^4$ can be found, \eg in \cite{Kol} and \cite{Ro}.

\begin{thm}
  \label{th:onV3}
If $X\subset\PP^4$ is a smooth irregular surface contained in a cubic
hypersurface $V_3$, then $X$ is an elliptic scroll of degree $5$.
Moreover, a general cubic hypersurface containing $X$ is a Segre cubic
and its ten singular points lie on $X$.
\end{thm}

As we have already explained, our main tool is the exact sequence
\eqref{eq:koszulForS0} which, for $m=3$, takes the form
\begin{equation}
  \label{eq:koszulForS0V3}
  0\lra \Oo_X(K_X +2H+Z_1) \lra \Nn_X \lra \Jj_{Z_0}(3H-Z_1) \lra0.
\end{equation}
The corresponding cohomological sequence \eqref{eq:koszulForS-H1}
controlling the irregularity of $X$ becomes
\begin{equation}
  \label{eq:cohkoszulForS0V3}
  H^1(\Oo_X(2K_X +2H+Z_1)) \lra H^1(\Nn_X(K_X)
  \lra H^1(\Jj_{Z_0}(K_X +3H -Z_1)).
\end{equation}
To analyse the cohomology group on the left in the above sequence we
will need the following.

\begin{pro}
  \label{p:h12KPlus2H}
If $H^1(\Oo_X(2K_X+2H))\neq0$, then $X$ is an elliptic scroll of
degree $5$.
\end{pro}

\proof
Assume that $H^1(\Oo_X(2K_X+2H))$ does not vanish.  Using the Serre
duality and \cite{Re}, we deduce that the line bundle 
$\Oo_X(K_X+2H)$ is base point free but not big.  Since $\Oo_X(K_X+2H)$
is not trivial%
\footnote{Otherwise $\Oo_X(K_X)=\Oo_X(-2H)$ and hence $q(X)=0$.}, we
deduce that the morphism defined by $\Oo_X(K_X+2H)$ maps $X$ onto a
curve.  In other words, there is a morphism $\phi:X\to B$ with
connected fibres onto a smooth curve $B$, and a base point free line
bundle $\Oo_B(D)$ on $B$ such that $\Oo_X(K_X+2H)=\phi^\ast\Oo_B(D)$.
This implies that
\begin{equation*}
  \label{eq:KHF}
  K_X + 2H = \deg(D) F,
\end{equation*}
where $F$ is the class of a general fibre of $\phi$. Taking the
intersection with $F$ on both sides of the above identity implies that
$H\cdot F=1$.  This means that the fibres of $\phi$ are lines and that
$X$ is a minimal ruled surface embedded into $\PP^4$ by $\Oo_X(H)$ as
a scroll with irregularity $q=q(X)=g(B)$, the genus of $B$. It is well
known that the only irregular scroll in $\PP^4$ is an elliptic scroll
of degree $5$, see Lemma~\ref{l:scrollsInP4}.
\qed

Next we turn to the group on the right of the sequence
\eqref{eq:cohkoszulForS0V3}.

\begin{lem}
  \label{l:rightsidenoncone}
If $V_3$ is not a cone and $Z_1$, the divisorial part in
\eqref{eq:cohkoszulForS0V3}, is non-zero, then
$H^1(\Jj_{Z_0}(K_X+3H-Z_1))=0$.
\end{lem}

\proof
Assume that $H^1(\Jj_{Z_0}(K_X+3H-Z_1))\neq 0$.  From the
identification
\[
  H^1(\Jj_{Z_0}(K_X+3H-Z_1))^\ast
  \cong \Ext^1(\Jj_{Z_0}(3H-Z_1),\Oo_X),
\]
the supposed nonvanishing is interpreted as a nontrivial extension
\begin{equation}
  \label{extseqV3}
  0 \lra  \Oo_X \lra \Ee \lra \Jj_{Z_0}(3H-Z_1) \lra 0.
\end{equation}
Tensoring it with $\Oo_X(-H)$, we obtain
\begin{equation}
  \label{identEandJZ0}
  H^0(\Ee(-H)) \cong H^0(\Jj_{Z_0}(2H-Z_1).
\end{equation}
This and Lemma \ref{l:glGenAndConeCondition} imply 
\begin{equation}
  \label{atleast5}
  h^0(\Ee(-H)) \geq 5. 
\end{equation}
This is going to play the role of a destabilizing condition for $\Ee$.
Namely, let $t$ be a nonzero global section of $\Ee(-H)$. It gives
rise to the exact sequence
\begin{equation}
  \label{KoszEs}
  0 \lra \Oo_X(A)\lra \Ee(-H) \lra \Jj_{Z'}(H-A-Z_1) \lra 0,
\end{equation}
where $A$ is the divisorial part of the zero locus of $t$ and $Z'$ is
the $0$-dimensional residual part of $(t=0)$. We proceed now to
analyse the above sequence according to the dimension of the linear
system $|A|$.

If $h^0(\Oo_X(A))\leq 2$, then \eqref{atleast5} implies
$h^0(\Jj_{Z'}(H-A-Z_1))\geq 3$.  Since $Z_1\neq 0$, it follows that
$Z_1$ is a line and $A=0$.  But then,
$h^0(\Jj_{Z'}(H-A-Z_1))=h^0(\Jj_{Z'}(H-Z_1))=3$ and the sequence
\eqref{KoszEs}, together with the estimate \eqref{atleast5}, imply
\[
  5 \leq h^0(\Ee(-H))
  \leq h^0(\Oo_X) + h^0(\Jj_{Z'}(H-Z_1)) = 1+3 = 4,
\]
an obvious contradiction.  Thus $h^0(\Oo_X(A)) \geq 3$.

Combining the sequence \eqref{KoszEs} with the defining sequence
\eqref{extseqV3} tensored with $\Oo_X(-H)$, we obtain the nonzero
morphism
\[
  \Oo_X(A) \lra \Jj_{Z_0}(2H-Z_1).
\]
Since this morphism can not be an isomorphism\footnote{Otherwise the
  extension sequence \eqref{extseqV3} is trivial.}, it is given by a
nonzero section $e\in H^0(\Jj_{Z_0}(2H-Z_1-A))$ vanishing on the
nonzero divisor
\[
  E = (e=0) \in |2H-Z_1-A|.
\]
In particular, the image of the morphism
$H^0(\Oo_X(A))\stackrel{e}{\to} H^0(\Jj_{Z_0}(2H-Z_1))$ consists of
global sections of $\Jj_{Z_0}(2H-Z_1)) $ vanishing on $E$.  In view of
the global generation of $\Jj_{Z_0}(2H-Z_1)) $ outside $Z_0$, see
Lemma \ref{l:glGenAndConeCondition}, it follows that
$e\,H^0(\Oo_X(A))$ is a proper subspace of
$H^0(\Jj_{Z_0}(2H-Z_1))$. Combining this and the isomorphism
\eqref{identEandJZ0}, we deduce that the image of $H^0(\Oo_X(A))$
under the monomorphism in \eqref{KoszEs} is a proper subspace of
$H^0(\Ee(-H))$. Hence $H^0(\Jj_{Z'}(H-Z_1-A))\neq 0$.

Let $D=H-Z_1-A$.  By the above, it is an effective divisor and
$h^0(\Oo_X(H-D)) =h^0(\Oo_X(Z_1 +A)) \geq 3$.  This tells us that
either $D=0$ or it is a line and the above inequality must be an
equality.  The first possibility is equivalent to $A=H-Z_1$ and
$\Jj_{Z'}(H-Z_1-A)=\Oo_X$. This implies $h^0(\Ee(-H))\leq 4$ which
contradicts \eqref{atleast5}.  If $D=H-Z_1-A$
is a line, then the estimate 
\[
  h^0(\Oo_X(H-Z_1-D) = h^0(\Oo_X(A) \geq 3 
\]
implies that $Z_1+D$ is a line.  But since $Z_1\neq 0$, this is
impossible.
\qed

The above lemma implies the following.

\begin{lem}
  \label{H1Z1noncone}
If $V_3$ is not a cone and the divisorial part $Z_1$ in
\eqref{eq:cohkoszulForS0V3} is nonzero, then $X$ is a regular surface.
\end{lem}

\proof
From Lemma \ref{l:rightsidenoncone} and \eqref{eq:cohkoszulForS0V3} it
follows that the irregularity of $X$ is controlled by the group
$H^1(\Oo_X(2K_X+2H+Z_1)$.  This group is related to the group
$H^1(\Oo_X(2K_X+2H)$ considered in Proposition \ref{p:h12KPlus2H} via
the obvious exact sequence
\begin{equation}
  \label{cohZ1}
  H^1(\Oo_X(2K_X+2H) \lra H^1(\Oo_X(2K_X+2H+Z_1)
  \lra  H^1(\Oo_{Z_1}(2K_X+2H+Z_1)),
\end{equation}
the understanding of which requires to have a good grasp of $Z_1$.
This is the case, since we know that $Z_1$ is contained in the
$1$-dimensional part of the singular locus $\Sing(V_3)$.  For a cubic
hypersurface, which is not a cone over a cubic surface, the
$1$-dimensional part of its singular locus is known to be
\begin{enumerate}[label={\rm (\roman*)}]
  \listspace
\item
  a line, 
\item
  a conic (possibly singular),
\item
  a rational normal curve of degree $4$ in $\PP^4$.
\end{enumerate}
Thus $Z_1$ is one of the above possibilities and we analyse each of 
them separately.

\paragraph{Case (i) --- $Z_1$ is a line}
In the exact sequence \eqref{cohZ1}, the degree of the line bundle
appearing in the group $H^1(\Oo_{Z_1}(2K_X+2H+Z_1))$ satisfies
\[
  Z_1\cdot(2K_X+2H+Z_1) = Z_1\cdot K_X + 2 + Z_1\cdot(K_X+Z_1)
  = Z_1\cdot K_X +2 -2 = Z_1\cdot K_X = -Z_1^2 -2.
\]
If $Z_1^2\leq-1$, it follows that $h^1(\Oo_{Z_1}(2K_X+2H+Z_1))=0$.
Combining this and \eqref{cohZ1},  we obtain
\[
  h^1(\Oo_X(2K_X+2H+Z_1)) \leq h^1(\Oo_X(2K_X+2H)
\]
Moreover, from Proposition \ref{p:h12KPlus2H}, we know
that $h^1(\Oo_X(2K_X+2H)\neq 0$ only if $X$ is an elliptic scroll.
But the only lines on such a surface $X$ are the rulings.  Therefore
$Z_1^2\leq-1$ implies that $h^1(\Oo_X(2K_X+2H)$ and hence
$h^1(\Oo_X(2K_X+2H+Z_1))$ are both equal to zero.  Thus we may
assume $Z_1^2=0$, since $X$ is obviously regular if $Z_1^2>0$.
Hence $X$, if irregular, is a scroll.  The proof of Proposition
\ref{p:h12KPlus2H} tells us that it must be an elliptic scroll of
degree $d=5$.  In addition, computing the second Chern number of
$\Nn_X$ from the exact sequence \eqref{eq:koszulForS0V3}, we obtain
\[
  25 = d^2
  = (K_X+2H +Z_1)(3H-Z_1) + \deg(Z_0)
  = 3d +3 +\deg(Z_0)
  = 18+\deg(Z_0)
\]
and hence 
\begin{equation}
  \label{Z0line}
  \deg(Z_0) = 7. 
\end{equation}
We decompose the $0$-dimensional subscheme $Z_0$, the zero locus of
the section $s_0$ in \eqref{eq:koszulForS0V3}, into two parts,
\begin{equation}
  \label{Z0linedecomp}
  Z_0 = Z^1_0 +Z'_0,
\end{equation}
where $Z^1_0$ is the part of $Z_0$ lying on the line $Z_1$ and the
residual subscheme $Z'_0$ is supported on singular points of $V_3$
belonging to $X\smallsetminus Z_1$.

To estimate the degree of $Z^1_0$ we use the exact sequence
\[
  0 \lra \Oo_X(-kZ_1) \lra \Jj_{Z^1_0} \lra \Oo_{kZ_1}(-Z^1_0) \lra 0,
\] 
where $k \geq 1$ is the smallest multiple of $Z_1$ containing the
subscheme $Z^1_0$.  Tensoring with $\Oo_X(2H-Z_1)$ and using the fact
that the linear system $|\Jj_{Z^1_0}(2H-Z_1)|$ has at most a
$0$-dimensional base locus, see Lemma \ref{l:glGenAndConeCondition},
we deduce that $h^0(\Oo_{kZ_1}(-Z^1_0)\otimes\Oo_X(2H-Z_1))>0$. From
this and $Z_1^2=0$, we obtain
\[
  0\leq \deg (2H|_{Z_1}-Z^1_0) =2-\deg(Z^1_0)
\]
or, equivalently, $\deg(Z^1_0)\leq 2$.  This together
with \eqref{Z0line} tells us that the part $Z'_0$ of the decomposition
\eqref{Z0linedecomp} has degree at least $5$.  In particular, $V_3$
must have singular points lying on $X$ and outside the line $Z_1$.
However, this is impossible.  Indeed, let $p$ be a singular point of
$V_3$ lying on $X\smallsetminus Z_1$.  Then the plane $P$ spanned by
$Z_1$ and $p$ must be contained in $V_3$.  Now consider the pencil of
hyperplanes $\{V(t) \mid t\in\PP^1\}$ containing $P$.  Each $V(t)$
intersects $V_3$ along $P\bigcup Q_t$, where $Q_t$ is a quadric
surface in $V(t)$ passing through $Z_1$ and the point $p$.  The
hyperplane sections $H_t=V(t)\cdot X$ are reducible and have the form
\[
  H_t = V(t)\cdot X = Z_1 + \Gamma_t,
\]
where $\Gamma_t$ is the divisor on $X$ residual to $Z_1$.
Furthermore, $\Gamma_t$ has degree $4$ and is contained in $Q_t$.  One
sees immediately that $\Gamma_t$ can not be irreducible, since
otherwise it is a smooth elliptic curve in $Q_t$ and its intersection
with the ruling $Z_1$ must be $2$, but on $X$ the curve $\Gamma_t $ is
a section of $X$ and $\Gamma_t\cdot Z_1=1$.  Hence
$\Gamma_t=\Gamma_0+L_t$, where $\Gamma_0$ is a smooth plane cubic, the
fixed part of the pencil $\{H_t\}$, and $L_t$ is a ruling of $X$. This
means that the rulings $L_t$ of $X$ are rationally equivalent and this
is clearly impossible.  This completes the proof of the lemma in the 
case (i).

\paragraph{Case (ii) --- $Z_1$ is a conic}
If $Z_1$ is a smooth conic, then the argument is as in the previous
case: we compute the degree 
\[
  Z_1\cdot(2K_X+2H+Z_1)
  = Z_1\cdot K_X+4 +Z_1\cdot(K_X+Z_1)
= Z_1\cdot K_X+4 -2 =Z_1\cdot K_X+2=-Z_1^2
\]
and obtain $h^1(\Oo_{Z_1}(2K_X+2H+Z_1))=0$, unless $X$ is regular.
Hence
\[
  h^1(\Oo_X(2K_X+2H+Z_1)) \leq h^1(\Oo_X(2K_X+2H),
\]
with the conclusion, in view of Proposition \ref{p:h12KPlus2H}, that
the irregular surface $X$ must be an elliptic scroll.  But such a
surface can not contain conics.

If $Z_1$ is singular, then $Z_1=L_1+L_2$, the sum of two lines $L_i$,
$i=1,2$.  The cohomology group $H^1(\Oo_{Z_1}(2K_X+2H+Z_1))$ fits into
the following exact sequence
\begin{equation}
  \label{conicred}
  H^1(\Oo_{L_1}(2K_X+2H+L_1))
  \to H^1(\Oo_{Z_1}(2K_X+2H+Z_1))
  \to H^1(\Oo_{L_2}(2K_X+2H+Z_1))
\end{equation}
and we continue as in the case (i).  Namely, we compute the degree of
the line bundle $\Oo_{L_2}(2K_X+2H+Z_1)$:
\[
  L_2\cdot (2K_X+2H+Z_1)
  = L_2\cdot K_X+L_2\cdot L_1
  = -L_2^2 + L_2\cdot L_1 - 2.
\] 
This implies that $h^1(\Oo_{L_2}(2K_X+2H+Z_1))=0$, unless $L_1=L_2$
and $L_1^2 = 0$. The latter condition means that $X$ is a scroll and,
by the proof of Proposition \ref{p:h12KPlus2H}, see also
Lemma~\ref{l:scrollsInP4}, it must be an elliptic scroll of degree
$d=5$.  Computing the second Chern number of $\Nn_X$ from
\eqref{eq:koszulForS0V3}, we obtain
\begin{multline*}
  25 = d^2 = (K_X+2H +Z_1)(3H-Z_1) \\
  = (K_X+2H +2L_1)(3H-2L_1) +\deg(Z_0)
  = 3d +6+\deg(Z_0) = 21+\deg(Z_0).
\end{multline*}
Hence $\deg(Z_0)=4$ and from here on we repeat the argument form the
proof of case (i).

Thus we must have $H^1(\Oo_{L_2}(2K_X+2H+Z_1)) =0$.  This and
\eqref{conicred} imply
\[
h^1(\Oo_{Z_1}(2K_X+2H+Z_1)) \leq h^1(\Oo_{L_1}(2K_X+2H+L_1))
\]
which puts us back into the situation of the case of the line in (i).

\paragraph{Case (iii) --- $Z_1$ is a rational normal curve of degree $4$}
This case is very special since a rational normal curve $C$ of degree
$4$ that lies in the singular locus of $V_3$ forces the cubic
hypersurface $V_3$ to be the secant variety of $C$.  Hence
$\Sing(V_3)=C$ and the Jacobian ideal ${\III}_{V_3}$ is equal to the
ideal sheaf $\Jj_{C/\PP^4}$ of $C$ in $\PP^4$.  From this, the exact
sequence \eqref{eq:koszulForS0V3}---for a smooth surface $X$ contained
in $V_3$ and passing through $C$---becomes
\begin{equation*}
  \label{RNC}
  0 \lra \Oo_X(2H+K_X+C) \lra \Nn_X \lra \Oo_X(3H-C) \lra 0, 
\end{equation*}
It is now easy to see that $X$ is regular.  Indeed, the cohomological
sequence controlling the irregularity $q$ of $X$,
\[
 H^1(\Oo_X(2H+2K_X+C)) \lra H^1(\Nn_X(K_X)) \lra H^1(\Oo_X(3H-C+K_X)),
\]
tells us that
\[
  q = h^1(\Nn_X(K_X)) \leq h^1(\Oo_X(2H+2K_X+C)), 
\]
since the line bundle $\Oo_X(3H-C)$ is ample.

The group  $H^1(\Oo_X(2H+2K_X+C))$ is computed by the exact sequence
\[
  H^1(\Oo_X(2H+2K_X))
  \lra H^1(\Oo_X(2H+2K_X+C))
  \lra H^1(\Oo_C(2H+2K_X+C)).
\]
The group on the right is zero unless $C^2 \geq 6$ and this of course
is impossible on an irregular surface.  If $h^1(\Oo_C(2H+2K_X+C))=0$,
then $h^1(\Oo_X(2H+2K_X+C))\leq h^1(\Oo_X(2H+2K_X))$ and Proposition
\ref{p:h12KPlus2H} tells us that if $h^1(\Oo_X(2H+2K_X))\neq0$, then
$X$ must be a an elliptic scroll.  But of course $C$ can not lie on
such a surface. Therefore, $h^1(\Oo_X(2H+2K_X))$ and hence
$h^1(\Oo_X(2H+2K_X+C))$ are both equal to zero.  This completes the
proof of the lemma.
\qed

We know now that if $V_3$ is not a cone and contains an irregular
surface $X$, then the exact sequence \eqref{eq:koszulForS0V3} must
have the form
\begin{equation}
  \label{eq:koszulForS0V30}
  0 \lra \Oo_X(K_X+2H) \lra \Nn_X \lra \Jj_Z(3H)\lra 0,
\end{equation}
where $Z$ is a $0$-dimensional scheme, the intersection of $X$ with
the singular locus of $V_3$.  Hence, the cohomology sequence
controlling the irregularity of $X$ becomes
\[
  H^1(\Oo_X(2K_X+2H) \lra H^1(\Nn_X(K_X)) \lra H^1(\Jj_Z(K_X+3H)).
\]

\begin{lem}
  \label{l:cohJZV3} 
If $X$ is irregular, then $H^1(\Jj_Z(K_X+3H))=0$.
\end{lem}

\proof
Assume $H^1(\Jj_Z(K_X+3H))\neq0$.  Using our approach, we interpret
this nonvanishing as a nontrivial extension sequence
\begin{equation}
  \label{extZV3}
  0 \lra \Oo_X \lra \Ee \lra \Jj_Z(3H) \lra 0.
\end{equation}
Tensoring it with $\Oo_X(-H)$, we obtain
\begin{equation*}
  \label{inentEJZV3}
  H^0(\Ee(-H)) \cong H^0(\Jj_Z(2H)).
\end{equation*}
Combining this and Lemma \ref{l:glGenAndConeCondition} gives
\begin{equation}
  \label{atleast5Z}
  h^0(\Ee(-H)) \geq 5.
\end{equation}
From here on we proceed as in the proof of Lemma
\ref{l:rightsidenoncone} to obtain a `destabilizing' sequence
\begin{equation}
  \label{EedestseqV3}
  0 \lra \Oo_X(A) \lra \Ee(-H) \lra \Jj_{Z'}(H-A) \lra 0.
\end{equation}

\smallskip
\paragraph{Case $h^0(\Oo_X(A)) \leq 2$} 
In this case \eqref{atleast5Z} implies $h^0(\Jj_{Z'}(H-A)) \geq 3$ and
we claim that $A=0$. Indeed, if $A\neq 0$, then the previous
inequality tells us that $A$ is a line, $\Jj_{Z'}(H-A)=\Oo_X(H-A)$,
and $h^0(\Oo_X(H-A))=3$.  From this, the exact sequence
\eqref{EedestseqV3}, and the inequality \eqref{atleast5Z}, we obtain
\[
  5 \leq h^0(\Oo_X(A)) + h^0(\Oo_X(H-A)) = h^0(\Oo_X(A)) +3.
\]
Equivalently, $h^0(\Oo_X(A))\geq 2$, \ie a line on $X$ moves in a
linear system.  But this is impossible on an irregular surface.

Once we know that $A=0$, the exact sequence \eqref{EedestseqV3}
becomes
\begin{equation}
  \label{EedestseqV3Aeq0}
  0 \lra \Oo_X \lra \Ee(-H) \lra \Jj_{Z'}(H) \lra 0.
\end{equation}
Together with \eqref{atleast5Z} gives
\[
  h^0(\Jj_{Z'}(H)) \geq 4.
\]
Hence either $Z'=p$, a single point, or $Z'=0$.  The first possibility
is interpreted as $p$ being a base point of $\Oo_X(K_X+H)$. Since
$d=H^2 \geq 5$, we can apply \cite[Theorem 1]{Re} to deduce that $X$
must be a scroll.  Furthermore, by the proof of Proposition
\ref{p:h12KPlus2H}, $X$ is an elliptic scroll of degree $d=5$.  This
together with \eqref{eq:koszulForS0V30} tells us that
\[
  25 = d^2 = c_2 (\Nn_X) = \deg(Z) +3d = \deg(Z) +15
\]
or, equivalently, that $\deg(Z)=10$.  However, from
\eqref{EedestseqV3Aeq0} and \eqref{extZV3} it follows that
\[
1 =\deg(Z') =c_2 (\Ee(-H)) =\deg(Z) -2H^2 =10-2\cdot5 =0
\]
which is absurd.

We turn now to the second possibility, $Z'=0$.  The exact sequence
\eqref{EedestseqV3Aeq0} takes the form
\[
  0 \lra \Oo_X \lra \Ee(-H) \lra \Oo_X(H) \lra 0
\]
and this implies that $\Ee(-H)\cong\Oo_X(H)\oplus\Oo_X$ since
$\Ext^1(\Oo_X(H),\Oo_X)\cong H^1(\Oo_X(-H))=0$.
Geometrically, this means that $Z$ is a complete intersection of two
effective divisors $H_1$ and $H_2$ in $|H|$ and $|2H|$ respectively.
In particular, we obtain
\begin{equation}
  \label{sectionsJZ2H}
  H^0(\Jj_Z(2H))
  = \{ h_1 h + \gamma h_2 \mid h\in H^0(\Oo_X(H)),\,\gamma\in\CC \},
\end{equation}
where $h_i$ is a section corresponding to $H_i$, for $i=1,2$.  Now,
from the proof of Lemma \ref{l:glGenAndConeCondition}, we recall that
$H^0(\Jj_Z(2H))$ contains the $5$-dimensional subspace
\[
  H^0(\frak{J}_{V_3}(2))|_X
  = \{ \partial_v (f)|_X \mid v\in H^0(\Oo_X(H))^\ast \}
\]
spanned by the restriction to $X$ of the partial derivatives of $f$, a
homogeneous polynomial defining $V_3$.  From this and the description
of $H^0(\Jj_Z(2H))$ in \eqref{sectionsJZ2H} it follows that the
intersection
\[
  H^0(\frak{J}_{V_3}(2))|_X \, \bigcap \,\, h_1 H^0(\Oo_X(H))
\]
is a $4$-dimensional vector space.  This implies that we can choose
homogeneous coordinate functions $X_i$, $i=0,\ldots,4$, in $\PP^4$ so
that the partial derivatives have the form
\[
  \frac{\partial f}{\partial{X_i}}
  = \htilde_1 T_j  \textq{for} 0\leq j\leq3
\]
where the $T_j$'s are some linear forms and on $H^0(\Oo_X(H))^\ast$
and $\htilde_1$ is the linear form corresponding to $h_1$.  It follows
that the remaining partial derivative 
$\frac{\partial f}{\partial{X_4}}$ defines a quadric hypersurface $Y$
which intersects the hyperplane $(\htilde_1=0)$ along a quadric
surface $Q$ contained in the singular locus of the cubic hypersurface
$V_3$.  This means that the secant variety of $Q$ is contained in
$V_3$.  Since the latter is irreducible, it follows that $Q$ is a
double plane.  But such a plane must intersect $X$ along a
$1$-dimensional subscheme which is contrary to our assumption.

\paragraph{Case $h^0(\Oo_X(A)) \geq 3$} 
We go back to the destabilizing sequence \eqref{EedestseqV3}.  We
argue as in the proof of Lemma \ref{l:rightsidenoncone} to deduce that
$H^0(\Oo_X(A))\hookrightarrow H^0(\Ee(-H))$ is a proper subspace.
Hence the divisor $D=H-A$ is effective.  This gives
\[
  3 \leq h^0(\Oo_X(A)) = h^0(\Oo_X(H-D)).
\]
Thus $D$ is either a line and the inequality above must be equality,
or $D=0$.  The first possibility implies
\[
  5 \leq h^0(\Ee(-H)) 
  \leq h^0(\Oo_X(A)) + h^0(\Jj_{Z'}(D))
  = 3 +h^0(\Jj_{Z'}(D))
\]
with the conclusion that the line $D$ moves in a linear system  on $X$
contradicting the hypothesis $X$ irregular.  The second possibility,
$D=0$, implies $A=H$.  In this case the destabilizing sequence
\eqref{EedestseqV3} takes the form
\[
  0 \lra \Oo_X(H) \lra \Ee(-H) \lra \Oo_X \lra 0.
\]
Since the epimorphism above induces a surjection on the level of
global sections, we deduce again that 
$\Ee(-H)\cong \Oo_X(H)\oplus\Oo_X$, a situation discarded in the
first part of the proof.
\qed

Next we investigate the possibility of $V_3$ being a cone over a cubic
surface in $\PP^3$.

\begin{lem}
  \label{l:V3cone} 
If the cubic hypersurface $V_3$ is a cone, then it contains no smooth
irregular surface.
\end{lem}

\proof 
We begin with the case of $V_3$ being a cone with vertex at a point
$x_0$ over a cubic surface $S$ which is not a cone.  We go back to the
sequence \eqref{eq:koszulForS0V3} and adapt our arguments thereafter
to the case at hand.

\paragraph{Claim}
If $Z_1 \neq 0$, then $H^1(\Jj_{Z_0}(K_X+3H-Z_1))=0$.

\medskip

We proceed as in the proof of Lemma \ref{l:rightsidenoncone} by
studying the extension sequence 
\[
  0 \lra \Oo_X \lra \Ee \lra \Jj_{Z_0}(3H-Z_1) \lra 0.
\]
As before we have the isomorphism 
$H^0(\Ee(-H))\cong H^0(\Jj_{Z_0}(2H-Z_1))$ and from the proof of
Lemma \ref{l:glGenAndConeCondition}, it follows that
\[
  h^0(\Jj_{Z_0}(2H-Z_1)) \geq 4. 
\]
This inequality and the isomorphism just above it give
\begin{equation*}
  \label{E-HsectionsCone}
  h^0(\Ee(-H)) \geq 4.
\end{equation*}
This leads again to the exact sequence \eqref{KoszEs}
\[
  0 \lra \Oo_X(A) \lra \Ee(-H) \lra \Jj_{Z'}(H-Z_1-A) \lra 0,
\]
where $A$ is an effective divisor and $H^0(\Jj_{Z'}(H-Z_1-A))\neq 0$.
Arguing as in the proof of Lemma \ref{l:rightsidenoncone} and using
the assumption that the base of the cone $V_3$ is not a cone, we
obtain
\begin{equation}
  \label{dimA}
  h^0(\Oo_X(A)) = h^0(\Jj_{Z'}(H-Z_1-A)) = 2. 
\end{equation}

Set $A_0$ (resp. $A'$) the fixed (resp. moving) part of $|A|$ and
consider $\Oo_X(H-Z_1-A_0)$.  It is easy to see that
$h^0(\Jj_{Z'}(H-Z_1-A_0))=3$.  From this and \eqref{dimA} it follows
that $A_0=0$, $Z_1$ is a line, $Z'=0$, and the linear map
\[
  H^0(\Oo_X(A))\otimes H^0(\Oo_X(H-Z_1-A)) \lra H^0(\Oo_X(H-Z_1))
\]
has a nontrivial kernel.  Hence there is a
basis $\{x,x'\}$ of $ H^0(\Oo_X(A))$ and nonzero elements 
$y,y'\in H^0(\Oo_X(H-Z_1-A))$ such that
\[
  xy -x'y' = 0 
\]
in $H^0(\Oo_X(H-Z_1))$.  Furthermore, $xy-x'y'$ can be viewed as the
restriction of an element from $\Sym^2 H^0(\Oo_X(H))$.  Since no
quadric hypersurface contains $X$ the above equality implies that
$y=\lambda x'$ and $y'=\lambda x$, for some nonzero $\lambda\in\CC$.
In particular, we obtain $\Oo_X(A)=\Oo_X(H-Z_1-A)$ or, equivalently,
\begin{equation}
  \label{2A}
  \Oo_X(H-Z_1) = \Oo_X(2A).
\end{equation}
Since $h^0(\Oo_X(2A)) =h^0(\Oo_X(H-Z_1)) =3$, we deduce 
$\Sym^2 H^0(\Oo_X(A))\cong H^0(\Oo_X(2A))$. Hence $|A|$ is base point
free and therefore, $A^2=0$.  This and \eqref{2A} imply that the
linear system $|2A|=|H-Z_1|$ is composed with a pencil, \ie the
morphism defined by $|H-Z_1|$ factors as follows,
\[
  \phi_{|H-Z_1|}: X \stackrel{|A|}{\lra} \PP^1
  \stackrel{|\Oo_{\PP^1}(2)|}{\lra} \PP^2.
\]
Equivalently, the projection from the line $Z_1$ maps $X$ onto a
conic.  But this means that $X$ is contained in a quadric
hypersurface, contrary to the hypothesis that $3$ is the least degree
of a hypersurface containing $X$.  This completes the proof of the
claim.

\medskip

Next we investigate the cohomology group $H^1(\Oo_X(2K_X +2H +Z_1))$
in the exact sequence \eqref{eq:cohkoszulForS0V3}.  For this we need
to control the divisor $Z_1$.  This divisor depends on the
singularities of the cubic surface $S$, the base of the cone
$V_3$.  The following two possibilities may occur:
\begin{enumerate}[label={\rm \arabic*)}]
  \listspace
\item
  $S$ has isolated singular points and then $Z_1$ is composed of one
  or two rulings of the cone $V_3$;
\item
  the singular locus $L$ of $S$ is a line and then $Z_1$ is a divisor
  contained in the plane $P$ spanned by $L$ and the vertex $x_0$ of
  the cone $V_3$. 
\end{enumerate}

The first possibility is dealt with in the same way as in the proof of
Lemma \ref{H1Z1noncone}, so we turn to the second possibility.  A
hyperplane $V$ passing through $P$ intersects $V_3$ along the
decomposable surface
\[
  V\cdot V_3 = 2P +P_V,
\]
where $P_V$ is the residual plane.  Hence the divisor $H_V =V\cdot X$
has the form 
\[
  H_V =2Z_1 + F_V,
\]
where $F_V$ is the divisor residual to $2Z_1$ and contained in the plane
$P_V$.  As in the proof of Lemma \ref{H1Z1noncone}, we relate 
$H^1(\Oo_X(2K_X+2H+Z_1))$ to $H^1(\Oo_X(2K_X+2H))$ via the sequence
\begin{equation}
  \label{seqH1Z1cone}
  H^1(\Oo_X(2K_X +2H)) \lra H^1(\Oo_X(2K_X +2H +Z_1))
  \lra H^1(\Oo_{Z_1} (2K_X +2H +Z_1))
\end{equation}
and investigate the cohomology group on the right of this sequence.

Since the divisor $Z_1=P\cdot X$ is the scheme-theoretic intersection
of a plane with the surface $X$, a result of Ellia and Folegatti,
\cite{EliFo}, implies that $Z_1$ is a reduced divisor.  Hence, for
an irreducible component $C$ of $Z_1$, we have
\[
  h^1(\Oo_{C}(2K_X +2H +Z_1))
  = h^1(\omega_C\otimes\Oo_{C}(K_X +2H+Z^C_1)
  = h^0(\Oo_{C}(-(K_X +2H)-Z^C_1),
\]
where $\omega_C$ is the dualizing sheaf of $C$, $Z^C_1$ the component
of $Z_1$ complementary to $C$ and the second equality is the Serre
duality on $C$.  Since $\Oo_X(K_X+2H)$ is base point free, see
\cite{Re}, the above identity tells us that
$h^1(\Oo_{C}(2K_X+2H+Z_1))=0$ unless $Z_1=C$ and 
$\Oo_{C}(K_X +2H)=\Oo_C$. But then $H\cdot C=1$ and $C^2=0$.  Hence
$X$ is a scroll and by Lemma \ref{l:scrollsInP4}, it is an
elliptic scroll of degree $d=5$.  Since this possibility was ruled out
in the proof of Lemma \ref{H1Z1noncone}, we obtain 
$H^1(\Oo_{Z_1} (2K_X +2H +Z_1))=0$.  This together with
\eqref{seqH1Z1cone} imply that the nonvanishing of 
$H^1(\Oo_X(2K_X +2H +Z_1))$ can only occur if 
$H^1(\Oo_X(2K_X +2H))\neq 0$.  This, in view of Proposition
\ref{p:h12KPlus2H}, implies that $X$ is an elliptic scroll of degree
$d=5$ and $Z_1$ is a ruling of $X$.  Thus we are back in the situation
ruled out in the proof of Lemma \ref{H1Z1noncone}.

\medskip

Our considerations are now reduced to the case when the divisorial
part $Z_1$ in the exact sequence \eqref{eq:koszulForS0V3} is
zero.  Hence that sequence has the form 
\[
0\lra \Oo_X(K_X +2H) \lra \Nn_X \lra \Jj_{Z_0} (3H)\lra 0
\]
and, as before, the irregularity of $X$ is controlled by the groups
$H^1(\Oo_X(2K_X +2H))$ and $H^1(\Jj_{Z_0}(3H+K_X))$.  

According to Proposition \ref{p:h12KPlus2H}, the nonvanishing of
$H^1(\Oo_X(2K_X +2H))$ may occur only if $X$ is an elliptic scroll of
degree $d=5$.  This degree is a numerical obstacle for $X$ to be
contained in a cubic cone.  Indeed, take a general hyperplane $V_1$ in
$\PP^4$ complementary to $x_0$, the vertex of the cone $V_3$, and
consider the projection of $X$ from $x_0$ into $V_1$. This gives rise
to a (rational) map
\begin{equation}
  \label{mapfx0}
  f_{x_0} : X --\rightarrow V_1 
\end{equation} 
onto a cubic surface $S:=V_1 \bigcap V_3$.  The degree of this map is
as follows: 
\begin{equation}
  \label{degfx0}
  \deg(f_{x_0}) =
  \begin{cases}
    \frac{d}{3} &\text{if }\, x_0\notin X,\\
    \frac{d-1}{3} &\text{if }\, x_0 \in X.
  \end{cases}
\end{equation}
For $d=5$ none of the above values is an integer.  Hence the vanishing
of $H^1(\Oo_X(2K_X +2H))$.

\paragraph{Claim}
$H^1(\Jj_{Z_0}(3H+K_X))=0$.

\medskip

The proof of this assertion follows the same pattern as the one of
Lemma \ref{l:cohJZV3}.  Namely, the nonvanishing of 
$H^1(\Jj_{Z_0}(3H+K_X))$ is interpreted as a nontrivial extension
\begin{equation}
  \label{extEconeV3}
  0 \lra \Oo_X \lra \Ee \lra \Jj_{Z_0}(3H) \lra 0. 
\end{equation}
This implies 
\begin{equation}
  \label{E-Hcone}
  h^0(\Ee(-H)) = h^0(\Jj_{Z_0}(2H)) \geq 4,
\end{equation}
where the inequality comes from the assumption that the surface $S$,
the base of the cone $V_3$, is not a cone; see the proof of Lemma
\ref{l:glGenAndConeCondition} for details.  In particular, as in the
proof of Lemma \ref{l:cohJZV3}, we have an exact sequence
\begin{equation}
  \label{destEconeV3}
  0 \lra \Oo_X(A) \lra \Ee(-H) \lra \Jj_{Z'}(H-A) \lra 0,
\end{equation}
where $A$ is an effective divisor and $h^0(\Jj_{Z'}(H-A))\geq 1$.
Putting this sequence together with the sequence \eqref{extEconeV3},
we obtain the relation
\begin{equation}
  \label{Z0Z'rel}
  \deg(Z_0)-2d = \deg(Z') + A\cdot (H-A).
\end{equation}
We wish to understand the geometric ingredients---the divisor $A$ and
the $0$-dimensional subscheme $Z'$---involved in this relation.

To begin with, we observe that $A$ is nonzero.  Indeed, if $A=0$, then
the relation \eqref{Z0Z'rel} reads
\begin{equation}
  \label{Z0Z'}
  \deg(Z')=\deg(Z_0)-2d.
\end{equation}
To make use of this relation, as well as of \eqref{Z0Z'rel} later on,
we give an upper bound for $\deg(Z_0)$.  Namely, from the fact that
the divisorial part $Z_1$ is zero, it follows that the cubic surface
$S$, the base of the cone $V_3$, has only isolated singularities and
that the map $f_{x_0}$ in \eqref{mapfx0} is finite outside $x_0$.
Furthermore, the subscheme $Z_0$ is the scheme-theoretic intersection
of $X$ with the rulings of the cone $V_3$ over the singular locus
$\Sing(S)$ of $S$.  Hence we get the upper bound
\[
  \deg(Z_0)
  \leq \deg(f_{x_0})\,\deg(\Sing(S))
  \leq \frac{d}{3}\,\deg(\Sing(S)),
\]
where the last inequality comes from \eqref{degfx0}.  In addition,
from the classification of normal cubic surfaces that are not cones,
see \cite[p. 448]{Dol}, it follows that $\deg(\Sing(S))\leq 6$.
Substituting this estimate into the above inequality, we obtain
\begin{equation}
  \label{degZ0coneV3}
  \deg(Z_0) \leq 2d.
\end{equation}
This and the relation \eqref{Z0Z'} imply $\deg(Z')=0$.  Hence the
exact sequence \eqref{destEconeV3} takes the form
\[
  0 \lra \Oo_X \lra \Ee(-H) \lra \Oo_{X}(H) \lra 0,
\]
situation we have already considered in the proof of Lemma
\ref{l:cohJZV3}.

We know now that the divisor $A$ in \eqref{destEconeV3} must be
nonzero.  Assume $h^0(\Oo_X(A))=1$.  From this and \eqref{E-Hcone} it
follows that
\[
  h^0(\Jj_{Z'}(H-A))\geq 3
\]
implying that $A$ is a line, the above inequality is an equality, and
$Z'=0$.  This and \eqref{Z0Z'rel} imply
\[
  \deg(Z_0)-2d = A\cdot (H-A) =1-A^2,
\]
which, together with the upper bound \eqref{degZ0coneV3}, implies
$A^2\geq 1$.  Hence $X$ is regular (one sees easily that $X$ is a
rational surface).  Thus we may assume $h^0(\Oo_X(A))\geq 2$. In fact
the equality must hold, since, if $h^0(\Oo_X(A))\geq 3$, then the
divisor $L=H-A$ is a line and the relation \eqref{Z0Z'rel} reads
\[
  \deg(Z_0)-2d = \deg(Z')+(H-L)\cdot L
  = \deg(Z')+1-L^2.
\]
This and the upper bound \eqref{degZ0coneV3} imply $L^2\geq 1$,
 the situation we have already discarded.

From $h^0(\Oo_X(A))=2$ and the inequality \eqref{E-Hcone} we deduce
\[
  h^0(\Jj_{Z'}(H-A)) \geq 2.
\]
Furthermore, if the inequality is strict, we arrive again at the
situation ruled out previously: $A$ is a line moving in a linear
system.  Thus
\[
  h^0(\Oo_X(A)) = h^0(\Jj_{Z'}(H-A)) = 2,
\]
a situation we have already encountered in the first part of the
proof.  Arguing as there, we show that the linear system $|A|$ has at
most a $0$-dimensional base locus.  Hence $A\cdot (H-A) \geq 0$. This,
the relation \eqref{Z0Z'rel}, and the upper bound \eqref{degZ0coneV3}
imply that all inequalities involved must be equalities.  In
particular, 
\[
  A\cdot (H-A)=0.
\]
Since both $A$ and $(H-A)$ are effective, nonzero divisors that add up
to $H$, a very ample divisor, the above identity is impossible.  This
completes the proof of the lemma in the case the surface $S$, the base
of the cone $V_3$, is not a cone.

\medskip

We now turn to the case when $S$ is a cone with vertex $x_S$ over a
plane cubic curve $C$.  If $C$ is smooth, then $V_3$ is a singular
scroll ruled by the planes $P_t=\Span(t\bigcup L)$, for $t\in C$, with
the singular locus $\Sing(V_3)$ being the line 
$L=\Span(x_0\bigcup x_S)$.  This is analogous to the situation we
arrived at in the proof of Theorem \ref{th:q=0V2}.  Hence the
sequence \eqref{eq:koszulForS0V3} takes the form%
\footnote{We consider here only the case $L\subset X$, since the other
  case, $L\not\subset X$, is treated in exactly the same way as in the
  proof of Theorem~\ref{th:q=0V2}.} 
\[
  0 \lra \Oo_X(K_X +2H +L) \lra \Nn_X \lra \Oo_X(3H -L) \lra 0.
\] 
It follows that the irregularity of $X$ is controlled by the
cohomology group $H^1(\Oo_X(2K_X +2H +L))$.  Its nonvanishing, as we
have seen on several occasions, may occur only if $X$ is an elliptic
scroll of degree $d=5$ and $L$ is a ruling of the scroll.  But
computing the second Chern classes from the above sequence, we obtain
\[
  25 = d^2 = (K_X +2H +L)\cdot(3H -L) = 18,
\]
an obvious contradiction.

If $C$ is singular, let $c_0$ be its (unique) singular point.  The
plane $P_{c_0}=\Span(c_0 \bigcup L)$ is the singular plane of $V_3$.
Furthermore, the planes $P_t=\Span(t\bigcup L)$ give rise to a
rational family of curves
\[
  \{F_t = P_t \cdot X \mid t\in C\}.
\] 
Let $F$ be the divisor class of this family.  Then by taking a general
hyperplane in $\PP^4$ passing through the line $L$, we have
\[
  \Oo_X(H) = \Oo_X(3F).
\]
From $h^0(\Oo_X(F))=h^0(\Oo_X(H-2F))$ it follows that
$h^0(\Oo_X(F))=2$.  Since the linear map
\[
  \Sym^2(H^0(\Oo_X(F))) \lra H^0(\Oo_X(2F))
\]
is injective, we deduce
\[
  3 \leq h^0(\Oo_X(2F)) = h^0(\Oo_X(H-F))
\]
and hence $F$ is a line.  Since $F$ moves in a pencil the surface $X$
must be rational.  This completes the proof of the lemma.
\qed

We know now that an irregular $X$ contained in a cubic hypersurface
must be an elliptic scroll of degree $5$. Furthermore, we have the
following.

\begin{lem}
  \label{l:singV3forscroll} 
Let $X$ be an elliptic scroll of degree $d=5$ in $\PP^4$ and let $V_3$
be a cubic hypersurface containing it. Then the scheme
$Z=X\bigcap\Sing(V_3)$ is $0$-dimensional of degree $10$.
\end{lem}

\proof 
From the proof of Lemma \ref{l:V3cone} we know that $V_3$ is not a
cone.  Combining this with Lemma \ref{H1Z1noncone}, we deduce that the
scheme $Z=X\bigcap\Sing(V_3)$ is $0$-dimensional. Hence the sequence
\eqref{eq:koszulForS0V3} has the form
\begin{equation*}
  \label{seqZscroll}
  0\lra \Oo_X(K_X +2H) \lra \Nn_X \lra \Jj_Z(3H) \lra 0.
\end{equation*}
Computing the second Chern classes from this sequence, we obtain
\[
\deg(Z) =\deg \left(X\bigcap\Sing(V_3)\right) =10.
\]
\qed

To make the paper self-contained, we include the following result
which is a well-known part of the classification of surfaces in
$\PP^4$, see \cite{Lant}.

\begin{lem}
  \label{l:scrollsInP4}
The only irrational scrolls in $\PP^4$ are elliptic scrolls of degree
$5$.
\end{lem}

\proof
Let $X$ be a $\PP^1$-bundle over a smooth connected curve $B$ of genus
$q \geq 1$ and let $\Oo_X(H)$ be a very ample line bundle on $X$
defining an embedding of $X$ into $\PP^4$ as a scroll.  In particular,
a smooth divisor in the linear system $|H|$ is isomorphic to $B$.  This
together with the adjunction formula implies
\[
  H\cdot K_X+ H^2 = H\cdot K_X+ d = 2q-2.
\]
Substituting this and the Chern numbers $\chi(\Oo_X)=1-q$,
$K_X^2=8(1-q)$ into the double point formula, we obtain
\begin{equation}
  \label{eq:dpfForScrolls}
  d^2-5d = 6(q-1).
\end{equation}
The main point of the argument is to compare the genus $q$ of $B$ in
the above formula with the Castelnuovo upper bound on the genus of a
smooth curve in a projective space.  Namely, the scroll 
$X\subset \PP^4$ is interpreted as an embedding
\[
  \phi: B \lra \Gr(1,\PP^4)
\]
into the Grassmannian $\Gr(1,\PP^4)$ of lines in $\PP^4$.  Setting
$\Gg$ to be the pullback under $\phi$ of the universal subbundle of
the Grassmannian, we identify $X$ with the projectivization
$\PP(\Gg)$. Then $\Oo_X(H)$, the line bundle embedding $X$ into
$\PP^4$, is such that the direct image 
$\pi_\ast(\Oo_X(H))\cong\Gg^\ast$.  In particular,
\begin{equation*}
  \label{Pluckd}
  \deg(c_1(\Gg^\ast)) = d.
\end{equation*}
Composing $\phi$ with the Pl\"ucker embedding of $\Gr(1,\PP^4)$ gives
the embedding
\begin{equation*}
  \label{BinP9}
  \psi: B \hookrightarrow \PP^9
\end{equation*}
realized by the subsystem of $|\Bigwedge^2\Gg^\ast|$ corresponding
to the image of the obvious homomorphism 
\[
  \rho: \Bigwedge^2 H^0(\Gg^\ast) \lra H^0(\Bigwedge^2\Gg^\ast).
\]
We claim that $\ker(\rho)$ has dimension at least $5$.  Indeed, assume
\[
  \dim(\ker(\rho)) \leq 4.
\]
Then the image of $\rho$ has dimension $N+1\geq 6$ and $\psi$ embeds
$B$ into $\PP^N$ as a nondegenerate curve of degree $d$.  The
Castelnuovo bound on the genus $q$ in $\PP^N$ gives
\[
  q-1 \leq sd -\frac{s(s+1)(N-1)}{2} -N-1,
\] 
where $d=s(N-1)+r$ for an integer $0\leq r<N-1$.  Rewriting the
expression on the right as a function of $d$, $r$, and $N$, we obtain
\[
  q-1
  \leq \frac{1}{2(N-1)}\,(d^2 -r^2) - \frac{1}{2}\,d
   + \frac{r}{2} - (N+1),
\]
In our situation $N\in\{5,\ldots,9\}$, hence we can consider the
larger bound
\[
  q-1 \leq \frac{1}{8}\,d^2 - \frac{1}{2}\,d
  - \frac{r^2}{16}+\frac{r}{2}-6.
\]
Putting it together with \eqref{eq:dpfForScrolls} gives
\[
  d^2 -5d
  \leq \frac{3}{4}\,d^2-3d-\frac{3r^2}{8}+3r-36.
\]
This can be rewritten in the form
\[
  \bigg(\frac{1}{2}\,d-2 \bigg)^2
  = \frac{1}{4}\,d^2 -2d + 4
  \leq -\frac{3r^2}{8}+3r-32 \leq -26
\]
which is false.

We know now that the kernel of $\rho$ is at least of dimension $5$.
Geometrically this means that the projectivized subspace
$\PP(\ker(\rho))$ intersects the Grassmannian variety
$\Gr(2,H^0(\Gg))\cong\Gr(1,\PP^4)$ of decomposable tensors in
$\PP(\Bigwedge^2 H^0(\Gg^\ast)) \cong \PP^9$ along a subscheme of
dimension at least $1$.  Each decomposable tensor $g\wedge g'$ in this
intersection, viewed as a section of $\Bigwedge^2\Gg^\ast$, is zero.
Equivalently, the two sections $g$ and $g'$ correspond, under the
isomorphism $H^0(\Gg^\ast)\cong H^0(\Oo_X(H))$, to two hyperplanes
$H_g$ and $H_{g'}$ in $\PP^4$ such that the plane 
$P=H_g\bigcap H_{g'}$ intersects the scroll $X$ along a curve
$\Gamma_{g,g'}$ which is a section of the structure projection
$\pi:X\to B$.  The above shows that there is a family
$\{\Gamma_t\}_{t\in T}$ of such sections parametrized by an
irreducible curve $T$.

Observe that every two planes $P_t$ and $P_{t'}$, for $t\neq t'$,
intersect at a single point.  This implies that
$\Gamma_t\cdot\Gamma_{t'}\leq 1$.  We easily check that
$\Gamma_t\cdot\Gamma_{t'}= 1$.  In particular, setting $\Gamma$ to be
the class of $\{\Gamma_t \}_{t\in T}$ in the N\'eron-Severi group of
$X$, we obtain
\[
  \Gamma^2=1.
\]
We write $H=\Gamma+aF$, where $F$ is the class of a ruling of $X$.
Then $d=H^2= 2a+1$ and the degree of $\Gamma$ is 
$d_{\Gamma}=H\cdot \Gamma=d-a =2a+1-a=a+1$.  Since $\Gamma$ is a plane
curve, we have
\[
  2(q-1)=(a+1)(a-2).
\]
This, the identity $d=2a+1$ and \eqref{eq:dpfForScrolls} imply
\[
  3(a+1)(a-2) = 6(q-1) = d^2-5d = (2a+1)^2 -5(2a+1) = (2a+1)(2a-4)
\]
or, equivalently,
\[
  (a-2)(a-1) = 0.
\]
This leads to two solutions $d=5$ and $3$, which are, respectively, an
elliptic scroll of degree $5$ and a rational scroll of degree $3$.
\qed

\section{On the elliptic scroll of degree $5$ and the Segre cubic}
\label{s:ellscrollSegre}

In the previous section we have characterized an elliptic scroll $X$
of degree $5$ in $\PP^4$ as being the only irregular surface lying on
a cubic hypersurface. Such scrolls are notorious and they have been
subject to extensive study; see, \eg \cite{Hu, ADHPR,ADHPR2} and the 
references therein.  The main objective of this section is to 
(re)establish a relationship between two entities related to the 
embedding of $X$ in $\PP^4$:
\begin{itemize}
  \listspace
\item 
  the space $I_X(3)$ of cubic hypersurfaces in $\PP^4$ containing $X$
\item
  the space of global sections $H^0(\Nn^\ast_X(3H))$ of the conormal
  bundle $\Nn^\ast_X$ of $X$ in $\PP^4$ twisted by $\Oo_X(3H)$, where
  $\Oo_X(H)$ is a line bundle realizing the embedding of $X$ in
  $\PP^4$.
\end{itemize}

Using the above notation, we formulate the main result of this section.

\begin{thm}
  \label{th:cubicsSections}
Let $X$ be an elliptic scroll of degree $5$.
  \begin{enumerate}[label={\rm \arabic*)}]
    \listspace
\item
  The sheaf $\Nn^\ast_X(3H)$ is a rank $2$ vector bundle generated by
  its global sections, with Chern invariants
  $c_1(\Nn^\ast_X(3H))=H-K_X$ and $c_2(\Nn^\ast_X(3H))=10$. 
\item
  There is a natural isomorphism
  $H^0(\Nn^\ast_X(3H))\cong I_X(3)\cong \CC^5$.
\item
  Every nonzero global section $s$ of $\Nn^\ast_X(3H)$ has
  $0$-dimensional zero locus of degree $10$.  Under the above
  correspondence, the scheme of zeros $Z_s=(s=0)$ is the
  scheme-theoretic intersection of $X$ with $\Sing(V_3(s))$, the
  singular locus of the cubic hypersurface $V_3(s)\in|I_X(3)|$
  corresponding to $s$.  In particular, every global section $s$ with
  $Z_s=(s=0)$ consisting of ten distinct points, corresponds to a
  Segre cubic $V_3(s)$, whose set of nodes $\Sing(V_3(s))=Z_s$.  
\end{enumerate}       
\end{thm}

\proof
The assertion about the Chern invariants is obvious.  Of course the
relation between global sections of $\Nn^\ast_X(3H)$ and cubic
hypersurfaces through $X$ has been at the origin of our considerations
in Section \ref{s:V3WithIsolatedSingularities} and stems from the
identification $\Nn^\ast_X=\Jj_X/\Jj^2_X$, where $\Jj_X$ is the ideal
sheaf of $X$ in $\PP^4$.  From the exact sequence
\[
  0 \lra \Jj^2_X(3) \lra \Jj_X(3) \lra \Nn^\ast_X(3H) \lra 0
\]
follows the inclusion\footnote{The inclusion comes from
$H^0(\Jj^2_X(3))=0$, since we know that $X$ is not contained in a
quadric hypersurface, see Theorem \ref{th:q=0V2}.}
\begin{equation}
  \label{I3inject}
  0 \lra I_X(3) = H^0(\Jj_X(3)) \lra H^0(\Nn^\ast_X(3H)).
\end{equation}
From the exact sequence
$0 \to \Jj_X \to \Oo_{\PP^4} \to \Oo_X \to 0$
tensored with  $\Oo_{\PP^4}(3)$ we obtain the estimate
\begin{equation}
  \label{I3lb}
  h^0(\Jj_X(3)) \geq h^0(\Oo_{\PP^4}(3))-h^0(\Oo_X(3H))
  = \binom{7}{3} -\frac{1}{2}\,(9H^2 -3H\cdot K_X)
  = 35-30 =5,
\end{equation}
where we used the easily verified property
\[
  H^i(\Oo_X(kH)) = 0  \textq{for all} k>0 \textq{and} i=1,2.
\]
The above estimate is actually an equality.\footnote{One easily
verifies that $X\subset\PP^4$ is a projectively normal
embedding.  Since we do not use this aspect anywhere, the above
mentioned equality is established differently in the proof of Lemma
\ref{l:aboutConormal3H}.} This as well as the isomorphism in
\eqref{I3inject}, and hence the assertion 2) of the theorem, follow
from the next lemma.

\begin{lem}
  \label{l:aboutConormal3H}
The vector bundle $\Nn_X^\ast(3H)$ is globally generated and
$h^0(\Nn_X^\ast(3H))=5$. 
\end{lem}

\proof
Let $\pi:X\to E$ be the structure morphism, \ie $E$ is an elliptic
curve and $\pi$ is a $\PP^1$-fibration over $E$.  The line bundle
$\Oo_X(H)$ defining the embedding of $X$ into $\PP^4$ as a scroll has
degree $1$ on the fibres of $\pi$, \ie 
\[
  \Oo_X(H)\otimes\Oo_F = \Oo_F(1) \cong \Oo_{\PP^1}(1)
\]
on every fibre $F$ of $\pi$. 

We wish to understand the restriction of $\Nn_X^\ast(3H)$ to a fibre
$F$.  Since $\det(\Nn_X)=\Oo_X(5H+K_X)$, we have
\[
  \det(\Nn_X)\otimes\Oo_F = \Oo_F(5H+K_X) = \Oo_F(3).
\]
From this it follows that
\begin{equation}
  \label{NonF}
  \Nn_X^\ast(3H)\otimes\Oo_F
  = \Nn_X^\ast\otimes\Oo_F(3)
  = \Nn_X^\ast\otimes\det(\Nn_X)\otimes\Oo_F
  \cong \Nn_X\otimes\Oo_F
  \cong \Oo_F(1)\oplus\Oo_F(2),
\end{equation}
where the last isomorphism follows from the global generation of
$\Nn_X(-H)\otimes\Oo_F=\Nn_X\otimes\Oo_F(-1)$.  In particular, the
restriction $\Nn_X^\ast(3H)\otimes\Oo_F$ is globally generated on
every fibre $F$ of $\pi$.  So to obtain the global generation of
$\Nn_X^\ast(3H)$, it is sufficient to show the surjectivity
\[
  H^0(\Nn_X^\ast(3H)) 
  \lra H^0(\Nn_X^\ast(3H)\otimes\Oo_F) 
  = H^0(\Oo_F(1)\oplus\Oo_F(2))
\]
for every fibre $F$ of $\pi$.  For this we use the inclusion
\eqref{I3inject} to obtain the composition
\vspace{7ex}
\begin{equation}
  \label{diagtofibre}
  \begin{tikzpicture}[overlay,every node/.style={draw=none},
    ->,inner sep=1.1ex]
    \matrix [draw=none,row sep=3.5ex,column sep=4.2ex]
    {
      && \node (02) {$0$}; \\
      && \node (12) {$I_X(3)$}; \\
      \node (20) {$0$};
      & \node (21) {$H^0(\Nn_X^\ast(3H -F))$};
      & \node (22) {$H^0(\Nn_X^\ast(3H))$};
      & \node (23) {$H^0(\Nn_X^\ast(3H)\otimes\Oo_F)$};
      & \node (24) {$0$}; \\
    };
    \path
      (02) edge (12)
      (12) edge (22)
      (12) edge (23)
      (20) edge (21)
      (21) edge (22)
      (22) edge (23)
      (23) edge (24)
    ;
  \end{tikzpicture}
  \vspace{8.5ex}
\end{equation}
In addition, from the proof of Theorem \ref{th:onV3}, we know that the
sections of $\Nn_X^\ast(3H)$ coming from $I_X(3)$ have $0$-dimensional
schemes of zeros.  Hence, the image of $I_X(3)$ in
$H^0(\Nn_X^\ast(3H))$ is complementary to the kernel of the
restriction homomorphism
\[
  H^0(\Nn_X^\ast(3H)) \lra H^0(\Nn_X^\ast(3H)\otimes\Oo_F) 
\]
in \eqref{diagtofibre}.  From this it follows that the slanted arrow
in \eqref{diagtofibre} is injective.  Since the target of that arrow
is a space of dimension $5$, see \eqref{NonF}, we deduce the
inequality $\dim (I_X(3))\leq 5$.  This and the estimate \eqref{I3lb}
imply
\begin{equation}
  \label{dimI3}
  \dim(I_X(3)) = 5.
\end{equation} 
Therefore, the slanted arrow in \eqref{diagtofibre} is an isomorphism,
hence the global generation of $\Nn_X^\ast(3H)$.

We now turn to the assertion $h^0(\Nn_X^\ast(3H))=5$.  From
\eqref{diagtofibre}, we already know that $h^0(\Nn_X^\ast(3H))\geq 5$.
Let us assume that the inequality is strict.  From the Koszul sequence
\begin{equation}
  \label{KoszulZsscroll}
  0 \lra \Oo_X
  \stackrel{s}{\lra} \Nn_X^\ast(3H)
  \stackrel{\wedge s}{\lra} \Jj_{Z_s}(H-K_X) \lra 0
\end{equation}
of a general global section $s$ of $\Nn_X^\ast(3H)$, we have
\begin{equation}
  \label{eq:intermediateBound}
  h^0(\Jj_{Z_s}(H-K_X)) \geq h^0(\Nn_X^\ast(3H))-1 \geq 5.
\end{equation}
Furthermore, considering another general global section of
$\Nn_X^\ast(3H)$, we obtain the smooth curve $\Gamma=(\gamma=0)$,
where $\gamma$ is the section of $\det(\Nn_X^\ast(3H))=\Oo_X(H-K_X)$
corresponding to $s\wedge s'$ under the natural homomorphism
\[
  \Bigwedge^2 H^0(\Nn_X^\ast(3H))
  \lra H^0(\det(\Nn_X^\ast(3H))) = H^0(\Oo_X(H-K_X)).
\]  
This gives rise to the exact sequence
\begin{equation*}
  \label{ZsonGamma}
  0 \lra \Oo_X \stackrel{\gamma}{\lra} \Jj_{Z_s}(H-K_X)
  \lra \Oo_\Gamma((H-K_X)|_\Gamma-{Z_s}) \lra 0.
\end{equation*}
From this and \eqref{eq:intermediateBound} we obtain
\begin{equation}
  \label{h0atleast4}
  h^0(\Oo_\Gamma((H-K_X)|_\Gamma-{Z_s}))
  \geq h^0(\Jj_{Z_s}(H-K_X))-1 \geq 5-1 = 4.
\end{equation}
On the other hand, the degree of $\Oo_\Gamma((H-K_X)|_\Gamma-{Z_s})$
is
\[
  (H-K_X)\cdot\Gamma-\deg({Z_s}) 
  = (H-K_X)^2 - c_2(\Nn_X^\ast(3H))
  = 15-10 = 5,
\]
while the genus of $\Gamma$, by the adjunction formula, is $6$. This
implies that $\Oo_\Gamma((H-K_X)|_\Gamma-{Z_s})$ is special and by the
Clifford inequality
\[
  h^0(\Oo_\Gamma((H-K_X)|_\Gamma-Z_s))
  \leq \frac{\deg((H-K_X)|_\Gamma-Z_s)}{2}+1
  = \frac{5}{2} +1.
\] 
Combining this inequality with \eqref{h0atleast4}, we obtain
\[
  4 \leq h^0(\Oo_\Gamma((H-K_X)|_\Gamma-Z)) 
  \leq \frac{5}{2} +1,
\]
an obvious contradiction.
\qed

\likeproof[End of the proof of Theorem \ref{th:cubicsSections}]
Since the spaces $H^0(\Nn^\ast_X(3H))$ and $I_X(3)$ are both
$5$-dimensional, see \eqref{dimI3} for the latter, an immediate
corollary of Lemma \ref{l:aboutConormal3H} is that the inclusion
\eqref{I3inject} is an isomorphism.

Once the isomorphism in Theorem \ref{th:cubicsSections}, 2) is
established, we also deduce that every nonzero global section of
$\Nn^\ast_X(3H)$ has $0$-dimensional scheme of zeros, since this is
now equivalent to the property of cubics in $I_X(3)$ to produce
sections of $\Nn^\ast_X(3H)$ with $0$-dimensional zero locus, see
Lemma \ref{l:singV3forscroll}.  This proves the first assertion of
Theorem \ref{th:cubicsSections}, 3).

The degree of the scheme of zeros of a nonzero section of
$\Nn^\ast_X(3H)$ is the value of the second Chern class (identified
with its degree) $c_2 (\Nn^\ast_X(3H))$ and this value is $10$ by the
part 1) of the theorem.

Next we turn to the assertion that $Z_s$, the zero scheme of a nonzero
global section $s$ of $\Nn^\ast_X(3H)$, is the scheme theoretic
intersection of $X$ with the singular locus of the cubic hypersurface
$V_3 (s)$ corresponding to $s$ under the isomorphism in part 2) of the
theorem.  From the general discussion about the relation of the ideal
sheaf $\Jj_{Z_s}$ of $Z_s$ and the restriction of the Jacobian ideal
$\III_{V_3 (s)}$ to $X$ we know that
\[
  \Jj_{Z_s} \supset \III_{V_3 (s)}\otimes\Oo_X,
\]
\ie $Z_s$ is contained in the scheme theoretic intersection of $X$
with the singular locus of the cubic hypersurface $V_3 (s)$.  To see
the equality, it is enough to show that the generators of 
$\III_{V_3 (s)} $, the partial derivatives of a homogeneous polynomial
defining ${V_3 (s)}$, restricted to $X$, generate the sheaf
$\Jj_{Z_s}(2H)$.  This follows immediately from the epimorphism
\[
  \Nn_X(-H) \lra \Jj_{Z_s}(2H)
\]
coming from the Koszul sequence \eqref{KoszulZsscroll} tensored with
$\Oo_X(H+K_X)$ and the surjective morphism
\[
  H^0(\Oo_X(H))^\ast\otimes\Oo_X \lra \Nn_X(-H).
\]
The resulting composition
\[
  H^0(\Oo_X(H))^\ast\otimes\Oo_X \lra \Jj_{Z_s}(2H)
\]
is surjective and is described explicitly by the partial derivatives
of a homogeneous polynomial defining ${V_3(s)}$, see the proof of
Lemma \ref{l:glGenAndConeCondition} for details.  Hence the asserted
equality
\begin{equation}
  \label{JZeqJacV3}
  \Jj_{Z_s} = \III_{V_3(s)}\otimes\Oo_X.
\end{equation}

We are left with the last assertion of 3) of the theorem, stating that
sections of $\Nn^\ast_X(3H)$ with simple zeros correspond to Segre
cubics in $|I_X(3)|$.  Indeed, let $s$ be a global section of
$\Nn^\ast_X(3H)$ with $Z_s =(s=0)$ consisting of ten distinct points.
From the equality \eqref{JZeqJacV3}, we deduce that the singular locus
$\Sing(V_3(s))$ of the cubic $V_3(s)$ contains ten distinct
points. It will be enough to show that $\Sing(V_3(s))$ is
$0$-dimensional, since then \eqref{JZeqJacV3} tells us that the
singular locus $\Sing(V_3(s))=Z_s$ and it is composed of ten ordinary
double points.  It is well known that such a cubic hypersurface is a
Segre cubic (see \cite{Dol3} for an inspiring introduction to the
subject). 

Let us check now that $\Sing(V_3(s))$ is $0$-dimensional.  By Lemma
\ref{l:V3cone}, $V_3(s)$ is not a cone.  Then the possibilities for
the one dimensional part of $\Sing(V_3(s))$ are a line, a conic
(possibly singular), or a rational normal curve of degree $4$ in
$\PP^4$.

If a conic $C$ is a component of $\Sing(V_3(s))$, then the plane $P$
spanned by $C$ is contained in $V_3$. We examine the pencil of
hyperplanes $V_t$ in $\PP^4$ passing through $P$. The intersection
$V_t \cdot V_3(s)$ is reducible
\[
  V_t \cdot V_3(s)=P\cup Q_t ,\quad\forall t,
\]
where $Q_t$ is a quadric surface residual to the plane $P$ such that
$Q_t \cap P=C$.  The hyperplane section
\begin{equation}
  \label{HtreducibleV3scroll}
  H_t =V_t \cdot X =B+\Gamma_t
\end{equation}
is also reducible, where $B=P\cdot X$ is the $1$-dimensional part of
the base locus of the pencil $\{H_t\}$.  Being a plane divisor, $B$
can be either a ruling of $X$ or its plane cubic section.  The latter
case implies $B\cdot C=6$.  So the part $Z_C$ of the scheme $Z_s$
contained in the intersection $B\cap C$ has degree at least $6$.  On
the other hand the degree of the subscheme $Z_B$ of $Z_s$ contained in
$B$ is at most $B\cdot(H-K_X)=4$ (we use here the fact that $Z_s$ is
contained in an irreducible divisor in the linear system $|H-K_X|$).
Hence $B$ must be a ruling of $X$.  From this and
\eqref{HtreducibleV3scroll} it follows that for a general $t$, the
curve $\Gamma_t$ is an elliptic curve of degree $4$ contained in
$Q_t$.  In particular, the scheme-theoretic intersection
$C\cap\Gamma_t$ is a $0$-dimensional scheme of degree $4$.
Furthermore, since $C$ is not contained in $X$, this scheme must be
contained in the base locus of the pencil $\{\Gamma_t\}$.  Thus we
obtain
\[
  4 \leq \Gamma^2_t = (H-B)^2 = 3
\]
a contradiction.

If a rational normal curve $R$ of degree $4$ is in $\Sing(V_3(s))$,
then $\Sing(V_3(s))=R$ and $Z_s=X\cdot R$.  It follows that
\begin{equation}
  \label{2HthroughZR}
  h^0(\Jj_{Z_s}(2H)) \geq h^0(\Jj_{R/\PP^4}(2)) = 6,
\end{equation}
where $\Jj_{R/\PP^4}$ is the ideal sheaf of $R$ in $\PP^4$.  Since
$Z_s$ lies on a smooth curve $\Gamma\in|H-K_X|$, we can calculate
$h^0(\Jj_{Z_s}(2H))$ from the exact sequence
\[
  0 \lra \Oo_X(-\Gamma) \lra \Jj_{Z_s}
  \lra \Oo_{\Gamma}(-Z_s) \lra 0
\]
tensored with $\Oo_X(2H)$. This gives
\[
  h^0(\Jj_{Z_s}(2H))
  = h^0(\Oo_{\Gamma}(2H|_{\Gamma}-Z_s))
  = 5+h^1(\Oo_{\Gamma}(2H|_{\Gamma}-Z_s))
  = 5+h^0(\Oo_{\Gamma}(Z_s-H|_{\Gamma}).
\]
This and \eqref{2HthroughZR} imply
$h^0(\Oo_{\Gamma}(Z_s-H|_{\Gamma})=1$ or, equivalently,
$\Oo_{\Gamma}(Z_s)=\Oo_{\Gamma}(H)$ and this contradicts Corollary
\ref{cor:HonGammanotZs} below.  

We turn now to the remaining case: the $1$-dimensional locus of
$\Sing(V_3(s))$ is a line $L$.  We divide the scheme $Z_s$ into two
parts,
\[
  Z_s = Z_L +Z',
\]
where $Z_L$ is the part of $Z_s$ contained in $L$ and 
$Z'=Z_s \setminus Z_L$ is the residual part.  It is easy to see that
$\deg(Z_L)\leq 3$.  But then $Z'$ consists of at least $7$ distinct
points.  For every point $z'\in Z'$ the plane $P_{z'}=\Span(z'\cup L)$
is contained in $V_3(s)$.  Furthermore, $P_{z'_1} \neq P_{z'_2}$ for
any pair of distinct points $z'_1, z'_2 \in Z'$ since otherwise
the line $L'=\Span\,\{z'_1,z'_2 \}$ intersects $L$ and hence, is a
component of the singular locus of $V_3(s)$.

Every plane $P_{z'}$ intersects $X$ along a plane curve.  This curve
is either a ruling of $X$ or its plane cubic section. Assume there is
$z'_0 \in Z'$ such that that $P_{z'_0}\cdot X=\Gamma$ is a plane
cubic section.  Then $\Gamma \cap L=Z_L$ and all other planes
$P_{z'}$, $z'\neq z'_0$, intersect $X$ along a ruling.  Since those
rulings must pass through one of the points of $Z_L$, the number of
such planes is at most $3$.  This makes the degree of $Z'$ at most
$4$.  This is contrary to the estimate $\deg(Z')\geq 7$.  Hence every
plane $P_{z'}$ intersects $X$ along a ruling.  But then 
$3 \geq \deg (Z_L) \geq \deg(Z') \geq 7$ which is impossible.  This
completes the proof of the theorem.
\qed

\begin{rem} 
  \label{r:degcubics} 
From the proof of Theorem \ref{th:cubicsSections}, 3), it follows that
if a cubic hypersurface $V_3(s)$ contains $X$ and has $1$-dimensional
singular locus, then its $1$-dimensional part must be a single
line. Furthermore, if this possibility occurs, the global section $s$
of $\Nn^\ast(3H)$ corresponding to $V_3(s)$ under the isomorphism
in Theorem \ref{th:cubicsSections}, 2), must have multiple zeros.  In
the appendix, see \eqref{Velcubic} and the discussion preceding it, we
give an explicit geometric construction of a general cubic in
$|I_X(3)|$, singular along a line.  Hence, that line is precisely the
$1$-dimensional part of the singular locus of such a cubic.  In
addition, we show that the isomorphism in Theorem
\ref{th:cubicsSections}, 2), matches precisely the global sections of
$\Nn^\ast(3H)$ having multiple zeros with the cubics in $|I_X(3)|$
having a (unique) line in their singular locus, see Proposition
\ref{p:degcubics}.
\end{rem} 

In the course of the proof of Theorem \ref{th:cubicsSections}, 3), we
used the fact that on a smooth curve $\Gamma \in |H-K_X |$ containing
$Z_s =(s=0)$, the zero scheme of a nonzero global section $s$ of
$\Nn^\ast (3H)$, the line bundles $\Oo_{\Gamma}(Z_s)$ and
$\Oo_{\Gamma}(H)$ are not isomorphic; see Corollary
\ref{cor:HonGammanotZs}.  This is a part of the proof of the
following.

\begin{lem}
\label{h0JZs2Hscroll}
For every nonzero global section $s\in H^0(\Nn^\ast_X(3H))$ with
$Z_s=(s=0)$ one has $h^0(\Jj_{Z_s}(2H))=5$ and 
$h^1(\Jj_{Z_s}(2H))=0$.
\end{lem}

\proof
By Theorem \ref{th:cubicsSections}, 3), the subscheme $Z_s=(s=0)$ is
$0$-dimensional.  Hence the Koszul sequence of $s$ 
\[
  0 \lra \Oo_X \stackrel{s}\lra \Nn^\ast_X(3H)\stackrel{s\wedge}
  \lra \Jj_{Z_s}(H-K_X) \lra 0
\]
is exact.  Tensoring it with $\Oo_X(K_X+H)$ and using the
identification $\Nn^\ast_X(K_X+4H)\cong\Nn_X(-H)$ we obtain
\[
  0 \lra \Oo_X(K_X+H)\stackrel{s}\lra \Nn_X(-H)\stackrel{s\wedge}
  \lra \Jj_{Z_s}(2H) \lra 0.
\]
From this it follows that $h^0(\Jj_{Z_s}(2H))=h^0(\Nn_X(-H))$ is
independent of the choice of $s$.  So to compute $h^0(\Jj_{Z_s}(2H))$
we choose $s$ with simple zeros and $s'\in H^0(\Nn^\ast_X(3H))$ so
that the curve $\Gamma=(s\wedge s'=0)\in|H-K_X|$ is smooth.  The curve
$\Gamma$ passes through $Z_s$ and gives the following exact sequence
\begin{equation}
  \label{GammaZs}
  0 \lra \Oo_X(-\Gamma)\stackrel{\gamma}
  \lra \Jj_{Z_s} \lra \Oo_{\Gamma}(-Z_s) \lra 0,
\end{equation}
where $\gamma$ is the global section of $\Oo_{X}(H-K_X)$ corresponding
to $s\wedge s'$ under the natural homomorphism 
$\Bigwedge^2 H^0(\Nn^\ast_X(3H))\to H^0(\det(\Nn^\ast_X(3H)))=
H^0(\Oo_{X}(H-K_X))$. 
Tensoring the above sequence with $\Oo_X(2H)$, we deduce
\[
  h^0(\Jj_{Z_s}(2H)) = h^0(\Oo_{\Gamma}(2H|_{\Gamma}-Z_s)) 
  = 5 +h^1(\Oo_{\Gamma}(2H|_{\Gamma}-Z_s))
  = 5 +h^0(\Oo_{\Gamma}(Z_s-H|_{\Gamma})),
\]
where the second equality is the Riemann-Roch for
$\Oo_{\Gamma}(2H|_{\Gamma}-Z_s)$ and the third one is the Serre
duality. Thus the first assertion of the lemma is equivalent to
\begin{equation}
  \label{HonGammaandZs}
  \Oo_{\Gamma}(Z_s) \neq \Oo_{\Gamma}(H).
\end{equation}
Assume the contrary.  Then the exact sequence \eqref{GammaZs} tensored
with $\Oo_X(H-K_X)$ takes the form
\[
  0 \lra \Oo_X \stackrel{\gamma}\lra \Jj_{Z_s}(H-K_X) 
  \lra \Oo_{\Gamma}(-K_X)\lra 0.
\]
This implies
\begin{equation}
  \label{estim-KonGamma}
  h^0(\Oo_{\Gamma}(-K_X)) \geq h^0(\Jj_{Z_s}(H-K_X))-1
  \geq h^0(\Nn^\ast_X(3H))-2 = 3.
\end{equation}
On the other hand we have
\[
  0 \lra \Oo_X(-K_X-\Gamma) \lra \Oo_X(-K_X) 
  \lra \Oo_{\Gamma}(-K_X)\lra 0.
\]
Since $\Oo_X(-K_X -\Gamma)=\Oo_X(-H)$, the above implies
$H^0(\Oo_{\Gamma}(-K_X))= H^0(\Oo_X(-K_X))$.  However, we know that
the last space is $1$-dimensional.  Hence $h^0(\Oo_{\Gamma}(-K_X))=1$
contradicting the estimate in \eqref{estim-KonGamma}.

The second assertion of the lemma about the vanishing of
$H^1(\Jj_{Z_s}(2H))$ follows immediately from the first assertion and
the Riemann-Roch for $\Jj_{Z_s}(2H)$.
\qed

\begin{cor}
  \label{cor:HonGammanotZs}
Let $s$ be a nonzero global section of $\Nn^\ast_X(3H)$ whose zero
locus $Z_s=(s=0)$ is contained in a smooth curve 
$\Gamma \in |H-K_X|$. Then the line bundle 
$\Oo_{\Gamma}(H|_{\Gamma }-Z_s)\neq\Oo_{\Gamma}$.
\end{cor}

\proof
The assertion is a restatement of the identity \eqref{HonGammaandZs}
proved in the previous lemma.
\qed

\section{Irregular surfaces on hypersurfaces of degree 
$4$ with non-degenerate isolated singularities} 

In this section we consider irregular surfaces $X\subset\PP^4$
contained in a hypersurface of degree $4$ and not in one of a smaller
degree. Our main result is as follows.

\begin{thm}
  \label{th:irrsurfV4}
Let $X\subset\PP^4$ be a smooth surface with $m_X=4$ and assume $X$ to
be contained in a quartic hypersurface $V_4$ with at most ordinary
double points.  Then $X$ is regular, with the possible exception of $X$
being a degree $8$ elliptic conic bundle with $H\cdot K_X=0$ and
$K^2_X=-8$.  If such a situation occurs, then $X$ must pass through
precisely $32$ singular points of $V_4$.
\end{thm}

\noindent
The exceptional possibility in the above theorem is an elliptic conic
bundle discovered about $20$ years ago by Abo, Decker, and Sasakura by
using a certain vector bundle of rank $5$ on $\PP^4$, see
\cite{AbDeSa}.  Shortly afterwards, Ranestad in \cite{Ra}, gave a
geometric construction of a general elliptic conic bundle in $\PP^4$
as the image of an elliptic scroll in $\PP^4$ under a certain Cremona
transformation.  In the sequel we refer to those elliptic conic
bundles as ADSR elliptic conic bundle.

From the first work cited above one knows that the space of quartics
$H^0(\Jj_X(4))$ containing $X$ is of dimension $6$.  However, at the
time of writing this paper, we do not know if there are quartics in
$H^0(\Jj_X(4))$ with only ordinary double points.

\medskip

Our proof of Theorem \ref{th:irrsurfV4} follows the same line of
thinking as in the case of surfaces contained in hypersurfaces of
degree $3$.  Namely, we assume $X\subset\PP^4$ to be an irregular
surfaces with $m_X =4$ and lying on a quartic hypersurface $V_4$ with
only ordinary double points.  Our general situation recorded by the
sequence \eqref{eq:koszulForS} takes the form
\begin{equation}
  \label{eq:koszulForSV4}
  0 \lra \Oo_X(K_X+H) \lra \Nn_X \lra \Jj_{Z}(4H) \lra 0,
\end{equation}
where $Z$ is the $0$-dimensional subscheme of $X$ supported on the
singular locus of $V_4$ and defined, at each point $p$ of the support
of $Z$, by the restriction to $X$ of the Jacobian ideal
$\III_{V_4,p}$.  In particular, we have\footnote{The identity
\eqref{JZJacV4} is valid as long as a section defining the Koszul
sequence \eqref{eq:koszulForSV4} has a $0$-dimensional scheme of
zeros.}
\begin{equation}
  \label{JZJacV4}
  \Jj_{Z} = \III_{V_4}\otimes\Oo_X,
\end{equation}
where $\III_{V_4}$ denotes the sheaf of the Jacobian ideal of $V_4$.  

Expressing the second Chern class of $\Nn_X$ from the exact sequence
\eqref{eq:koszulForSV4} provides a new ``double point formula''
\begin{equation}
  \label{eq:ndp4}
  d^2 = \deg(Z) + 4H\cdot(H+K_X) = \deg(Z)+8(g-1),
\end{equation}
where $g=g(H)$ is the geometric genus of a general hyperplane section.
In the sequel we refer to this identity as the ndp formula.

The cohomological sequence \eqref{eq:koszulForS-H1} controlling the
irregularity of $X$ becomes 
\begin{equation}
  \label{eq:koszulForS-H1_V4}
  H^1(\Oo_X(2K_X+H)) \lra \Nn_X(K_X) \lra H^1(\Jj_{Z}(K_X+4H)). 
\end{equation}
Therefore, an understanding of the irregularity is reduced to the
study of the cohomology groups $H^1(\Oo_X(2K_X+H))$ and
$H^1(\Jj_{Z}(4H+K_X))$.  In particular, we need to control the scheme
$Z$ which, in view of the identity \eqref{JZJacV4}, comes down to
controlling the singular locus $\Sing(V_4)$ of $V_4$.  Under our
assumption on the isolated singularities of $V_4$, the singular locus
$\Sing(V_4)$ is the set of ordinary double points of $V_4$ and we can
quote a result of A.~Varchenko, \cite{Va}, for the estimate 
$\deg(\Sing(V_4))\leq 45$.  This together with \eqref{JZJacV4} gives
\begin{equation}
  \label{eq:upperBoundZ0V4}
  \deg(Z) \leq \deg(\Sing(V_4)) \leq 45.
\end{equation}

With this estimate of $\deg(Z)$ recorded, we turn now to the study of
the cohomology groups $H^1(\Oo_X(2K_X+H))$ and
$H^1(\Jj_{Z}(K_X+4H))$ in \eqref{eq:koszulForS-H1_V4}.

\subsection{The study of $H^1(\Oo_X(2K_X+H))$}

By the Serre duality, $H^1(\Oo_X(2K_X+H))^\ast=H^1(\Oo_X(-(K_X+H)))$.
Thus the question of the (non)vanishing of this group comes down to
understanding the geometric properties of the divisor $K_X+H$.  The
next lemma is an easy consequence of \cite{Re}.

\begin{lem}
  \label{l:basept-scrollbasepf-conicf}
Let $X$ be an irregular surface and $H$ be a very ample divisor on
$X$.  Then the following assertions hold.
\begin{enumerate}[label={\rm\arabic*)}]
  \listspace
\item
  $\Oo_X(K_X+H)$ has base points if and only if $X$ is a ruled surface
  and its embedding by $\Oo_X(H)$ is a scroll.
\item
  If $\Oo_X(K_X+H)$ is base point free and $H^1(\Oo_X(-(K_X+H)))\neq0$,
  then $X$ is a birationally ruled surface embedded by $\Oo_X(H)$ as a
  conic bundle over a smooth curve $B$ of genus $q=q(X)$.
\end{enumerate}
\end{lem}

\proof
The assumption that $X$ is irregular implies $H^2 \geq 5$.  Then, by
\cite[Theorem 1,(i)]{Re}, a base point of $\Oo_X(K_X+H)$ gives rise to
an effective divisor $D\subset X$ passing through a base point of
$|K_X+H|$ such that $H\cdot D=1$ and $D^2=0$.  It follows that $D$ is
a line in the embedding given by $\Oo_X(H)$.  Hence the Albanese map
$a:X\to\Alb(X)$ must contract $D$ to a point.  The fact that $D^2=0$
implies that the map $a$ factors through a smooth curve
$B\subset\Alb(X)$ of genus $q=q(X)$ and $a:X\to B$ is a
$\PP^1$-fibration, with $D$ one of the fibres.  Thus $X$ is a ruled
surface embedded by $\Oo_X(H)$ as a scroll.  The assertion in the
other direction is obvious.

We turn now to the assertion 2) of the lemma.  The hypotheses imply
that $\Oo_X(K_X+H)$ is nef but not big, \ie that $(K_X+H)^2=0$.  Since
$\Oo_X(K_X+H) \neq \Oo_X$ (otherwise $\Oo_X(K_X)=\Oo_X(-H)$ and hence
$q=h^1(\Oo_X)=h^1(\Oo_X(K_X))=h^1(\Oo_X(-H))=0$), the linear system
$|K_X+H|$ induces a morphism whose image is a curve.  More precisely,
there is a morphism
\begin{equation}
  \label{pi-conicfibration}
  \pi: X \lra B
\end{equation}
onto a smooth curve $B$ with connected fibres and a base point free
line bundle $\Oo_B(D)$ on $B$ such that 
$\Oo(K_X+H)=\pi^\ast(\Oo_B(D))$.  From this we obtain the relation
\begin{equation*}
  \label{K+H=f}
  K_X+H = \deg(D)F
\end{equation*}
in $\NS(X)$, where $F$ stands for the class of a fibre of $\pi$.
Intersecting with $F$ the above identity, we deduce that $H\cdot F=2$,
\ie that $\pi$ in \eqref{pi-conicfibration} is a conic fibration.
Thus a general fibre of $\pi$ is a smooth conic and there is at most a
finite number of singular fibres.  {\it A priori}, a singular fibre is
either the union of two lines intersecting transversely or a double
line.  The latter, however, is impossible since $F=2L$ with $L$ a line
leads to $K_X\cdot L=L^2=-1$ which contradicts $0=F^2=4L^2$.
\qed

We apply the above result to a surface $X$ subject to the hypotheses of 
Theorem \ref{th:irrsurfV4}

\begin{pro}
  \label{p:2KPlusH}
Let $X$ and $V_4$ be as in Theorem \ref{th:irrsurfV4} and assume $X$
to be irregular.  Then $\Oo_X(K_X+H)$ is base point free.
Furthermore, $\Oo_X(K_X+H)$ is big and hence $H^1(\Oo_X(2K_X+H))=0$
with the possible exception of $X$ being an ADSR elliptic conic bundle.
If such a situation occurs, then $X$ must pass through precisely $32$
singular points of $V_4$.
\end{pro}

\proof
By Serre duality $H^1(\Oo_X(2K_X+H))\cong H^1(\Oo_X(-(K_X+H)))^\ast$.
According to Lemma \ref{l:basept-scrollbasepf-conicf}, the latter
group is nonzero if $X$ is embedded, either as a scroll, or as a conic
bundle.  The first possibility implies that $X$ is an elliptic scroll
of degree $5$, see Lemma \ref{l:scrollsInP4}.  But such a scroll, as we
have seen in the previous section, is contained in hypersurfaces of
degree $3$, hence it can not occur here.  We turn to the second
possibility: $X$ is birational to a ruled surface embedded into
$\PP^4$ by $\Oo_X(H)$ as a conic bundle.  More precisely, from the
proof of Lemma \ref{l:basept-scrollbasepf-conicf}, 2), the line bundle
$\Oo (K_X+H)$ induces a morphism $\pi:X\to B$ onto a smooth curve $B$
of genus $q=q(X)$, the irregularity of $X$, such that the $H$-degree
of the fibres of $\pi$ is $2$.  A general fibre of $\pi$ is a smooth
plane conic, while the singular fibres are reduced singular conics.

If $X$ is not minimal, then $\pi$ factors trough a minimal model of
$X$, call it $X'$, and gives the diagram
\vspace{-.5ex}
\[
  \begin{tikzpicture}[every node/.style={draw=none},
    ->,inner sep=1.1ex]
    \matrix [draw=none,row sep=4ex,column sep=4.2ex]
    {
      \node (11) {$X$};
      && \node (13) {$X'$}; \\
      & \node (22) {$B$}; \\
    };
    \path
      (11) edge node[above=-.2ex,scale=.75] {$\sigma$} (13)  
      (11) edge node[below left=-.1ex,scale=.75] {$\pi$} (22)
      (13) edge node[below right=-.4ex,scale=.75] {$\pi'$} (22)
    ;
  \end{tikzpicture}
  \vspace{-2ex}
\]
where $\sigma$ is the blow-down of a collection of $(-1)$-curves on
$X$ and $X'$ is a ruled surface over $B$ with $\pi'$ its structure
morphism.  The collection of exceptional curves is a choice of one out
of the two irreducible components of each reducible fibre of
$\pi$---these are the only $(-1)$-curves of $X$.  This implies, in
particular, that $X$ is the blow-up of $X'$ at distinct points.  Let
$\delta$ be the number of blown-up points and let $\{C_j\}$ be the
collection of the exceptional $(-1)$-curves on $X$.  Then we write
\begin{equation*}
  \label{eq:strictTransform}
  K_X = \sigma^\ast (K_{X'})+\Delta,
\end{equation*}
where $K_{X'}$ is the canonical divisor of $X'$ and 
$\Delta =\sum_{j=1}^\delta C_j$ is the sum of the exceptional
$(-1)$-curves blown-down by $\sigma$.  This implies
\begin{equation}
  \label{Chernscrollbl}
  K^2_X = K^2_{X'} -\delta = -8(q-1) -\delta
  \textq{and}
  c_2 = -4(q-1) +\delta.
\end{equation}
Since $\Oo_X(H+K_X)$ is composed of a pencil, we also have
$(H+K_X)^2 =0$, hence
\[
  K^2_X =-d-2H\cdot K_X.
\]
From this and the first identity in \eqref{Chernscrollbl}, it follows
that
\begin{equation}
  \label{eq:HKconicbundle}
  H\cdot K_X = 4(q-1)-\frac{d-\delta}{2}.
\end{equation}
Substituting this and the Chern numbers computed in
\eqref{Chernscrollbl} into the double point formula, we deduce the
identity
\begin{equation}
  \label{dp-conicbundle}
  d^2-\frac{15}{2}\,d-\frac{1}{2}\,\delta = 16(q-1).
\end{equation}

Since $X$ is a conic bundle, we can associate to $\pi:X\to B$ the
embedding
\begin{equation}
  \label{phiP2bundle}
  \phi: B \lra \Gr(2,\PP^4)
\end{equation}
of the base curve $B$ into the Grassmannian $\Gr(2,\PP^4)$ of planes
in $\PP^4$, where $\phi$ sends a point $b\in B$ to the plane $P_b$
spanned by the conic $F_b =\pi^{-1}(b)$, the fibre of $\pi$.  Let
$\Uu$ be the pullback under $\phi$ of the universal subbundle of
$\Gr(2,\PP^4)$.  It is a rank $3$ bundle on $B$ and $\phi$ induces the
morphism
\begin{equation}
  \label{eq:phitilde_P2Bundle}
  \phitilde: \PP(\Uu) \lra \PP^4
\end{equation}
defined on its projectivization $\PP(\Uu)$.  The image of $\phitilde$
is a $3$-fold which, set-theoretically, is the union of the family of
planes $\{P_b \}_{b\in B}$.  In particular, the line bundle 
$\Oo_{\PP(\Uu)}(1):=\phitilde^\ast(\Oo_{\PP^4}(1)) $ is determined by
the identity
\begin{equation}
  \label{eq:UAndH}
  \Uu^\ast =\rho_\ast(\Oo_{\PP(\Uu)}(1)) = \pi_\ast(\Oo_X(H)),
\end{equation}
where $\rho:\PP(\Uu)\to B$ is the structure projection.  We will need
to know the degree of $\Uu^\ast$.

\paragraph{Claim}
$\deg(\Uu^\ast)=\dfrac{3d-\delta}{4}$.

\medskip

To justify the claim, we start by computing the holomorphic Euler
characteristic of $\Uu^\ast$ in \eqref{eq:UAndH}:
\[
  \chi(\Uu^\ast)
  = \chi(\Oo_X(H)) = \frac{H^2-H\cdot K_X}{2} + \chi(\Oo_X)
  = \frac{d-H\cdot K_X}{2}-(q-1).
\] 
On the other hand, the Riemann-Roch for $\Uu^\ast$ on $B$ gives 
$\chi(\Uu^\ast)=\deg (\Uu^\ast)-3(q-1)$.  Putting the two
expressions for $\chi(\Uu^\ast)$ together, we deduce
\[
  \deg(\Uu^\ast) = \frac{d-H\cdot K_X}{2}+2(q-1).
\]
This and the expression for $H\cdot K_X$ in \eqref{eq:HKconicbundle}
imply the equality of the claim.

\medskip
 
Set 
\begin{equation}
  \label{eq:d'}
  d':=\frac{3d-\delta}{4}.
\end{equation}
The geometric meaning of $d'$ is two-fold:
\begin{enumerate}[label={\rm \arabic*)}]
  \listspace
\item
  If $\phi': B\hookrightarrow \PP^9$ is the composition of $\phi$ with
  the Pl\"ucker embedding of $\Gr(2,\PP^4)$ then $d'$ is the degree of
  the image $B'=\phi'(B)$. 
\item 
  If $V$ is the image of $\phitilde$, then $d'=\deg(V)$.
\end{enumerate}
It will be more convenient at this point to express the double point
formula \eqref{dp-conicbundle} in terms of $d$ and $d'$:
\begin{equation}
  \label{dpddprime}
  d^2 -9d +2d' =16(q-1).
\end{equation}
We also bring in the ndp formula \eqref{eq:ndp4} to obtain an upper
bound for $d'$
\[
\begin{split}
  d^2 &= 4(H^2 +H\cdot K_X)+\deg(Z) \\
  &= 4\bigg(d+4(q-1)-\frac{d-\delta}{2} \bigg) +\deg(Z) \\
  &= 16(q-1)+8(d-d')+\deg(Z),
\end{split}
\]
where the second equality uses \eqref{eq:HKconicbundle} and the last
one \eqref{eq:d'}.  This expression for $d^2$ and the identity
\eqref{dpddprime} give
\begin{equation}
  \label{Zdd'}
  \deg(Z)=d+6d'.
\end{equation}
From this, the upper bound $\deg(Z)\leq 45$ in
\eqref{eq:upperBoundZ0V4}, and $d\geq 5$, we deduce that $d'\leq 6$.
This upper bound tells us that the curve $B'$, the image of $\phi$,
spans a subspace $\PP^N$ of dimension $3\leq N\leq 5$.  The value
$N=2$ is excluded; indeed, if $B'$ spans a plane, then that plane
either intersects $\Gr(2,\PP^4)$ along $B'$ or it is contained in the
Grassmannian $\Gr(2,\PP^4)$.  The first possibility means that $B'$ is
a conic and hence, $q=0$, while the second possibility tells us that
all planes $P_b$, $b\in B$, intersect along a line, call it $l$; but
then $l\cong\PP^1\subset X$ is a multi-section of $\pi:X\to B$
and this forces $q$ to be zero again.  Furthermore, if $N=5$, then
$d'= 6$ and $B$ must be an elliptic curve, \ie $q=1$.  Substituting
these values in \eqref{dpddprime}, we obtain
\[
  d^2 -9d+12=0
\]
which has no integer solutions.  Thus $ N=3$ or $4$.

If $N=4$, then $d'=5$ or $6$.  The first possibility implies again
that $q=1$ and the formula \eqref{dpddprime} becomes $d^2-9d+10=0$
with no integer solution.  The second possibility, $d'=6$, implies
that $q=1,2$. The first value has been ruled out in the discussion of
the case $N=5$.  As for the second, the formula \eqref{dpddprime}
becomes $d^2-9d-4=0$ with no integer solutions.

Thus $N=3$ is the only admissible value, while $d'=4,5$ or $6$.  For
$d'=6$, the Castelnuovo upper bound on genus gives $q\leq 4$. Only
$q=4$ is compatible with \eqref{dpddprime}, leading to the equation
$d^2 -9d-36=0$.  Hence $d=12$. Substituting into \eqref{Zdd'}, we
obtain the contradiction
\[
  45 \geq \deg(Z) = 12+36 = 48.
\]
For $d'=5$, the Castelnuovo upper bound implies $q=1,2$.  The first
value was already discarded in the case $N=4$, while the second value
substituted into \eqref{dpddprime} gives $d^2-9d-6=0$ with no
integral solutions.

Thus we are left with $d'=4$ and hence $q=1$.  These values
substituted in \eqref{dpddprime} yield $d^2-9d+8=0$ with the integer
solution $d=8$. We now go to \eqref{eq:d'} to deduce
\[
K^2_X =-\delta =-8.
\]
This together with \eqref{eq:HKconicbundle} imply $H\cdot K_X=0$.
Thus $X$ is an ADSR elliptic conic bundle.  Furthermore, from the
formula \eqref{Zdd'} it follows that $\deg(Z)=32$.  This completes the
proof of the proposition.
\qed

As we have mentioned in the discussion following the statement of
Theorem~\ref{th:irrsurfV4}, we do not know if an ADSR elliptic conic
bundle is contained in a quartic hypersurface with isolated ordinary
double points only.  On the other hand, such a surface, by its very
definition, is contained in a distinguished quartic whose singular
locus is $2$-dimensional.  This implicitly appears in the works
\cite{AbDeSa} and \cite{Ra}.  In the proof of
Proposition~\ref{p:2KPlusH} this distinguished quartic is $V$, the
image of the morphism $\widetilde{\phi}:\PP(\Uu)\to\PP^4$ in
\eqref{eq:phitilde_P2Bundle}.  The following statement summarizes the
properties of the vector bundle $\Uu$ and its relation to the geometry
of $V$.

\begin{pro}
  \label{p:U}
Let \/$\Uu^\ast$ be the vector bundle defined in \eqref{eq:UAndH} and
let $V$ be the image of the morphism
$\widetilde{\phi}:\PP(\Uu)\to\PP^4$ in \eqref{eq:phitilde_P2Bundle}.
Then
\begin{enumerate}[label={\rm \arabic*)}]
  \listspace
\item 
$\Uu^\ast$ has the form
\begin{equation}
  \label{Udirsum}
  \Uu^\ast=\Oo_B \oplus \Ff^\ast,
\end{equation}
where $\Ff^\ast$ is a rank $2$ bundle on $B$ fitting into the exact
sequence 
\begin{equation}
  \label{Ffilt}
  0\lra \Oo_B(D)\lra \Ff^\ast \lra \Oo_B(D') \lra 0,
\end{equation}
with $\Oo_B(D)$ and $\Oo_B(D')$ being line bundles of degree $2$.
\item
$V\subset\PP^4$ is a hypersurface of degree $4$.  It is a cone with
vertex $[v]$, the image of the section of \,$\PP(\Uu)$ corresponding
to the trivial summand in the direct sum decomposition
\eqref{Udirsum}.
\item
The summand $\Ff^\ast$ in \eqref{Udirsum} determines a distinguished
divisor \,$\PP(\Ff)$ of \,$\PP(\Uu)$.  Its image $S$ under
$\widetilde{\phi}$ is a base of the cone $V$.  In particular, $S$ is a
quartic surface with the singular locus $\Sing(S)$ consisting of
either one or two (skew) lines in $\PP^3$.  The latter possibility
occurs when the exact sequence \eqref{Ffilt} splits.
\item
The singular locus $\Sing(V)$ of $V$ is the cone over $\Sing(S)$ with
vertex at $[v]$.  In particular, set-theoretically, it is composed of
either one or two planes depending on whether or not the sequence
\eqref{Ffilt} is non split.
\end{enumerate} 
\end{pro}

\proof
Consider the pullback under the morphism $\phi$ in \eqref{phiP2bundle}
of the dual of the universal sequence on $\Gr(2,\PP^4)$
\begin{equation}
  \label{dualunivseq}
  0 \lra \Ff \lra H^0(\Uu^\ast) \otimes \Oo_B
  \lra \Uu^\ast \lra 0.
\end{equation}
This implies that $H^0(\Ff)=0$.  The subbundle $\Ff\subset
H^0(\Uu^\ast)\otimes\Oo_B$ defines the morphism 
\[
\phi\dual : B\lra  \Gr(1,(\PP^4)\,\dual\,), 
\]
dual of the morphism $\phi$ in \eqref{phiP2bundle}, where
$(\PP^4)\dual =\PP(H^0(\Uu^\ast))$.  Composing $\phi\dual$ with the
Pl\"{u}cker embedding 
$\Gr(1,(\PP^4)\dual)\subset\PP(\Bigwedge^2H^0(\Uu^\ast))\cong(\PP^9)\dual$, 
we obtain the embedding
\[
  \psi: B\hookrightarrow (\PP^9)\dual.
\]
The image of $\psi$, as the image of the embedding $\phi'$ in the
proof of Proposition \ref{p:2KPlusH}, spans a $\PP^3$. Hence the image
of the linear map
\[
  w: \Bigwedge^2 H^0(\Uu^\ast)^\ast \lra H^0(\det(\Ff^\ast))
\]
defining the morphism $\psi$ is $4$-dimensional, while $\ker(w)$ is a
$6$-dimensional subspace of $\Bigwedge^2 H^0(\Uu^\ast)^\ast$.  This
means that $\PP(\ker(w))$ intersects the Grassmann variety of
decomposable tensors in $\PP(\Bigwedge^2 H^0(\Uu^\ast)^\ast)$ along a
subscheme of dimension at least $1$ implying that the linear map
\begin{equation}
  \label{nontrivker}
  H^0(\Uu^\ast)^\ast \lra H^0(\Ff^\ast)
\end{equation}
has a non-trivial kernel.  Indeed, let $l\wedge l'$ be a nonzero
decomposable tensor in $\Bigwedge^2 H^0(\Uu^\ast)^\ast$ lying in the
kernel of $w$.  We may assume that the pencil $\Span(l,l')$ injects
into $H^0(\Ff^\ast)$ under the map in \eqref{nontrivker}, since
otherwise we are done.  Thus we can think of $l$ and $l'$ as two
linearly independent global sections of $\Ff^\ast$ which are
proportional, \ie  $l\wedge l'$ is zero as a section of
$\det(\Ff^\ast)$.  Hence the Koszul sequence of one of these sections
gives rise to an exact sequence
\begin{equation}
  \label{seqFdual}
  0\lra \Oo_B(D)\lra \Ff^\ast \lra \Oo_B(D')\lra 0,
\end{equation}
where $h^0(\Oo_B(D))\geq 2$. Hence $\deg(D) \geq 2$.  Furthermore,
since $\deg(\Ff^\ast) =\deg(\Uu^\ast) =4$ and the quotient $\Oo_B(D')$
must be generated by its global sections and nontrivial (the latter
comes from $H^0(\Ff)=0$), we deduce $\deg(D) =\deg(D')= 2$.  This and
the exact sequence \eqref{seqFdual} imply $h^0(\Ff^\ast)=4$.  Since
$H^0(\Uu^\ast)^\ast \cong H^0(\Oo_X(H)$ is $5$-dimensional, we deduce
that the kernel in \eqref{nontrivker} is nontrivial.

Considering the dual of \eqref{dualunivseq}, we see that the kernel in
\eqref{nontrivker} is $H^0(\Uu)$.  Now, from $H^0(\Uu) \neq 0$ and the
global generation of $\Uu^\ast$, we deduce the direct sum
decomposition
\begin{equation}
  \label{dirsum1}
  \Uu^\ast \cong \Oo_B \oplus \Gg.
\end{equation}
This together with \eqref{dualunivseq} implies that the direct summand
$\Gg$ fits into the exact sequence
\begin{equation}
  \label{FandG}
  0 \lra \Ff \lra H^0(\Gg)\otimes \Oo_B \lra \Gg \lra 0.
\end{equation}
It remains to identify $\Gg$ with $\Ff^\ast$.  To this end, we go back
to the morphism $\phi': B\to\PP^9$ in the proof of Proposition
\ref{p:2KPlusH} and recall that its image spans a $\PP^3$.  This means
that the linear map
\[
  \Bigwedge^3 H^0(\Uu^\ast) \lra H^0(\det(\Uu^\ast))
\]
has the image of dimension $4$ or, equivalently, the kernel of
dimension $6$.  This together with the decomposition in
\eqref{dirsum1} implies that the linear map
\[
  \Bigwedge^2 H^0(\Gg) \lra H^0(\det(\Gg))
\]
has the kernel, call it $W$, of dimension at least $2$.  Combining
this and the second exterior power of \eqref{FandG} gives
$h^0(\Ff\otimes\Gg)=h^0(\End(\Ff^\ast,\Gg))\geq 2$.  It follows that a
general morphism $\Ff^\ast\to\Gg$ is an isomorphism.  Thus 
$\Gg\cong\Ff^\ast$. Substituting into \eqref{dirsum1}, we deduce the
decomposition asserted in \eqref{Udirsum}.  Furthermore, we have
$h^0(\End(\Ff^\ast,\Ff^\ast))\geq 2$. A nontrivial endomorphism of
$\Ff^\ast$ gives rise to the exact sequence \eqref{Ffilt}. 

The remaining assertions of the proposition are obvious geometric
analogues of the properties of \,$\Uu^\ast$ (resp. $\Ff^\ast$) in 1).
\qed

\begin{rem}
  \label{r:rulingDegreeV4} 
  Let $X$ be a smooth surface in $\PP^4$ such that
  \vspace{-1ex}
\begin{description}
  \listspace
\item[---]
  $X$ is birational to an irregular ruled surface,
\item[---]
  $m_X=4$ and $X$ is contained in a quartic hypersurface with only
  nondegenerate isolated singularities, 
\item[---]
  the degree of $X$ in $\PP^4$ is not $8$.
\end{description}
\vspace{-1ex}
Then from the proof of Proposition~\ref{p:2KPlusH} it follows that the
fibres of $X$ are embedded as curves of degree at least $3$.
\end{rem}

\subsection{The study of $H^1(\Jj_{Z}(K_X+4H))$}

Our result here is as follows.

\begin{pro}
  \label{p:h1JZV4}
Let $X$ and $V_4$ be as in Theorem \ref{th:irrsurfV4} and assume $X$
to be irregular. Then $H^1(\Jj_{Z}(K_X+4H))=0$ with a possible
exception of $X$ being an ADSR elliptic conic bundle.  If such a
situation occurs, then $h^1(\Jj_{Z}(K_X+4H)=1$ and 
$Z=X\bigcap\Sing(V_4)$ is a subset of $32$ nodes of $V_4$.
\end{pro}

We assume the nonvanishing of the cohomology group
$H^1(\Jj_{Z}(K_X+4H))$.  The identification
\begin{equation*}
  \label{Z0Ext1}
  H^1(\Jj_{Z}(K_X+4H))^\ast = \Ext^1(\Jj_{Z}(4H),\Oo_X)
\end{equation*}
provided by the Serre duality, gives rise to a nontrivial extension
sequence
\begin{equation}
  \label{eq:theExtension}
  0 \lra \Oo_X \lra \Ee \lra \Jj_{Z}(4H) \lra 0.
\end{equation}
We may assume the sheaf $\Ee$ in the middle of this sequence to be
locally free.  The Bogomolov semistability condition for this sheaf
reads
\begin{equation}
  \label{eq:BogomForEV4}  
  0 \leq 4c_2(\Ee)-c_1^2(\Ee) = 4\deg(Z)-16H^2 = 4(\deg(Z)-4d).
\end{equation}
From this and the upper bound $\deg(Z)\leq 45$, it follows that $\Ee$
is Bogomolov unstable provided $d\geq 12$.  For this reason, the proof
of Proposition \ref{p:h1JZV4} is naturally divided into two parts:
\begin{itemize}
  \listspace
\item
  the first part rules out the case $d\geq 12$ by examining the
  geometric consequences of the Bogomolov instability of $\Ee$;
\item
  the second part deals with the remaining values $5\leq d\leq11$
  for the degree of $X$.
\end{itemize}

\paragraph{First part: $d\geq12$}
We begin by recording some geometric consequences of the Bogomolov
instability condition $\deg(Z)<4d$ for the sheaf $\Ee$ in
\eqref{eq:theExtension}.

Let $\Oo_X(A)$ be the maximal Bogomolov destabilizing subsheaf of
$\Ee$. Combining the inclusion $\Oo_X(A)\to\Ee$ with the defining
extension sequence \eqref{eq:theExtension} gives the diagram
\vspace{15ex}
\begin{equation}
  \label{d:EunstV4}
  \begin{tikzpicture}[overlay,every node/.style={draw=none},
    ->,inner sep=1.1ex]
    \matrix [draw=none,row sep=3.5ex,column sep=4ex]
    {
      && \node (02) {$0$}; \\
      && \node (12) {$\Oo_X(A)$}; \\
      \node (20) {$0$};
      & \node (21) {$\Oo_X$};
      & \node (22) {$\Ee$};
      & \node (23) {$\Jj_Z(4H)$};
      & \node (24) {$0$}; \\
      && \node (32) {$\Jj_{Z'}(E)$}; \\      
      && \node (42) {$0$}; \\      
    };
    \path
      (02) edge (12) 
      (12) edge (22)
      (22) edge (32)
      (32) edge (42)

      (12) edge (23)
      (21.-47) edge (32)
      
      (20) edge (21)
      (21) edge (22)
      (22) edge (23)
      (23) edge (24)
    ;
  \end{tikzpicture}
  \vspace{16.5ex}
\end{equation}
where $Z'$ is a $0$-dimensional subscheme of $X$ and $\Jj_{Z'}$ is its
ideal sheaf.  The Bogomolov destabilizing condition tells us that the
divisor
\[
  B := A-2H
\]
is in $N^+(X)$, the positive cone of $X$.  In particular,
\begin{equation}
  \label{AEB}
  A = 2H+B \textq{and}  E = 2H-B.
\end{equation}
The first formula implies that the slanted arrows in the above diagram
are nonzero.  From this it follows that $E$ is an effective, nonzero
divisor and that those arrows are defined by multiplication by a
section $e\in H^0(\Oo_X(E))$ corresponding to $E$.  In particular,
from the upper (resp. lower) slanted arrow, we deduce that $Z$
(resp. $Z'$) is contained in $E$.  On the other hand we also know that
the linear system $|\Jj_{Z}(3H)|$ is base point free outside $Z$.
Hence there is a reduced irreducible divisor $C\in|3H|$ containing
$Z$.  Therefore, $Z \subset C\cdot E$ and thus subject to the estimate
\begin{equation}
  \label{Z0compiV4ineq}
  \deg(Z) \leq 3H\cdot E.
\end{equation} 
Computing the second Chern number $c_2(\Ee)=\deg(Z)$ from the
vertical sequence in \eqref{d:EunstV4} gives
\begin{equation}
  \label{Z0compiV4}
  \deg(Z)
  = A\cdot E +\deg(Z')
  = (4H-E)\cdot E +\deg(Z')
  = 4H\cdot E -E^2 +\deg(Z').
\end{equation}
Together with the inequality \eqref{Z0compiV4ineq}, this expression
implies 
\[
  E^2 \geq H\cdot E +\deg(Z').
\]  
This inequality acquires more geometry by observing that the divisor
$E-H$ is effective\footnote{This is seen by tensoring
  \eqref{d:EunstV4} with $\Oo_X(-H)$ and showing that
  $H^0(\Jj_{Z'}(E-H))\neq0$; the argument is exactly the same as in the
  proof of Lemma \ref{l:rightsidenoncone}.}, 
allowing us to write $E=H+R$, where $R$ is an effective divisor.
Combining this and the second formula in \eqref{AEB}, we have
\begin{equation}
  \label{HBR}
  H = B+R.
\end{equation}
In addition, by substituting the formulas from \eqref{AEB} into
\eqref{Z0compiV4}, we obtain the identity
\begin{equation}
  \label{eq:ZFuncitonOfBAndZPrime}
  \deg(Z) = A\cdot E +\deg(Z')
  = (2H+B)\cdot(2H-B)+\deg(Z')
  = 4d-B^2 +\deg(Z'),
\end{equation}
which, together with the bound $\deg(Z)\leq45$, enables us
to control the degree $d$ and the intersection number $H\cdot R$.

\begin{lem}
  \label{l:HRleq2}
If $d\geq12$, then $H\cdot R\leq2$, \ie $R$ is empty, a line, or a
conic.
\end{lem}

\proof
Assume that $H\cdot R\geq 3$.  By the Hodge index we have
\[
  B^2 \leq \frac{(H\cdot B)^2}{d} = \frac{(H\cdot(H-R))^2}{d}
  = \frac{(d-H\cdot R)^2}{d} \leq \frac{(d-3)^2}{d}
\]
where the last inequality uses the fact that 
$\frac{(d-H\cdot R)^2}{d}$ is a decreasing function of $H\cdot R$ on
the interval $[3,d]$.  Substituting this upper bound for $B^2$ in
\eqref{eq:ZFuncitonOfBAndZPrime} gives the estimate
\begin{equation}
  \label{Z0d}
  45 \geq \deg(Z) = 4d-B^2 +\deg(Z')
  \geq 4d-\frac{(d-3)^2}{d} +\deg(Z')
  = 3d +6 -\frac{9}{d} +\deg(Z').
\end{equation}
From this it follows that the only admissible values for $d$ are $12$
and $13$.  Furthermore, the ndp formula \eqref{eq:ndp4} yields the
divisibility condition
\begin{equation}
  \label{div8}
  \deg(Z) \equiv d^2 \,(\text{mod}\,\, 8).
\end{equation}
For $d=13$, it implies $\deg(Z)\leq41$, hence \eqref{Z0d} becomes
\[
  41 \geq \deg(Z) \geq 3\cdot13+6-\frac{9}{13} +\deg(Z') > 44
\]
which is absurd.  The same argument rules out the other admissible
value as well.
\qed

Next we examine the three possibilities for $R$ provided by the
previous lemma.

If $R=0$, then $B=H$ and $Z'=0$ and the formula
\eqref{eq:ZFuncitonOfBAndZPrime} becomes $45\geq \deg(Z)=3d$.  Hence
$12\leq d\leq 15$.  But none of these values satisfies the
divisibility condition \eqref{div8}.

If $R$ is a line, then $B=H-R$ and $B^2=d-2+R^2$.  Substituting the
self-intersection number into
\eqref{eq:ZFuncitonOfBAndZPrime} gives
\[
  45 \geq \deg(Z) = 3d +2 -R^2 +\deg(Z').
\]
Hence $d=12,\,13$ or $14$.  The last two values together with the
divisibility condition in \eqref{div8} imply that $\deg(Z)=41$ and
$44$ respectively.  But both values force $R^2=\deg(Z')=0$.  Thus an
irregular $X$ is a scroll and this is ruled by Remark
\ref{r:rulingDegreeV4}.  Therefore, we are left with $d=12$ and
$\deg(Z)=40$.  However, the inequality in \eqref{Z0compiV4ineq} now
reads
\[
  40 \leq 3H\cdot E = 3H\cdot(H+R) = 36+3 = 39
\]
which is absurd.

If $R$ is a conic, then $H\cdot B=10$ and $B^2=d-4+R^2$.  Substituting
into \eqref{eq:ZFuncitonOfBAndZPrime} gives
\[
 45 \geq \deg(Z) = 3d +4 -R^2 +\deg(Z').
\]
This together with the divisibility condition \eqref{div8} implies
$\deg(Z)=40$, $d=12$, and $\deg(Z')=R^2=0$.  The last equality
together with $H\cdot R=2$ tells us that an irregular surface $X$ is a
conic bundle and this is impossible by Remark~\ref{r:rulingDegreeV4}.
This completes the treatment of the case $d\geq 12$ and thus proves
the vanishing of $H^1(\Jj_{Z}(K_X+4H))$ for these values of $d$.

\paragraph{Second part: $5\leq d\leq11$}
Though the argument comes down to a case by case consideration, there
is a basic feature that is common to all of them.  This aspect will be
explained next.  We begin by recalling that the values of the degree
$d$ control the values of $\deg(Z)$ via the divisibility condition
\eqref{div8}, $\deg(Z)\equiv d^2\,(\hspace{-1.25ex}\mod 8)$.  Hence we
obtain
\begin{equation}
  \label{mod8restriction}
  \deg(Z) =
  \begin{cases}
    8k+1, &\textq{if} d\in\{5,7,9,11\} \\
    8k+4, &\textq{if} d\in \{6,10\} \\
    8k,   &\textq{if} d =8
  \end{cases}
\end{equation}
where $1\leq k\leq 5$, in view of the upper bound
$\deg(Z)\leq 45$.  Thus for every value of $d$ we obtain a list of
admissible values for $\deg(Z)$.  These values are divided into two
types according to whether the sheaf $\Ee$ in the extension sequence
\eqref{eq:theExtension} is semistable or unstable in the sense of
Bogomolov.

When $\Ee$ is unstable, all the possibilities are discarded by
exploiting the decomposition $H=B+R$ in \eqref{HBR} and by observing
that the divisor $R$ has to move in a linear system in order for the
hypersurface $V_4$ to have isolated singularities.

In the semistable case there are two basic ingredients:
\begin{itemize}
  \listspace
\item[$-$]
  we check that $X$ is birational to a ruled surface of irregularity
  $q$,
\item[$-$]
  the two conditions, $X$ birational to a ruled surface and
  $\Oo_X(K_X+H)$  base point free and big, exclude all but one
  possibility: $X$ is an ADSR elliptic conic bundle, see Proposition
  \ref{p:2KPlusH}. 
\end{itemize}

With these general guidelines in mind, we proceed with the second
part of the proof according to the possible values of the degree,
$5\leq d\leq11$. 

\medskip 

$\bullet$ 
The case $d=11$. According to \eqref{mod8restriction}, we
have $\deg(Z)=8k+1 \leq 41$.  This upper bound insures that the sheaf
$\Ee$ in the middle of the extension sequence \eqref{eq:theExtension}
is still Bogomolov unstable, see \eqref{eq:BogomForEV4}.  Thus the
argument used in the case $d\geq12$ applies and we obtain the
decomposition
\[
H=B+R
\]
 as in
\eqref{HBR}. This implies 
\[
  H\cdot B = H\cdot(H-R) = d-H\cdot R \leq d
\] 
 and, hence, by the Hodge index, $B^2\leq d$.  Substituting into
\eqref{eq:ZFuncitonOfBAndZPrime} gives
\begin{equation}
  \label{degZlb3d}
  \deg(Z) = 4d-B^2 +\deg(Z') \geq 3d+\deg(Z') \geq 3d = 33.
\end{equation}
Thus $\deg(Z)=33$ or $41$.  Furthermore, for the first value all the
above inequalities must be equalities and hence $R$ and $Z$ must both
be zero.  It follows that $E=H$ and that the vertical sequence in
\eqref{d:EunstV4} takes the form
\begin{equation}
  \label{destseq1133}
  0 \lra \Oo_X(3H) \lra \Ee \lra \Oo_X(H) \lra 0.
\end{equation} 
We have already encountered a similar situation in the proof of Lemma
\ref{l:cohJZV3} and we use the same argument here.  Namely, recall the
identification
\begin{equation}
  \label{identE-H1133}
  H^0(\Ee(-H)) \cong H^0(\Jj_Z(3H))
\end{equation} 
resulting from the defining extension sequence \eqref{eq:theExtension}
tensored with $\Oo_X(-H)$.  On the other hand, the destabilizing
sequence \eqref{destseq1133} gives
\[
  0 \lra H^0(\Oo_X(2H)) \lra H^0(\Ee(-H)) \lra H^0(\Oo_X). 
\]
This together with the fact that the above inclusion
$H^0(\Oo_X(2H))\hookrightarrow H^0(\Ee(-H))$ must be proper, implies
that the arrow on the right must be onto.  Hence 
$H^0(\Oo_X(2H)) \hookrightarrow H^0(\Ee(-H))$ is a codimension $1$
subspace of $H^0(\Ee(-H))$. This and the isomorphism
\eqref{identE-H1133} tell us that the subspace $e\,H^0(\Oo_X(2H))$,
where $e\in H^0(\Oo_X(H))$ is a section defining the divisor $E$, is a
codimension $1$ subspace of $H^0(\Jj_Z(3H))$.

Next recall that the space $H^0(\Jj_Z(3H))$ contains the
$5$-dimensional subspace $W$, spanned by the restrictions to $X$ of
the partial derivatives of a homogeneous polynomial, call it $f$,
defining the quartic hypersurface $V_4$ that contains $X$.  The above
discussion implies that the subspace $e\,H^0(\Oo_X(2H))$ intersect $W$
along a subspace of codimension at most $1$.  This means that we can
choose homogeneous coordinate functions $X_i$, $i=0,\ldots,4$, on
$\PP^4$ so that
\[
  \frac{\partial f}{\partial X_i}
  = h\gamma_i, \textq{for} i=0,\ldots,3,
\]
where $\gamma_i \in H^0(\Oo_{\PP^4}(2))$ and 
$h \in H^0(\Oo_{\PP^4}(1))$ is the linear form corresponding to the
section $e$ under the identification 
$H^0(\Oo_{\PP^4}(1)) \cong H^0(\Oo_X(H))$.  From this it follows that
$V_4$ is singular along the subvariety 
$\big(h=\frac{\partial f}{\partial X_4}=0 \big)$ which contradicts the
assumption that $V_4$ has only isolated singularities.

Notice that the above argument remains valid as long as the space
$H^0(\Jj_{Z'}(R))$---where the cokernel of 
$H^0(\Oo_X(A-H)\to H^0(\Ee(-H))$ lives---is $1$-dimensional.  This
will be our tool to rule out all the remaining cases whenever the
sheaf $\Ee$ in \eqref{eq:theExtension} is Bogomolov unstable.  We give
a full account of this for $\deg(Z)=41$.

From the identity
\[
  B^2 = 4d-\deg(Z) +\deg(Z')
  = 44-41 +\deg(Z') = 3+\deg(Z') 
\]
we deduce the inequality $B^2\geq 3$.  This and the Hodge index give
$H\cdot B\geq 6$ or, equivalently,
\[
  H\cdot R = H\cdot (H-B) = d-H\cdot B \leq 5.
\]
From the inequality \eqref{Z0compiV4} we also obtain the lower bound
\[
  H\cdot R \geq \frac{\deg(Z)}{3} -d = \frac{41}{3}-11,
\]
\ie $H\cdot R \geq 3$, and proceed according to the possible values
of 
\[
  H\cdot R = 3, 4, \text{ or } 5.
\]

As we have said above, the main idea is, as in the case of
$\deg(Z)=33$, to show that the space $H^0(\Jj_{Z'}(R))$ is
$1$-dimensional.  For this we will also need the formula
\begin{equation}
  \label{ZprimeR1141}
  \deg(Z')
  = B^2-3 = (H-R)^2-3 = d-3-2H\cdot R+R^2 = 8-2H\cdot R+R^2
\end{equation}
which relates $\deg(Z')$ to the intersection numbers $H\cdot R$ and
$R^2$. 

\medskip

1) $H\cdot R=3$.  We wish to analyse the possibility
$h^0(\Oo_X(R))\geq 2$.  Let $R_0$ (resp. $R'$) be the fixed
(resp. moving) part of the linear system $|R|$.  If $R_0\neq 0$, then
$H\cdot R'\leq 2$ and hence, by Hodge index, $R'^2\leq 0$. Since the
linear system $|R'|$ has at most a $0$-dimensional base locus, we
deduce the equality $R'^2=0$.  But this means that $X$ is either a
scroll or a conic bundle and, in view of Remark
\ref{r:rulingDegreeV4}, neither possibility is allowed.  Thus we may
assume that the linear system $|R|$ has at most a $0$-dimensional base
locus.  This and the Hodge index imply $R^2=0$.  Hence $|R|$ is base
point free.  At the same time, observe that
$h^0(\Oo_X(B))=h^0(\Oo_X(H-R))\geq 1$, \ie $B$ is effective and
$h^0(\Oo_X(H-B))=h^0(\Oo_X(R))\geq 2$.  Since $B$ can not be a line,
we deduce that if $R$ moves on $X$, then $h^0(\Oo_X(R))=2$. From the
formula \eqref{ZprimeR1141}, we also find that $\deg(Z')=2$.  Hence
\[
  h^0(\Jj_{Z'}(R)) \leq h^0(\Oo_X(R))-1 = 1
\]
and from here on we conclude as in the case $\deg(Z)=33$.

\medskip

2) $H\cdot R=4$.  The argument follows the same pattern as in the
previous case.  Write $R=R_0 +R'$ as above. If $R_0 \neq 0$, then
$H\cdot R' \leq 3$ and we deduce that $R'^2 = 0$.  Hence, $|R'|$ is a
base point free pencil inducing the morphism
\begin{equation}
  \label{f1141}
  f: X \lra \PP^1.
\end{equation}
The fibres of $f$ must be connected (otherwise we are back to the
situation of $X$ being a scroll or a conic bundle) with $H\cdot R'=3$
(resp. $H\cdot R_0=1$) and hence, $f$ is an elliptic fibration%
\footnote{The other possibility is that $f$ has rational fibres, but
  then $X$ is rational.}  with fibres embedded by $\Oo_X(H)$ as
plane curves of degree $3$. The union of the planes spanned by the
fibres of $f$ form a hypersurface, call it $V$. The degree of this
hypersurface is subject to the inequality
\[
  \deg(V) \leq \deg(f_\ast(\Oo_X(H)))
\]
Setting $f_\ast(\Oo_X(H))=\bigoplus^3_{i=1}\Oo_{\PP^1}(a_i)$, we
obtain
\[
  \deg(V)
  \leq \sum^3_{i=1}a_i
  = h^0(f_\ast(\Oo_X(H)))-3
  = h^0(\Oo_X(H))-3 = 5-3 = 2,
\]
contradicting the assumption on the smallest degree of a hypersurface
containing $X$.

We have just shown that the linear system $|R|$ has at most a
$0$-dimensional base locus.  By the Hodge index $R^2\leq 1$.  Hence
either $R^2=0$ or $R^2=1$.

If $R^2=0$, then $|R|$ is base point free and induces a morphism as in
\eqref{f1141}.  We may again assume that the fibres of this morphism
are connected, since otherwise $X$ is either a scroll or a conic
fibration.  In view of the degree $H\cdot R=4$ of $R$, the
morphism $f:X\to\PP^1$ is either an elliptic fibration or a fibration
by plane curves of degree $4$.  Hence
 \begin{equation}
  \label{KdotR1141}
  K_X\cdot R = 0 \textq{or} K_X\cdot R = 4.
\end{equation}
The ndp formula tells us that $K_X\cdot H=9$ and together with
\eqref{KdotR1141} implies
\[
  K_X\cdot B = K_X\cdot(H-R) = 9 \textq{or} 5.
\] 
Since $B^2 =(H-R)^2 =11-2H\cdot R =11-2\cdot 4=3$, we obtain 
\begin{equation}
  \label{B1141}
  B^2+K_X\cdot B = 12 \textq{or} 8
\end{equation}
which is the degree of the dualizing sheaf of $B$.  But $B$ is as a
divisor of degree $d_B=H\cdot B=7$ contained in a plane%
\footnote{This follows from $h^0(\Oo_x(H-B))=h^0(\Oo_X(R))=2$.}.  
Hence the degree of its dualizing sheaf verifies 
\[
  d_B(d_B-3) = 7\cdot4 = 28,
\]
which does not match any of the values in \eqref{B1141}.

If $R^2=1$, then $|R|$ must have a unique base point and we may assume
that a general member of the linear system is either an elliptic curve
or a smooth plane curve of degree $4$.  These possibilities lead to
$K_X\cdot R =-1$ and $K_X\cdot R=3$ respectively.  Together with the
identity $K_X\cdot H=9$ (the ndp formula), they imply
\[
  K_X\cdot B = 10 \textq{(resp. $=6$).}
\]
But $B^2=(H-R)^2 =4$, hence the degree of the dualizing sheaf of $B$
verifies 
\[
  B^2 + K_X\cdot B = 14 \textq{(resp. $=10$)}
\]
and, as in the case $R^2=0$, neither value agrees with $d_B(d_B-3)=28$.

\medskip

3) $H\cdot R =5$. Then $H\cdot B=6$, $B^2=3$, and $R^2=2$.  Writing
$R=R_0+R'$ as before, we see that $R_0\neq 0$ reduces to the cases
previously considered.  So we may assume that the linear system $|R|$
has at most a $0$ dimensional base locus.  Furthermore, from $R^2=2$
it follows that a general member of $|R|$ is irreducible.  If, in
addition, a general curve in the linear system is not contained in a
hyperplane, then by the Castelnuovo upper bound on the genus of curves
in $\PP^4$, a general member of $|R|$ is a smooth curve of genus $1$.
This and the adjunction for $R$ give $K_X\cdot R=-R^2=-2$ implying
that $X$ is birational to a ruled surface with $q=1$.  From the ndp
formula in \eqref{eq:ndp4}, we have $H\cdot K_X=9$ and obtain 
$K_X^2 =-17$.

Next we consider the line bundle $\Oo_X(K_X+2R)$.  From the
Riemann-Roch,
\[
 h^0(\Oo_X(K_X+2R)) = (K_X+2R)\cdot R = R^2 = 2.
\]  
On the other hand, $(K_X+2R)^2 =K^2_X+4K_X\cdot R +4R^2=K^2_X=-17$.
Hence the linear system $|K_X+2R|$ has a fixed part $F\neq 0$.  In
particular, there is a dense open subset $F'$ of $F$ such that through
every point $x\in F'$ passes an irreducible curve $R_x\in|R|$.  But
$x$ is a base point of $\Oo_X(K_X+2R)$ and by \cite[Theorem 1]{Re},
there is an irreducible curve $C \subset X$ passing through $x$ with
the property $C\cdot R_x =C\cdot R =0$.  Hence the curves $C$ and
$R_x$ have a common component.  Since both curves are irreducible, we
obtain $C=R_x$ and, therefore, $0=C\cdot R_x=R^2_x =R^2=2$, an obvious
contradiction.

Now we know that a general curve in the linear system $|R|$ is
contained in a hyperplane, hence that
$H^0(\Oo_X(B))=H^0(\Oo_X(H-R))\neq 0$.  This tells us that $B$ is an
effective divisor.  Furthermore, $h^0(\Oo_X(H-B))=h^0(\Oo_X(R))\geq 2$
tells us that equality must hold, \ie $B$ is a plane curve subject to 
$d_B=H\cdot B=6$ and $B^2=3$.  Hence
\[
  B^2 +K_X\cdot B = d_B(d_B-3) = 18.
\]
It follows that $K_X\cdot B=18-B^2=15$.  But then  
\[
  15 = K_X\cdot B = K_X\cdot (H-R)
  = K_X\cdot H - K_X\cdot R \leq 9-(-2-R^2) = 9+4 = 13
\]
gives an obvious contradiction.  This completes the case $d=11$. 

\medskip

$\bullet$ The case $d=10$.  
 As in the case $d=11$, we obtain two
possible values for $\deg(Z)$:  $36$ and $44$.  For the first value,
the sheaf $\Ee$ in \eqref{eq:theExtension} is Bogomolov unstable.  The
considerations are the same as in the case $d=11$ and we leave the
details to the reader.  For the latter, \ie if $\deg(Z)=44$, then we
follow the plan outlined in the paragraph prior to the case $d=11$.

We begin by showing that $X$ is birational to a ruled surface.  The
ndp formula reads
\[
  H^2 + H\cdot K_X = 14,
\]
hence $H\cdot K_X=4$.  Combining this with the Hodge Index, gives
$K^2_X \leq 1$.  Substituting this upper bound into the double point
formula yields
\[
  -20 = 2K^2_X-12\chi \leq 2-12\chi.
\]
Hence $\chi\leq 1$ and we obtain
\begin{equation}
  \label{K210}
  K^2_X = -10+6\chi \leq -4.
\end{equation}
We are now ready to show that $X$ is birational to a ruled surface.
Since we are assuming that $X$ is irregular, it is enough to show that
the Kodaira dimension of $X$ is negative.  Assume the opposite and let
$X_0$ be the minimal model of $X$ with $\sigma:X\to X_0$ the
sequence of blow down maps.  We write
\[
K_X =\sigma^\ast K_0 +\Delta,
\]
where $K_0$ is the canonical divisor of $X_0$ and $\Delta$ is the
exceptional divisor composed of the blown-down curves.  Since
$K^2_0\geq 0$, the estimate in \eqref{K210} tells us that $\sigma$ is
the composition of at least four blow-downs.  It follows that $\Delta$
has at least four irreducible components, \ie that 
$H\cdot\Delta\geq 4$.  This gives
\begin{equation}
  \label{HK010}
  4 = H\cdot K_X
  = H\cdot\sigma^\ast K_0 + H\cdot \Delta
  \geq H\cdot\sigma^\ast K_0 +4
\end{equation}
or, equivalently, $H\cdot\sigma^\ast K_0\leq 0$.  Since 
$\sigma^\ast K_0$ is nef, it follows that $H\cdot\sigma^\ast K_0=0$
and $X$ is of Kodaira dimension zero.  By the Enriques-Kodaira
classification, we must have $\chi=0$.  Returning to the equality in
\eqref{K210}, we obtain $K^2_X=-10$.  But then $\Delta$ has at least
ten irreducible components and the estimate in \eqref{HK010} gives
$H\cdot\sigma^\ast K_0\leq -6$ which is impossible.  Hence $X$ is
birational to a ruled surface and $\chi(\Oo_X)=1-q$, which,
substituted into the equality \eqref{K210}, gives
\begin{equation}
  \label{eq:K2_d=10}
  K^2_X = -10+6\chi(\Oo_X) = -10-6(q-1).
\end{equation}

We pursue by studying the adjoint divisor $K_X+H$.  We have, using
\eqref{eq:K2_d=10}, 
\[
  (K_X+H)^2 = K^2_X +2H\cdot K_X +H^2
  = -10-6(q-1) +2\cdot 4 + 10 = 8-6(q-1)
\]
and recalling that $\Oo_X(K_X+H)$ is base point free, see Proposition
\ref{p:2KPlusH}, we deduce that either $q=1$ or $q=2$.

If $q=2$, then $(K_X+H)^2=2$.  By Riemann-Roch, we have 
\begin{equation*}
  h^0(\Oo_X(K_X+H)) = \frac{(K_X+H)\cdot H}{2}+\chi(\Oo_X) = 7-1 = 6.
\end{equation*}
But, according to Proposition \ref{p:2KPlusH}, the linear system
$|K_X+H|$ is base point free and hence defines a morphism $X\to\PP^5$
whose image is a surface of degree at most $(K_X+H)^2=2$ and this is
impossible.

If $q=1$, then $K^2_X =-10$ and $(K_X+H)^2=8$.  The latter and
Proposition \ref{p:2KPlusH} tell us that $|K_X+H|$ is base point free
and defines a morphism
\begin{equation}
  \label{adjmor1044}
  f: X \lra \PP^6
\end{equation}
which must be birational onto its image.  We also record the
projection morphism
\begin{equation}
  \label{projmor1044}
  \pi: X \lra B
\end{equation}
onto an elliptic curve $B$ with rational fibres.  In particular, we
wish to understand the degree of the fibres of $\pi$ with respect to
$H$.  For this we look at the double adjoint line bundle
$\Oo_X(2K_X+H)$.  By Riemann-Roch
\[
  h^0(\Oo_X(2K_X+H)) = \frac{(2K_X+H)\cdot(K_X+H)}{2}
  = \frac{K_X\cdot (K_X+H)+(K_X+H)^2}{2} = 1.
\]
Let $D$ be the divisor defined by a nonzero section of
$\Oo_X(2K_X+H)$.  It is a non-zero divisor, since
$D^2=(2K_X+H)^2=-24$.  Furthermore, since $(K_X+H)\cdot D=2$, and
since $|K_X+H|$ is base point free, an irreducible component $C$ of
$D$ verifies one of the following possibilities:
\begin{equation}
  \label{irrcomD1044}
  \textq{\small (i)} (K_X+H)\cdot C = 0,
  \quad
  \textq{\small (ii)} (K_X+H)\cdot C = 1,
  \quad
  \textq{\small (iii)} (K_X+H)\cdot C = 2.
\end{equation}

The curves $C$ of type (i) are precisely the curves contracted by the
morphism $f$ in \eqref{adjmor1044}.  From $(K_X+H)\cdot C=0$ it
follows that $C$ must be a smooth rational curve with 
$C^2 =H\cdot C -2$.  This implies that $H\cdot C=1$ or $2$, since 
$C^2 \leq 0$.  The second value means that $X$ is a conic bundle which
is impossible by Remark \ref{r:rulingDegreeV4}.  Thus the curves of
type (i) are $(-1)$-curves on $X$, contracted by $\Oo_X(K_X+H)$, and
mapped onto lines by $\Oo_X(H)$ in the embedding $X\subset \PP^4$.  In
particular, all these irreducible components of $D$ are contained in
the fibres of the projection $\pi:X\to B$.

The curves $C$ of type (ii) are smooth rational curves (they are lines
with respect to the morphism $f$) with $C^2=H\cdot C-3$.  Hence
$H\cdot C=1,2$ or $3$.  The last value is impossible, since otherwise
$3$ is the $H$-degree of all fibres of $\pi$ and they all have
intersection $-1$ with $D$ which is impossible.  Thus the curves of
type (ii) are smooth rational curves on $X$ of self-intersection,
either $-1$, or $-2$.  Moreover, they are conics and lines
respectively in the embedding by $\Oo_X(H)$ and are contained in the
fibres of the projection $\pi: X\to B$.  Also observe that there are
at most $2$ distinct such curves in $D$.

We claim that there is no curve of type (iii).  Indeed, all points of
such a curve $C$ are fixed points of $\Oo_X(2K_X+H)$ and according to
\cite[Theorem 1]{Re}, through every point $x\in C$ passes an
irreducible curve of type (i) or (ii). This is impossible since those
curves are rigid.

From the above analysis of the irreducible components of $D$ it
follows that 
\[
  D = C_1+C_2+ D_0,
\]
where $C_i$, $i=1,2,$ are the curves (not necessarily distinct) of
type (ii) and $D_0$ is the residual part of $D$ composed of curves of
type (i).  In addition, \cite[Theorem 1]{Re} stipulates the existence
of a unique divisor $E_x$ passing through every point $x\in C_i$, 
subject to either
\begin{equation}
  \label{Ex1044_a}\tag{a}
  (K_X+H)\cdot E_x = 0  \textq{and} E_x^2 = -1
\end{equation}
or
\begin{equation}
  \label{Ex1044_b}\tag{b}
  (K_X+H)\cdot E_x = 1  \textq{and} E_x^2 = 0
\end{equation}
We know that the irreducible components of $E_x$ in \eqref{Ex1044_a}
are the $(-1)$-curves of type (i) and no such curve passes through a
general point of $C_i$.  So for all points $x$ in the complement of
some finite subset of $C_i$, the corresponding divisor $E_x$ passing
through $x$ is of type \eqref{Ex1044_b}.  This implies that $E_x$ must
be the fibre of $\pi$ containing $C_i$.  Thus if $F$ is the class of a
fibre of the projection $\pi$ in \eqref{projmor1044}, then
\[
 (K_X+H)\cdot F  =1.
\]
or, equivalently, $H\cdot F =3$.  But then
\[
  D\cdot F =(2K_X+H)\cdot F = 2K_X\cdot F + H\cdot F
  = -4+3 = -1
\]
contradicting the fact that $D$ is effective.  This completes the
proof of the case $d=10$.

\medskip

$\bullet$ The case $d=9$. The possible values for $\deg(Z)$ are $41$
and $33$, where the first (resp. second) value is `stable'
(resp. `unstable').  We treat only the stable possibility, \ie
\[
  \deg(Z)=41.
\]
From the ndp formula we obtain $H\cdot K_X=1$.  This and the Hodge
index imply $K^2_X\leq 0$. From the double point formula we obtain
\[
  0 \geq K^2_X = 6\chi-7
\]
and hence $\chi \leq 1$. Arguing as in the case $d=10$, $\deg(Z)=44$,
we deduce that $X$ has negative Kodaira dimension or, equivalently,
that $X$ is birational to a ruled surface of irregularity $q$.  Thus
$\chi=1-q$ and $K^2_X=-7-6(q-1)$.

We proceed with the study of $\Oo_X(K_X+H)$.  The self-intersection
\[
  (K_X+H)^2 =11+K^2_X =4-6(q-1)
\]
and Proposition \ref{p:2KPlusH} tell us that the only possibility for
the irregularity is $q=1$.  Hence $(K_X+H)^2=4$ and 
$K^2_X =-7$. Computing the genus of a smooth curve in the adjoint
linear system $|K_X+H|$ gives
\[
(K_X+H)^2 +(K_X+H)\cdot K_X =4+1-7=-2,
\]
\ie such a curve is rational.  But then $X$ can not be irregular.

\medskip

$\bullet$ The case $d=8$.  
The possible values for $\deg(Z)$ are 
\[
\deg(Z) =40, \,32, \text{ and } 24,
\] 
where the last one is the only unstable value.  It can be treated
easily as in the case $d=11$, $\deg(Z)=33$, so we turn towards the two
remaining stable values.

From the ndp formula we obtain
\[
  H\cdot K_X = 8-2k,
\]
where $k=4$ or $5$.  The negativity of the Kodaira dimension follows
readily.  Hence $\chi=1-q$ and the double point formula gives
\[
  K_X^2 =-8-6(q-1)-5(4-k),  \textq{where} k=4, \,5.
\]
We now compute the self-intersection
\[
  (K_X+H)^2 =-6(q-1) -(4-k),  \textq{where} k=4, \,5.
\]
This and Proposition \ref{p:2KPlusH}, according to which
$\Oo_X(K_X+H)$ is base point free, give the only possibility $q=1$,
$k=4$, $(K_X+H)^2=0$.  This corresponds to the exceptional possibility
of an ADSR elliptic conic bundle.

$\bullet$ The case $d\leq7$.
All those values are incompatible with $\Oo_X(K_X+H)$ being base point
free and big. Indeed, the last condition stipulates $(K_X+H)^2>0$,
which we rewrite as
\[
  K_X^2 > -d-2H\cdot K_X.
\]
We substitute this into the double point formula and obtain
\[
  d^2-8d+12\chi > H\cdot K_X.
\]
On the other hand, expressing $H\cdot K_X$ from the ndp formula
\eqref{eq:ndp4} gives 
\begin{equation}
  \label{eq:expressionForHK}
  H\cdot K_X = \frac{1}{4}\,(d^2-4d-\deg(Z)).
\end{equation}
Putting it together with the previous inequality, we have
\begin{equation}
  \label{eq:dSmallerThan9}
  \frac{3}{4}\,d^2-7d+\frac{1}{4}\,\deg(Z) +12\chi > 0.
\end{equation}
The left hand side of this inequality is an increasing function of $d$
for $d\geq5$.  So, for $d\leq7$, we have
$\frac{3}{4}\,d^2-7d\leq-\frac{49}{4}$.  Substituting into
\eqref{eq:dSmallerThan9} and using $\deg(Z)\leq 45$ give
\[
  12\chi > \frac{49}{4}-\frac{1}{4}\,\deg(Z) \geq 1.
\]
Hence  
\begin{equation}
  \label{chiatleast1}
  \chi\geq1.
\end{equation} 
Going back to \eqref{eq:expressionForHK} we obtain
\begin{equation}
  \label{HK67}
  H\cdot K_X = \frac{1}{4}\,(d^2-4d-\deg(Z))
  \leq \frac{1}{4}\,(7^2-4\cdot 7 -\deg Z)
  = \frac{1}{4}(21-\deg Z),
\end{equation}
for $d\leq 7$. In particular, if $d\leq 7$ and $\deg(Z)>21$, then
$H\cdot K_X$ is negative.  Hence the Kodaira dimension of $X$ is
negative.  This together with the estimate \eqref{chiatleast1} gives
$1\leq \chi =1-q$ meaning that $q=0$.

If $\deg(Z)\leq 21$, then $\deg(Z)\leq 21<4d$, provided 
$6\leq d\leq7$.  Hence the values of $\deg(Z)\leq 21$ are unstable and
one has the estimate
\[
  \deg(Z) > 3d,\,\, \mbox{for $d=6,7$,}
\]
coming from the first inequality in \eqref{degZlb3d} and the
divisibility condition \eqref{div8}.  This implies the only
possibility: $d=6$ and $\deg(Z)=20$.  Substituting these values into
the first equality in \eqref{HK67} gives
\[
  H\cdot K_X = \frac{1}{4}\,(d^2-4d-\deg(Z))
  = -2.
\]
Hence, as before, $q=0$.

For the remaining value $d=5$, the identity \eqref{eq:expressionForHK}
and the divisibility condition \eqref{div8} give
\[
 H\cdot K_X = \frac{1}{4}\,(5-\deg(Z)) \textq{and} \deg(Z) =8k+1.
\]
Therefore, as before, the Kodaira dimension of $X$ is negative, except
possibly in the case $\deg(Z)=1$.  But then the 
subscheme $Z$ is a single point, call it $p$, and the group
$H^1(\Jj_p(K_X+4H))=0$, since by \cite[Theorem 1]{Re}, the line bundle
$\Oo_X(K_X+4H)$ is very ample. This completes the study of cases
$5\leq d \leq 11$.

\section{Albanese dimension}

In this section we consider the Albanese dimension of surfaces lying
on small degree hypersurfaces.  Since we have already seen that there
is no irregular surface lying on a quadric hypersurface and there is
only the elliptic scroll on a cubic hypersurface, our examination will
focus on hypersurfaces of degree $4$ and $5$.  From now on we set
$m=m_X=4$ or $5$---the smallest degree of a hypersurface containing a
smooth surface $X\subset\PP^4$---and we ask the following questions:
\begin{enumerate}[label={\rm \arabic*)}]
  \listspace
\item
  Can $X$ be of Albanese dimension $2$?
\item
  If the answer to question 1) is affirmative, can the irregularity of
  $X$ be bounded? 
\end{enumerate}

To investigate these questions we assume $X$ of Albanese dimension $2$
and consider the extension sequence \eqref{eq:extSeq0} associated to
the cohomology class $\xi\in H^1(\Theta_X(-K_X-(5-m)H))$ arising from a
hypersurface of degree $m$ containing $X$ (see Lemma \ref{l:EP-coh}).
Observe: the assumption that the Albanese dimension of $X$ is $2$
insures that $X$ has non-negative Kodaira dimension and hence the
class $\xi=\delta_X(s)$ in Lemma \ref{l:EP-coh} is nonzero.  In other
words the sequence \eqref{eq:extSeq0} is non split.

\subsection{The Albanese dimension of $X\subset \PP^4$
  with $m_X =4$}

This subsection is devoted to the proof of the following.

\begin{thm}
  \label{th:albaneseDimension}
If $X\subset\PP^4$ is a smooth surface with $m_X=4$, then the
Albanese dimension of $X$ is at most $1$.
\end{thm}

\proof
Assume that $X$ is subject to the hypothesis of the theorem and that
it has the Albanese dimension $2$.  Since $m=m_X=4$, the exact sequence
\eqref{eq:extSeq0} becomes 
\begin{equation}
  \label{eq:extSeq-m=4}
  0 \lra \Oo_X(-K_X-H) \lra \Tt_\xi \lra \Omega_X \lra 0.
\end{equation}
The assumption on the Albanese dimension means that $\Omega_X$ is
generically generated by global sections.  In particular, $q(X)\geq2$
and $K_X$ is effective.  The last property implies that the
homomorphism
\begin{equation}
  \label{eq:extSeqOnGlobalSections}
  H^0(\Tt_\xi) \lra H^0(\Omega_X)
\end{equation}
induced by the epimorphism in \eqref{eq:extSeq-m=4} is an isomorphism:
the injectivity follows from the obvious vanishing of
$H^0(\Oo_X(-K_X-H))$ and the surjectivity is insured by the following
lemma.

\begin{lem}
  \label{l:KPlusHNef}
If $K_X$ is effective, then the divisor $K_X+H$ is nef and big.
In particular, $H^1(\Oo_X(-K_X-H))=0$.
\end{lem}

\likeproof[Proof of Lemma \ref{l:KPlusHNef}]
Let $C$ be an effective curve such that $(K_X+H)\cdot C\leq 0$.  It
follows that $K_X\cdot C \leq -H\cdot C<0$.  But $K_X$ is effective,
hence $C^2<0$ forcing $C$ to be a $(-1)$-curve, \ie 
$ C^2=K_X\cdot C=-1$.  Since $H\cdot C \geq 1$, we deduce
$(K_X+H)\cdot C\geq 0$ and hence $(K_X+H)\cdot C=0$.  Thus $K_X+H$ is
nef.  This implies the non-negativity $(K_X+H)\cdot K_X \geq 0$.
Hence
\[
  (K_X+H)^2 = (K_X+H)\cdot K_X+(K_X+H)\cdot H 
  \geq (K_X+H)\cdot H \geq H^2 >0.
\] 
\qed

Once the isomorphism \eqref{eq:extSeqOnGlobalSections} has been
established, we define the subsheaf $\Ff \subset\Tt_\xi$ as the
saturation of the subsheaf generated by the global sections of
$\Tt_\xi$.  In particular, the inclusion $\Ff\subset\Tt_\xi$ composed
with the epimorphism in \eqref{eq:extSeq-m=4} gives the morphism
$\phi:\Ff\to\Omega_X$ which is generically surjective.  Hence the rank
of $\Ff$ is at least $2$.  On the other hand \eqref{eq:extSeq-m=4}
tells us that the determinant of $\Tt_\xi$ is $\Oo_X(-H)$.  Hence
$\Ff$ must be of rank $2$ and Lemma~\ref{l:mainTechnical} applies to
give us the decomposition
\begin{equation}
  \label{eq:decompotisionOfK}
   K_X = L+E
\end{equation}
where $L=c_1(\Ff)$ and $E=c_1(\coker(\phi))$ is an effective nonzero
divisor. Furthermore, $L$ is effective as well, since by definition
$H^0(\Ff)\cong H^0(\Omega_X)$ and $\Ff$ is generically generated by
its global sections.  This and Lemma~\ref{l:mainTechnical}, 3), imply
that $L=0$.  Hence $\Ff=\Oo_X\oplus\Oo_X$ and therefore, $q(X)=2$. In
addition, the trivial subsheaf $\Ff$ provides
\begin{equation}
  \label{seq-destV4}
  0 \lra \Oo_X \oplus \Oo_X \lra \Tt_\xi \lra \Jj_Z(-H) \lra 0,
\end{equation}
a destabilizing sequence for the sheaf $\Tt_{\xi}$.  This sequence
will play an important role in a later part of the argument.

In fact, more can be extracted from part 3) of Lemma
\ref{l:mainTechnical}.  Namely, since $L=0$ (and $m_X=4$) that
assertion also tells us that we have a unique nonzero global section
$\tau$ of $\Theta_{\PP^4}\otimes\Oo_X(-H)$ such that its image in
$H^0(\Nn_X(-H))$ is the section $e\cdot s$, where $s$ is the global
section of $\Nn_X(-K_X-H)$ defined by a quartic hypersurface
containing $X$ and $e$ is a global section of $\Oo_X(K_X)$
corresponding to the divisor $E$ in \eqref{eq:decompotisionOfK}.
Furthermore, putting the normal sequence of $X$ and the Euler sequence
of $\Theta_{\PP^4}\otimes\Oo_X$ together and tensoring everything with
$\Oo_X(-H)$, we obtain the diagram
\vspace{16ex}
\begin{equation}
  \label{d:normalEuler}
  \begin{tikzpicture}[overlay,every node/.style={draw=none},
    ->,inner sep=1.1ex]
    \matrix [draw=none,row sep=3.5ex,column sep=4.2ex]
    {
      && \node (-12) {$0$}; \\
      && \node (02) {$\Oo_X(-H)$}; \\
      && \node (12) {$H^0(\Oo_X(H))^\ast\otimes\Oo_X$}; \\
      \node (20) {$0$};
      & \node (21) {$\Theta_X(-H)$};
      & \node (22) {$\Theta_{\PP^4}\otimes\Oo_X(-H)$};
      & \node (23) {$\Nn_X(-H)$};
      & \node (24) {$0$}; \\
      && \node (32) {$0$}; \\      
    };
    \path
      (-12) edge (02) 
      (02) edge (12) 
      (12) edge (22)
      (22) edge (32)

      (12) edge node[above right=-.3ex,scale=.75] {$\eta$} (23)
      
      (20) edge (21)
      (21) edge (22)
      (22) edge (23)
      (23) edge (24)
    ;
  \end{tikzpicture}
  \vspace{16.5ex}
\end{equation}
from which it follows that the global section $\tau$ of
$\Theta_{\PP^4}\otimes\Oo_X(-H)$ constructed above comes from a
unique element $v\in H^0(\Oo_X(H))^\ast$.  This vector has the
following geometric meaning.

\begin{lem}
  \label{l:coneOverS}
Let $V_4$ be a quartic hypersurface containing $X$ and let $[v]$ be
the point of $\PP(H^0(\Oo_X(H))^\ast)$ corresponding to the vector
$v$ above.  Then $V_4$ is a cone with the vertex $[v]$ over a quartic
surface $S$ in $\PP^3$.
\end{lem}

\likeproof[Proof of Lemma \ref{l:coneOverS}]
To establish the asserted relation of $V_4$ with the vector $v$ we
recall that $V_4$ gives rise to the nonzero global section $s$ of
$\Nn_X(-K_X-H)$.  Let $\Sigma$ be the $1$-dimensional part of the zero
locus of this section and denote by $\sigma$ a section of
$\Oo_X(\Sigma)$ corresponding to $\Sigma$.  Then $s=\sigma s'$, where
$s'\in H^0(\Nn_X(-K_X-H-\Sigma))$ and has a $0$-dimensional zero locus
$Z_{s'}$.  Therefore the Koszul sequence of $s'$
\[
  0 \lra \Oo_X \stackrel{s'}{\lra} \Nn_X(-K_X-H-\Sigma)
  \stackrel{s'\wedge}{\lra} \Jj_{Z_{s'}}(3H-K_X-2\Sigma) \lra 0
\]
is exact. Tensoring it with $\Oo_X(K_X+\Sigma)$, we obtain
\begin{equation*}
  \label{Koszul4}
  0 \lra \Oo_X(K_X+\Sigma) \stackrel{s'}{\lra} \Nn_X(-H)
  \stackrel{s'\wedge}{\lra} \Jj_{Z_{s'}}(3H-\Sigma) \lra 0.
\end{equation*}
From this sequence, it follows that the global section 
$es=e\sigma s'\in H^0(\Nn_X(-H))$, which interests us, lies in the
kernel of the homomorphism 
$H^0(\Nn_X(-H))\to H^0(\Jj_{Z_{s'}}(3H-\Sigma))$.  Furthermore, by the
definition of $v\in H^0(\Oo_X(H))^\ast$, we have
\begin{equation*}
  \label{s-v}
  es = \eta(v),
\end{equation*}
where $\eta$ is the slanted arrow in the diagram
\eqref{d:normalEuler}.  On the other hand, the composition
\[
  H^0(\Oo_X(H))^\ast\otimes\Oo_X \stackrel{\eta}{\lra} \Nn_X(-H) 
  \stackrel{s'\wedge\cdot}{\lra} \Jj_{Z_{s'}}(3H-\Sigma)
  \stackrel{\sigma}{\hookrightarrow} \Jj_{Z_{s'}}(3H)
\]
is given by the partial derivatives of a polynomial defining $V_4$.
More precisely, if $F$ is a homogeneous polynomial defining $V_4$,
then the composition $\sigma\circ (s'\wedge\cdot)\circ\eta$ sends
every $w\in H^0(\Oo_X(H))^\ast$ to $\partial_w(F)|_X $, the derivative
of $F$ in the direction of $w$ restricted to $X$.  Evaluating the
composition on the vector $v$, we obtain
\[
  \partial_v(F)|_X = \sigma\cdot ((s'\wedge\cdot)\circ\eta)(v)
  = \sigma (s'\wedge\eta(v)) = \sigma (s'\wedge(es)) 
  = \sigma  s'\wedge(e\sigma s') = 0.
\]
Since $X$ is contained in no hypersurface of degree $3$, the above
implies that $\partial_v(F)=0$ in $\Sym^3 (H^0(\Oo_X(H)))$.
Equivalently, $F\in\Sym^4(\ker(v))$, \ie the quartic hypersurface
$V_4$ is the cone in $\PP(H^0(\Oo_X(H))^\ast)$ with vertex $[v]$ and
base the quartic surface $S$ defined by $F$ in 
$\PP(\ker(v)^\ast) \cong \PP^3$.
\qed

In the course of the proof of the preceding lemma we introduced the
divisorial part $\Sigma$ of the section $s\in H^0(\Nn_X(-K_X-H))$
defined by $V_4$.  Geometrically, $\Sigma$ is the $1$-dimensional part
of the singular locus of $V_4$ contained in $X$.  Our next result
shows

\begin{lem}
  \label{l:1dimPart}
  $\Sigma=0$.
\end{lem}

\likeproof[Proof of Lemma \ref{l:1dimPart}]
Assume $\Sigma$ is not zero.  Then we claim that $\Sigma$ is composed
of $(-1)$-curves which are lines with respect to $H$.  Indeed, let $C$
be a reduced, irreducible component of $\Sigma$ and let $\gamma$ be a
global section defining $C$.  Then $s=s'\gamma$, with 
$s'\in H^0(\Nn_X(-K_X-H-C))$. This section gives the cohomology class
$\delta_X(s')=\xi'\in H^1(\Theta_X(-K_X-H-C))$ which is related to
$\xi$ by the formula $\gamma\xi'=\xi$.  In particular, $\xi'$ is
nonzero and gives rise to a nontrivial extension
\begin{equation}
  \label{eq:extMinusC}
  0 \lra \Oo_X(-K_X-H-C)) \lra \Tt_{\xi'} \lra \Omega_X \lra 0.
\end{equation}
In contrast to the homomorphism induced by the extension
\eqref{eq:extSeq-m=4}, the homomorphism
\begin{equation}
  \label{eq:HZeroExtensionMinusC}
  H^0(\Tt_{\xi'}) \lra H^0(\Omega_X)
\end{equation}
induced by the epimorphism in \eqref{eq:extMinusC} fails to be an
isomorphism, since otherwise the argument of the first part of the
proof gives a decomposition $K_X=L'+E'$, with $L'$ effective, while
Lemma \ref{l:mainTechnical}, 3), applied to ${\xi'=\delta_X(s')}$
yields $-(L'+C)\cdot H \geq 0$, an obvious contradiction.

The failure of the homomorphism \eqref{eq:HZeroExtensionMinusC} to be
an isomorphism implies that $H^1(\Oo_X(-K_X-H-C))\neq0$.  But from
the exact sequence
\[
  0 \lra \Oo_X(-K_X-H-C) \lra \Oo_X(-K_X-H) \lra \Oo_C(-K_X-H)
  \lra 0
\]
and Lemma~\ref{l:KPlusHNef}, it follows
\[
  H^0(\Oo_C(-K_X-H)) \cong H^1(\Oo_X(-K_X-H-C)) \neq 0.
\]
Hence $(K_X+H)\cdot C\leq0$.  In the proof of Lemma~\ref{l:KPlusHNef}
such a curve $C$ is identified as a $(-1)$-curves with 
$H\cdot C=1$.

Next we wish to identify the configuration of the lines composing
$\Sigma$ inside the quartic cone $V_4$.  We begin by observing that no
line in $\Sigma$ can be a ruling of the cone, since otherwise such a
ruling $L$ must be contained in the singular locus of $V_4$ and hence
connect the vertex $[v]$ to a singular point of a base $S$ of the
cone; but then the linear system $|H-L|$ must have $L$ in its base
locus which is absurd since $|H-L|$ is base point free---the linear
system $|H-L|$ corresponds to the projection of $X$ from the line $L$
and this is a morphism $X\to\PP^2$.

Once we know that the lines composing $\Sigma$ are not rulings of
$V_4$, we deduce that their projection from the vertex $[v]$ are lines
contained in the singular locus of a base $S$ of the cone.  Thus each
line $C$ of $\Sigma$ is contained in the plane $P_C$ spanned by $C$
and $[v]$.  The plane $P_C$ is contained in the singular locus of
$V_4$.  In addition, since the lines composing $\Sigma$ are
disjoint\footnote{If there are two intersecting lines $C$ and
$C^{\prime}$ in $\Sigma$, then $(C+C')^2 =0$ and this contradicts the
assumption that $X$ is of Albanese dimension $2$.}, $C$ is the only
component of $\Sigma$ contained in $P_C$.  From this it follows that
the one dimensional part of $X\cap P_C$ is the line $C$.  This means
that the pencil of hyperplanes $\{V_t\}_{t\in\PP^1}$ in $\PP^4$
cutting out the plane $P_C$, intersects $X$ along reducible divisors
$H_t=V_t\cdot X$ containing the line $C$ with multiplicity at least
$2$.  From this it follows that the plane $P_C$ is the embedded
tangent plane of $X$ for every $x\in C$.  But this contradicts the fact
that the Gauss map of $X$ is finite, see \cite{Za}.
\qed

The last lemma implies that the quartic surface $S$ in Lemma
\ref{l:coneOverS} has at most isolated singularities and the
hypersurface $V_4$ is the cone over $S$ with vertex $[v]$.  Our
surface $X$ is a smooth divisor in the cone $V_4$ and the projection
from the vertex $[v]$ defines a morphism
\[
  p_v: X \lra S
\]
onto a normal quartic surface $S\subset\PP^3$.  In particular, the
degree of $p_v$ is given by the formula
\begin{equation*}
  \label{degpv}
  \deg(p_v) =
  \begin{cases}
    \frac{d}{4},  & \text{if } [v]\notin X \\
    \frac{d-1}{4},& \text{if } [v]\in X.
  \end{cases}
\end{equation*}
Since the singular locus of $V_4$ is the union of its rulings joining
$[v]$ to the points of $\Sing(S)$, the singular locus of $S$, and by
Lemma \ref{l:1dimPart} none of these is contained in $X$, we deduce
that the zero-locus $Z_0=(s=0)$ of the section 
$s\in H^0(\Nn_X(-K_X-H))$ defined by $V_4$ is $0$-dimensional
of degree
\begin{equation*}
  \label{degZ0-V4-initial}
  \deg(Z_0) =
  \begin{cases}
    \deg(\Sing(S))\cdot \frac{d}{4}, & \textq{if} [v]\notin X,\\
    \deg(\Sing(S))\cdot \frac{d-1}{4} +1, & \textq{if} [v]\in X.
  \end{cases}
\end{equation*}
The above formulas and the well known fact that 
$\deg(\Sing(S))\leq 16$ imply
\begin{equation}
  \label{degZ0-V4}
  \deg(Z_0) \leq 4d.
\end{equation}

Next we relate the normal sequence of $X$ with the extension sequence
\eqref{eq:extSeq-m=4} defined by 
$\xi=\delta_X(s)\in H^1(\Theta_X(-K_X-H))$.  For this we write the
section $s$ as a morphism of sheaves
\begin{equation}
  \label{s-morphV4}
  \Oo_X(K_X+H) \stackrel{s}{\lra} \Nn_X.
\end{equation}
Applying $\Ext^1(\bullet,\Theta_X)$, we have
\begin{equation}
  \label{ext-hom}
  \Ext^1(\Nn_X,\Theta_X)
  \stackrel{s}{\lra} \Ext^1(\Oo_X(K_X+H),\Theta_X).
\end{equation}
In particular, the extension class $\nnn\in\Ext^1(\Nn_X,\Theta_X)$
corresponding to the normal sequence
\[
  0 \lra \Theta_X \lra \Theta_{\PP^4}|_X \lra \Nn_X \lra 0
\]
goes under the homomorphism in \eqref{ext-hom} to the extension class
$\nnn\cdot s$.  But the coboundary map 
$\delta_X: H^0(\Nn_X(-K_X-H)\to H^1(\Theta_X(-K_X-H))$ is
precisely the cup-product with $\nnn$.  Hence 
$\nnn\cdot s =\delta_X(s)=\xi$.  This means that the morphism in
\eqref{s-morphV4} extends to a morphism of extensions
\[
  \begin{tikzpicture}[every node/.style={draw=none},
    ->,inner sep=1.1ex]
    \matrix [draw=none,row sep=4.5ex,column sep=4.2ex]
    {
      \node (10) {$0$};
      & \node (11) {$\Theta_X$};
      & \node (12) {$\Tt^\ast_\xi$};
      & \node (13) {$\Oo_X(K_X+H)$};
      & \node (14) {$0$}; \\
      \node (20) {$0$};
      & \node (21) {$\Theta_X$};
      & \node (22) {$\Theta_{\PP^4}\otimes\Oo_X$};
      & \node (23) {$\Nn_X$};
      & \node (24) {$0$}; \\
    };
    \path
      (10) edge (11) 
      (11) edge (12) 
      (12) edge (13)
      (13) edge (14)
      (11) edge[-,double distance=.33ex] (21)
      (12) edge (22)
      (13) edge node[right=-.2ex,scale=.75] {$s$} (23)
      
      (20) edge (21)
      (21) edge (22)
      (22) edge (23)
      (23) edge (24)
    ;
  \end{tikzpicture}
  \vspace{-1ex}
\]
where the sequence on the top is the dual of \eqref{eq:extSeq-m=4}.
This can be completed to the following commutative diagram
\vspace{16.5ex}
\begin{equation}
  \label{d:twoextV4}
  \begin{tikzpicture}[overlay,every node/.style={draw=none},
    ->,inner sep=1.1ex]
    \matrix [draw=none,row sep=3.7ex,column sep=4.2ex]
    {
      && \node (02) {$0$};
      & \node (03) {$0$}; \\
      \node (10) {$0$};
      & \node (11) {$\Theta_X$};
      & \node (12) {$\Tt^\ast_\xi$};
      & \node (13) {$\Oo_X(K_X+H)$};
      & \node (14) {$0$}; \\
      \node (20) {$0$};
      & \node (21) {$\Theta_X$};
      & \node (22) {$\Theta_{\PP^4}\otimes\Oo_X$};
      & \node (23) {$\Nn_X$};
      & \node (24) {$0$}; \\
      && \node (32) {$\Jj_{Z_0}(4H)$};
      & \node (33) {$\Jj_{Z_0}(4H)$}; \\
      && \node (42) {$0$};
      & \node (43) {$0$}; \\
    };
    \path
      (10) edge (11) 
      (11) edge (12) 
      (12) edge (13)
      (13) edge (14)
      (02) edge (12)
      (03) edge (13)
      (11) edge[-,double distance=.33ex] (21)
      (12) edge (22)
      (13) edge node[right=-.2ex,scale=.75] {$s$} (23)
      (22) edge (32)
      (23) edge node[right=-.2ex,scale=.75] {$s\wedge$} (33)
      (32) edge (42)
      (33) edge (43)
      
      (20) edge (21)
      (21) edge (22)
      (22) edge (23)
      (23) edge (24)
      (32) edge[-,double distance=.33ex] (33)
    ;
  \end{tikzpicture}
  \vspace{17.5ex}
\end{equation}
This diagram will enable us to control the subscheme $Z$ in the
destabilizing sequence \eqref{seq-destV4} of $\Tt_\xi $.

\begin{lem}
  \label{Z-V4}
The subscheme $Z$ is either $[v]$ or empty. In particular, 
$\deg(Z)\leq 1$.
\end{lem}

\likeproof[Proof of Lemma \ref{Z-V4}]
Dualizing \eqref{seq-destV4} and tensoring it with $\Oo_X(-H)$, we
obtain a section $t\in H^0(\Tt^\ast_\xi(-H))$ whose zero
locus is $Z$.  On the other hand the middle column of the diagram
\eqref{d:twoextV4} tensored with $\Oo_X(-H)$ tells us that $t$ can be
identified with the section 
$\tau\in H^0(\Theta_{\PP^4}\otimes\Oo_X(-H))$ corresponding to the
vector $v$ defining the vertex $[v]$ of $V_4$.  In particular, the
scheme $Z_{\tau}=(\tau=0)\supset Z$.  To understand $Z_{\tau}$ we
consider the contraction with $\tau$
\[
  \contraction \tau : 
  H^0(\Omega_{\PP^4}(2)\otimes\Oo_X) \lra H^0(\Oo_X(H)).
\]
This implies that the image of $\contraction\tau$ is contained in 
$H^0(\Jj_{Z}(H))$.
But 
\[
  H^0(\Omega_{\PP^4}(2)\otimes\Oo_X) \cong \Bigwedge^2H^0(\Oo_X(H))
\]
and under the identification of $\tau$ with $v$, the above
contraction becomes simply the contraction with 
$v \in H^0(\Oo_X(H))^\ast$,
\[
  \contraction v: \Bigwedge^2H^0(\Oo_X(H)) \lra H^0(\Oo_X(H)).
\]
This implies 
\[
  \im(\contraction\tau) = \im(\contraction v) = v^{\perp}, 
\]
\ie $\im(\contraction\tau)$ is identified with the space of
hyperplanes in $\PP^4$ passing through $[v]$.  Hence, 
$Z\subset Z_\tau$ is either $\emptyset$ or $[v]$.
\qed

We now have everything to rule out the existence of $X$.  Indeed, from
the destabilizing sequence of $\Tt_{\xi}$ we have
\[
  \deg(Z) = c_2(\Tt_{\xi}) = c_2(\Tt^\ast_{\xi}).
\]
This and the middle column in \eqref{d:twoextV4} give
\[
  10d = c_2(\Theta_{\PP^4}\otimes\Oo_X)
  = c_2(\Tt^\ast_{\xi}) +\deg(Z_0) +4d
  = \deg(Z) +\deg(Z_0) +4d.
\]
From this and Lemma \ref{Z-V4} it follows
\[
  \deg(Z_0) \geq 6d-1.
\]
Putting together this inequality and the upper bound for $\deg(Z_0)$
in \eqref{degZ0-V4}, we obtain $2d\leq 1$ which is clearly
impossible.
\qed

\subsection{The surfaces $X \subset \PP^4$ with $m_X =5$
  and of Albanese dimension $2$}

We now turn to the consideration of surfaces $X\subset\PP^4$ of
Albanese dimension $2$ and $m_X=5$.  To apply our method in this case
we need the additional assumption of $X$ being minimal.  With this in
mind we proceed as in the case $m_X=4$.  Namely, let $V_5$ be a quintic
hypersurface containing $X$ and let $s$ be a nonzero global section of
$\Nn_X(-K_X)$ defined by $V_5$.  That section is used to obtain the
cohomology class $\xi=\delta_X(s)$, which by Lemma \ref{l:EP-coh} is
nonzero. The class $\xi$ is interpreted as a nontrivial extension
\begin{equation}
  \label{extSeqV5}
    0 \lra \Oo_X(-K_X) \lra \Tt_\xi \lra \Omega_X \lra 0
\end{equation}
which we use to gain an insight into the geometry of $X$ and $V_5$.

\begin{lem}
  \label{l:KX-5}
$X$ is either of general type or an abelian surface of degree $d=10$.
\end{lem}

\proof
The assumption that $X$ is of Albanese dimension $2$ implies that
$K_X$ is effective. The minimality of $X$ ensures that $K_X$ is nef.
Hence either $K^2_X>0$ and $X$ is of general type, or $K^2_X=0$.  In
the latter case, we deduce that $K_X=0$, and by the Enriques-Kodaira
classification, $X$ is an abelian surface.  From this and the double
point formula, we obtain that $X$ is of degree $d=10$.
\qed

Abelian surfaces of degree $10$ are of course Horrocks-Mumford
surfaces and it is well-known that they are contained in a quintic
hypersurface and not in one of a smaller degree.  From now on we
assume that $X$ is of general type.  This assumption implies that the
homomorphism
\begin{equation*}
  \label{hom-sectionsV5}
  H^0(\Tt_\xi) \lra H^0(\Omega_X)
\end{equation*}
induced by the epimorphism in \eqref{extSeqV5} is an isomorphism.
This allows us to define the saturation $\Gg$ of the subsheaf of
$\Tt_\xi$ generated by its global sections.

\begin{lem}
  \label{l:concerningG}
The following possibilities may arise:
\begin{enumerate}[label={\rm \arabic*)}]
  \listspace
\item 
  The rank of $\Gg$ is $3$ and then $\Gg=\Tt_\xi\cong\Oo^{\oplus
    3}_X$.  In particular, the irregularity $q=3$, $\Omega_X$ is
  generated by its global sections and $K^2_X = c_2$.
\item
  The rank of $\Gg$ is $2$. 
\end{enumerate}
\end{lem}

\proof
From the definition of $\Gg$ and the assumption that $X$ is of
Albanese dimension $2$, it follows that the rank of $\Gg$ is at least
$2$.

If $\rank(\Gg)=3$, then one can choose a subspace $V\subset H^0(\Gg)$
of dimension $3$ so that the evaluation morphism $V\otimes\Oo_X\to\Gg$
is generically an isomorphism.  This, followed by the inclusion
$\Gg\hookrightarrow\Tt_\xi$, gives a morphism
\[
  V\otimes\Oo_X \lra \Tt_\xi
\]
which is generically an isomorphism.  Since $\det(\Tt_\xi)=\Oo_X$, the
above morphism must be an isomorphism and we obtain the
identifications
\begin{equation*}
  \label{Txi=trivial}
 V\otimes\Oo_X \cong \Gg = \Tt_\xi.
\end{equation*}
This implies
\[
  H^0(\Omega_X) 
  \cong H^0(\Gg)
  \cong H^0(V\otimes\Oo_X)=V
  \cong \CC^3.
\]
Hence $q=3$ and the exact sequence
 \eqref{extSeqV5} takes the form
\[
  0 \lra \Oo_X(-K_X) \lra V\otimes\Oo_X \lra \Omega_X \lra 0.
\]
In particular, the above sequence implies that the cotangent bundle
$\Omega_X$ is generated by its global sections and the Chern numbers
of $X$ are subject to $K^2_X =c_2$.
\qed

We now investigate the case $\rank(\Gg)=2$.

\begin{lem}
\label{l:Grk2}
If $\rank(\Gg)=2$, then $q=h^0(\Gg)=2$ or $3$.  Furthermore, if $q=3$,
then the following possibilities may arise: 
\begin{enumerate}[label={\rm \arabic*)}]
  \listspace
\item 
  $\det(\Gg)=\Oo_X(H)$.
\item
  $h^0(\det(\Gg))=2$ and $X$ admits a fibration $p:X\to B$ onto a
  smooth curve $B$ of genus $g_B=2$ such that 
  \[
    \det(\Gg) = \Oo_X(p^\ast K_B +R),
  \]
  where $K_B$ is the canonical divisor of $B$ and $R$ is the fixed part
  of the linear system $|\det(\Gg)|$.
\end{enumerate}
\end{lem}

\proof
The condition $\rank(\Gg)=2$ allows us to apply Lemma
\ref{l:mainTechnical} and deduce:
\begin{enumerate}[label={\rm \roman*)}]
  \listspace
\item 
  the decomposition of the canonical divisor
  \[
    K_X = L +E,
  \]
  where $L=c_1(\Gg)$ and $E$ is an effective nonzero divisor supported
  on the cokernel of the morphism $\phi:\Gg\to\Omega_X$,
\item
  a nonzero global section 
  $\tau\in H^0(\Theta_{\PP^4}\otimes\Oo_X(-L))$.
\end{enumerate}

To analyse the situation further, we take a nonzero global section $g$
of $\Gg$ and write the exact sequence
\begin{equation}
  \label{KoszGgV5}
  0 \lra \Oo_X(F) \lra \Gg \lra \Jj_A(L-F) \lra 0,
\end{equation}
associated to $g$.  In this sequence, $F$ is the divisorial part of
the scheme $Z_g=(g=0)$, $A$ is the $0$-dimensional subscheme obtained
after dividing $g$ by a section of $\Oo_X(F)$ corresponding to $F$,
and $\Jj_A$ is the ideal sheaf of $A$.  The assumption that the
Albanese dimension of $X$ is $2$ implies that the homomorphism
\[
  H^0(\Gg) \lra H^0(\Jj_A(L-F))
\]
induced by the epimorphism in \eqref{KoszGgV5} is nonzero.  In
particular,
\[
  h^0(\Oo_X(L-F)) \geq h^0(\Gg) - h^0(\Oo_X(F)) 
  = q - h^0(\Oo_X(F)) \geq 1 
\]
and we have the estimate
\begin{multline}
  \label{h0(L)-q}
  h^0(\Oo_X(L)) 
  \geq h^0(\Oo_X(F)) + h^0(\Oo_X(L-F)) - 1 \\
  \geq h^0(\Oo_X(F)) + q - h^0(\Oo_X(F)) - 1
  = q-1,
\end{multline}
implying $h^0(\Oo_X(L))\geq 2$ if $q\geq 3$. 

From now on we assume $q\geq 3$ and write
\begin{equation}
  \label{LMR}
  L = M+R,
\end{equation}
where $|M|$ and $R$ are the moving and the fixed part of $|L|$
respectively.  We divide our considerations according to the dimension
$h^0(\Oo_X(L))=h^0(\Oo_X(M))$.

\paragraph{Case $h^0(\Oo_X(L))=h^0(\Oo_X(M))\geq 3$} 
The nonzero global section 
$\tau\in H^0(\Theta_{\PP^4}\otimes\Oo_X(-L))$, see the beginning of
the proof, gives rise to a nonzero global section of
$\Theta_{\PP^4}\otimes\Oo_X(-M))$.  Since $M$ is nef and
big\footnote{$|M|$ is not composed of a pencil, since otherwise the
image of the Albanese map is $1$-dimensional.}, the Euler sequence for
$\Theta_{\PP^4}\otimes\Oo_X$ tensored with $\Oo_X(-M)$ implies
$H^0(\Oo_X(H-M))\neq 0$.

We claim that $M=H$.  Indeed, assuming $\Gamma=H-M\neq0$, we have
\[
  H^0(\Oo_X(H-\Gamma)) = H^0(\Oo_X(M)) 
  \cong H^0(\Oo_X(L)).
\]
This implies $h^0(\Oo_X(H-\Gamma)) =h^0(\Oo_X(L))\geq 3$.  Hence
$\Gamma$ must be a line and the last inequality must be an equality.
Since
\[
  H = M+\Gamma = L -R +\Gamma = K_X+\Gamma -E -R,
\]
taking the intersection with $\Gamma$ on both sides, we deduce
\[
  1 = H\cdot\Gamma 
  = (K_X+\Gamma)\cdot\Gamma - (E+R)\cdot\Gamma
  = -2-(E+R)\cdot \Gamma.
\]
From this, it follows
 \[
R\cdot\Gamma=-3-E\cdot\Gamma.
\]  
In addition, if $\Gamma$ is not in $E$, then $E\cdot\Gamma\geq 0$, and
if $\Gamma$ is an irreducible component of $E$, then
$E\cdot\Gamma\geq-2$, see \cite{NaRe}. This and the above identity
imply $R\cdot\Gamma\leq-1$.  Hence $\Gamma$ is a component of $R$ and
we can rewrite the relation \eqref{LMR} as follows
\[
  L = M+R = M+\Gamma+R' = H+R',
\]
where $R'$ is an effective divisor.  But this gives the contradiction
\[
  3 = h^0(\Oo_X(L)) \geq h^0(\Oo_X(H)) = 5. 
\]

Once we know that $M=H$, 
the identity \eqref{LMR} reads $L=H+R$ and we claim that $R=0$.
Indeed, if this is not the case, we take $C$, a reduced, irreducible
 component of $R$ and use the fact that
$H^0(\Theta_{\PP^4}\otimes\Oo_X(-H-C))\neq 0$. This together with
the Euler sequence implies $H^1(\Oo_X(-H-C))\neq 0$.  But from the
exact sequence
\[
  0 \lra \Oo_X(-C) \lra \Oo_X \lra \Oo_C \lra 0
\]
tensored with $\Oo_X(-H)$, we see that 
$H^1(\Oo_X(-H-C))\cong H^0(\Oo_C(-H))=0$.

Next we show that $q=3$.  For this we go back to the exact sequence
\eqref{KoszGgV5} which now takes the form
\begin{equation}
  \label{KoszGgV5two}
  0 \lra \Oo_X(F) \lra \Gg \lra \Jj_A(H-F) \lra 0.
\end{equation}
Furthermore, we may assume that there are two linearly independent
global sections $g$ and $g'$ of $\Gg$ which are proportional, \ie the
exterior product $g\wedge g'$ is zero, viewed as a section of
$\det(\Gg)$ (otherwise we are done by \cite[Lemma 5.4]{NaRe}).  With
such a choice of $g\in H^0(\Gg)$ in constructing the exact sequence
\eqref{KoszGgV5two}, we obtain that the line bundle $\Oo_X(F)$ in that
sequence has $h^0(\Oo_X(F))\geq 2$. From the isomorphism 
$H^0(\Gg) \cong H^0(\Omega_X)$ it also follows that the sections $g$,
$g'$ correspond to two linearly independent holomorphic $1$-forms,
call them $\omega$, $\omega'$, subject to $\omega \wedge \omega'=0$ as
a section of $\Oo_X(K_X)$.  By Castelnuovo--de Franchis theorem, it
follows that $X$ admits a morphism
\begin{equation}
  \label{pXtoBV5}
  p:X\lra B
\end{equation}
onto a smooth curve $B$ of genus $g_B\geq 2$ such that
$\omega=p^\ast(\eta)$ and $\omega'=p^\ast(\eta')$, where 
$\eta,\eta' \in H^0(\Oo_B(K_B))$.  Hence
\begin{equation}
  \label{pKBandF}
  p^\ast(\Oo_B(K_B)) =\Oo_X(F).
\end{equation}

Furthermore, from the assumption on the Albanese dimension of $X$ we
know that $q>g_B$.  Hence
\[
  H^0(\Jj_A(H-F)) \supset \coker(H^0(\Oo_X(F)) \lra H^0(\Gg)) \neq 0.
\]
In particular, $\Gamma=H-F$ is an effective nonzero divisor. From this
it follows that
\[
  h^0(\Oo_X(F)) = h^0(\Oo_X(H-\Gamma)) \leq 3,
\]
with the equality holding if and only if $\Gamma$ is a line in
$\PP^4$.  We claim that the equality is impossible. Indeed, if
$\Gamma$ is a line, then
\[
  |F| = |H-\Gamma|
\]
is base point free and hence, by \eqref{pKBandF}, the divisor $F$ is
composed of the fibres of the morphism $p$ in \eqref{pXtoBV5}.  But
$\Gamma$ is a rational curve and hence, must be contained in a fibre
of that morphism.  Then
\[
  0 = F\cdot\Gamma = (H-\Gamma)\cdot\Gamma = 1-\Gamma^2.   
\]
Since $\Gamma^2<0$, the above identity is an obvious contradiction.  

Thus $h^0(\Oo_X(F))\leq 2$. This and the hypothesis 
$h^0(\Oo_X(F))\geq 2$ imply
\[
h^0(\Oo_X(F))=2.
\]
which, together with \eqref{pKBandF} leads to $g_B=2$.

Next we claim that $h^0(\Oo_X(H-F))=1$.  Indeed, assume
$h^0(\Oo_X(H-F))\geq 2$.  This means that the divisors of the linear
system $|F| =|p^\ast (K_B)|$ are contained in a plane in $\PP^4$.  But
a general divisor of $|F|$ contains two disjoint irreducible curves
and this can not happen for plane curves.

\paragraph{The case $h^0(\Oo_X(L))=h^0(\Oo_X(M)) =2$} 
In this case the estimate \eqref{h0(L)-q} tells us that $q\leq 3$ and
hence, $q=3$, in view of the assumption $q\geq 3$.  Furthermore,
 the natural homomorphism
\[
 \Bigwedge^2 H^0(\Gg) \lra H^0(\det(\Gg)) =  H^0(\Oo_X(L))
\]
has a nontrivial kernel and this implies, as in the previous case,
that the sheaf $\Oo_X(F)$ in \eqref{KoszGgV5} satisfies
$h^0(\Oo_X(F))\geq 2$.  By the hypothesis on the Albanese dimension,
the sheaf $\Jj_A (L-F)$ in that exact sequence must have 
$h^0(\Jj_A (L-F))\geq 1$.  This and the first inequality in
\eqref{h0(L)-q} tells us that $h^0(\Oo_X(F))\leq h^0(\Oo_X(L))=2$.
Hence $h^0(\Oo_X(L)) =h^0(\Oo_X(F))=2$ and $h^0(\Oo (L-F)) =1$.

From the isomorphism
$H^0(\Gg) \cong H^0(\Omega_X)$
we also conclude that $H^0(\Oo_X(F))$ defines a two dimensional
subspace of $H^0(\Omega_X)$ contained in the kernel of the natural
homomorphism
\[
  \Bigwedge^2 H^0(\Omega_X) \lra H^0(\det(\Omega_X))
  = H^0(\Oo_X(K_X)).
\]
By Castelnuovo--de Franchis theorem, this means that $X$ admits a
morphism $p:X\to B$ onto a smooth curve $B$ of genus 
$g_B\geq 2$. This together with the hypothesis on
the Albanese dimension of $X$ imply 
\[
 2\leq g_B < q =3.
\]
Hence $g_B=2$ and $\Oo_X(F) =p^\ast (\Oo_B(K_B))$.  In addition, from
$h^0(\Oo_X(F))=h^0(\Oo_X(L))=2$, we deduce the formula
\[
  L = p^\ast K_B + R,
\]
where $R$ is the fixed part of $|L|$.
\qed

We summarize the above discussion in the following statement.

\begin{thm}
  \label{th:V5Alb2}
If $X\subset \PP^4$ is a minimal surface with $m_X=5$ and of 
Albanese dimension $2$, then its irregularity $q(X)=2$ or $3$.  
Furthermore, if $q(X)=3$, then one of the following possibilities
may occur: 
\begin{enumerate}[label={\rm \arabic*)}]
  \listspace
\item
  The cotangent bundle $\Omega_X$ is generated by its global sections
  and $K^2_X=c_2$,
\item
  The canonical divisor $K_X$ admits the decomposition $K_X =L+E$, with
  $L$ and $E$ effective nonzero divisors.  The decomposition is subject
  to the following properties:
  \begin{enumerate}[label={\rm \alph*)}]
  \item
    The divisor $L$ is either equal to $H$, or $X$ admits a fibration
    $p:X\to B$ over a smooth curve $B$ of genus $2$ and
    $L=p^\ast (K_B)+R$, where $K_B$ is the canonical divisor of $B$
    and $R$ is the fixed part of the linear system $|L|$. 
  \item
    $L\cdot H \leq H^2 =d$ and $K^2_X-c_2 \leq \frac{2}{3}\,L^2$.
  \end{enumerate}
\end{enumerate}
\end{thm}

\proof
The assertions 1) and 2), a), are Lemma \ref{l:concerningG} and Lemma 
\ref{l:Grk2} respectively.  The assertion b) is Lemma 
\ref{l:mainTechnical}, 3) and c) is as in Theorem \ref{th:m5}, ii).
\qed

\appendix

\section{The projective bundle $\PP(\Nn_X(-3H))$ and the 
  embedding $X\subset \PP^4$} 

In this appendix we return to an elliptic scroll $X$ of degree $5$ in
$\PP^4$. In Theorem \ref{th:cubicsSections} we established an
isomorphism between the space of global sections of $\Nn^\ast_X(3H)$
and the space of cubic hypersurfaces $I_X(3)$ containing $X$, see
Section \ref{s:ellscrollSegre} for notation.  The subscheme 
$Z_s =(s=0)$ of zeros of a nonzero global section $s$ of
$\Nn^\ast_X(3H)$ is identified with the scheme-theoretic intersection
of $X$ with the singular locus $\Sing(V_3(s))$ of the cubic
hypersurface $V_3(s)$ corresponding to $s$ under the isomorphism
\begin{equation}
  \label{Ndualsectionscubics}
  H^0(\Nn^\ast_X(3H)) \cong I_X(3).
\end{equation} 
This, in particular, allows for a purely geometric way to recover $X$
from the space $I_X(3)$---a geometric counterpart of a well known
algebraic fact that the homogeneous ideal of $X$ is generated in
degree $3$.

More conceptually, the isomorphism \eqref{Ndualsectionscubics}
suggests a sort of `duality' between $\Nn^\ast_X(3H)$ and $X$ embedded
in $\PP^4$ by $\Oo_X(H)$.  This is the main theme of this section.  Of
course, one aspect of the above mentioned duality is well-known---the
famous quadro-cubic transformation of Cremona relating a normal
elliptic quintic curve in $\PP^4$ with an elliptic scroll (in another
copy of $\PP^4$) see \cite{ADHPR,ADHPR2}.  So we do not claim any
novelty in the results exposed here.  However, placing the vector
bundle $\Nn^\ast_X(3H)$ in the center of the study, revisiting various
aspects of the geometry of the scroll $X$ and of the Segre cubics
containing it via the properties of $\Nn^\ast_X(3H)$, seem to be new
and fruitful.

\medskip

Set $\Ee=\Nn_X(-3H)$. The projectivization $Y=\PP(\Ee)$ with the
structure projection 
\begin{equation}
  \label{eq:pYToX}
  p: Y=\PP(\Ee) \lra X
\end{equation}
is equipped with the line bundle $\Oo_Y(1)$ chosen so that the direct
image
\[
  p_\ast\Oo_Y(1) = \Ee^\ast = \Nn^\ast_X(3H).
\]
By Lemma \ref{l:aboutConormal3H}, the vector bundle $\Ee^\ast$ is
globally generated. Hence $\Oo_Y(1)$ is globally generated and defines
a morphism
\begin{equation}
  \label{phiYtoPW}
  \phi: Y=\PP(\Ee) \lra \PP(H^0(\Ee^\ast)^\ast) 
\end{equation}
where the target space is, in view of Theorem \ref{th:cubicsSections},
2), a $4$-dimensional projective space.  We let 
\[
  V := H^0(\Oo_X(H)) \textq{and}
  W := H^0(\Ee^\ast) = H^0(\Nn^\ast_X(3H))
\]
and want to keep a clear distinction between the following two
geometric incarnations of the projective space $\PP^4$:
\begin{itemize}
  \listspace
\item
  the projective space $\PP(V^\ast)$, where $X$ lives,
\item
  the projective space $\PP(W^\ast)$, the target of the morphism
  $\phi$ defined in \eqref{phiYtoPW}. 
\end{itemize}
One of our goals here is to describe how to go between these two
spaces.

\subsection{On the geometry of the scroll $X$}

To describe the geometry of the morphism $\phi$ (resp. of the embedding
$X\subset \PP(V^\ast)$) we recall that $X$ is a $\PP^1$-bundle over an
elliptic curve which will be denoted by $E$ and we let
\begin{equation}
  \label{piXtoE}
  \pi: X \lra E
\end{equation}
be the structure projection.  It is well known that $X$ can be
identified with $\Sym^2(E)$, the second symmetric power of $E$.  Then
the structure projection becomes the Abel-Jacobi map which takes a
subscheme $D\subset E$ of degree $2$, viewed as a point of
$\Sym^2(E)$, to (the isomorphism class of) the line bundle $\Oo_E(D)$,
viewed as a point of $E$.  We are making here a (non-canonical)
identification of $E$ with $\Jac_2(E)$, the variety of isomorphism
classes of line bundles of degree $2$ on $E$.  From this it follows
that $X$ admits the obvious double covering
\begin{equation}
  \label{sdoublecover}
  \tau : E\times E \lra \Sym^2(E)=X
\end{equation}
which sends a point $(e,e') \in E\times E$ to $\tau(e,e')$, the
subscheme of $E$ of degree $2$ supported on $e$ and $e'$.  This equips
$X$ with a distinguished family of curves 
\begin{equation}
  \label{curvesGp}
  \Gamma_e:= \tau(\nu_1^{-1}(e)) = \tau(\nu_2^{-1}(e)),
\end{equation}
where $e\in E$ and $\nu_i$ denotes the projection of $E\times E$ on the
$i$-th factor, for $i=1,2$.  These curves will play an important role
in the sequel, so we mention some of their properties, immediate
consequences of the definition \eqref{curvesGp}.
\begin{equation}
  \label{eq:GammaFamily}
  \begin{aligned}
    i) & \quad\text{For every $e\in E$, the curve $\Gamma_e$ is a
      section of the projection $\pi$ in \eqref{piXtoE}.}\\
    ii) & \quad \text{$\Gamma_e\cdot \Gamma_{e'} =1$, for every $e,e'
      \in E$.} \\
    iii) & \quad h^0(\Oo_X(\Gamma_e))=1.
  \end{aligned}
\end{equation} 
The curves $\{\Gamma_e\}_{\ e\in E}$ are used to define embeddings of
$X$ into $\PP^4$.

\begin{pro}
  \label{p:embedGpplusD} 
Let $\Oo_E(D)$ be a line bundle of degree $2$ on $E$.  Then for every
point $e\in E$ the line bundle $\Oo_X(\Gamma_e +\pi^\ast(D))$ is very
ample and it defines an embedding of $X$ into $\PP^4$ as a scroll of
degree $5$.
\end{pro}

\proof
Let $l$ be the class of a fibre of $\pi:X\to E$ in the N\'eron-Severi
group of $X$.  The canonical divisor has the form 
$K_X=-2\Gamma_e+l$ and we set
\begin{equation}
  \label{Hadjform}
  H := \Gamma_e +\pi^\ast(D) = K_X+(H-K_X),
\end{equation}
where $H-K_X=3\Gamma_e+l$ is nef and big.  The very ampleness of $H$
now follows easily from \cite[Theorem 1, 2)]{Re}. 

From the formula \eqref{Hadjform}, we see that 
$h^i(\Oo_X(\Gamma_e+\pi^\ast(D))=0$, for $i=1,2$.  Hence by
Riemann-Roch we obtain
\[
  h^0(\Oo_X(\Gamma_e +\pi^\ast(D))
  = \chi(\Oo_X(\Gamma_e +\pi^\ast(D))
  = \frac{1}{2}\left((\Gamma_e +2l)^2
    -(\Gamma_e +2l)\cdot(-2\Gamma_e +l)\right)
  = 5.
\]
Thus $\Oo_X(H)$ embeds $X$ into $\PP^4$.  Moreover, since 
$H\cdot l= (\Gamma_e +2l)\cdot l =1$, it follows that $\Oo_X(H)$
embeds $X$ as a scroll of degree
\[
  H^2 = (\Gamma_e +2l)^2 = \Gamma^2_e +4\Gamma_e\cdot l = 1+4
  = 5.
\]
\qed

From now on we fix a point $o\in E$ and a line bundle $\Oo_E(D)$ of
degree $2$ on $E$.  According to the previous proposition, this gives
an embedding of $X$ into $\PP(V^\ast)$ defined by $\Oo_X(H)$, where
$H:=\Gamma_o+\pi^\ast(D)$ and $V=H^0(\Oo_X(H))$.

\begin{pro}
  \label{p:embedGpplusD2}
Under the embedding of $X$ defined by $\Oo_X(H)$, the curves
$\{\Gamma_e\}_{e\in E}$ are embedded as plane curves of degree $3$.
Conversely, a plane cubic curve contained in $X$ is one of the curves
of the family $\{\Gamma_e\}_{e\in E}$.
\end{pro}

\proof
From
$\Oo_X(H-\Gamma_e)=\Oo_X(\Gamma_o +\pi^\ast(D)-\Gamma_e)=\pi^\ast
\Oo_E(D+(o-e))$ it follows that 
\[
  h^0(\Oo_X(H-\Gamma_e))
  = h^0(\pi^\ast \Oo_E(D+(o-e))) = h^0(\Oo_E(D+(o-e)))
  = 2.
\]
This means that under the embedding by $\Oo_X(H)$, the image of 
$\Gamma_e$ is contained in a plane.  Its degree is computed by the
intersection number 
\[
  H\cdot \Gamma_e = (\Gamma_o+\pi^\ast(D))\cdot \Gamma_e = 1+2 = 3,
\]
where the properties {\it i)} and {\it ii)} in \eqref{eq:GammaFamily}
are used.

Conversely, let $C$ be a plane cubic contained in $X$.  Then
\[
  h^0(\Oo_X(H-C)) = 2 \textq{and} H\cdot(H-C) = 2.
\]
Hence, the effective divisor $H-C$ is composed of rational
curves.  Since the only rational curves on $X$ are the rulings (the
fibres of $\pi$) of $X$, we deduce the linear equivalence 
\[
  H-C \sim \pi^\ast(D'),
\]
for some effective divisor $D'$ of degree $2$ on $E$.  Hence
\[
  C \sim H-\pi^\ast(D') = \Gamma_o +\pi^\ast(D-D') = \Gamma_e,
\]
where $e=o+e'$ and $e'$ is the point of $E$ corresponding to
$\Oo_E(D-D')$ under the identification $E=\Jac_0(E)$.  Since
$h^0(\Oo_X(\Gamma_e))=1$, see the property {\it iii)} in
\eqref{eq:GammaFamily}, it follows that $C=\Gamma_e$.
\qed

\begin{rem}
Viewing $X$ as a subvariety of $\PP(V^\ast)$, one obtains the family
$\{\Gamma_e\}_{e\in E}$ as follows: consider two rulings $l_a$ and
$l_{a'}$ of $X$ lying over the points $a$ and $a'$ of $E$ and consider
the hyperplane $V_{a,a'}$ in $\PP(V^\ast)$ spanned by the rulings
$l_a$ and $l_{a'}$.  The hyperplane section 
$ H_{a,a'}=V_{a,a'}\bigcap X$ is a reducible divisor in $|H|$ of the
form
\[
  H_{a,a'} = C +l_a +l_{a'},
\]
where $C$ is the component of $H_{a,a'}$ of degree $H\cdot C=3$,
residual to $l_a +l_{a'}$.  Hence it must be irreducible and a section
of the projection $\pi:X\to E$.  Therefore, $C$ is a plane cubic
contained in $X$ and, in view of Proposition \ref{p:embedGpplusD2}, it
must be a curve in $\{\Gamma_e\}_{e\in E}$.  The first part of the
proof of Proposition \ref{p:embedGpplusD2} shows that every curve in
$\{\Gamma_e\}_{e \in E}$ is obtained in this way.
\end{rem}

The family of curves $\{\Gamma_e \}_{e\in E}$ gives rise to the family
of planes $\{\Pi_e \}_{e\in E}$, where $\Pi_e$ is the plane spanned by
$\Gamma_e$.  To understand the properties of this family of planes we
have the following.

\begin{pro}
  \label{p:planebundleinPV}
Let $\Ff:= \nu_{1\ast}(\tau^\ast\Oo_X(H))$.  Then $\Ff$ is a rank $3$
vector bundle on $E$ of degree $5$, generated by its global sections
with $h^0(\Ff)=5$.
\end{pro}

\proof
By definition, $\rank(\Ff\otimes\Oo_{E,e})=h^0(\Oo_{\Gamma_e}(H))=3$
for every $e\in E$.  Hence $\Ff$ is locally free of rank $3$.  By
Riemann-Roch
\[
  \deg(\Ff) = \chi(\Ff)
  = \chi({\nu_1}_\ast(\tau^\ast \Oo_X(H)))
  =\chi(\tau^\ast \Oo_X(H)) = H^2 = 5.
\]
Furthermore, 
$H^0(\Ff)\cong H^0(\tau^\ast\Oo_X(H))\cong H^0(\Oo_X(H))$.  Hence
$h^0(\Ff)=h^0(\Oo_X(H))=5$. 

The global generation of $\Ff$ follows from the identifications
\[
  H^0(\Ff) \cong H^0(\Oo_X(H))
  \lra H^0(\Oo_{\Gamma_e}(H)) \cong H^0(\Ff\otimes\Oo_{E,e})
\]
and the fact that, for every $e\in E$, the above restriction
homomorphism is surjective.
\qed

Set $T:=\PP(\Ff^\ast)$ to be the projectivization of the dual of $\Ff$
and $\rho: T\to E$ the structure projection.  Let $\Oo_T(1)$ be chosen
so that $\rho_\ast \Oo_T(1)=\Ff$.  In view of Proposition
\ref{p:planebundleinPV}, the line bundle $\Oo_T(1)$ is generated by
its global sections and hence gives a morphism
\[
  \psi: T \lra \PP(V^\ast)
\]
By definition, the fibres of $T$ are mapped by $\psi$ to the family of
planes $\{\Pi_e\}_{e\in E}$. In particular, we obtain

\begin{pro}
  \label{p:psimap}
The image $T'$ of $\psi$ is a hypersurface of degree $5$ containing
$X$ in its singular locus.
\end{pro}

\proof
From the definition of $\psi$ it follows
\[
  \deg(\psi)\,\deg(T')= c_1^3(\Oo_T(1)) = \deg(\Ff) =5,
\]
where the last equality comes from Proposition
\ref{p:planebundleinPV}.  The above implies $\deg(\psi)=1$ and
$\deg(T')=5$.

The hypersurface $T'$ obviously contains $X$.  To see that $X$ is
contained in its singular locus, we observe that the planes $P_e$ and
$P_{e'}$, for $e\neq e'\in E$, intersect at a single point.  Indeed,
otherwise the hyperplane spanned by $P_e\bigcup P_{e'}$ gives rise to
a divisor $D\in|H|$ of the form
\[
  D = \Gamma_e +\Gamma_{e'} +R
\]
where $R$ is some effective divisor.  The divisor $D$ has a wrong
degree on the rulings of $X$. 

The point of intersection $P_e\bigcap P_{e'}$ is obviously a singular
point of $T'$.  Furthermore, this point is the point of intersection
$\Gamma_e\cdot\Gamma_{e'}$ and hence belongs to $X$. Varying $e$ and
$e'$, the points $\Gamma_e\cdot\Gamma_{e'}$ form a Zariski dense open
subset of $X$.  Therefore, $X$ is contained in the singular locus of
$T'$.
\qed

Besides being the union of the family of planes $\{\Pi_e\}_{e\in E}$, 
the variety $T'$ in Proposition \ref{p:psimap} can be also characterized 
as follows.

\begin{pro}
  \label{p:trisecantstoX}
The variety $T'$ in Proposition \ref{p:psimap} is the union of the
proper trisecant lines of $X\subset\PP(V^\ast)$.
\end{pro}

\proof
Let $\Pi^{\dual}_e $ be the dual of the plane $\Pi_e$. Then the family
of the dual planes $\{\Pi^{\dual}_e\}_{e\in E}$ gives rise to a family
of proper trisecants of $X$ in $\PP(V^\ast) $.  We need to check that
every proper trisecant of $X$ belongs to this family.  Let $l$ be a
proper trisecant line of $X$, \ie $Z_l :=l\cdot X$ is a \linebreak
$0$-dimensional subscheme of $X$ of degree at least $3$.  Then 
\[
  h^1(\Jj_{Z_l}(H)) = \deg(Z_l) -2 \geq 3-2 =1.
\]
This and the Serre duality 
$\Ext^1(\Jj_{Z_l}(H),\Oo_X(K_X))\cong H^1(\Jj_{Z_l}(H))^\ast\neq0$ 
gives rise to a nontrivial extension sequence
\[
  0 \lra \Oo_X(K_X) \lra \Ee \lra \Jj_{Z_l}(H) \lra 0.
\]
The sheaf $\Ee$ in the middle of that sequence is torsion free with
$\det(\Ee) =\Oo_X(H+K_X)$ and 
$h^0(\Ee) \geq h^0(\Jj_{Z_l}(H))-h^1(\Oo_X(K_X)) =3-1=2$. Hence the
saturation of the subsheaf generated by global section of $\Ee$ is a
torsion free subsheaf of $\Ee$ of rank one. Hence it has the form
$\Jj_{Z_1}(A)$, where $A$ (resp. $Z_1$) is an effective divisor
(resp. $0$-dimensional subscheme) on $X$, and gives rise to the
following destabilizing sequence of $\Ee$
\[
0 \lra \Jj_{Z_1}(A) \lra \Ee \lra \Jj_{Z_2}(K_X+H-A) \lra 0.
\]
Combining this together with the defining extension sequence implies
that there is an effective nonzero divisor $\Gamma\in |H-A|$ passing
through $Z_l$. From $h^0(\Oo_X(H-\Gamma)) =h^0(\Oo_X(A))\geq 2$ and
the fact that $\Gamma$ can not be a line, it follows that the equality
must hold and therefore, $\Gamma$ is a plane curve. Hence $\Gamma$ is
a plane cubic section of $X$ and, in view of Proposition
\ref{p:embedGpplusD2}, $\Gamma=\Gamma_e$, for some $e\in E$. This
implies $Z_l =l\cdot \Gamma_e \subset \Pi_e$ and hence the line $l$
correspond to a point in the dual plane $\Pi^{\dual}_e$.
\qed

\subsection{From the morphism $\phi$ in \eqref{phiYtoPW} to
  $X\subset \PP(V^\ast)$} 

The above discussion gives a clear geometric picture on the side of
the embedding $X\subset \PP(V^\ast)$. We now turn to the other side,
the projectivization $\PP(\Ee)$ and the morphism
\[ 
\phi: Y =\PP(\Ee) \lra \PP(W^\ast)
\]
introduced in \eqref{phiYtoPW}, where $W=H^0(\Ee^\ast)$.

\begin{pro}
  \label{p:descriptionOfY0}
For $\phi: Y=\PP(\Ee)\to\PP(W^\ast)$, the following statements hold.
\begin{enumerate}[label={\rm \arabic*)}]
\item 
The $3$-fold $Y$ contains a distinguished divisor
$\Sigma\cong E\times E$ and
\[
  p|_{\Sigma} : E\times E \cong \Sigma \lra X \cong \Sym^2(E),
\]
the restriction to $\Sigma$ of the structure projection $p$ in
\eqref{eq:pYToX}, is the double covering $\tau$ defined in
\eqref{sdoublecover}.
\item
Let $\Delta_E$ denote the diagonal of $E\times E$. Then the image
$\phi(E\times E)=\phi(\Delta_E)$ is an embedding of $E$ as an elliptic
normal curve of degree $5$.
\item
The image $Y_0=\phi(Y)\subset \PP(W^\ast)$ is a quintic hypersurface
which is the secant variety of the elliptic normal curve 
$\phi(\Delta_E)$.
\item
The composition $\pi\circ p: Y\to E$ is a fibration whose fibres
$(\pi\circ p)^{-1}(e) =S_e$ are isomorphic to the Hirzebruch surface
$\Sigma_1 =\PP(\Oo_{\PP^1}\oplus\Oo_{\PP^1}(-1))$.  For every 
$e\in E$, the restriction $\phi|_{F_e} :S_e\to\PP(W^\ast)$ is an
embedding.  The image $S'_e :=\phi(S_e)$ is a rational cubic scroll
containing the elliptic curve $\phi(\Delta_E)$.
This curve is the image under $\phi$ of
the intersection $\Sigma \cdot S_e$.  The intersection, viewed as a
divisor in $ S_e$, has the form $2L_e +3f$, where $L_e$ is a unique
section of $S_e$ with $L^2_e=-1$ and $f$ is the class of a fibre of
$p$.  In particular, the image $L'_e:=\phi(L_e)$ is a line
intersecting $\phi(\Delta_E)$ transversely at a single point.
\item 
The projection $p: Y\to X$ admits a distinguished section 
$\gamma: X\to Y$.  The image of the composition 
$\phi\circ\gamma: X\to\PP(W^\ast)$ is a singular scroll of degree
$15$ containing $\phi(\Delta_E)$.
\end{enumerate}
\end{pro}

\noindent
For the convenience of the reader, the following diagram summarizes
all the morphisms appearing in (the proof of) the proposition and we
suggest to consult this diagram while advancing through the proof.
\[
  \begin{tikzpicture}[every node/.style={draw=none},
    ->,inner sep=1.1ex]
    \matrix [draw=none,row sep=4ex,column sep=4.2ex]
    {
      & \node (01) {$\Delta_E$}; \\
      & \node (11) {$E\times E$};
      & \node (12) {$\PP(\Ee)=Y$};
      & \node (13) {$\PP(W^\ast)$}; \\
      \node (20) {$\PP(V^\ast)$};
      & \node (21) {$X$}; \\
      & \node (31) {$E$}; \\
    };
    \path
      (11) edge[right hook->] node[above=-.1ex,scale=.75] {$i$} (12) 
      (12) edge node[below=-.1ex,scale=.75] {$\phi$} (13)
      (01) edge[right hook->] (11)
      (11) edge node[left=-.2ex,scale=.75] {$\tau$} (21)
      (21) edge node[left=-.2ex,scale=.75] {$\pi$} (31)
     
      (01) edge[bend left=15pt]
        node[above=.5ex,scale=.75] {$\mu=\phi|_{\Delta_E}$} (13)
      (12) edge node[below right=-.3ex,scale=.75] {$p$} (21)
      (21) edge[left hook->] (20)
    ;
  \end{tikzpicture}
\]

\proof
We begin by showing that the natural double covering
\[
  \tau : E\times E \lra \Sym^2(E) \cong X
\]
admits a lifting to $Y=\PP(\Ee)$, \ie that there is a morphism 
$i: E\times E\to Y$ which fits into the following commutative diagram
\vspace{5ex}
\begin{equation}
  \label{d:liftE2}
    \begin{tikzpicture}[overlay,every node/.style={draw=none},
    ->,inner sep=1.1ex]
    \matrix [draw=none,row sep=4ex,column sep=4.2ex]
    {
      & \node (12) {$Y=\PP(\Ee)$}; \\
      \node (21) {$E\times E$};
      & \node (22) {$X$}; \\
    };
    \path 
      (21) edge[right hook->] node[above left=-.3ex,scale=.75] {$i$} (12)
      (12) edge node[right=-.3ex,scale=.75] {$p$} (22)
      (21) edge node[above=-.2ex,scale=.75] {$\tau$} (22)
    ;
  \end{tikzpicture}
  \vspace{5.5ex}
\end{equation}
It is well known that this amounts to having a line bundle $\Mm$ on
$E\times E$ together with a surjective morphism 
\[
  \tau^\ast(\Ee^\ast) \lra \Mm,
\]
see \eg \cite{Ha}, Ch. II, Proposition 7.12.  In order to construct 
$\Mm$ and the morphism, we investigate the
restriction of $\Ee$ to the curves $\{\Gamma_e\}_{e\in E}$ in
\eqref{curvesGp}.

\begin{lem}
  \label{l:ErestGe}
For every $e\in E$, the restriction of $\Ee$ to $\Gamma_e$ has the
form
\[
  \Ee\otimes\Oo_{\Gamma_e} = \Nn_X(-3H)\otimes\Oo_{\Gamma_e} 
  = \Oo_{\Gamma_e}\oplus\Oo_{\Gamma_e}(K_X-H).
\]
\end{lem}

\likeproof[Proof of Lemma \ref{l:ErestGe}]
Set $\Gamma=\Gamma_e$.  In view of Proposition \ref{p:embedGpplusD2},
the curve $\Gamma$ is embedded into $\PP(V^\ast)$ as a plane curve of
degree $3$.  Hence its normal bundle $\Nn_{\Gamma/\PP(V^\ast)}$ in
$\PP(V^\ast)$ has the form
\[
  \Nn_{\Gamma/\PP^4}
  \cong \Oo_\Gamma(H)\oplus\Oo_\Gamma(H)\oplus\Oo_\Gamma(3H).
\]
Using this decomposition in the short exact sequence of normal bundles
of the inclusions $\Gamma\subset X\subset\PP(V^\ast)$, we obtain
\[
  0 \lra \Oo_\Gamma(\Gamma)
  \lra \Oo_\Gamma(H)\oplus\Oo_\Gamma(H)\oplus\Oo_\Gamma(3H)
  \lra \Nn_X\otimes\Oo_{\Gamma} \lra 0.
\]
This gives a nonzero morphism
$\Oo_\Gamma(3H)\to\Nn_X\otimes\Oo_{\Gamma}$ or, equivalently, a
nonzero global section, call it $\eta$, of
$\Nn_X(-3H)\otimes\Oo_{\Gamma}$.  This together with the fact that the
dual vector bundle $\Nn^\ast_X(3H)\otimes\Oo_{\Gamma}$ is globally
generated, see Theorem \ref{th:cubicsSections}, implies that $\eta$ is
nowhere vanishing and gives a splitting
\[
  \Ee\otimes\Oo_{\Gamma} = \Nn_X(-3H)\otimes\Oo_{\Gamma}
  = \Oo_{\Gamma}\oplus\Oo_{\Gamma}(K_X-H)
\]
as asserted.
\qed

We now take the pullback $\tau^\ast(\Ee)$ and recall that the curve
$\Gamma_e$, for every $e\in E$, is the image under $\tau$ of the fibre
$\nu_1^{-1}(e)$ of the projection $\nu_1 :E\times E\to E$.  This and
the decomposition of Lemma \ref{l:ErestGe} imply that the direct image
${\nu_1}_\ast (\tau^\ast(\Ee))$ is a line bundle, call it $\Ll$.  This
can be interpreted as a nowhere vanishing section of
${\nu_1}_\ast(\nu_1^\ast(\Ll^{-1})\otimes \tau^\ast(\Ee))$ which in
turn gives a nowhere vanishing global section of
$\nu^\ast_1(\Ll^{-1})\otimes \tau^\ast(\Ee)$ or, equivalently, the
inclusion of bundles
\[
  \nu^\ast_1(\Ll) \hookrightarrow \tau^\ast(\Ee).
\]
Dualizing gives the sought after epimorphism
\begin{equation}
  \label{surjfori}
  \tau^\ast(\Ee^\ast) \lra \nu^\ast_1(\Ll^{-1}) \lra 0
\end{equation}
and hence the diagram \eqref{d:liftE2}.  Furthermore, under the
morphism $i:E\times E\to Y$ associated to the above epimorphism, the
line bundle $\nu^\ast_1(\Ll^{-1})$ becomes the pullback 
$i^\ast(\Oo_Y (1))$.  Hence the composition morphism
\[
  \phi \circ i: E\times E \lra \PP(W^\ast)
\]
factors through the morphism 
\begin{equation}
  \label{mumorph}
  \mu: E \lra \PP(W^\ast),
\end{equation}
defined by $\Ll^{-1}$.  Since the projection $\nu_1$ restricted to the
diagonal $\Delta_E$ induces the isomorphism of $\Delta_E$ with the
base $E$ (of the first projection), we deduce
\[
  (\phi\circ i)(E\times E) = (\phi\circ i)(\Delta_E).
\]

Geometrically, the morphism $i$ in \eqref{d:liftE2}, corresponding to
the epimorphism \eqref{surjfori}, maps the fibre 
$pr^{-1}_1(e) \cong \Gamma_e$ to the section of the ruled surface
$\PP(\Ee\otimes\Oo_{\Gamma_e})$ defined by the trivial summand of the
decomposition 
$\Ee\otimes\Oo_{\Gamma_e}=\Oo_{\Gamma_e}\oplus\Oo_{\Gamma_e}(K_X-H)$
described in Lemma \ref{l:ErestGe}. Under the morphism $\phi$ that
section is contracted to the point $\mu (e)$, while the ruled surface
$\PP(\Ee\otimes\Oo_{\Gamma_e})$ is mapped by $\phi$ to the cone over a
normal quartic elliptic curve%
\footnote{The base of the cone is the image of $\Gamma_e$ by
  $\Oo_{\Gamma_e}(H-K_X)$ which is a line bundle of degree
  $(H-K_X)\cdot \Gamma_e =3+1=4$.}
with vertex at $\mu (e)$.  From this
and the properties of $\Nn^\ast(3H)$ we can draw several conclusions.
\begin{enumerate}[label={\rm \alph*)}]
  \listspace
\item
The image $\mu(E)$ of the morphism $\mu$ in \eqref{mumorph} can not be
contained in a plane.  Indeed, otherwise a plane containing $\mu(E)$
gives rise to a pair of linearly independent global sections of
$\Ee^\ast=\Nn^\ast(3H)$ which are proportional, but such sections
must have $1$-dimensional zero locus and this, according to Theorem
\ref{th:cubicsSections}, 3), is impossible.
\item
The degree of $\mu$ onto its image is $1$.  Indeed, otherwise take two
general points $x$ and $x'$ of $\mu(E)$ and let $e_1$ and $e_2$
(resp. $e'_1$ and $e'_2$) be two distinct points lying in the preimage
$\mu^{-1}(x)$ (resp. $\mu^{-1}(x')$); now take a plane $P$ passing
through the points $x$ and $x'$ chosen on $\mu(E)$.  Then all rulings
of the cones $\phi(\PP(\Ee\otimes\Gamma_{e_i}))$ and
$\phi(\PP(\Ee\otimes\Gamma_{e'_i}))$, for $i=1,2$, intersect $P$.
From the identification
\[
  W = H^0(\Nn^\ast(3H)) \cong H^0(\Oo_Y(1)),
\]
we deduce that two hyperplanes in $\PP(W^\ast)$ cutting out the plane
$P$ give rise to two linearly independent global sections, say $s$ and
$s'$, of $\Nn^\ast(3H)$ such that $\gamma =s\wedge s'$, viewed as a
global section of $\det(\Nn^\ast(3H))=\Oo_X(H-K_X)$.  The section
$\gamma$ vanishes along the divisor
\[
  \Gamma_{e_1}+\Gamma_{e_2}+\Gamma_{e'_1}+\Gamma_{e'_2} = 4\Gamma_o
\]
which is impossible, since the divisor $H-K_X-4\Gamma_o=-\Gamma_o+l$
can not be numerically equivalent to an effective divisor.
\end{enumerate}

Let $E'$ be a Zariski dense open subset of $E$, where $\mu$ is an
embedding.  Take two distinct points $e_1$ and $e_2$ on $E'$. The
corresponding curves $\Gamma_{e_1}$ and $\Gamma_{e_2}$ intersect at a
point, call it $x$.  From the properties established above, we know
that the fibre $L_x$ of the projection $p:Y\to X$ is mapped by $\phi$
into the line joining $\mu(e_1)$ and $\mu(e_2)$.  Hence the secant
variety of $\mu(E')$ is contained in the image $Y_0=\phi(Y)$.  This
implies, in view of the irreducibility of $Y_0$, that the curve
$\mu(E)$ spans $\PP(W^\ast)$.  Furthermore, for $e\in E'$, the union
of the secant lines of $\mu(E)$ passing through $\mu(e)$ must be the
cone $\phi(\PP(\Ee\otimes\Gamma_{e}))$ with the vertex $\mu(e)$ and
base a smooth elliptic curve of degree $4$.  Hence this base curve is
the image of the projection of $\mu(E)$ from the point $\mu(e)$.  This
implies that $\mu$ must be an embedding, $\mu(E)$ is an elliptic
normal curve of degree $5$, and the secant variety of $\mu(E)$ is
contained in $Y_0$.  Since both are irreducible they must coincide.
The assertion that the degree of $Y_0$ is $5$ can now be deduced
either from the fact that the secant variety of an elliptic normal
curve in $\PP^4$ has degree $5$ or from the direct computation
\[
  \deg(Y_0) = c_1^3(\Oo_Y (1))) =  c^2_1(\Ee^\ast) -c_2(\Ee^\ast)=
  (H-K_X)^2 -10=15-10=5.
\]
The second equality uses the Chern invariants of 
$\Ee^\ast =\Nn^\ast_X(3H)$ computed in Theorem \ref{th:cubicsSections}.

The only remaining statement of Proposition \ref{p:descriptionOfY0},
1)\,-\,3), to prove is that $i:E\times E\to Y$ in \eqref{d:liftE2} is an
embedding.  From the work we have already done, this is immediate,
since the image of $i$ is a bi-section with respect to the projection
$p:Y\lra X$ and it is an isomorphism outside the diagonal $\Delta_E$.
On the diagonal $\Delta_E$, we know that the composition $\phi\circ i$
is an embedding.  Hence $i$ is an embedding everywhere.

Turning to the part 4) of Proposition \ref{p:descriptionOfY0}, we
observe that the fibre $(\pi\circ p)^{-1}(e)=S_e$ over $e\in E$ is
$\PP(\Ee\otimes\Oo_{F_e} )$, where $F_e=\pi^{-1}(e)$. This together
with 
\[
  \Ee\otimes\Oo_{F_e} = \Nn_X(-3H)\otimes\Oo_{F_e}
  \cong \Oo_{\PP^1}(-1)\oplus\Oo_{\PP^1}(-2) 
\]
(see \eqref{NonF} for the last isomorphism) proves that $S_e$ is
isomorphic to the Hirzebruch surface $\Sigma_1$.  Furthermore, the
line bundle $\Oo_Y(1)$ restricted to $S_e$ is
$\Oo_{S_e}(1):=\Oo_Y(1)\otimes\Oo_{S_e}$ which is very ample, since
$p_\ast(\Oo_{S_e}(1))=\Nn^\ast_X(3H)\otimes\Oo_{F_e}\cong
\Oo_{\PP^1}(1)\oplus\Oo_{\PP^1}(2)$ 
is obviously very ample.  The morphism $\phi|_{S_e}$ is defined by
$\Oo_{S_e}(1)$ and hence it embeds $S_e$ into $\PP(W^\ast)$ as a
scroll of degree 
$\deg(\Nn^\ast_X(3H)\otimes\Oo_{F_e})=\deg(\Oo_{\PP^1}(1)
\oplus\Oo_{\PP^1}(2))=3$, \ie the image $S'_e =\phi (S_e)$ is a
rational normal scroll of degree $3$ in $\PP(W^\ast)$. Setting $L_e$
to be the minimal section of $S_e$, we deduce
\[
  L_e = c_1 (\Oo_{S_e}(1))-2f,
\]
where $f$ is the class of a fibre of $p$.  In particular, $L_e$ is a
$(-1)$-curve on $S_e$ and its image $L'_e =\phi(L_e)$ is a line in
$\PP(W^\ast)$.  The rulings of $S'_e$ cut out on the elliptic curve
$\phi(\Delta_E)$ the pencil of degree $2$ corresponding to the fibre
$F_e\subset X\cong Sym^2(E)$.  Hence $\phi(\Delta_E)$ is a divisor on
$S'_e$ of the form $2L'_e +bl$ for some integer $b$, and where $l$,
the image of $f$, is the class of a ruling of $S'_e$.  The integer is
determined from the equation
\[
  5 = \deg(\phi(\Delta_E))
  = c_1(\Oo_{S'_e}(1))\cdot(2L'_e +bl) = 2+b.
\]
Hence $\phi(\Delta_E)=2L'_e+3l$.  From this it follows that 
\[
  L'_e \cdot \phi(\Delta_E) = L'_e \cdot (2L'_e +3l)
  = -2+3 = 1.
\]
Thus $ L'_e$ intersects $\phi (\Delta_E)$ transversely at a single
point. 

The minimal section $L_e$ of $S_e$ is a distinguished section of $p$
over the fibre $F_e =\pi^{-1}(e)$.  Varying $e$ in $E$ gives rise to
the section $\gamma$ of $p$ whose existence is claimed in the part 5)
of the proposition.  To be more formal, we seek a sub-linebundle of
$\Ee=\Nn(-3H)$ which coincides with $\Oo_{F_e}(-1)$ on every fibre
$F_e$ of $\pi$.  The construction is the same as in the proof of part
1) of the proposition.  Namely, we know that
$\Nn(-2H)\otimes\Oo_{F_e}\cong\Oo_{\PP^1}(-1)\oplus\Oo_{\PP^1}$.
Hence, the direct image $\pi_\ast(\Nn(-2H))$ is a line bundle on $E$
which will be denoted $\Oo_E (D)$.  This gives rise to the exact
sequence
\[
  0 \lra \pi^\ast(\Oo_E (D))
  \lra \Nn(-2H) \lra \Oo_X(K_X+H-\pi^\ast(D)) \lra 0
\]
and hence, the sought after sub-linebundle of $\Nn(-3H)$ is 
$\Oo_X(-H +\pi^\ast D)$.  Furthermore, the second Chern class
computation from the above exact sequence yields $\deg(D)=-5$. In
addition, by construction, the section $\gamma: X\lra Y$ corresponding
to the subbundle $\Oo_X(-H +\pi^\ast D)\subset\Nn(-3H)$ has the
property $\Oo_X(H -\pi^\ast D)=\gamma^\ast(\Oo_Y(1))$.  Hence the
image $X'=(\phi\circ\gamma)(X)$ has degree
\[
\deg(X') =(H-\pi^\ast D)^2 =(H+5f)^2=5+10=15.
\]
\qed

\begin{rem*}
  \label{r:aboutquadrics}
All the facts in Proposition \ref{p:descriptionOfY0}, with the
possible exception of identifying the $\PP^1$-bundle $\PP(\Ee)$, have
been proved in \cite{ADHPR}.  It seems to us that taking the vector
bundle $\Ee^\ast=\Nn^\ast_X(3H)$ on the scroll $X$ as the starting
point makes its relation with the associated elliptic quintic curve
more natural.

The property of $\PP(\Ee)$ being a fibration by Hirzebruch surfaces
$\Sigma_1$, which comes almost for free in our exposition, is all one
needs to establish the relationship with the space of quadrics passing
through the elliptic normal quintic curve $\phi(\Delta_E)$ and
eventually, show that $\PP(\Ee)$ and 
$\PP(\Nn_{\phi(\Delta_E)}\otimes\Oo_{\PP(W^\ast)}(-2))$
are related by a birational morphism%
\footnote{The morphism is the relative blow-down of $\PP(\Ee)$ along
  the image of the section $\gamma$ in Proposition
  \ref{p:descriptionOfY0}, 5).},
where $\Nn_{\phi(\Delta_E)}$ is the normal bundle of
$\phi(\Delta_E)$ in $\PP(W^\ast)$.  See \cite{ADHPR} for details. 
\end{rem*}

With the above description of the morphism $\phi:Y\to\PP(W^\ast)$ we
can now explain how to go from $\phi$ back to the embedding 
$X\subset\PP(V^\ast)$.  For this, we use the following notation.  
\begin{itemize}
  \listspace
\item 
  $\langle z,z'\rangle$ is the line through the distinct points
  $z,z'\in\PP^4$
\item
  $Y_x=p\inv(x)$, the fibre of $p$ over $x\in X$. 
\end{itemize}

\begin{pro}
  \label{p:HZH}
The isomorphism
\begin{equation}
  \label{eq:isoH0EdualW}
  W= H^0(\Nn^\ast(3H)) = H^0(\Ee^\ast) \cong H^0(\Oo_Y(1)),
\end{equation}
establishes the following geometric correspondence between the
hyperplanes in $\PP(W^\ast)$ and $0$-cycles on $X$.

Let $M$ be a hyperplane in $\PP(W^\ast)$ intersecting the elliptic
normal curve $\phi(\Delta_E)$ along the divisor 
$D_M=M\cdot\phi(\Delta_E)$ consisting of $5$ distinct points%
\footnote{In what follows we abuse the notation by tacitly identifying
  $\phi(\Delta_E)$ with $E$.}.
Then the configuration of the ten secant lines of $\phi(\Delta_E)$
\[
  \sum_{e\neq e'\in D_M} \langle e,e'\rangle,
\]
gives rise to the $0$-cycle
\begin{equation}
  \label{eq:ZH}
  Z_{M} = \sum_{e\neq e'\in D_M}  p_\ast(Y_{\tau(e,e')})
\end{equation}
on $X$, where $\phi(Y_{\tau(e,e')})=\langle e,e'\rangle$.  That
$0$-cycle is the scheme of zeros of a unique projective section
$[s_M]\in \PP(H^0(\Nn^\ast(3H)))$ corresponding to $M$ under the
isomorphism \eqref{eq:isoH0EdualW}.
\end{pro}

\proof
Let $[s_M]\in \PP(H^0(\Nn^\ast(3H)))$ be the projective section
corresponding to a hyperplane $M$. The the zero locus 
$(s_M =0) \subset X$ parametrizes the fibres of the projection
$p:Y\to X$ which are mapped onto the secant lines of $\phi(\Delta_E)$
contained in $M$ and those are precisely the ones appearing on the
right side of the equality \eqref{eq:ZH}.
\qed

We can now reconstruct\footnote{We are grateful to I. Dolgachev for
  pointing out to us the construction that follows; according to him,
  it was Segre's way to see an elliptic scroll inside a Segre cubic.} 
$X\subset \PP(V^\ast)$ from the correspondence
\[
  \PP(W) \ni M \mapsto Z_{M}
  = \sum_{e\neq e'\in D_M} p_\ast(Y_{\tau(e,e')})
\]
described in Proposition \ref{p:HZH}.  The $0$-cycle $Z_M$ admits a
distinguished decomposition into $5$ subcycles
\begin{equation}
  \label{ZHdecomp}
  Z_{M} = \bigcup_{e\in D_M} Z^{e}_{M}
  \textq{where}
  Z^{e}_{M} = \sum_{e'\in D_M\smallsetminus\{e\}} p_\ast(Y_{\tau(e,e')}),
\end{equation}
\ie $Z^{e}_{M}$ parametrizes the fibres of $p:Y\to X$ mapped by
$\phi$ onto the rulings of the cone
$\phi(\PP(\Ee\otimes\Oo_{\Gamma_e}))$ that are contained in $M$.
Hence we have
\begin{equation*}
  \label{ZHdecomp1}
  Z_M^e \subset \Gamma_e  \textq{for every} e\in D_M.
\end{equation*}
The configuration of $5$ curves $\{\Gamma_e\}_{e\in D_M}$, on the side
of $\PP(V^\ast)$, gives rise to $5$ planes $\{\Pi_e\}_{e\in D_M}$,
where $\Pi_e$ is the span of $\Gamma_e$ in the embedding
$X\subset\PP(V^\ast)$.  Now one recovers the Segre cubic $V_3(s_M)$,
the cubic hypersurface corresponding to $s_M$ under the isomorphism 
$H^0(\Nn^\ast(3H))\cong I_X(3)$ in Theorem \ref{th:cubicsSections}, as
the union of the lines in $\PP(V^\ast)$ intersecting (any) four of the
five planes of the collection $\{\Pi_e \}_{e\in D_M}$.  The $0$-cycle
$Z_{M}=(s_M =0)=\Sing(V_3 (s_M))$ is seen in $\PP(V^\ast)$ as the
cycle
\[
Z_{M} = \sum_{e\neq e' \in D_M} \Pi_e \bigcap \Pi_{e'}.
\]
composed of points of pairwise intersections of the planes 
$\{\Pi_e \}_{e\in D_M}$.  As $M$ varies in the complement of the dual
variety of the elliptic normal curve $\phi(\Delta_E)$, the zero cycles
$Z_{M}$ sweep out the scroll $X$.

\medskip

Let $U$ be the Zariski dense open subset of $\PP(W)$ parametrizing
regular global sections of $\Nn^\ast(3H)$. From the construction
above $U$ coincides with the set of hyperplanes in $\PP(W^\ast)$
intersecting the elliptic normal quintic curve 
$\phi(\Delta_E)\subset \PP(W^\ast)$ along five distinct points.  Thus
the complement
\begin{equation*}
  \label{dualElcurve}
  T_{\phi(\Delta_E)} := \PP(W^\ast) \setminus U
\end{equation*}
is the dual variety of $\phi(\Delta_E)$, \ie the closed points of
$T_{\phi(\Delta_E)}$ parametrize the hyperplanes in $\PP(W^\ast)$
which are tangent to $\phi(\Delta_E)$ at some point.  Furthermore, the
construction described above and the identification
\begin{equation}
  \label{PWcubicsApp}
  \PP(W) \cong |I_X(3)|
\end{equation}
in Theorem \ref{th:cubicsSections}, 3), imply that the Zariski open
subset $U$, under the above identification, corresponds to Segre
cubics in $\PP(V^\ast)$ containing $X$.  Thus the dual variety
$T_{\phi(\Delta_E)}$, via the identification \eqref{PWcubicsApp},
parametrizes cubic hypersurfaces in $|I_X(3)|$ with degenerate
singular locus.  To see this degeneracy we take a hyperplane $M$ in
$\PP(V^\ast)$ which has a contact of order $2$ with $\phi(\Delta_E)$
at a point $e_0$ and is transversal to $\phi(\Delta_E)$ at the
remaining three points which we denote $e_i$, $i=1,2,3$.  Thus the
divisor $D_M =M\cdot\phi(\Delta_E)$ has the form
\[
  D_M = e_0 +t_{e_0} +e_1+e_2+e_3,
\]
where $t_{e_0}$ denotes a tangent vector of $\phi(\Delta_E)$ at
$e_0$. The corresponding configuration of secant lines of
$\phi(\Delta_E)$ contained in $M$ consists of:
\begin{itemize}
  \listspace
\item
  $\langle e_0,e_i\rangle$, $i=1,2,3$, with multiplicity $2$, 
\item
  the embedded tangent line $\langle e_0,t_{e_0}\rangle$ of
  $\phi(\Delta_E)$ at $e_0$,
\item
  three lines $\langle e_i,e_j\rangle$, $1\leq i<j\leq 3$. 
\end{itemize}
This accounts for ten secant lines counted with multiplicities. Since
the scheme of zeros $Z_M$ of the projective section $[s_M] \in \PP(W)$
corresponding $M$ under the isomorphism in \eqref{eq:isoH0EdualW} must
parametrize the secants of $\phi(\Delta_E)$ contained in $M$ and
$\deg(Z_M) =10$, we deduce that no other secant of $\phi(\Delta_E)$ is
contained in $M$ and the $0$-cycle $Z_M$ has the form
\begin{equation}
  \label{ZMdegen}
  Z_M
  = 2\sum^3_{i=1}\tau(e_0,e_i) + \tau(e_0,e_0)
  + \sum_{1\leq i<j\leq 3}\tau(e_i,e_j),
\end{equation}
where $\tau:E\times E\to X$ is the double covering in
\eqref{p:psimap}.
 
Now we move on the side of the embedding $X\subset\PP(V^\ast)$.
Denote by $\Gamma_i$ the curves of the family $\{\Gamma_b\}_{b\in E}$
in \eqref{curvesGp} corresponding to the points $e_i$, with
$i=0,\ldots,3$, (we are making here the obvious identification of $E$
and $\phi(\Delta_E)$) and let 
$\{\Pi_i=\Span(\Gamma_i)\}_{0\leq i\leq3}$ be the corresponding
configuration of planes.  By construction, each curve $\Gamma_i$ is
identified with the projection of $\phi(\Delta_E)$ from the point
$e_i$ and comes along with a distinguished divisor $D_i$ in the linear
system $|\Oo_{\Gamma_i}(H-K_X)|=|\Oo_{\Gamma_i}(H+\Gamma_i)|$:
\[
  \begin{aligned}
    D_0 &= e_0 +e_1+e_2+e_3 \\
    D_i &= e_0 +t_{e_0} + \sum_{j\neq i}e_j  \textq{for} i=1,2,3.
  \end{aligned}
\]
In particular, $e_0=\Gamma_0\cdot\Gamma_0$ and 
$e_1+e_2+e_3 =D_0-e_0 \sim H\cdot\Gamma_0$.  This means that the plane
cubic section $\Gamma_0\subset \Pi_0$ comes along with the line
$l_0=\Span\{e_1,e_2,e_3\}$ spanned by the three colinear points
$e_i$ ($i=1,2,3$).
 
With the four distinct planes $\{\Pi_i\}_{0\leq i\leq3}$ the
 construction of Segre 
associates
\begin{equation}
  \label{Velcubic}
  V(e_0,l_0):=
  \text{the union of lines in $\PP(V^\ast)$ intersecting the
    four planes $\{\Pi_i\}_{0\leq i\leq3}$.}
\end{equation}
This is a cubic hypersurface containing $X$ and its singular locus
$\Sing(V(e_0,l_0))$, according to Theorem \ref{th:cubicsSections}, 3),
intersects $X$ along the $0$-cycle $Z_M$ in \eqref{ZMdegen}.  In
particular, $V(e_0,l_0)$ is singular at three colinear points
$e_1,e_2,e_3$ and hence, it must be singular along the line $l_0$.
This together with Remark \ref{r:degcubics} implies that the
$1$-dimensional part of $\Sing(V(e_0,l_0))$ is the line $l_0$.  In
addition, the formulas $Z_M =\Sing(V(e_0,l_0))\cdot X$ and
\eqref{ZMdegen} tell us that $V(e_0,l_0)$ is singular at four distinct
points, $\tau(e_0,e_0)$, $\tau(e_1,e_2)$, $\tau(e_1,e_3)$ and
$\tau(e_2,e_3)$, lying outside of $l_0$.  Since this is the maximal
possible number of isolated singularities for a cubic hypersurface in
$\PP^4$ with precisely one line as the $1$-dimensional part of its
singular locus, we deduce
\[
  \Sing(V(e_0,l_0)) = l_0
  \,\bigcup\, \{\tau(e_0,e_0),\tau(e_1,e_2),\tau(e_1,e_3),\tau(e_2,e_3)\}. 
\]

By the above construction, a Zariski dense open subset of the dual
variety $T_{\phi(\Delta_E)}$ corresponds to the open part of
`degenerate' cubics in $|I_X(3)|$, \ie cubics having precisely a
line as the $1$-dimensional part of their singular locus.  Since
$T_{\phi(\Delta_E)}$ is irreducible, we deduce the following.

\begin{pro}
  \label{p:degcubics}
Under the isomorphism $\PP(W)\cong|I_X(3)|$ in \eqref{PWcubicsApp},
the dual variety $T_{\phi(\Delta_E)}$ corresponds to the hypersurface
$\mathsf{C}_1$ in $|I_X(3)|$ parametrizing the cubic hypersurfaces
containing $X$ and having precisely a line as the $1$-dimensional part
of their singular locus. Furthermore, the singular lines of the cubics
in $\mathsf{C}_1$ are precisely the trisecant lines of 
$X\subset \PP(V^\ast)$. In particular, $T_{\phi(\Delta_E)}$ and
$\mathsf{C}_1$ are both isomorphic to the variety of trisecants of
$X$.
\end{pro}
 

\begin{rem*}
The desingularization of the dual variety $T_{\phi(\Delta_E)}$ is the
projectivization $\PP(\Ss)$ of
$\Ss:=N_{\phi(\Delta_E)}^\ast\otimes\Oo_{\PP(W^\ast)}(1)$, the twisted
conormal bundle of $\phi(\Delta_E)\subset\PP(W^\ast)$, while the
desingularization of the variety of trisecants of
$X\subset\PP(V^\ast)$ is $\PP(\Ff)$, where $\Ff$ is the bundle defined
in Proposition \ref{p:planebundleinPV}.  So the last assertion of
Proposition \ref{p:degcubics} implies an isomorphism between
$\PP(\Ss)$ and $\PP(\Ff)$, hence an identification
$\Ss\cong\Ff\otimes\Ll$, for some line bundle $\Ll$ on $E$.  Comparing
the degrees on both sides, we obtain $\deg(\Ll)=-5$.  Thus
\[
  \Ll = \Oo_{\phi(\Delta_E)}(-1)\otimes\Ll',
\]
where $\Ll'\in\Pic^0(E)$ and
$\Oo_{\phi(\Delta_E)}(1)=\Oo_{\PP(W^\ast)}(1)\otimes\Oo_{\phi(\Delta_E)}$.
Thus, the above identification takes the
form 
\[
  N_{\phi(\Delta_E)}^\ast\otimes\Oo_{\PP(W^\ast)}(2)
  \cong \Ff\otimes\Ll'.
\]
This isomorphism can be viewed as the vector bundle version of the
quadro-cubic Cremona transformation, see \cite{ADHPR} for another
proof of this isomorphism.
\end{rem*}

\subsection{The decomposition \eqref{ZHdecomp} and the
  configuration $(10_4 , 15_6)$ of Segre} 

It is well-known that the ten ordinary double points of a Segre cubic
hypersurface in $\PP^4$ have remarkable combinatorial properties:
there are fifteen planes, each spanned precisely by four singular
points, and through every singular point pass precisely six of the
fifteen planes. This is called a $(10_4,15_6)$ configuration.  In this
subsection we revisit this configuration in the light of the
correspondence between the embedding $X\subset\PP(V^\ast)$ and the
geometry of the twisted conormal bundle $\Nn_X^\ast(3H)$ encapsulated
in the properties of the morphism 
$\phi: Y=\PP(\Nn_X(-3H))\to\PP(W^\ast)$, see Proposition
\ref{p:descriptionOfY0}.  We keep the notation of the previous
subsection unless stated otherwise.

On the side of the vector bundle $\Nn_X^\ast(3H)$ we take a global
section $s$ whose zero locus $Z_s=(s=0)$ consists of ten distinct
points on $X$.  According to Theorem \ref{th:cubicsSections}, the
section $s$ corresponds to a Segre cubic denoted
$V_3(s)\in\PP(I_X(3))$, and under this correspondence
\[
  \Sing(V_3(s))=Z_s,
\]
where $\Sing(V_3(s))$ stands for the singular locus of $V_3(s)$.  So
from the Segre's works we know that $Z_s$ is a $(10_4,15_6)$
configuration.  We wish to be more specific and identify the fifteen
planes as well as the $4$-point subcycles of $Z_s$ spanning these
planes in the light of geometry of the morphism $\phi$ and of the
embedding $X\subset\PP(V^\ast)$.

We know from the previous section that $\phi$ gives rise to a
distinguished decomposition of $Z_s$ into subcycles of degree
$4$. Namely, the isomorphism
\[
  H^0(\Nn^\ast(3H)) \cong H^0(\Oo_Y(1)) = W
\]
identifies $s$ with the hyperplane $M_s$ in $\PP(W^\ast)$.  Denote by
$D_s$ the divisor cut out by $M_s$ on the elliptic normal curve
$\phi(\Delta_E)$, \ie
\[
  D_s: = M_s \bigcap \phi(\Delta_E).
\]
Then \eqref{ZHdecomp} provides the decomposition
\begin{equation}
  \label{Zsdecomp}
  Z_s = \bigcup_{e\in D_s} Z^e_s,
\end{equation}
where, to shorten the notation, we write $Z^e_s$ instead of
$Z^e_{M_s}$.  We know that each $Z_s^e$ is a subcycle of degree $4$
in $Z_s$ lying on the plane elliptic curve
$\Gamma_e\subset\PP(V^\ast)$ and no other point of $Z_s$ is contained
in $\Gamma_e$.  In particular, the plane $\Pi_e$ spanned by $\Gamma_e$
contains $Z_s^e$ and no other point of $Z_s$.  Furthermore, since no
three points of $Z_s$ can be
colinear\footnote{Otherwise the line through the three colinear points
  of $Z_s$ is a singular line of the cubic $V_3(s)$ which is
impossible.}, the plane $\Pi_e$ is spanned by $Z_s^e$.  This gives
us a collection of five planes $\{\Pi_e \}_{e\in D_s}$ with the
property that for every $e\neq e'$ the intersection
\[
 \Pi_e \bigcap \Pi_{e'} = \Gamma_e\cdot\Gamma_{e'}
\] 
is a single point, see the proof of Proposition \ref{p:psimap}.
Furthermore, the $0$-cycle $Z_s$ is seen in the embedding 
$X\subset\PP(V^\ast)$ as the cycle of points in $\PP(V^\ast)$ formed
by the above pairwise intersections
\begin{equation}
  \label{ZsinPV}
  Z_s=\sum_{e\neq e' \in D_s} \Pi_e \bigcap \Pi_{e'}.
\end{equation}
In the sequel we set 
\begin{equation}
  \label{edote'}
  e\cdot e':= \Pi_e \bigcap \Pi_{e'} \textq{for} e\neq e' \in D_s.
\end{equation}
With this notation the formula \eqref{ZsinPV} reads
\[
  Z_s=\sum_{e\neq e' \in D_s}e\cdot e',
\]
while, for every $e\in D_s$, the subcycle $Z^e_s$ in the decomposition
\eqref{Zsdecomp} take the form
\[
  Z^e_s = \Pi_e \bigcap Z_s = \sum_{ e' \neq e \in D_s}e\cdot e'.
\]

We now give a description of the remaining ten planes of $(10_4,15_6)$
configuration.

\begin{lem}
  \label{planesee'}
For every $e\neq e' \in D_s$, let
\begin{equation*}
  \label{Zee'}
  Z^{e\cdot e'}_s
  = Z_s \smallsetminus \left(Z^e_s\bigcup Z^{e'}_s\right) + e\cdot e'.
\end{equation*}
Then $Z^{e\cdot e'}_s$ is a $4$-degree subcycle of $Z_s$ that spans
in $\PP(V^\ast)$ a plane denoted by $\Pi_{e\cdot e'}$.
\end{lem}

\proof  
The fact that $Z^{e\cdot e'}_s$ consists of four points is obvious.
To see that these points span a plane in $\PP(V^\ast)$ we use the
identification of $Z_s$ with the singular locus $\Sing(V_3(s))$ of the
Segre cubic $V_3(s)$, see Theorem \ref{th:cubicsSections}.  We fix the
point $e\cdot e'=\Pi_e\bigcap\Pi_{e'} $ and consider the projection
from this point onto a complementary $\PP^3$ which intersects $V_3(s)$
transversely along a smooth cubic surface, call it $F$.  Then the
remaining nine points of $Z_s=\Sing(V_3(s))$ project to the nine nodes
of a $(3,3)$-divisor
\[
  A = F \bigcap Q,
\]
where $Q$ is a smooth quadric in $\PP^3$, the image of the tangent
cone of $V_3(s)$ at the point $e\cdot e'$.  The nine nodes force the
divisor $A$ to be completely reducible, \ie $A$ has the form
\[
  A = \sum^3_{i=1}f_i +\sum^3_{i=1}g_i,
\]
where $f_i$ (resp. $g_i$), for $i=1,2,3$, are three disjoint rulings
of $Q$ such that the $0$-cycle
\[
  Z'_s = \sum_{1\leq i,j\leq 3} f_i\bigcap g_j
\]
is the image of $Z_s\smallsetminus\{e\cdot e'\}$ under the projection
from $e\cdot e'$, see \cite{Dol} for an excellent account of the above
construction.
 
Under the projection from $e\cdot e'$, the planes $\Pi_{e}$ and
$\Pi_{e'}$ go to two skew rulings, say $f_1$ and $f_2$.  Hence the
points $Z^e_s$ and $Z^{e'}_s$ are mapped to the points
$\sum_{j=1}^3f_1\bigcap g_j$ and $\sum_{j=1}^3f_2\bigcap g_j$
respectively.  Therefore, the remaining points
$Z_s\smallsetminus(Z^e_s \bigcup Z^{e'}_s)$ of
$Z_s\smallsetminus\{e\cdot e'\} $ must go to 
$\sum_{j=1}^3f_3\bigcap g_j$.  Hence the three points 
$Z_s \smallsetminus (Z^e_s\bigcup Z^{e'}_s)$ together with $e\cdot e'$
lie in the plane spanned by the line $f_3$ and the point $e\cdot e'$.
Using again the fact that no three points of $Z_s$ are colinear, we
deduce the assertion of the claim.
\qed

The five planes $\{\Pi_e\}_{e\in D_s}$ together with the planes
$\Pi_{e\cdot e'}$, ($e\neq e' \in D_s$), in Claim \ref{planesee'}
account for fifteen planes in the $(10_6,15_4)$ configuration.  This
collection of planes will be denoted by
\begin{equation}
  \label{planecollec}
  \PPP_s := \{\Pi_e, \Pi_{e\cdot e'} \mid e\in D_s, e\neq e' \in D_s \}.
\end{equation}
The following is obvious from the construction.

\begin{lem}
  \label{cl:interofplanes} 
Every point $e\cdot e'\in Z_s$ lies precisely on the following six
planes of the collection $\PPP_s$,
\[
  \PPP^{e\cdot e'}_s
  = \{ \Pi_e, \Pi_{ e'}, \Pi_{e\cdot e'}, \Pi_{e''\cdot e'''}
  \mid e''\neq e''' \in D_s \smallsetminus \{e,e'\} \}.
\]
Furthermore, each subset of the partition 
$\PPP^{e\cdot e'}_s={}^1\PPP^{e\cdot e'}_s\bigcup{}^2\PPP^{e\cdot
  e'}_s$, where
\[
  {}^1\PPP^{e\cdot e'}_s = \{\Pi_e, \Pi_{ e'}, \Pi_{e\cdot e'} \}
  \textq{and}
  {}^2\PPP^{e\cdot e'}_s =
  \{ \Pi_{e''\cdot e'''}
  \mid e''\neq e''' \in D_s\smallsetminus\{e,e'\} \},
\]
consists of three planes of $\PPP^{e\cdot e'}_s$ which intersect
precisely at $e\cdot e'$, while the planes taken from different
subsets intersect along a line.  More precisely, if 
$D_s=\{e,e',e'',e''',c\}$, then
\[
  \Pi_e \cap  \Pi_{e''\cdot e'''} = \langle e\cdot e', e\cdot c \rangle, 
  \quad
  \Pi_{e'} \bigcap  \Pi_{e''\cdot e'''}
  = \langle e\cdot e', e'\cdot c \rangle,
  \quad
  \Pi_{e\cdot e'} \bigcap  \Pi_{e''\cdot e'''}
  = \langle e\cdot e', e''\cdot e''' \rangle,
\]
where for two distinct points $x, y$ in a projective space, 
$\langle x,y\rangle$ denotes the line spanned by those points.
\end{lem}

\begin{rem}
The above two claims give a precise recipe of how to recover the
collection $\PPP_s$ of the fifteen planes of the $(10_4, 15_6)$
configuration from the cycle $Z_s$ of ten points on 
$X\subset V_3(s)$.  It should be pointed out that the collection
$\PPP_s$ has an extra feature of being `polarized' into two types of
planes: $\Pi_e$, $e\in D_s$, and $\Pi_{e\cdot e'}$, for 
$e\neq e'\in D_s$. The presence of this polarization is due of course
to the fact that $Z_s$ is not just the singular locus of a Segre cubic
$V_3(s)$, but also lies on the scroll $X$ inside of $V_3(s)$. To be
even more precise, the above polarization emerges from the fact that
$Z_s$ is the zero locus of a section of a vector bundle on $X$.
\end{rem}

The two types of planes in $\PPP_s$ play different roles with
respect to the embedding $X\subset \PP(V^\ast)$.  By construction, the
planes $\{\Pi_e \}_{e\in D_s}$ are distinguished by the property that
the intersection $X \bigcap\Pi_e$ is the plane elliptic curve
$\Gamma_e$, for every $e\in D_s$.  The following lemma gives a
similar characterization of the planes $\Pi_{e\cdot e'}$.

\begin{lem}
  \label{planePee'}
For every pair $e\neq e' \in D_s$, the plane $\Pi_{e\cdot e'}$
intersects the scroll $X$ along the subcycle $Z_s^{e\cdot e'}$ and the
ruling $l_{e\cdot e'}$ of $X$ passing through the point $e\cdot e'$.
\end{lem}

\proof
Let $\{c,c',c''\}$ be the complement $D_s\smallsetminus\{e,e'\}$.
From Claim \ref{planesee'} it follows that the subcycle
\[
  Z_s^{e\cdot e'} = c\cdot c' + c\cdot c'' +c'\cdot c'' +e\cdot e'
\]
is contained in the intersection $X\bigcap\Pi_{e\cdot e'}$.  Observe
that the line $L=\langle c\cdot c',c\cdot c''\rangle$ is contained
in the plane $\Pi_c$ and hence it intersects the curve $\Gamma_c$ along
three points (the degree $3$ divisor)
\[
  T = L\cdot \Gamma_c =c\cdot c'+ c\cdot c'' +t.
\]
The same holds for the line 
$L' =\langle c\cdot c', c'\cdot c'' \rangle$ (resp.  
$L'' =\langle c\cdot c'', c'\cdot c'' \rangle$) and leads to 
\[
  T'= L'\cdot \Gamma_c = c\cdot c'+ c'\cdot c'' +t'
  \qquad(\text{resp. } T'' = c\cdot c'' + c'\cdot c'' +t'').
\]
It follows that the plane $\Pi_{e\cdot e'}$ intersects $X$ along seven
points, $Z^{e\cdot e'} +t+t'+t''$.  Since the degree of $X$ is $5$, we
deduce that the intersection of $\Pi_{e\cdot e'}$ and $X$ is not
proper, \ie that the intersection has a $1$-dimensional component,
call it $F$.  The scheme $F$ can be either a plane cubic or a ruling
of $X$.  According to Proposition \ref{p:embedGpplusD2}, the first
possibility means that $F$ is one of the curves 
$\{\Gamma_b\}_{b\in E}$ and this is clearly impossible.  Hence $F$ is
a ruling of $X$.

To identify this ruling we observe that it must meet all curves
$\{\Gamma_b\}_{b\in E}$.  In particular, it must intersect
$\Gamma_e$. Hence $F$ must pass through the intersection 
$\Pi_{e\cdot e'} \bigcap \Pi_{e}$ and this, in view of Lemma
\ref{cl:interofplanes}, is the point ${e\cdot e'}$.  Hence $F$ is the
ruling of $X$ passing through the point ${e\cdot e'}$ as asserted.
\qed

\subsection{Toward a categorification of the configuration
  $(10_4 , 15_6)$ of Segre}

Conceptually, the whole approach of our paper can be termed as a
representation of various geometric or cohomological entities attached
to a surface in $\PP^4$ in the category of complexes of coherent
sheaves on that surface.

In this subsection we apply this approach to the $(10_4,15_6)$
configuration of Segre considered in the previous section. Namely,
with our geometric set up of an elliptic scroll $X$ embedded in
$\PP(V^\ast)$, we have seen how the scheme of zeros $Z_s$ of a
regular\footnote{A {\it regular} global section is a global section with
  simple isolated zeros.} global 
section $s$ of $\Nn_X^\ast(3H)$ acquires the structure of the
$(10_4,15_6)$ configuration of Segre.  With our notation and results
from the previous subsection ,
\[
  Z_s = \sum_{e\neq e' \in D_s} e\cdot e',
\]
where $D_s$ is intrinsically definedby $s$; it is a set of five
distinct points on the elliptic curve $E$, the base of $X$.  See
\eqref{ZsinPV} and \eqref{edote'} for notation.  We have noticed that
$Z_s$ contains the distinguished subcycles of degree $4$,
\begin{equation}
  \label{subcycles}
  Z^{e}_s =\sum_{e'\in D_s\smallsetminus\{e\}} e\cdot e', \qquad
  Z^{e\cdot e'}_s
  = Z_s \smallsetminus {\left(Z^e_s \bigcup Z^{e'}_s\right)}+ e\cdot e',
\end{equation}
which have the geometric property of spanning the planes $\Pi_e$ and
$\Pi_{e\cdot e'}$ in $\PP(V^\ast)$.  These planes form the collection
$\PPP_s$ of fifteen planes in \eqref{planecollec}.  We suggest that
there is a lifting of $(Z_s,\PPP_s)$ to the category
$\mathfrak{Comp}(X)$ of (short) exact complexes of torsion free
sheaves on $X$ and hence, to the derived category $\DDD(X)$ of the
coherent sheaves on $X$.  Before we go on, let us be more precise
about our suggestion.

The first step of a categorification process is to assign a complex to
every subcycle in \eqref{subcycles}.  The second one is to turn
$(Z_s,\PPP_s)$ into a category and then to check that the morphisms of
that category go to morphisms of complexes.

The main result of this subsection is a realization of the first step.
As for the second, let us just indicate here how one could think of
$(Z_s,\PPP_s)$ as a category.  This can be achieved by turning
$\PPP_s$ into a graph, call it $\Cc(\PPP_s)$:
\begin{itemize}
  \listspace
\item[--]
  the vertices of $\Cc(\PPP_s)$ are the planes of the collection
  $\PPP_s$, 
\item[--]
  there is an edge between two vertices if and only if the
  corresponding planes intersect along a line. 
\end{itemize}
Of course, the graph $\Cc(\PPP_s)$ is obviously a category: the
objects are the vertices of $\Cc(\PPP_s)$ and the morphisms between
two objects, say from $\Pi$ to $\Pi'$, are the paths, composed of
edges of the graph, beginning at $\Pi$ and ending at $\Pi'$.  With
this understood, we define $\Cc(\PPP_s)$ to be the category of the
Segre configuration $(Z_s,\PPP_s)$ and we propose that there should be
geometrically interesting functor(s)
\begin{equation}
  \label{functor}
  \FFF: \Cc(\PPP_s) \lra \mathfrak{Comp}(X)
  \quad(\text{resp. } \DDD(X)).
\end{equation}

In the sequel, we construct such a functor $\FFF$ on the level of
objects.  This is essentially Serre construction and is based on the
following observation.

\begin{lem}
  \label{l:Zext}
Let $Z$ be one of the subcycles of degree $4$ in \eqref{subcycles}.
Then there is an extension sequence
\begin{equation}
  \label{Zext}
  0 \lra \Oo_X(K_X) \lra \Ff_Z \lra \Jj_Z(H) \lra 0
\end{equation}
intrinsically attached to $Z$, where $\Jj_Z$ is the ideal
sheaf of $Z$ in $X$ and $\Ff_Z$ is a locally free sheaf of rank $2$
with Chern invariants
\[
  c_1(\Ff_Z) = K_X+H  \textq{and} c_2(\Ff_Z) = -1.
\]
\end{lem}

\proof
The geometric condition of $Z$ spanning a plane in $\PP(V^\ast)$ is
translated, via the exact sequence
\[
  0 \lra \Jj_Z(H) \lra \Oo_X(H)\lra \Oo_Z(H) \lra 0,
\]
to the cohomological condition $h^1(\Jj_Z(H))=1$.  This and the Serre
duality
\[
  H^1(\Jj_Z(H)))^\ast \cong \Ext^1(\Jj_Z(H),\Oo_X(K_X)).
\]
imply that there is an extension as in \eqref{Zext} and such an
extension is unique, up to the $\CC^{\times}$-action of scaling the
morphisms in that sequence.  Furthermore, since for any proper subscheme
$Z'\subset Z$ the cohomology $H^1(\Jj_{Z'}(H))=0$, it follows by a
lemma of Serre, \cite[Lemma 5.1.2]{OScS}, that the sheaf $\Ff_Z$ in
\eqref{Zext} is locally free.  Its invariants are immediately deduced
from \eqref{Zext}.
\qed

Next we investigate the vector bundle $\Ff_Z$ in  \eqref{Zext}.

\begin{lem}
  \label{l:FZbundle}
The vector bundle $\Ff_Z$ in \eqref{Zext} is $H$-unstable.  More
precisely, there is an effective nonzero divisor $A_Z$ on $X$ such
that $\Ff_Z$ fits into the short exact sequence
\begin{equation}
  \label{FZdestseq}
  0 \lra \Oo_X(A_Z) \lra \Ff_Z \lra \Jj_{Z'}(K_X+H-A_Z) \lra 0,
\end{equation}
where $Z'$ is a $0$-dimensional subscheme of $X$ and $\Jj_{Z'}$ is its
ideal sheaf. Furthermore, $\Oo_X(A_Z)$ is the $H$-maximal
destabilizing subsheaf of $\Ff_Z$.
\end{lem}

\proof
From \eqref{Zext} it follows that
\[
  h^0(\Ff_Z) \geq h^0(\Jj_Z(H)) - h^1(\Oo_X(K_X)) = 2-1 = 1.
\]
Hence $\Ff_Z$ has a nonzero global section, call it $f$.  This and the
Chern invariant $c_2 (\Ff_Z)=-1$ computed in Lemma \ref{l:Zext}, imply
that the subscheme of zeros of $f$ must have a divisorial part which
is the divisor $A_Z$ of the lemma.  Hence $\Ff_Z(-A_Z)$ has a nonzero
global section $f'$ whose zero locus is $0$-dimensional. The asserted
sequence \eqref{FZdestseq} is the Koszul sequence of $f'$ tensored
with $\Oo_X(A_Z)$.

From $H\cdot A_Z > 0 =H\cdot (K_X+H) =H\cdot c_1(\Ff_Z)$ it follows
that $\Oo_X(A_Z)$ is $H$-destabilizing.  This together with the fact
that the quotient sheaf in \eqref{FZdestseq} is torsion free insures
the maximality of $\Oo_X(A_Z)$.
\qed

Next we show how the destabilizing sequence \eqref{FZdestseq}
distinguishes between the two types of subcycles in
\eqref{subcycles}.  For this we put the defining extension sequence
\eqref{Zext} together with the destabilizing one to obtain the
diagram.
\vspace{15ex}
\begin{equation}
  \label{d:Zextdest}
  \begin{tikzpicture}[overlay,every node/.style={draw=none},
    ->,inner sep=1.1ex]
    \matrix [draw=none,row sep=3.5ex,column sep=4.2ex]
    {
      && \node (02) {$0$}; \\
      && \node (12) {$\Oo_X(A_Z)$}; \\
      \node (20) {$0$};
      & \node (21) {$\Oo_X(K_X)$};
      & \node (22) {$\Ff_Z$};
      & \node (23) {$\Jj_Z(H)$};
      & \node (24) {$0$}; \\
      && \node (32) {$\Jj_{Z'}(K_X+H-A_Z)$}; \\      
      && \node (42) {$0$}; \\      
    };
    \path
      (02) edge (12) 
      (12) edge (22)
      (22) edge (32)
      (32) edge (42)

      (12) edge (23)
      (21) edge (32)
      
      (20) edge (21)
      (21) edge (22)
      (22) edge (23)
      (23) edge (24)
    ;
  \end{tikzpicture}
  \vspace{15ex}
\end{equation}
The morphisms defined by the slanted arrows are nonzero and hence give
rise to a nonzero effective divisor $B_Z \in|\Jj_Z(H-A_Z)|$
(resp. $|\Jj_{Z'}(H-A_Z)|$).  In particular, we obtain the
decomposition
\begin{equation*} 
  \label{HAZBZ}
  H=A_Z +B_Z
\end{equation*}
and claim the following.

\begin{lem}
  \label{l:AZ}
Let $\pi:X\to E$ be the structure projection of $X$ onto the elliptic
curve $E$.  Then $A_Z =\pi^\ast(\aaa)$, where $\aaa$ is a divisor of
degree either $1$ or $2$ on $E$.
\end{lem}

\proof
Observe that $h^0(\Oo_X(A_Z)) \leq 2$, since otherwise
\[
  3\leq h^0(\Oo_X(A_Z)) =h^0(\Oo_X(H-B_Z))
\]
implies that $B_Z$ is a line containing $Z$ which is impossible since
$Z$ spans a plane.  On the other hand, the Riemann-Roch for
$\Oo_X(A_Z)$ gives
\begin{equation} 
  \label{h0AZ}
  h^0(\Oo_X(A_Z)) = \frac{1}{2}(A^2_Z -A_Z\cdot K_X) +h^1(\Oo_X(A_Z)).
\end{equation}
From the vertical sequence in \eqref{d:Zextdest} we obtain
\[
  -1 =c_2(\Ff_Z) = A_Z\cdot (K_X+H-A_Z) +\deg(Z').
\]
This together with \eqref{h0AZ} implies
\[
  h^0(\Oo_X(A_Z))
  = \frac{1}{2}(A_Z\cdot H +1 +\deg(Z'))+h^1(\Oo_X(A_Z)).
\]
This identity and the above upper bound $h^0(\Oo_X(A_Z))\leq 2$ give
the following possibilities: 
\begin{enumerate}[label={\rm \arabic*)}]
  \listspace
\item
  $h^0(\Oo_X(A_Z)) =1$ and $A_Z\cdot H =1$, $\deg(Z')=0$,
\item
  $h^0(\Oo_X(A_Z)) =2$ and $A_Z\cdot H =2$, $\deg(Z')=1$,
\item
  $h^0(\Oo_X(A_Z)) =2$ and $A_Z\cdot H =3$, $\deg(Z')=0$.
\end{enumerate}
The third one can not hold, since it implies that $B_Z$ must consist
of two rulings of $X$ that contain $Z$.  This forces the two rulings
to be contained in the plane spanned by $Z$.

The possibility 1) (resp. 2)) implies that $A_Z$ is a (resp. the union
of two) ruling(s) of $X$.  Hence the assertion of the lemma.
\qed

Before we proceed, let us recall that the embedding 
$X\subset \PP(V^\ast)$ is defined by $\Oo_X(H)$, where 
\begin{equation}
  \label{formHrecall}
  H=\Gamma_o +\pi^\ast(D),
\end{equation}
with $o$ and $D$ being respectively a point and a divisor of degree $2$
on $E$.  With this in mind, we can now identify all the ingredients
involved in the diagram \eqref{d:Zextdest} for each type of subcycle
in \eqref{subcycles}.

\begin{pro}
Let $Z$ be one of the degree $4$ subcycles of $Z_s$ appearing in
\eqref{subcycles}.
\begin{enumerate}[label={\rm \arabic*)}]
  \listspace
\item 
  If $Z=Z^e_s$, $e\in D_s$, then the destabilizing sequence
  \eqref{FZdestseq} has the form
  \begin{equation*}
    \label{FZedestseq}
    0 \lra \Oo_X(\pi^\ast(\aaa_e)) \lra \Ff_{Z^e_s}
    \lra \Jj_{z_e}(K_X+\Gamma_e) \lra 0,
  \end{equation*}
  where $z_e$ is the point $\Gamma_e \cdot{\Gamma_e}$ and
  $\aaa_e$ is a divisor of degree $2$ on $E$ determined by the 
  linear equivalence
  \[
    \aaa_e \sim D + o -e
  \]
  with $D$ as in \eqref{formHrecall}.
  \item
  If $Z=Z^{e\cdot e'}_s$, $e\neq e'\in D_s$, then the destabilizing
  sequence \eqref{FZdestseq} has the form
  \begin{equation*}
    \label{FZee'destseq}
    0 \lra \Oo_X(l_{x}) \lra \Ff_{Z^{e\cdot e'}_s}
    \lra \Oo_X(K_X+H-l_{x}) \lra 0, 
  \end{equation*}
  where $l_{x}$ is the ruling of $X$ passing through a point 
  $x\in Z_s^{e\cdot e'}$.  Furthermore, there is a unique divisor
  $R_{e\cdot e'} \in |H-l_{x} |$ passing through $Z_s^{e\cdot e'} $ and
  subject to one of the following properties:
  \begin{itemize}
    \listspace
  \item
    either $e\cdot e'$ is a unique point of $Z_s^{e\cdot e'} $ lying on
    the ruling $l_{e\cdot e'}$ and then $x=e\cdot e'$ and $R_{e\cdot e'}$
    is a smooth elliptic curve of degree $4$ and $Z_s^{e\cdot e'}$ is its
    plane section,
  \item
    or $x\neq e\cdot e'$, then the divisor 
    $R_{e\cdot e'}=l_{e\cdot e'}+\Gamma_c$, for some 
    $c \in D_s \setminus\{e,e'\}$, the ruling $l_{e\cdot e'} $ passes
    through two points $e\cdot e', c'\cdot c''$ of $Z_s^{e\cdot e'}$,
    where $\{c',c''\} = D_s \setminus\{e,e',c\}$, and the point 
    $x=c\cdot c'$ or $c\cdot c''$.
  \end{itemize}
\end{enumerate}
\end{pro}

\proof
To prove 1), set $A_e:=H-\Gamma_e$ and observe that it has the form
$A_e =\pi^\ast(\aaa_e)$,
where $\aaa_e$ is a divisor of degree $2$ on $E$.  Since the subcycle
$Z=Z_s^e$ lies on the curve $\Gamma_e$, we have
\[
  H^0(\Jj_{Z_s^e}(H-\pi^\ast(\aaa_e)))
  = H^0(\Jj_{Z_s^e}(\Gamma_e)) \neq 0.
\] 
On the other hand, the defining sequence \eqref{Zext} for $Z=Z_s^e$,
tensored with $\Oo_X(-\pi^\ast(\aaa_e))$, gives
\[
  0 \lra H^0(\Ff_{Z_s^e}(-\pi^\ast(\aaa_e)))
  \lra H^0(\Jj_{Z_s^e}(\Gamma_e))
  \lra H^1(\Oo_X(K_X-\pi^\ast(\aaa_e))).
\]
Since by Serre duality 
$H^1(\Oo_X(K_X-\pi^\ast(\aaa_e))) =H^1(\Oo_X(\pi^\ast(\aaa_e)))^\ast=0$, 
we deduce
\[
  H^0(\Ff_{Z_s^e}(-\pi^\ast(\aaa_e)))
  \cong H^0(\Jj_{Z_s^e}(\Gamma_e)) \cong \CC.
\]
Hence $\Ff_{Z_s^e}(-\pi^\ast(\aaa_e))$ has, up to a nonzero scalar
multiple, a unique nonzero global section.  Furthermore, the scheme of
zeros of this section is obviously $0$-dimensional.  Its Koszul
sequence tensored with $\Oo_X(\pi^\ast(\aaa_e))$ gives
\[
  0 \lra \Oo_X(\pi^\ast(\aaa_e)) \lra \Ff_{Z_s^e}
  \lra \Jj_{Z'}(K_X+\Gamma_e) \lra 0,
\]
where $Z'$ is a single point, see the proof of Lemma \ref{l:AZ}.  Call
this point $z_e$.  From \eqref{d:Zextdest} it follows that 
$z_e \in \Gamma_e$.  To identify it, we restrict the diagram
\eqref{d:Zextdest} to the curve $\Gamma_e$ and obtain the identity
\begin{equation}
  \label{tauident}
  \Oo_{\Gamma_e}(\aaa_e + z_e)
  = \Oo_X(K_X)\otimes\Oo_{\Gamma_e}(Z_s^e),
\end{equation}
where we tacitly use the identification of $\Gamma_e$ with
$E$.  Furthermore, we have 
\[
  \Oo_{\Gamma_e}(\aaa_e) =\Oo_{\Gamma_e}(H-\Gamma_e)
  \textq{and}
  \Oo_{\Gamma_e}(Z_s^e) =\Oo_{\Gamma_e}(H-K_X),
\]
where the first equality is the definition of the divisor $\aaa_e$
and the second comes from realizing $Z_s^e$ as the complete
intersection of $\Gamma_e$ with a smooth curve in $|H-K_X|$ containing
$Z_s$.  Substituting into \eqref{tauident}, we obtain 
\[
  \Oo_{\Gamma_e}(z_e)\otimes\Oo_{\Gamma_e}(-\Gamma_e)
  = \Oo_{\Gamma_e}
\] 
or, equivalently, $z_e\sim\Gamma_e \cdot{\Gamma_e}$. This together
with $h^0( \Oo_{\Gamma_e}(z_e))=1$ imply the equality 
$z_e =\Gamma_e \cdot{\Gamma_e}$.

The linear equivalence asserted in 1) of the proposition follows
from writing
\[
\pi^\ast(\aaa_e) +\Gamma_e \sim H \sim \Gamma_o +\pi^\ast(D),
\]
where the last equivalence is \eqref{formHrecall}.  Hence 
$\aaa_e \sim D +o-e$ as divisors of $E$.

We turn now to the part 2) of the proposition.  We know that the
subcycle $Z_s^{e\cdot e'}$ does not lie on any of the curves of the
family $\{\Gamma_{b} \}_{b\in E}$.  This together with the first part
of the proof and Lemma \ref{l:AZ} implies that the destabilizing
sequence \eqref{FZdestseq} for $\Ff_{Z_s^{e\cdot e'}}$ has the form
\[
  0 \lra \Oo_X(l_x) \lra \Ff_{Z_s^{e\cdot e'}}
  \lra \Oo_X(K_X+H-l_x) \lra 0,
\]
where $l_x$ is the ruling passing through some point $x\in X$. The slanted arrow in the
upper right corner of \eqref{d:Zextdest} tells us that there is an effective divisor, call
it $R_{e\cdot e'}$, in $|H-l_x |$ passing through 
$Z_s^{e\cdot e'}$. The uniqueness of this divisor follows from $h^0(\Jj_{Z_s^{e\cdot e'}}(H-l_x)) =h^0(\Ff_{Z_s^{e\cdot e'}}(-l_x)) =1$, where the first equality comes from the horizontal sequence in \eqref{d:Zextdest} and the second one from the destabilizing sequence above.

 It
remains to analyse the properties of the divisor $R_{e\cdot e'}$ as well as the position of the point $x$. From the previous paragraph, we know already that
 $R_{e\cdot e'}$ is a unique effective divisor
in  $|H-l_x |$ passing through 
$Z_s^{e\cdot e'}$. From Lemma \ref{planePee'} we also know that any
divisor in $|H|$ passing through $Z_s^{e\cdot e'}$ must also contain
the ruling $l_{e\cdot e'}$.  Hence, if $l_x\neq l_{e\cdot e'}$, then
the divisor $R_{e\cdot e'}$ has the form
\[
  R_{e\cdot e'} = l_{e\cdot e'} +\Gamma_b,
\]
for some $b\in E$, and it must contain $Z_s^{e\cdot e'}$.  But
$l_{e\cdot e'}$ is allowed to contain at most two points of
$Z_s^{e\cdot e'}$, since no three points in $Z_s^{e\cdot e'}$ are
colinear.  Therefore, $\Gamma_b$ must contain a subscheme
$Z_{\Gamma_b}\subset Z_s^{e\cdot e'}$ consisting of at least two
points.  Restricting now the diagram \eqref{d:Zextdest} to $\Gamma_b$,
we obtain
\vspace{15.5ex}
\begin{equation*}
  \label{d:Zee'extdestGb}
  \begin{tikzpicture}[overlay,every node/.style={draw=none},
    ->,inner sep=1.1ex]
    \matrix [draw=none,row sep=3.5ex,column sep=3.5ex]
    {
      && \node (02) {$0$}; \\
      && \node (12) {$\Oo_{\Gamma_b }(l_x)$}; \\
      \node (20) {$0$};
      & \node (21) {$\Oo_X(K_X)\otimes\Oo_{\Gamma_b }(Z_{\Gamma_b})$};
      & \node (22) {$\Ff_Z\otimes\Oo_{\Gamma_b}$};
      & \node (23) {$\Oo_X(H)\otimes\Oo_{\Gamma_b}(-Z_{\Gamma_b})$};
      & \node (24) {$0$}; \\
      && \node (32) {$\Jj_{Z'}(K_X+H-A_Z)$}; \\      
      && \node (42) {$0$}; \\      
    };
    \path
      (02) edge (12) 
      (12) edge (22)
      (22) edge (32)
      (32) edge (42)

      (12) edge (23)
      (21) edge (32)
      
      (20) edge (21)
      (21) edge (22)
      (22) edge (23)
      (23) edge (24)
    ;
  \end{tikzpicture}
  \vspace{17ex}
\end{equation*}
where the slanted arrows now are zero morphisms.  Hence
\[
  \Oo_{\Gamma_b }(l_x)
  = \Oo_X(K_X)\otimes\Oo_{\Gamma_b }(Z_{\Gamma_b}).
\]
This implies
\[
  l_x|_{\Gamma_b}
  \sim Z_{\Gamma_b} +K_X|_{\Gamma_b}
  = Z_{\Gamma_b} -\Gamma_b|_{\Gamma_b}.
\]
In particular, $\deg(Z_{\Gamma_b})=2$. Hence, the remaining two point
of $Z^{e\cdot e'}$ must lie on $l_{e\cdot e'}$. 

To go further, recall that
\[
  Z^{e\cdot e'} =e\cdot e' +c\cdot c' +c\cdot c'' +c'\cdot c'',
\]
where $\{c,c',c''\}=D_s \smallsetminus \{e,e'\}$. Assume now that the
ruling $l_{e\cdot e'}$ passes through $c'\cdot c''$ (in addition to
the point $e\cdot e'$).  Then $\Gamma_b$ passes through $c\cdot c'$
and $c\cdot c''$ and hence $\Gamma_b$ must be $\Gamma_c$.

The curve $\Gamma_c$ intersects the plane $\Pi_{e\cdot e'}$, the span
of $Z^{e\cdot e'}$, along the points situated on the line 
$\langle c\cdot c',c\cdot c'' \rangle$, the span of $c\cdot c'$ and
$c\cdot c''$.  Namely, we have
\[
  \Gamma_c \bigcap \Pi_{e\cdot e'}
  = \Gamma_c\bigcap \langle c\cdot c',c\cdot c'' \rangle
  = c\cdot c' +c\cdot c'' +l_{e\cdot e'}\cdot \Gamma_c.
\]
Since the ruling $l_x$ must pass through one of these points, we
deduce that $x$ is either $c\cdot c'$ or $c\cdot c''$ as asserted.

We now assume that $e\cdot e'$ is the only point of $Z_s^{e\cdot e'}$
lying on the ruling $l_{e\cdot e'}$.  Then the preceding argument tells
us that $l_x =l_{e\cdot e'}$ and 
$R_{e\cdot e'} \in |H-l_{e\cdot e'}|$ is a unique divisor passing through 
$Z_s^{e\cdot e'}$.  Let us assume it to be reducible.  Then it has the
form
\[
  R_{e\cdot e'} = l_y + \Gamma_b.
\]
Running the argument involving the previous diagram, we deduce that
the ruling $l_y$ must contain two points of $Z_s^{e\cdot e'}$.  Hence
$l_y$ is contained in the plane $\Pi_{e\cdot e'}$ and then it
coincides with $l_{e\cdot e'}$.  This contradicts the assumption that
$l_{e\cdot e'}$ contains only one point of the subcycle 
$Z_s^{e\cdot e'}$.  Hence $R_{e\cdot e'}$ is irreducible and is a
smooth section of $\pi:X\to E$. Its degree 
$H\cdot R_{e\cdot e'}=H\cdot (H-l_{e\cdot e'})=4$. Since the
intersection 
$R_{e\cdot e'}\bigcap\Pi_{e\cdot e'}\supset Z_s^{e\cdot e'}$, we
deduce the equality
\[
  R_{e\cdot e'} \bigcap \Pi_{e\cdot e'} = Z_s^{e\cdot e'} .
\]
\qed

Now we define a functor $\FFF$ in \eqref{functor} on the level of
objects by
\[
\begin{aligned}
  &\FFF(\Pi_e)
  = \{ 0 \lra \Oo_X(K_X) \lra \Ff_{Z_s^e} \lra \Jj_{Z_s^e}(H) \lra 0 \},
  \textq{for every} e\in D_s, \\
  &\FFF(\Pi_{e\cdot e'})
  = \{ 0 \lra \Oo_X(K_X) \lra \Ff_{Z_s^{e\cdot e'}}
  \lra \Jj_{Z^{e\cdot e'}}(H) \lra 0 \},
  \textq{for every}  e\neq e'\in D_s.
\end{aligned}
\]
The further study of $\FFF$ and related topics will be considered
elsewhere.

\small{
\pagestyle{plain} 

}

\bigskip

\begin{tabular}[r]{p{22ex}rp{3ex}r}\small
  &  Daniel Naie && Igor Reider \\
  & Universit\'e d'Angers && Universit\'e d'Angers \\
  & 2, Bd Lavoisier && 2, Bd Lavoisier \\
  & France && France \\
  &daniel.naie@univ-angers.fr && igor.reider@univ-angers.fr
\end{tabular}

\end{document}